\documentclass[11pt]{amsart}

\pdfoutput=1 % for arXiv submission

\usepackage{pictures}
\usepackage{fntsymbols}
\usepackage{ourdiagrams} % changed from diagrams to ourdiagrams for arXiv submission

\makeatletter

\def\DeclareStructure#1#2{%
	\expandafter\def\csname Declare#1\endcsname##1{%
		\expandafter\def\csname structure:#1:\detokenize{##1}\endcsname
	}%
	\protected\def#2##1{%
		\ifcsname structure:#1:\detokenize{##1}\endcsname
			\message{^J===CATS===^JAsked for ##1 -- exists!^J}%
			\csname structure:#1:\detokenize{##1}\expandafter\endcsname
		\else
			\message{^J===CATS===^JAsked for ##1 -- using generic!^J}%
			\def\next{\csname structure-generic:#1\endcsname{##1}}%
			\expandafter\next
		\fi
	}%
	\@ifnextchar:%
		{\expandafter\def\csname structure-generic:#1\expandafter\endcsname\@gobble}%
		{\expandafter\let\csname structure-generic:#1\expandafter\endcsname\@firstofone}%
}

\DeclareStructure{Set}{\set}{\mathrm}
\DeclareStructure{Category}{\cat}:#1{\cat@generate{#1}}
\DeclareStructure{Bicategory}{\ccat}{\mathbf}

\newtoks\cat@name
\newtoks\cat@upper

\def\cat@generate#1{\begingroup
	\countdef\parser@state=0
	\countdef\parser@newstate=1
	\parser@state=0		% 0: active or macro; 1: lowercase or symbol; 2: uppercase
	\cat@name={}%        % a token list with the generated name
	\cat@upper={}%		% a token list with uppercase characters
	\cat@parse#1\relax
	\edef\next{\unexpanded{\DeclareCategory{#1}}{\the\cat@name}}\next
	\expandafter
\endgroup\the\cat@name}

\def\cat@parse#1{%
	\ifx\relax#1\relax
		\parser@newstate=4
	\else\ifcat\noexpand#1\relax\relax
		\parser@newstate=1
	\else\ifnum`#1>`Z\relax
		\parser@newstate=3
	\else\ifnum`#1<`A\relax
		\parser@newstate=3
	\else
		\parser@newstate=2
	\fi\fi\fi\fi
	\ifnum\parser@state=\parser@newstate\else
		\ifcase\parser@state
		\or
			\ifnum\parser@newstate<4
				\edef\next{\cat@name={\the\cat@name\/}}\next
			\fi
		\or
			\edef\next{\cat@name={%
				\the\cat@name\noexpand\mathcal{\the\cat@upper}%
			}}\next
			\cat@upper{}%
		\fi
		\parser@state\parser@newstate
	\fi
	\ifnum\parser@state=2
		\edef\next{\cat@upper={\the\cat@upper#1}}\next
	\else
		\edef\next{\cat@name={\the\cat@name\noexpand#1}}\next
	\fi
	\ifnum\parser@state<4	\expandafter\cat@parse\fi
}

\newcommand*{\Kar}{\mathrm{Kar}}

\DeclareCategory{RTfunc}{\mathcal{RT}}

\makeatother
\def\coev{\mathrm{coev}}
\def\ev{\mathrm{ev}}

\makeatletter

\begingroup
	\def\definering#1#2{\protected\gdef#1{\@ifstar{#2^*}{#2}}}%
	\definering\scalars\Bbbk
	\definering\Fld{\mathbb F}
	\definering\Z{\mathbb Z}
	\definering\Q{\mathbb Q}
	\definering\R{\mathbb R}
	\definering\C{\mathbb C}
	\definering\Zq{\mathbb Z[q^{\pm1}]}
	\definering\Zhtq{\mathbb Z[h,t,q^{\pm1}]}
\endgroup

\newcommand*\LieSL[1][2]{\mathfrak{sl}_{#1}}

\newcommand*{\Uqsl}[1][2]{\mathcal U_q(\mathfrak{sl}_{#1})}
\newcommand*{\Uqgl}[1][2]{\mathcal U_q(\mathfrak{gl}_{#1})}

\def\Disk{\@ifnextchar^{\mathbb D}{\mathbb D^2}}
\def\Ball{\@ifnextchar^{B}{B^3}}
\def\Ann{\mathbb A}

\def\mfld{\@ifnextchar^{M}{M^n}}%

\def\setcomponent@#1{%
	\@ifnextchar[{\setcomponent@weight#1}{\setcomponent@weight#1[]}
}
\def\setcomponent@weight#1[#2]#3{%
	\setcomponent@draw#1[#2]#3||\relax
}
\def\setcomponent@draw#1[#2]#3|#4|#5\relax{%
	\mathcal{#1}^{#4}_{#3}%
	\ifx\relax#2\relax\else(#2)\fi
}

\def\tangles@#1#2|#3|#4\relax{%
	{#1}^{#3}_{#2}%
}

\def\tangles#1{\tangles@{\mathcal T\!\mathit{an}}#1||\relax}
\def\ftangles#1{\tangles@{\mathcal{FT}\!\mathit{an}}#1||\relax}
\def\dtangles#1{\tangles@{\widetilde{\mathcal T}\!\mathit{an}}#1||\relax}

\def\gmatchings{\setcomponent@{GM}}

\def\revmatching#1{#1^!}

\DeclareMathOperator{\id}{id}		%	identity morphism
\DeclareMathOperator{\im}{im}		%	image
\def\tensor{\@ifstar\tensor@{\tensor@\otimes}}
\def\tensor@#1#2{\let\tensor@symb#1\tensor@@#2,\tensor@end,}
\def\tensor@@#1,{\tensor@expand#1\relax\tensor@@@}
\def\tensor@@@#1,{%
	\ifx\tensor@end#1\relax\else
		\tensor@symb\tensor@expand#1\relax
	\expandafter\tensor@@@\fi
}
\def\tensor@expand#1#2\relax{\ifx*#1#2\else#1_{#2}\fi}

\DeclareMathOperator\dcbtensor\heartsuit
\DeclareMathOperator\bdcbtensor\spadesuit

\def\utimes{\mytimes@\otimes}
\def\udtimes{\mytimes@{\hat\otimes}}

\def\mytimes@#1#2{%
	\mathchoice
		{\mytimes@@{#1}{#2}}%
		{\mytimes@@@{#1}{#2}\scriptstyle}%
		{\mytimes@@@{#1}{#2}\scriptscriptstyle}%
		{\mytimes@@@{#1}{#2}\scriptscriptstyle}%
}
\def\mytimes@clap#1#2{%
	\mathchoice
		{\mytimes@@clap{#1}{#2}}%
		{\mytimes@@@{#1}{#2}\scriptstyle}%
		{\mytimes@@@{#1}{#2}\scriptscriptstyle}%
		{\mytimes@@@{#1}{#2}\scriptscriptstyle}%
}

\def\mytimes@@#1#2{\mkern\thinmuskip\underset{\mkern-1mu#2}{#1}\mkern\thinmuskip}
\def\mytimes@@clap#1#2{\mkern\thinmuskip\underset{\mathclap{#2}}{#1}\mkern\thinmuskip}
\def\mytimes@@@#1#2#3{\mkern\thinmuskip{#1}_{\mkern-2mu\raisebox{-1pt}{$\m@th#3#2$}}\mkern\thinmuskip}

\def\shdw{\@ifnextchar[{\shdw@@}{\shdw@}}
\def\shdw@#1{\langle\!\langle#1\rangle\!\rangle}
\def\shdw@@[#1]#2{%
	\langle\!\langle#2\rangle\!\rangle
	\raisebox{-0.75ex}{$\scriptstyle\!#1$}%
}

\newcommand*{\TL}{\mathcal T\mkern-3mu\mathcal L}						% category
\newcommand*{\BN}{\cat{BN}}

\begingroup
	\def\defname#1{\expandafter\protected\expandafter\gdef\csname#1\endcsname}%
	\def\defmodulecats#1#2{%
		\defname{#1}{#2}%
		\defname{l#1}##1{%{##1\text{--}#2}%
			\raisebox{-0.35ex}{$\scriptstyle##1$}#2%
		}
		\defname{r#1}##1{%{#2\text{--}##1}%
			#2\raisebox{-0.35ex}{$\scriptstyle##1$}%
		}
		\defname{b#1}##1##2{%{##1\text{--}#2\text{--}##2}%
			\raisebox{-0.35ex}{$\scriptstyle##1$}%
			#2%
			\raisebox{-0.35ex}{$\scriptstyle##2$}%
		}
	}%
	\defmodulecats{Mod}{\cat{Mod}}
	\defmodulecats{PMod}{\cat{Proj}}
	\defmodulecats{ModG}{\mathrm{gr}\mkern-2mu\cat{Mod}}
	\defmodulecats{PModG}{\mathrm{gr}\cat{Proj}}
	\defname{Bimod}{\cat{Bimod}}
	\defname{gBimod}{\mathrm g\cat{Bimod}}
	\def\defbicats#1#2{%
		\defname{#1}{\ccat{#2}}%
		\defname{D#1}{\mathcal D(\ccat{#2})}%
		\defname{DD#1}{\mathcal D^-(\ccat{#2})}%
		\defname{g#1}{\@ifstar{\ccat{g#2}_0}{\ccat{g#2}}}%
		\defname{gD#1}{\mathcal D(\@ifstar{\ccat{g#2}_0)}{\ccat{g#2})}}%
		\defname{gDD#1}{\mathcal D^-(\@ifstar{\ccat{g#2}_0)}{\ccat{g#2})}}%

		\defname{e#1}{%
			\mathrlap{\hskip-0.15em\widetilde{\phantom{\rule{1.25em}{1.5ex}}}}%
			\ccat{#2}}%
		\defname{eD#1}{\mathcal D(%
			\mathrlap{\hskip-0.15em\widetilde{\phantom{\rule{1.25em}{1.5ex}}}}%
			\ccat{#2})}%
		\defname{eDD#1}{\mathcal D^-(%
			\mathrlap{\hskip-0.15em\widetilde{\phantom{\rule{1.25em}{1.5ex}}}}%
			\ccat{#2})}%
		\defname{eg#1}{%
			\mathrlap{\phantom{\ccat g}\hskip-0.15em\widetilde{\phantom{\rule{1.25em}{1.5ex}}}}%
			\@ifstar{\ccat{g#2}_0}{\ccat{g#2}}}%
		\defname{egD#1}{\mathcal D(%
			\mathrlap{\phantom{\ccat g}\hskip-0.15em\widetilde{\phantom{\rule{1.25em}{1.5ex}}}}%
			\@ifstar{\ccat{g#2}_0)}{\ccat{g#2})}}%
		\defname{egDD#1}{\mathcal D^-(%
			\mathrlap{\phantom{\ccat g}\hskip-0.15em\widetilde{\phantom{\rule{1.25em}{1.5ex}}}}%
			\@ifstar{\ccat{g#2}_0)}{\ccat{g#2})}}%
	}%
	\defbicats{BBimod}{Bimod}
	\defbicats{CCom}{Com}
	\defbicats{RRep}{Rep}
	\defbicats{DiagBimod}{DB}
\endgroup

\protected\def\Rep{\cat{Rep}}

\protected\def\Com{\mathrm{Com}}

\def\Com{\@ifnextchar^{\Com@}{\Com@@}}
\def\Com@^#1(#2){\mathit{Com}^{#1}(#2)}
\def\Com@@(#1){\mathit{Com}(#1)}

\def\HCom{\@ifnextchar^{\HCom@}{\HCom@@}}
\def\HCom@^#1(#2){\mathit{Com}^{#1}_{\!\raisebox{0.3ex}{$\scriptstyle/$}\mkern-2mu h}(#2)}
\def\HCom@@(#1){\mathit{Com}_{\!\raisebox{0.3ex}{$\scriptstyle/$}\mkern-2mu h}(#1)}

\def\SimplMod(#1){\cat{SMod}_{#1}}
\def\HoSimplMod(#1){\cat{HoSMod}_{#1}}

\newcommand*{\End}{\mathrm{End}}
\newcommand*{\Hom}{\mathrm{Hom}}

\newcommand*{\Id}{\mathrm{Id}}

\DeclareMathOperator\Tr{Tr}
\DeclareMathOperator\hTr{hTr}
\DeclareMathOperator\vTr{vTr}

\def\endofun{\Sigma}

\newcommand*{\FKh}{\mathcal{F}_{\!\textit{Kh}}}

\newcommand*{\Fweb}{\mathcal{F}_{\mkern-3mu\mathit w}}

\def\KhBracket{\@ifstar\KhBracketScaled\KhBracketSimple}
\def\KhCube{\@ifstar\KhCubeScaled\KhCubeSimple}

\newcommand*{\KhCubeScaled}[1]{\mathcal{I}\left(#1\right)}
\newcommand*{\KhCubeSimple}[1]{\mathcal{I}(#1)}
\newcommand*{\KhBracketScaled}[1]{\left\llbracket#1\right\rrbracket}
\newcommand*{\KhBracketSimple}[1]{\llbracket#1\rrbracket}

\def\wKhBracket{\@ifstar\wKhBracketScaled\wKhBracketSimple}
\def\wKhCube{\@ifstar\wKhCubeScaled\wKhCubeSimple}

\newcommand*{\wKhCubeScaled}[1]{\mathcal{I}_{\mathrm F}\left(#1\right)}
\newcommand*{\wKhCubeSimple}[1]{\mathcal{I}_{\mathrm F}(#1)}
\newcommand*{\wKhBracketScaled}[1]{\left\llbracket#1\right\rrbracket_{\mathrm F}}
\newcommand*{\wKhBracketSimple}[1]{\llbracket#1\rrbracket_{\mathrm F}}

\def\twosubs#1#2#3{%
	\ifx.#1%
		#3_{#2}%
	\else
		\raisebox{-0.85ex}{$\scriptstyle#1$}%
		#3%
		\raisebox{-0.85ex}{$\scriptstyle#2$}%
	\fi
}

\def\arcalg{\@ifnextchar|{\arcalg@comp}{\arcalg@whole}}
\def\arcalg@comp|#1#2|#3{\twosubs{#1}{#2}{\!\left(\arcalg@whole{#3}\right)\!}}
\def\arcalg@whole#1{H^{#1}}

\def\arcmod{\@ifnextchar|{\arcmod@comp}{\arcmod@whole}}
\def\arcmod@comp|#1#2|#3{\arcmod@whole{\ifx.#1\else\revmatching{#1}\fi#3#2}}
\def\arcmod@whole#1{\FKh(#1)}

\def\arcmap{\@ifnextchar|{\arcmap@comp}{\arcmap@whole}}
\def\arcmap@comp|#1#2|#3{\arcmap@whole{\ifx.#1\else\revmatching{#1}\fi#3#2}}
\def\arcmap@whole#1{\FKh(#1)}

\def\CKalg{\@ifnextchar|{\CKalg@comp}{\CKalg@whole}}
\def\CKalg@comp|#1#2|#3{\twosubs{#1}{#2}{\!\left(\CKalg@whole{#3}\right)\!}}
\def\CKalg@whole#1{A^{#1}}

\def\CKmod{\@ifnextchar|{\CKmod@comp}{\CKmod@whole}}
\def\CKmod@comp|#1#2|#3{\CKmod@whole{\ifx.#1\else\revmatching{#1}\fi#3#2}}
\def\CKmod@whole#1{C_{CK}(#1)}

\def\CKmap{\@ifnextchar|{\CKmap@comp}{\CKmap@whole}}
\def\CKmap@comp|#1#2|#3{\twosubs{#1}{#2}{(\CKmap@whole{#3})}}
\def\CKmap@whole#1{\@ifnextchar'{\varphi_{#1}}{\varphi_{#1}^{\mathstrut}}}

\def\arctimes#1{\mytimes@clap\otimes{\arcalg{#1}}}
\def\CKtimes#1{\mytimes@clap\otimes{\CKalg{#1}}}

\def\arcdtimes#1{\mytimes@clap{\hat\otimes}{\arcalg{#1}}}
\def\CKdtimes#1{\mytimes@clap{\hat\otimes}{\CKalg{#1}}}

\def\webtimes#1{\mytimes@clap\otimes{\webalg{#1}}}
\def\qtwebtimes#1{\mytimes@clap\otimes{\qtwebalg{#1}}}

\def\webalg{\@ifnextchar|{\webalg@comp}{\webalg@whole}}
\def\webalg@comp|#1#2|#3{\twosubs{#1}{#2}{\!\left(\webalg@whole{#3}\right)\!}}
\def\webalg@whole#1{\mathfrak{W}^{#1}}

\def\webmod{\@ifnextchar|{\webmod@comp}{\webmod@whole}}
\def\webmod@comp|#1#2|#3{\twosubs{#1}{#2}{\webmod@whole{#3}}}
\def\webmod@whole#1{\Fweb^\circ(#1)}

\def\webmap{\@ifnextchar|{\webmap@comp}{\webmap@whole}}
\def\webmap@comp|#1#2|#3{\twosubs{#1}{#2}{\webmap@whole{#3}}}
\def\webmap@whole#1{\Fweb^\circ(#1)}

\def\qtwebalg{\@ifnextchar|{\qtwebalg@comp}{\qtwebalg@whole}}
\def\qtwebalg@comp|#1#2|#3{\twosubs{#1}{#2}{\!\left(\qtwebalg@whole{#3}\right)\!}}
\def\qtwebalg@whole#1{\mathfrak{A}^{#1}}

\def\qtwebmod{\@ifnextchar|{\qtwebmod@comp}{\qtwebmod@whole}}
\def\qtwebmod@comp|#1#2|#3{\twosubs{#1}{#2}{\qtwebmod@whole{#3}}}
\def\qtwebmod@whole#1{C_{\mathfrak A}(#1)}

\def\qtwebmap{\@ifnextchar|{\qtwebmap@comp}{\qtwebmap@whole}}
\def\qtwebmap@comp|#1#2|#3{\twosubs{#1}{#2}{\qtwebmap@whole{#3}}}
\def\qtwebmap@whole#1{\@ifnextchar'{\varphi_{#1}}{\varphi_{#1}^{\mathstrut}}}

\def\Web{\@ifnextchar[\WebS\Web@}
\def\Web@{\cat{Web}}
\def\WebS[#1]{\Web@^{(#1)}}
\def\clWeb{\Web@^{\mathit{cl}}}

\def\Foam{\@ifnextchar[\FoamS\Foam@}
\def\Foam@{\cat{Foam}}
\def\FoamS[#1]{\Foam@^{(#1)}}
\def\clFoam{\Foam@^{\mathit{cl}}}

\def\cupfoam@#1{\@ifnextchar*{\cupfoam@get{#1}}{\@ifnextchar+{\cupfoam@getx{#1}}{\cupfoam@get{#1}}}}

\def\cupfoam@getx#1+#2{\cupfoam@get{#1#2}}
\def\cupfoam@get#1{\@ifnextchar[{\cupfoam@dots{#1}}{\cupfoam@nodots{#1}}}

\def\cupfoam@dots#1[#2]#3{#1(#3;#2)}
\def\cupfoam@nodots#1#2{#1(#2)}

\def\extpower{\@ifnextchar^\extpowersup\extpowernosup}%
\def\extpowersup^#1{\mbox{\Large $\wedge$}^{\mkern-6mu#1}}%
\def\extpowernosup{\mbox{\Large $\wedge$}}%

\def\Diff{\mathit{Diff}\@ifnextchar_{\!}{}}

\DeclareCategory{Web}{{\normalfont\textbf{Web}}}
\DeclareCategory{Foam}{{\normalfont\textbf{Foam}}}
\DeclareCategory{pFoam}{%
	\mskip-2mu{\mathrlap{\widetilde{\phantom{\normalfont\textbf{Fo}}}}}%
	\mskip 2mu{\normalfont\textbf{Foam}}}

\DeclareCategory{kMod}{{\normalfont\scalars{-}\textbf{Mod}}}

\def\flip#1{\@ifnextchar_{\mathrlap{\phantom{#1}^!}#1}{#1^!}}
\def\dotted#1{\@ifnextchar_{\mathrlap{\phantom{#1}^\bullet}#1}{#1^\bullet}}

\makeatother

\tikzset{%
	mid-arrow/.style 2 args={%
		decoration={%
			markings,
			mark=at position #1 with {\arrow{#2}}%
		},
		postaction={decorate}
	},
	->-/.style={mid-arrow={#1}{>}},           ->-/.default=0.5,
	-<-/.style={mid-arrow={#1}{<}},           -<-/.default=0.5,
	-^-/.style={mid-arrow={#1}{Bar[left]}},   -^-/.default=0.5,
	-v-/.style={mid-arrow={#1}{Bar[right]}},  -v-/.default=0.5,
}%

\colorlet{V1lineColor}{blue}
\colorlet{V2lineColor}{red!80!blue}
\colorlet{V0lineColor}{red!80!blue}
\colorlet{V2innerColor}{red!30}

\colorlet{BCshade}{black!15}

\tikzset{%
	V0/.style={%
		V0lineColor,
		line width=0.75pt,
		>={To[width=7pt, length=5pt]},
		dashed
	},
	V1/.style={%
		V1lineColor,
		line width=1.25pt,
		>={To[width=7pt, length=5pt]}
	},
	V2/.style={%
		V2innerColor,
		draw=V2lineColor,
		double=V2innerColor,
		double distance=1pt,
		line width=0.7pt,
		>={To[width=9pt]}
	},
	V1dot/.style={%
		V1lineColor,
		line width=1.25pt
	},
	V2dot/.style={%
		fill=V2innerColor,
		draw=V2lineColor,
		line width=0.75pt
	},
	BCedge/.style={%
		V1lineColor,
		line width=1.25pt,
		>={To[width=7pt, length=5pt]}
	},
	BCedge+/.style={%
		V2lineColor,
		line width=1.25pt,
		>={To[width=7pt, length=5pt]}
	}
}%

\colorlet{seamColor}{red!50!blue}
\colorlet{seamInnerColor}{seamColor!60!white}

\tikzset{%
	1facetFront/.style={blue!60!black, fill opacity=0.5, fill=blue!50},
	2facetFront/.style={red!80!blue,   fill opacity=0.5, fill=red!60},
	1facetBack/.style= {blue!60!black, fill opacity=0.5, fill=blue!30},
	2facetBack/.style= {red!80!blue,   fill opacity=0.5, fill=red!40},
	2facetInner/.style={red!80!blue,   fill opacity=0.5, fill=red!20},
	1facetLine/.style= {thin, draw=blue!60!black},
	2facetLine/.style= {thin, draw=red!80!blue},
	seam/.style={
		draw=seamColor,
		thin
	}
}

\tikzset{%
	web/bdry/.style={draw=black!80,thin},
	web/bg/.style={fill=black!10}
}

\tikzset{%
	arc hline/.style={
		color=black!60,
		thin
	},
	arc platform/.style={
		color=red,
		dashed,
		very thick
	},
	help line/.style={
		dashed,
		thin,
		color=black!60
	}
}%

\NewFntSymbol{PosCr}[r]{%
		\draw[->](0.5,-0.5) -- (-0.55,0.55);
		\draw[fntborder](-0.5,-0.5) -- (0.45,0.45);
		\draw[->](0.45,0.45) -- (0.55,0.55);
}%

\NewFntSymbol{NegCr}[r]{%
		\draw[->](-0.5,-0.5) -- (0.55,0.55);
		\draw[fntborder](0.5,-0.5) -- (-0.45,0.45);
		\draw[->](-0.45,0.45) -- (-0.55,0.55);
}%

\def\drawbdrypoint#1,{
	\ifx\relax#1\relax\else
		\webendpoint#1(\bdrypointx,0);
		\edef\bdrypointx{\numexpr\bdrypointx+1}%
		\expandafter\drawbdrypoint
	\fi}

\makeatletter
\def\webdrawdots#1#2(#3)#4{%
	\ifcase#1\or
		\fill[V1] (#3) circle[radius=3pt];
	\or
		\draw[V1,fill=white,line width=0.5pt] (#3) circle[radius=3pt];
	\fi
	\ifx#4;\expandafter\@gobble\else\expandafter\@firstofone\fi
	{\webdrawdots#2#4}%
}
\def\webendpoint{\@ifnextchar[\webendpoint@{\webendpoint@[above]}}
\def\webendpoint@[#1]#2#3(#4){%
	\filldraw[V#3dot] (#4) circle[radius=\the\dimexpr#3pt/2+1.5pt\relax]
		node[#1,text=V#3lineColor] {$\scriptstyle#2\vphantom+$}}
\makeatother

\DeclarePictureFamily{web}

\defwebpict{merge}(-1,-0.7)--(1,1) M=0.15 {%
	\draw[V2] (0,0) -- (0,1);
	\draw[V1] (0,0) arc[start angle=90, end angle=150, x radius=0.8, y radius=1.4];
	\draw[V1] (0,0) arc[start angle=90, end angle= 30, x radius=0.8, y radius=1.4]; 
}
\defwebpict{merge up}(-1,-0.7)--(1,1) M=0.15 {%
	\draw[V2,->-=0.7] (0,0) -- (0,1);
	\draw[V1,-<-=0.7] (0,0) arc[start angle=90, end angle=150, x radius=0.8, y radius=1.4];
	\draw[V1,-<-=0.7] (0,0) arc[start angle=90, end angle= 30, x radius=0.8, y radius=1.4];
}
\defwebpict{merge down}(-1,-1)--(1,0.7) M=-0.15 {%
	\draw[V2,->-=0.7] (0,0) -- (0,-1);
	\draw[V1,-<-=0.7] (0,0) arc[start angle=270, end angle=210, x radius=0.8, y radius=1.4];
	\draw[V1,-<-=0.7] (0,0) arc[start angle=270, end angle=330, x radius=0.8, y radius=1.4];
}

\defwebpict{split}(-1,-1)--(1,0.7) M=-0.15 {%
	\draw[V2] (0,0) -- (0,-1);
	\draw[V1] (0,0) arc[start angle=270, end angle=210, x radius=0.8, y radius=1.4];
	\draw[V1] (0,0) arc[start angle=270, end angle=330, x radius=0.8, y radius=1.4];
}
\defwebpict{split up}(-1,-1)--(1,0.7) M=-0.15 {%
	\draw[V2,-<-=0.7] (0,0) -- (0,-1);
	\draw[V1,->-=0.7] (0,0) arc[start angle=270, end angle=210, x radius=0.8, y radius=1.4];
	\draw[V1,->-=0.7] (0,0) arc[start angle=270, end angle=330, x radius=0.8, y radius=1.4];
}
\defwebpict{split down}(-1,-0.7)--(1,1) M=0.15 {%
	\draw[V2,-<-=0.7] (0,0) -- (0,1);
	\draw[V1,->-=0.7] (0,0) arc[start angle=90, end angle=150, x radius=0.8, y radius=1.4];
	\draw[V1,->-=0.7] (0,0) arc[start angle=90, end angle= 30, x radius=0.8, y radius=1.4];
}
 
\defwebpict{circle}(-0.5,-0.5)--(0.5,0.5) M=0 [#1]=[V1]{%
	\draw[#1] (0,0) circle[radius=0.4];
}

\defwebpict{vline}(-0.3,-1)--(0.3,1) M=0 [#1]=[V1] {%
	\draw[#1] (0,-1) -- (0,1);
}
\defwebpict{line}(-0.3,-0.5)--(0.3,0.5) M=0 [V#1#2]=[V1]{%
	\draw[V#1] (-0.2,-0.4) -- (0.2,0.4);
	\webdrawdots#20(0,0);
}

\defwebpict{T-shape}(-0.1,-0.8)--(0.8,0.8) M=0 {%
	\draw[V2] (0,-0.3) -- (0.7,-0.3)
	          (0, 0.3) -- (0.7, 0.3);
	\draw[V1] (0,-0.8) -- (0,0.8);
}

\defwebpict{line and cup}(-0.1,-0.8)--(0.8,0.8) M=0 {%
	\draw[V2] (0.7,0.3) arc[x radius=0.5, y radius=0.3,
	                        start angle=90, end angle=270];
	\draw[V1] (0,-0.8) -- (0,0.8);
}

\defwebpict{<bigon}(-0.5,-1)--(0.5,1) M=0 {%
	\draw[V2] (0.25pt,-0.5) .. controls (-0.5,-0.2) and (-0.5, 0.2) ..
	          (0.25pt, 0.5);
	\draw[V1] (0,-1) -- (0,-0.5)
	          (0, 0.5) -- (0, 1)
	          (0, 0.5) .. controls (0.5, 0.2) and (0.5,-0.2) .. (0,-0.5);
}

\defwebpict{bigon>}(-0.5,-1)--(0.5,1) M=0 {%
	\draw[V2] (-0.25pt,-0.5) .. controls ( 0.5,-0.2) and ( 0.5, 0.2) ..
	          (-0.25pt, 0.5);
	\draw[V1] (0,-1) -- (0,-0.5)
	          (0, 0.5) -- (0, 1)
	          (0, 0.5) .. controls (-0.5, 0.2) and (-0.5,-0.2) .. (0,-0.5);
}

\defwebpict{bubble}(-0.5,-1)--(0.5,1) M=0 {%
	\draw[V2] (0,-1) -- (0,-0.4) (0,0.4) -- (0,1);
	\draw[V1] (0,0) circle[radius=0.4];
}

\defwebpict{vert arcs}(-0.6,-0.6)--(0.6,0.6) M=0 {%
	\draw[V2,->-=0.85] (-0.5,-0.5) .. controls (-0.1,0) .. (-0.5, 0.5);
	\draw[V2,->-=0.85] ( 0.5, 0.5) .. controls ( 0.1,0) .. ( 0.5,-0.5);
}

\defwebpict{hor arcs}(-0.6,-0.6)--(0.6,0.6) M=0 {%
	\draw[V2,->-=0.85] (-0.5,-0.5) .. controls (0,-0.1) .. ( 0.5,-0.5);
	\draw[V2,->-=0.85] ( 0.5, 0.5) .. controls (0, 0.1) .. (-0.5, 0.5);
}

\defwebpict{circles}(-1.1,-0.5)--(1.1,0.5) M=0 [V#1#2]=[V1]{%
	\draw[V#1]
		(-0.6,-0.4) arc[radius=0.4, start angle=-90, end angle=90]
		( 0.6,-0.4) arc[radius=0.4, start angle=270, end angle=90];
	\draw[V#1lineColor,dotted]
		(-0.6,-0.4) arc[radius=0.4, start angle=270, end angle=90]
		( 0.6,-0.4) arc[radius=0.4, start angle=-90, end angle=90];
	\webdrawdots #200(-0.2,0)(0.2,0);
}

\defwebpict{merged}(-1.1,-0.5)--(1.1,0.5) M=0 [#1]=[V1] {%
	\draw[#1]
		(-0.6,-0.4) .. controls (-0.3,-0.4) and (-0.3,-0.2) ..
		( 0.0,-0.2) .. controls ( 0.3,-0.2) and ( 0.3,-0.4) .. (0.6,-0.4)
		(-0.6, 0.4) .. controls (-0.3, 0.4) and (-0.3, 0.2) ..
		( 0.0, 0.2) .. controls ( 0.3, 0.2) and ( 0.3, 0.4) .. (0.6, 0.4);
	\draw[#1lineColor, dotted]
		(-0.6,-0.4) arc[radius=0.4, start angle=270, end angle=90]
		( 0.6,-0.4) arc[radius=0.4, start angle=-90, end angle=90];
}

\defwebpict{inn bigon}(-0.4,-0.9)--(0.8,0.9) M=0 [#1,#2,#3]=[V0,V2,+]{%
	\fill[black!10]
		(0.5, 0.8) arc[radius=0.5, start angle=0, end angle=-90]
		           arc[radius=0.3, start angle=90, end angle=270]
					  arc[radius=0.5, start angle=90, end angle=0]
					  -| (0.7,0.8);
	\draw[V1] (0,-0.8) -- (0,0.8);
	\draw[#2] (0.5,-0.8) arc[radius=0.5, start angle= 0, end angle=90]
	          (0,  +0.3) arc[radius=0.5, start angle=-90, end angle=0];
	\ifx#3-\relax
		\draw[#2,<-] (0,-0.8) ++(30:0.5) -- ++(132:1pt);
	\else\ifx#3+\relax
		\draw[#2,<-] (0,0.8) ++(-30:0.5) -- ++(-132:1pt);
	\fi\fi
	\draw[#1] (0, 0.3) arc[radius=0.3, start angle=90, end angle=270];
}

\defwebpict{ext bigon}(-0.7,-0.9)--(0.8,0.9) M=0 [#1,#2,#3]=[V1,V0,+]{%
	\fill[black!10]
		(0.5, 0.8) arc[radius=0.5, start angle= 0, end angle=-90]
		           arc[radius=0.3, start angle=90, end angle=270]
					  arc[radius=0.5, start angle=90, end angle=0]
					  -| (-0.6,0.8);
	\draw[V1] (0,-0.8) -- (0,0.8);
	\draw[#2] (0,-0.3) arc[radius=0.5, start angle= 90, end angle=0]
	          (0, 0.3) arc[radius=0.5, start angle=-90, end angle=0];
	\ifx#3-\relax
		\draw[#2,<-] (0,-0.8) ++(30:0.5) -- ++(132:1pt);
	\else\ifx#3+\relax
		\draw[#2,<-] (0,0.8) ++(-30:0.5) -- ++(-132:1pt);
	\fi\fi
	\draw[#1] (0, 0.3) arc[radius=0.3, start angle= 90, end angle=270];
}

\defwebpict{ext cup}(-0.5,-0.9)--(0.8,0.9) M=0 [#1,#2]=[V2,+]{%
	\fill[black!10]
		(0.5, 0.8) .. controls (0.5, 0.3) and (0.3, 0.3) .. (0.3, 0)
	              .. controls (0.3,-0.3) and (0.5,-0.3) .. (0.5,-0.8) -| (-0.4, 0.8);
	\draw[V1] (0,-0.8) -- (0,0.8);
	\draw[#1] (0.5,0.8) .. controls (0.5, 0.3) and (0.3, 0.3) .. (0.3, 0)
	                    .. controls (0.3,-0.3) and (0.5,-0.3) .. (0.5,-0.8);
	\ifx#2+\relax
		\draw[#1,->] (0.3785, 0.2803) -- (0.4, 0.325);
	\else\ifx#2-\relax
		\draw[#1,->] (0.3785,-0.2803) -- (0.4,-0.325);
	\fi\fi
}

\defwebpict{inn cup}(-0.5,-0.9)--(0.8,0.9) M=0 [#1,#2]=[V2,+]{%
	\fill[black!10]
		(0.5, 0.8) .. controls (0.5, 0.3) and (0.3, 0.3) .. (0.3, 0)
	              .. controls (0.3,-0.3) and (0.5,-0.3) .. (0.5,-0.8) -| (0.7, 0.8);
	\draw[V1] (0,-0.8) -- (0,0.8);
	\draw[#1] (0.5,0.8) .. controls (0.5, 0.3) and (0.3, 0.3) .. (0.3, 0)
		                .. controls (0.3,-0.3) and (0.5,-0.3) .. (0.5,-0.8);
	\ifx#2+\relax
		\draw[#1,->] (0.3785, 0.2803) -- (0.4, 0.325);
	\else\ifx#2-\relax
		\draw[#1,->] (0.3785,-0.2803) -- (0.4,-0.325);
	\fi\fi
}

\defwebpict{nested}(-0.8,-0.8)--(0.8,0.8) M=0 [#1]=[V2]{%
	\draw[#1]
		(60:0.3) arc[radius=0.3, start angle=60, end angle=-60]
		(30:0.7) arc[radius=0.7, start angle=30, end angle=-30];
	\draw[#1lineColor,dotted]
		(60:0.3) arc[radius=0.3, start angle=60, end angle=300]
		(30:0.7) arc[radius=0.7, start angle=30, end angle=330];
}

\defwebpict{croissant}(-0.8,-0.8)--(0.8,0.8) M=0 [#1]=[V2]{%
	\draw[#1]
		( 60:0.3) .. controls ++(-60:0.3) and ++(-30:0.2) .. ( 30:0.7)
		(-60:0.3) .. controls ++( 60:0.3) and ++( 30:0.2) .. (-30:0.7);
	\draw[#1lineColor,dotted]
		(60:0.3) arc[radius=0.3, start angle=60, end angle=300]
		(30:0.7) arc[radius=0.7, start angle=30, end angle=330];
}

\defwebpict{arc res}(-0.5,-0.5)--(0.5,0.5) M=0 {%
	\draw[V1,->] (-0.5,-0.5) .. controls (-0.08,0) and (-0.08,0) .. (-0.5,0.5);
	\draw[V1,->] ( 0.5,-0.5) .. controls ( 0.08,0) and ( 0.08,0) .. ( 0.5,0.5);
}

\defwebpict{zip res}(-0.5,-0.5)--(0.5,0.5) M=0 {%
	\draw[V2] (0,-0.2) -- (0,0.2);
	\draw[V1] (-0.5,-0.5) .. controls (-0.3,-0.3) and (-0.2,-0.2) .. (0,-0.2)
	                      .. controls (0.2,-0.2) and (0.3,-0.3) .. (0.5,-0.5);
	\draw[V1,<->] (-0.5,0.5) .. controls (-0.3,0.3) and (-0.2,0.2) .. (0,0.2)
	                      .. controls (0.2,0.2) and (0.3,0.3) .. (0.5,0.5);
}

\defwebpict{blue-red-blue}(-0.6,-0.7)--(0.6,0.7) M=0 {%
	\fill[black!10] (0,-0.7) rectangle (0.6,0.7);
	\draw[V0] (0,-0.7) -- (0,0.7);
	\draw[V1] (-0.4,-0.7) .. controls (-0.2,0) .. (-0.4,0.7)
	          ( 0.4,-0.7) .. controls ( 0.2,0) .. ( 0.4,0.7);
}

\defwebpict{bigon-blue}(-0.6,-0.7)--(0.6,0.7) M=0 {%
	\fill[black!10] (0,-0.7) .. controls (0.5,0) .. (0,0.7) -| (0.6,-0.7);
	\draw[V0] (0,-0.7) .. controls (0.5,0) .. (0,0.7);
	\begin{scope}
		\clip (0.4,-0.7) .. controls (-0.1,0) .. (0.4,0.7) -- (0.6,0) -- cycle;
		\draw[V2] (0,-0.7) .. controls (0.5,0) .. (0,0.7);
	\end{scope}
	\draw[V1] (-0.4,-0.7) .. controls (-0.3,0) .. (-0.4,0.7)
	          ( 0.4,-0.7) .. controls (-0.1,0) .. ( 0.4,0.7);
}

\defwebpict{blue-bigon}(-0.6,-0.7)--(0.6,0.7) M=0 {%
	\fill[black!10] (0,-0.7) .. controls (-0.5,0) .. (0,0.7) -| (0.6,-0.7);
	\draw[V0] (0,-0.7) .. controls (-0.5,0) .. (0,0.7);
	\begin{scope}
		\clip (-0.4,-0.7) .. controls (0.1,0) .. (-0.4,0.7) -- (-0.6,0) -- cycle;
		\draw[V2] (0,-0.7) .. controls (-0.5,0) .. (0,0.7);
	\end{scope}
	\draw[V1] (-0.4,-0.7) .. controls (0.1,0) .. (-0.4,0.7)
	          ( 0.4,-0.7) .. controls (0.3,0) .. (0.4,0.7);
}

\defwebpict{zipped}(-0.6,-0.7)--(0.6,0.7) M=0 {%
	\fill[black!10] (0,-0.7) rectangle (0.6,0.7);
	\draw[V0] (0,-0.7) -- (0,0.7);
	\draw[V2] (0,-0.25) -- (0,0.25);
	\draw[V1] (-0.4, 0.7) .. controls (-0.2, 0.1) and (0.2, 0.1) .. (0.4, 0.7)
	          (-0.4,-0.7) .. controls (-0.2,-0.1) and (0.2,-0.1) .. (0.4,-0.7);
}

\defwebpict{joined circles}(-1.1,-0.5)--(1.1,0.5) M=0 [#1]=[]{%
	\draw[V0] (-0.6,0) -- (0.6,0);
	\draw[V2] (-0.2,0) -- (0.2,0);
	\draw[V1]
		(-0.6,-0.4) arc[radius=0.4, start angle=-90, end angle=90]
		( 0.6,-0.4) arc[radius=0.4, start angle=270, end angle=90];
	\draw[V1lineColor,dotted]
		(-0.6,-0.4) arc[radius=0.4, start angle=270, end angle=90]
		( 0.6,-0.4) arc[radius=0.4, start angle=-90, end angle=90];
	\webdrawdots #100(-0.3182,0.2828)(0.3182,0.2828);
}

\defwebpict{unjoined circles}(-1.1,-0.5)--(1.1,0.5) M=0 [#1]=[]{%
	\draw[V0] (0,-0.5) -- (0,0.5);
	\draw[V1]
		(-0.6,-0.4) arc[radius=0.4, start angle=-90, end angle=90]
		( 0.6,-0.4) arc[radius=0.4, start angle=270, end angle=90];
	\draw[V1lineColor,dotted]
		(-0.6,-0.4) arc[radius=0.4, start angle=270, end angle=90]
		( 0.6,-0.4) arc[radius=0.4, start angle=-90, end angle=90];
	\webdrawdots #100(-0.3182,0.2828)(0.3182,0.2828);
}

\defwebpict{joined merged}(-1.1,-0.5)--(1.1,0.5) M=0 {%
	\draw[V0] (0,-0.5) -- (0,0.5);
	\draw[V2] (0,-0.2) -- (0,0.2);
	\draw[V1]
		(-0.6,-0.4) .. controls (-0.3,-0.4) and (-0.3,-0.2) ..
		( 0.0,-0.2) .. controls ( 0.3,-0.2) and ( 0.3,-0.4) .. (0.6,-0.4)
		(-0.6, 0.4) .. controls (-0.3, 0.4) and (-0.3, 0.2) ..
		( 0.0, 0.2) .. controls ( 0.3, 0.2) and ( 0.3, 0.4) .. (0.6, 0.4);
	\draw[V1lineColor, dotted]
		(-0.6,-0.4) arc[radius=0.4, start angle=270, end angle=90]
		( 0.6,-0.4) arc[radius=0.4, start angle=-90, end angle=90];
}

\defwebpict{unjoined merged}(-1.1,-0.5)--(1.1,0.5) M=0 {%
	\draw[V0] (-0.6,0) -- (0.6,0);
	\draw[V1]
		(-0.6,-0.4) .. controls (-0.3,-0.4) and (-0.3,-0.2) ..
		( 0.0,-0.2) .. controls ( 0.3,-0.2) and ( 0.3,-0.4) .. (0.6,-0.4)
		(-0.6, 0.4) .. controls (-0.3, 0.4) and (-0.3, 0.2) ..
		( 0.0, 0.2) .. controls ( 0.3, 0.2) and ( 0.3, 0.4) .. (0.6, 0.4);
	\draw[V1lineColor, dotted]
		(-0.6,-0.4) arc[radius=0.4, start angle=270, end angle=90]
		( 0.6,-0.4) arc[radius=0.4, start angle=-90, end angle=90];
}
\NewFntSymbol{foam-1cup}(-0.5,-0.7)(0.5,0.2)<2.5>[r]{%
	\path[1facetBack] (-0.5,0)
			arc[start angle=180, end angle=0, x radius=0.5, y radius=0.15]
			.. controls (0.5,-1) and (-0.5,-1) .. (-0.5,0);
	\path[1facetFront] (-0.5,0)
			arc[start angle=180, end angle=360, x radius=0.5, y radius=0.15]
			.. controls (0.5,-1) and (-0.5,-1) .. (-0.5,0);
	\draw[1facetLine,line width=\fntlinewidth]
	    (0.0,0) ellipse[x radius=0.5, y radius=0.15]
			(0.5,0) .. controls (0.5,-1) and (-0.5,-1) .. (-0.5,0);
}
\NewFntSymbol{foam-1cup*}(-0.5,-0.7)(0.5,0.2)<2.5>[r]{%
	\path[1facetBack] (-0.5,0)
			arc[start angle=180, end angle=0, x radius=0.5, y radius=0.15]
			.. controls (0.5,-1) and (-0.5,-1) .. (-0.5,0);
	\path[1facetFront] (-0.5,0)
			arc[start angle=180, end angle=360, x radius=0.5, y radius=0.15]
			.. controls (0.5,-1) and (-0.5,-1) .. (-0.5,0);
	\draw[1facetLine,line width=\fntlinewidth]
	    (0.0,0) ellipse[x radius=0.5, y radius=0.15]
			(0.5,0) .. controls (0.5,-1) and (-0.5,-1) .. (-0.5,0);
	\fill[black] (0.2,-0.35) circle[radius=0.3ex];
}
\NewFntSymbol{foam-2cup}(-0.5,-0.7)(0.5,0.2)<2.5>[r]{%
	\path[2facetBack] (-0.5,0)
			arc[start angle=180, end angle=0, x radius=0.5, y radius=0.15]
			.. controls (0.5,-1) and (-0.5,-1) .. (-0.5,0);
	\path[2facetFront] (-0.5,0)
			arc[start angle=180, end angle=360, x radius=0.5, y radius=0.15]
			.. controls (0.5,-1) and (-0.5,-1) .. (-0.5,0);
	\draw[2facetLine,line width=\fntlinewidth]
	    (0.0,0) ellipse[x radius=0.5, y radius=0.15]
			(0.5,0) .. controls (0.5,-1) and (-0.5,-1) .. (-0.5,0);
}

\NewFntSymbol{foam-1cap}(-0.5,1)(0.5,-0.3)<2.5>[r]{%
	\path[1facetBack] (-0.5,0)
			arc[start angle=180, end angle=360, x radius=0.5, y radius=0.15]
			.. controls (0.5, 1) and (-0.5, 1) .. (-0.5,0);
	\draw[1facetLine,line width=0.5\fntlinewidth]
	    (0.5,0) arc[start angle=0, end angle=180, x radius=0.5, y radius=0.15];
	\path[1facetFront] (-0.5,0)
			arc[start angle=180, end angle=360, x radius=0.5, y radius=0.15]
			.. controls (0.5, 1) and (-0.5, 1) .. (-0.5,0);
	\draw[1facetLine,line width=\fntlinewidth]
			(0.5,0) .. controls (0.5,1) and (-0.5,1) .. (-0.5,0)
			arc[start angle=180, end angle=360, x radius=0.5, y radius=0.15] -- cycle;
}
\NewFntSymbol{foam-1cap*}(-0.5,1)(0.5,-0.3)<2.5>[r]{%
	\path[1facetBack] (-0.5,0)
			arc[start angle=180, end angle=360, x radius=0.5, y radius=0.15]
			.. controls (0.5, 1) and (-0.5, 1) .. (-0.5,0);
	\draw[1facetLine,line width=0.5\fntlinewidth]
	    (0.5,0) arc[start angle=0, end angle=180, x radius=0.5, y radius=0.15];
	\path[1facetFront] (-0.5,0)
			arc[start angle=180, end angle=360, x radius=0.5, y radius=0.15]
			.. controls (0.5, 1) and (-0.5, 1) .. (-0.5,0);
	\draw[1facetLine,line width=\fntlinewidth]
			(0.5,0) .. controls (0.5,1) and (-0.5,1) .. (-0.5,0)
			arc[start angle=180, end angle=360, x radius=0.5, y radius=0.15] -- cycle;
	\fill[black] (0.2,0.35) circle[radius=0.3ex];
}
\NewFntSymbol{foam-2cap}(-0.5,0.7)(0.5,-0.2)<2.5>[r]{%
	\path[2facetBack] (-0.5,0)
			arc[start angle=180, end angle=360, x radius=0.5, y radius=0.15]
			.. controls (0.5, 1) and (-0.5, 1) .. (-0.5,0);
	\draw[2facetLine,line width=0.5\fntlinewidth]
	    (0.5,0) arc[start angle=0, end angle=180, x radius=0.5, y radius=0.15];
	\path[2facetFront] (-0.5,0)
			arc[start angle=180, end angle=360, x radius=0.5, y radius=0.15]
			.. controls (0.5, 1) and (-0.5, 1) .. (-0.5,0);
	\draw[2facetLine,line width=\fntlinewidth]
			(0.5,0) .. controls (0.5,1) and (-0.5,1) .. (-0.5,0)
			arc[start angle=180, end angle=360, x radius=0.5, y radius=0.15] -- cycle;
}

\NewFntSymbol{foam-saddle-front}(-0.8,-0.2)(0.8,1.2)<2.5>[r]{%
	\path[2facetBack] (-0.5,0.2)
	    arc[start angle= 90, end angle=0, radius=0.2]
			arc[start angle=180, end angle=0,  x radius=0.3, y radius=0.5]
			arc[start angle=180, end angle=90, x radius=0.5, y radius=0.2]
			-- (0.8,1.2) -- (-0.5,1.2) -- cycle;
	\draw[2facetLine,line width=\fntlinewidth]
		(-0.5,0.8) -- (-0.5,1.2) -- (0.8,1.2) -- (0.8,0.2)
	    arc[start angle=90, end angle=127, x radius=0.5, y radius=0.2];
	\draw[2facetLine,line width=0.5\fntlinewidth]
		(-0.5,0.8) -- (-0.5,0.2) arc[start angle= 90, end angle=0, radius=0.2]
		( 0.3,0.0) arc[start angle=180, end angle=127, x radius=0.5, y radius=0.2];
	\draw[2facetFront,line width=\fntlinewidth] (-0.8,-0.2)
	    arc[start angle=270, end angle=360, x radius=0.5, y radius=0.2]
			arc[start angle=180, end angle=0,   x radius=0.3, y radius=0.5]
			arc[start angle=180, end angle=270, radius=0.2]
			-- (0.5,0.8) -- (-0.8,0.8) -- cycle;
}
\NewFntSymbol{foam-saddle-side}(-0.8,-1.2)(0.8,0.2)<2.5>[r]{%
	\path[2facetBack] (-0.5,0.2)
	    arc[start angle= 90, end angle=0, radius=0.2]
			arc[start angle=180, end angle=360, x radius=0.3, y radius=0.5]
			arc[start angle=180, end angle=90,  x radius=0.5, y radius=0.2]
			-- (0.8,-0.8) -- (-0.5,-0.8) -- cycle;
	\draw[2facetLine,line width=0.5\fntlinewidth]
		(0.5,-0.8) -- (-0.5,-0.8) -- (-0.5,-0.16);
	\path[2facetFront] (-0.8,-0.2)
	    arc[start angle=270, end angle=360, x radius=0.5, y radius=0.2]
			arc[start angle=180, end angle=360, x radius=0.3, y radius=0.5]
			arc[start angle=180, end angle=270, radius=0.2]
			-- (0.5,-1.2) -- (-0.8,-1.2) -- cycle;
	\draw[2facetLine, line width=\fntlinewidth] (-0.5,-0.16) -- (-0.5,0.2)
		arc[start angle=90,  end angle=0, radius=0.2]
		arc[start angle=360, end angle=270, x radius=0.5, y radius=0.2]
		-- (-0.8,-1.2) -- (0.5,-1.2) -- (0.5,-0.2)
		arc[start angle=270, end angle=180, radius=0.2]
		arc[start angle=180, end angle=90,  x radius=0.5, y radius=0.2]
		-- (0.8,-0.8) -- (0.5,-0.8)
		(-0.3,0) arc[start angle=180, end angle=360, x radius=0.3, y radius=0.5];
}

\NewFntSymbol{foam-detach-cup}(-0.4,-0.9)(1,0.9)<3>[r]{%
	\fill[1facetFront] (-0.3, 0.1) -- (0.3, 0.9) -- (0.3, -0.1) -- (-0.3, -0.9) -- cycle;
	\draw[1facetLine,line width=0.5\fntlinewidth]
		(-0.15,-0.7) -- (0.3,-0.1) -- (0.3,0.15);
	\draw[1facetLine,line width=\fntlinewidth]
		(-0.15,-0.7) -- (-0.3,-0.9) -- (-0.3,0.1) -- (0.3,0.9) -- (0.3,0.15);
	\fill[2facetBack] (0,0)
		arc[start angle=90, end angle=0, x radius=0.15, y radius=0.3] -| (0.9, 0.7)
		arc[start angle=90, end angle=180, x radius=0.5, y radius=0.2]
		arc[start angle=0, end angle=-90, x radius=0.4, y radius=0.5];
	\draw[2facetLine,line width=0.5\fntlinewidth] (0.15,-0.3) -- (0.6,-0.3);
	\draw[seam,line width=0.5\fntlinewidth]
		(0,0) arc[start angle=90, end angle=0, x radius=0.15, y radius=0.3];
	\fill[2facetFront] (0.4,0.5)
		arc[start angle=0, end angle=-90, x radius=0.4, y radius=0.5]
		arc[start angle=90, end angle=180, x radius=0.15, y radius=0.7]
		-| (0.6,0.3) arc[start angle=270, end angle=180, radius=0.2];
	\draw[2facetLine,line width=\fntlinewidth] (0.6,-0.3) -| (0.9,0.7)
		arc[start angle= 90, end angle=180, x radius=0.5, y radius=0.2]
		arc[start angle=180, end angle=270, radius=0.2] |- (-0.15,-0.7)
		(0.4,0.5) arc[start angle=0, end angle=-90, x radius=0.4, y radius=0.5];
	\draw[seam,line width=\fntlinewidth]
		(0,0) arc[start angle=90, end angle=180, x radius=0.15, y radius=0.7];
}
\NewFntSymbol{foam-attach-cup}(-0.4,-0.9)(1,0.9)<3>[r]{%
	\fill[1facetFront] (-0.3, 0.1) -- (0.3, 0.9) -- (0.3, -0.1) -- (-0.3, -0.9) -- cycle;
	\draw[1facetLine,line width=\fntlinewidth]
		(0.3,0.7) -- (0.3,0.9) -- (-0.3,0.1) -- (-0.3,-0.9) -- (0.27,-0.14);
	\draw[1facetLine,line width=0.5\fntlinewidth]
		(0.27,-0.14) -- (0.3,-0.1) -- (0.3,0.7);
	\fill[2facetBack] (0.15, 0.7) -| (0.9, -0.3)
		arc[start angle=90, end angle=180, x radius=0.5, y radius=0.2]
		arc[start angle=0, end angle=90, x radius=0.4, y radius=0.5]
		arc[start angle=270, end angle=360, x radius=0.15, y radius=0.7];
	\draw[2facetLine,line width=0.5\fntlinewidth]
		(0.9,-0.3) arc[start angle=90, end angle=180, x radius=0.5, y radius=0.2];
	\draw[seam,line width=0.5\fntlinewidth]
		(0.15,0.7) arc[start angle=360, end angle=270, x radius=0.15, y radius=0.7];		
	\draw[seam,line width=\fntlinewidth]
		(0.15,0.7) arc[start angle=360, end angle=325, x radius=0.15, y radius=0.7];
	\fill[2facetFront] (-0.15, 0.3) -| (0.6, -0.7)
		arc[start angle=270, end angle=180, radius=0.2]
		arc[start angle=0, end angle=90, x radius=0.4, y radius=0.5]
		arc[start angle=270, end angle=180, x radius=0.15, y radius=0.3];
	\draw[seam,line width=\fntlinewidth]
		(-0.15,0.3) arc[start angle=180, end angle=270, x radius=0.15, y radius=0.3];
	\draw[2facetLine,line width=\fntlinewidth]
		(0,0) arc[start angle=90, end angle=0, x radius=0.4, y radius=0.5]
		      arc[start angle=180, end angle=270, radius=0.2] |- (-0.15,0.3)
		(0.15,0.7) -| (0.9,-0.3)
			arc[start angle=90, end angle=125, x radius=0.5, y radius=0.2];
}

\DeclarePictureFamily{foam}

\def\foamdot(#1){\fill (#1) circle[radius=0.5ex]}
\def\foamdualdot(#1){\draw[line width=0.15ex] (#1) circle[radius=0.45ex]}

\def\foamorient(#1) Rx=#2 Ry=#3 S=#4 E=#5;{\begin{scope}[shift={(#1)}]
	\coordinate (arrowst) at (#5:#2 and #3);
	\ifnum #4 < #5
		\path (arrowst) ++(#5+90:#2 and #3) coordinate(arrowdir);
	\else
		\path (arrowst) ++(#5-90:#2 and #3) coordinate(arrowdir);
	\fi
	\draw[stealth-] ($(arrowst)!3pt!(arrowdir)$) -- (arrowst)
		arc[x radius=#2, y radius=#3, start angle=#5, end angle=#4];
\end{scope}}

\def\getlength#1(#2){%
  \path let \p{x}=(#2), \n{xlen}={veclen(\x{x},\y{x})}
  in \pgfextra{\xdef#1{\n{xlen}}}%
}

\def\foamarrow#1[#2]#3(#4){%
	\ifx\relax#3\relax
		\foamarrowx#1[#2]0({#4})%
	\else
		\foamarrowx#1[#2]#3({#4})%
	\fi\foamarrownext
}

\def\foamarrowx#1[#2]#3(#4){%
	\ifx<#1\relax
		\edef\foamarrownext{\unexpanded{\draw[#2, line width=0.75pt]}
			(#4) ++(#3:-2.5pt) arc[radius=5pt,
				start angle=\the\numexpr -90 + #3,
				end angle=\the\numexpr -25 + #3]
			(#4) ++(#3:-2.5pt) arc[radius=5pt,
				start angle=\the\numexpr 90 + #3,
				end angle=\the\numexpr 25 + #3]}%
	\else\ifx>#1\relax
		\edef\foamarrownext{\unexpanded{\draw[#2, line width=0.75pt]}
			(#4) ++(#3:2.5pt) arc[radius=5pt,
				start angle=\the\numexpr -90 + #3,
				end angle=\the\numexpr -155 + #3]
			(#4) ++(#3:2.5pt) arc[radius=5pt,
				start angle=\the\numexpr 90 + #3,
				end angle=\the\numexpr 155 + #3]}%
	\else\let\foamarrownext\relax\fi\fi
}

\def\foamdotsonsphere#1{%
	\ifcase 0#1\or
		\foamdot(0.2, 0.25);\or
		\foamdot(0.22, 0.25);\foamdot(0, 0.25);\or
		\foamdot(0.22, 0.25);\foamdot(0, 0.25);\foamdot(-0.22, 0.25);\or
		\foamdualdot(0.12,0.25);\or
		\foamdot(0.2,0.25);\foamdualdot(-0.05,0.25);
	\fi
}

\deffoampict{1sphere}(-0.6,-0.6)--(0.6, 0.6) M=0 [#1]=[] {%
	\fill[1facetBack] (0,0) circle[radius=0.5];
	\draw[1facetFront] (0,0) circle[radius=0.5];
	\draw[1facetLine] (-0.5,0)
		arc[start angle=180, end angle=360, x radius=0.5, y radius=0.15];
	\foamdotsonsphere{#1}
}

\deffoampict{2sphere}(-0.6,-0.6)--(0.6, 0.6) M=0 [#1]=[] {%
	\getlength\xradius({$(0.5,0) + (-1.7pt,0)$});
	\getlength\yradius({$(0,0.15) + (0,-1pt)$});
	\fill[2facetInner] (0,0) circle[radius=0.5];
	\draw[2facetFront](0,0) circle[radius=\xradius];
	\draw[2facetLine]
		(-\xradius,0)
			arc[start angle=180, end angle=360, x radius=\xradius, y radius=\yradius];
	\draw[2facetFront] (0,0) circle[radius=0.5];
	\draw[2facetLine]
		(-0.5 ,0) arc[start angle=180, end angle=360, x radius=0.5 , y radius=0.15];
	\foamdotsonsphere{#1}
}

\deffoampict{sphere and membrane}(-0.6,-0.6)--(0.6,0.6) M=0 [#1#2]=[>]{%
	\fill[1facetBack]  (0,0) circle[radius=0.5];
	\draw[2facetInner] (0,-0.85pt) ellipse[x radius=0.5, y radius=0.15];
	\draw[2facetFront] (0, 0.85pt) ellipse[x radius=0.5, y radius=0.15];
	\draw[1facetFront] (0,0) circle[radius=0.5];
	\draw[seam]
		(0.5,-0.85pt) arc[x radius=0.5, y radius=0.15, start angle=360, end angle=180]
		(0.5, 0.85pt) arc[x radius=0.5, y radius=0.15, start angle=360, end angle=180];
	\foamarrow#1[seam](0,-0.15);
	\foamdotsonsphere{#2}
}%

\deffoampict{sphere in plane}(-1.5,-0.6)--(1.5,0.6) M=0 [#1#2]=[>]{%
	\fill[1facetBack] (0.5,0)
			arc[radius=0.5, start angle=360, end angle=180]
			arc[x radius=0.5, y radius=0.15, start angle=180, end angle=360];
	\fill[1facetFront] (0.5,0)
			arc[radius=0.5, start angle=360, end angle=180]
			arc[x radius=0.5, y radius=0.15, start angle=180, end angle=360];
	\draw[1facetLine]
			(0.5,0) arc[radius=0.5, start angle=360, end angle=180];
	\fill[2facetInner] (0.5,-0.85pt)
			arc[x radius=0.5, y radius=0.15, start angle=0, end angle=360]
			-- ++(0.9,0.3) -- ++(-0.8,-0.6) -- ++(-2,0) -- ++(0.8, 0.6) -- ++(2,0)
			-- cycle;
	\draw[2facetLine] (0,-0.85pt)
		++(1.4,0.3) -- ++(-0.8,-0.6) -- ++(-2,0) -- ++(0.8, 0.6) -- ++(2,0) -- cycle;
	\fill[2facetFront] (0.5, 0.85pt)
			arc[x radius=0.5, y radius=0.15, start angle=0, end angle=360]
			-- ++(0.9,0.3) -- ++(-0.8,-0.6) -- ++(-2,0) -- ++(0.8, 0.6) -- ++(2,0)
			-- cycle;
	\draw[2facetLine] (0, 0.85pt)
		++(1.4,0.3) -- ++(-0.8,-0.6) -- ++(-2,0) -- ++(0.8, 0.6) -- ++(2,0) -- cycle;
	\draw[seam]
		(0,  -0.85pt) ellipse[x radius=0.5, y radius=0.15]
		(0.5, 0.85pt) arc[x radius=0.5, y radius=0.15, start angle=0, end angle=180];
	\fill[1facetBack] (0.5,0)
			arc[radius=0.5, start angle=0, end angle=180]
			arc[x radius=0.5, y radius=0.15, start angle=180, end angle=360];
	\fill[1facetFront] (0.5,0)
			arc[radius=0.5, start angle=0, end angle=180]
			arc[x radius=0.5, y radius=0.15, start angle=180, end angle=360];
	\draw[1facetLine] (0.5,0)
			arc[radius=0.5, start angle=0, end angle=180];
	\draw[seam]
		(0.5, 0.85pt) arc[x radius=0.5, y radius=0.15, start angle=360, end angle=180];
	\foamarrow#1[seam](0,-0.15);
	\foamarrow#1[2facetLine](-1,-0.3);
	\foamarrow#1[2facetLine]36(1.1,0.075);
	\foamdotsonsphere{#2}
}

\def\foamdotsonplane#1{%
	\ifcase 0#1\or
		\foamdot(0,0);\or
		\foamdot(-0.15,0);\foamdot(0.15,0);\or
		\foamdot(0,0);\foamdot(-0.25,0);\foamdot(0.25,0);\or
		\foamdualdot(0,0);\or
		\foamdualdot(-0.15,0);\foamdot(0.15,0);\or
		\foamdualdot(-0.15,0);\foamdualdot(0.15,0);
	\fi
}

\deffoampict{1plane}(-0.9,-0.5)--(0.9,0.5) M=0 [#1]=[]{%
	\draw[1facetFront]
		(0.85,0.4) -- ++(-0.45,-0.8) -- ++(-1.25,0) -- ++(0.45,0.8) -- cycle;
	\foamdotsonplane{#1}
}

\deffoampict{2plane}(-0.9,-0.5)--(0.9,0.5) M=0 [#1#2]=[x]{%
	\draw[2facetInner] ( 0.2pt,-0.85pt)
		++(0.85,0.4) -- ++(-0.45,-0.8) -- ++(-1.25,0) -- ++(0.45,0.8) -- cycle;
	\draw[2facetFront] (-0.2pt, 0.85pt)
		++(0.85,0.4) -- ++(-0.45,-0.8) -- ++(-1.25,0) -- ++(0.45,0.8) -- cycle;
	\foamarrow#1[2facetLine](-0.25,-0.4);
	\foamarrow#1[2facetLine]180(0.25,0.4);
	\foamdotsonplane{#2}
}

\deffoampict{1cup}(-0.7,-0.4)--(0.7,0.75) M=0.15 [#1#2]=[0]{%
	\fill[1facetBack] (-0.5,0.5)
		.. controls (-0.5,-0.5) and (0.5,-0.5) .. (0.5,0.5)
		arc[x radius=0.5, y radius=0.15, start angle=0, end angle=180];
	\fill[1facetFront] (-0.5,0.5)
		.. controls (-0.5,-0.5) and (0.5,-0.5) .. (0.5,0.5)
		arc[x radius=0.5, y radius=0.15, start angle=360, end angle=180];
	\draw[1facetLine] (-0.5,0.5) .. controls (-0.5,-0.5) and (0.5,-0.5) .. (0.5,0.5)
		(0,0.5) ellipse[x radius=0.5, y radius=0.15];
	\ifcase #1\or
		\foamdot(0.15,0);\or
		\foamdualdot(0.15,0);%
	\fi
	\foamarrow#2[1facetLine](0,0.35);
}

\deffoampict{2cup}(-0.7,-0.4)--(0.7,0.75) M=0.15 [#1#2]=[x]{%
	\getlength\innradx({$(0.5,0) + (-1.7pt,0)$});
	\getlength\innrady({$(0,0.15) + (0,-1pt)$});
	\fill[2facetInner] (-0.5,0.5)
		.. controls (-0.5,-0.5) and (0.5,-0.5) .. (0.5,0.5)
		arc[x radius=0.5, y radius=0.15, start angle=0, end angle=180];
	\fill[2facetBack] (-\innradx,0.5)
		.. controls (-\innradx,-\innradx) and (\innradx,-\innradx)
		.. (\innradx,0.5)
		arc[x radius=\innradx, y radius=\innrady, start angle=0, end angle=180];
	\draw[2facetLine] (-\innradx,0.5)
		.. controls (-\innradx,-\innradx) and (\innradx,-\innradx)
		.. (\innradx,0.5)
		(0,0.5) ellipse[x radius=\innradx, y radius=\innrady];
	\fill[2facetInner] (-\innradx,0.5)
		.. controls (-\innradx,-\innradx) and (\innradx,-\innradx)
		.. (\innradx,0.5)
		arc[x radius=\innradx, y radius=\innrady,
		    start angle=360, end angle=180];
	\fill[2facetFront] (-0.5,0.5)
		.. controls (-0.5,-0.5) and (0.5,-0.5) .. (0.5,0.5)
		arc[x radius=0.5, y radius=0.15, start angle=360, end angle=180];
	\draw[2facetLine] (-0.5,0.5) .. controls (-0.5,-0.5) and (0.5,-0.5) .. (0.5,0.5)
		(0, 0.5) ellipse[x radius=0.5, y radius=0.15];
	\foamarrow#1[2facetLine]({$(0,0.35) + (0,0.6pt)$});
	\ifcase 0#2\or
		\foamdot(0.15,0);\or
		\foamdualdot(0.15,0);%
	\fi
}

\deffoampict{1cup and disk}(-0.7,-0.8)--(0.7,0.65) B=-0.2 [#1#2]=[>]{%
	\fill[1facetBack] (-0.5,0.5) -- (-0.5,0)
		.. controls (-0.5,-1) and (0.5,-1) .. (0.5,0) -- (0.5,0.5)
		arc[x radius=0.5, y radius=0.15, start angle=0, end angle=180];
	\draw[2facetInner] (0,-0.85pt) ellipse[x radius=0.5, y radius=0.15];
	\draw[2facetFront] (0, 0.85pt) ellipse[x radius=0.5, y radius=0.15];
	\fill[1facetFront] (-0.5,0.5) -- (-0.5,0)
		.. controls (-0.5,-1) and (0.5,-1) .. (0.5,0) -- (0.5,0.5)
		arc[x radius=0.5, y radius=0.15, start angle=360, end angle=180];
	\draw[1facetLine] (-0.5,0.5) -- (-0.5,0)
		.. controls (-0.5,-1) and (0.5,-1) .. (0.5,0) -- (0.5,0.5)
		(0,0.5) ellipse[x radius=0.5, y radius=0.15];
	\draw[seam]
		(-0.5, 0.85pt) arc[x radius=0.5, y radius=0.15, start angle=180, end angle=360]
		(-0.5,-0.85pt) arc[x radius=0.5, y radius=0.15, start angle=180, end angle=360];
	\foamarrow#1[seam](0,-0.15);
	\ifcase 0#2\or
		\foamdot(0.15,-0.5);\or
		\foamdualdot(0.15,-0.5);%
	\fi
}

\deffoampict{1cap}(-0.7, 0.4)--(0.7,-0.75) M=0 [#1]=[]{%
	\fill[1facetBack] (-0.5,-0.5)
		.. controls (-0.5, 0.5) and (0.5, 0.5) .. (0.5,-0.5)
		arc[x radius=0.5, y radius=0.15, start angle=0, end angle=180];
	\draw[1facetLine] (0.5,-0.5)
		arc[x radius=0.5, y radius=0.15, start angle=0, end angle=180];
	\draw[1facetFront] (-0.5,-0.5)
		.. controls (-0.5, 0.5) and (0.5, 0.5) .. (0.5,-0.5)
		arc[x radius=0.5, y radius=0.15, start angle=360, end angle=180];
	\ifcase 0#1\or
		\foamdot(0.15,0);\or
		\foamdualdot(0.15,0);%
	\fi
}

\deffoampict{2cap}(-0.7, 0.4)--(0.7,-0.75) M=0 [#1#2]=[x]{%
	\getlength\innradx({$(0.5,0) + (-1.7pt,0)$});
	\getlength\innrady({$(0,0.15) + (0,-1pt)$});
	\fill[2facetInner] (-0.5,-0.5)
		.. controls (-0.5, 0.5) and (0.5, 0.5) .. (0.5,-0.5)
		arc[x radius=0.5, y radius=0.15, start angle=0, end angle=180];
	\draw[2facetLine] (0.5,-0.5)
		arc[x radius=0.5, y radius=0.15, start angle=0, end angle=180];
	\fill[2facetBack] (-\innradx,-0.5)
		.. controls (-\innradx, \innradx) and (\innradx, \innradx)
		.. (\innradx,-0.5)
		arc[x radius=\innradx, y radius=\innrady, start angle=0, end angle=180];
	\draw[2facetLine] (-\innradx,-0.5)
		.. controls (-\innradx, \innradx) and (\innradx, \innradx)
		.. (\innradx,-0.5)
		(0,-0.5) ellipse[x radius=\innradx, y radius=\innrady];
	\fill[2facetInner] (-\innradx,-0.5)
		.. controls (-\innradx, \innradx) and (\innradx, \innradx)
		.. (\innradx,-0.5)
		arc[x radius=\innradx, y radius=\innrady,
		    start angle=360, end angle=180];
	\draw[2facetFront] (-0.5,-0.5)
		.. controls (-0.5, 0.5) and (0.5, 0.5) .. (0.5,-0.5)
		arc[x radius=0.5, y radius=0.15, start angle=360, end angle=180];
	\foamarrow#1[2facetLine]({$(0,-0.65) + (0,0.6pt)$});
	\ifcase 0#2\or
		\foamdot(0.15,0);\or
		\foamdualdot(0.15,0);%
	\fi
}

\deffoampict{1cap and disk}(-0.7, 0.8)--(0.7,-0.65) M=0.1 [#1#2]=[>]{%
	\fill[1facetBack] (-0.5,-0.5) -- (-0.5,0)
		.. controls (-0.5, 1) and (0.5, 1) .. (0.5,0) -- (0.5,-0.5)
		arc[x radius=0.5, y radius=0.15, start angle=0, end angle=180];
	\draw[2facetInner] (0,-0.85pt) ellipse[x radius=0.5, y radius=0.15];
	\draw[2facetFront] (0, 0.85pt) ellipse[x radius=0.5, y radius=0.15];
	\draw (0.5,-0.5) arc[x radius=0.5, y radius=0.15, start angle=0, end angle=180];
	\draw[1facetFront] (-0.5,-0.5) -- (-0.5,0)
		.. controls (-0.5, 1) and (0.5, 1) .. (0.5,0) -- (0.5,-0.5)
		arc[x radius=0.5, y radius=0.15, start angle=360, end angle=180];
	\draw[seam]
		(-0.5, 0.85pt) arc[x radius=0.5, y radius=0.15, start angle=180, end angle=360]
		(-0.5,-0.85pt) arc[x radius=0.5, y radius=0.15, start angle=180, end angle=360];
	\foamarrow#1[seam](0,-0.15);
	\ifcase 0#2\or
		\foamdot(0.15, 0.5);\or
		\foamdualdot(0.15, 0.5);%
	\fi
}

\deffoampict{3facets}(-1,-0.5)--(0.65,0.5) M=0 [#1#2]=[.]{%
	\draw[1facetFront]
		(-0.5, 0) arc[radius=1, start angle=180, end angle=150] --
		++( 1, 0) arc[radius=1, start angle=150, end angle=210] --
		++(-1, 0) arc[radius=1, start angle=210, end angle=180];
	\draw[2facetInner] (0.5,-0.85pt) -- ++(-0.5,-0.3) -- ++(-1,0) -- ++(0.5,0.3);
	\draw[seam] (0.5,-0.85pt) -- ++(-1,0);
	\draw[2facetFront] (0.5, 0.85pt) -- ++(-0.5,-0.3) -- ++(-1,0) -- ++(0.5,0.3);
	\draw[seam] (0.5, 0.85pt) -- ++(-1,0);
	\ifx+#1%
		\ifcase 0#2
			\foamdot(0.25,0.25);\or
			\foamdot(0.25,0.25);\or
			\foamdualdot(0.25,0.25);
		\fi
	\else\ifx-#1%
		\ifcase 0#2
			\foamdot(0.35,-0.3);\or
			\foamdot(0.35,-0.3);\or
			\foamdualdot(0.35,-0.3);
		\fi
	\fi\fi
}

\deffoampict{1neck}(-0.6,-1.2)--(0.6,1.2) M=0 {%
	\fill[1facetBack] (-0.5,-1) .. controls (-0.3,0) .. (-0.5,1)
		arc[x radius=0.5, y radius=0.15, start angle=180, end angle=0]
		.. controls (0.3,0) .. (0.5,-1)
		arc[x radius=0.5, y radius=0.15, start angle=0, end angle=180];
	\draw[1facetLine] (-0.5,-1)
		arc[x radius=0.5, y radius=0.15, start angle=180, end angle=0];
	\fill[1facetFront](-0.5,-1) .. controls (-0.3,0) .. (-0.5,1)
		arc[x radius=0.5, y radius=0.15, start angle=180, end angle=360]
		.. controls (0.3,0) .. (0.5,-1)
		arc[x radius=0.5, y radius=0.15, start angle=360, end angle=180];
	\draw[1facetLine] (-0.5,1) .. controls (-0.3,0) .. (-0.5,-1)
		arc[x radius=0.5, y radius=0.15, start angle=180, end angle=360]
		.. controls (0.3,0) .. (0.5,1);
	\draw[1facetLine] (0, 1) ellipse[x radius=0.5, y radius=0.15];
}
\deffoampict{2neck}(-0.6,-1.2)--(0.6,1.2) M=0 {%
	\getlength\innradx({$(0.5,0) + (-1.7pt,0)$});
	\getlength\innrady({$(0,0.15) + (0,-1pt)$});
	\fill[2facetInner](-0.5,-1) .. controls ++(0.2,1) .. (-0.5,1)
		arc[x radius=0.5, y radius=0.15, start angle=180, end angle=0]
		.. controls ++(-0.2,-1) .. (0.5,-1)
		arc[x radius=0.5, y radius=0.15, start angle=0, end angle=180];
	\draw[2facetLine] (0.5,-1)
		arc[x radius=0.5, y radius=0.15, start angle=0, end angle=180];
	\fill[2facetBack](-\innradx,-1) .. controls ++(0.2,1) .. (-\innradx,1)
		arc[x radius=\innradx, y radius=\innrady, start angle=180,end angle=0]
		.. controls ++(-0.2,-1) .. (\innradx,-1)
		arc[x radius=\innradx, y radius=\innrady, start angle=0,end angle=180];
	\draw[2facetLine] (\innradx,-1)
		arc[x radius=\innradx, y radius=\innrady, start angle=0,end angle=180];
	\fill[2facetInner](-\innradx,-1) .. controls ++(0.2,1) .. (-\innradx,1)
		arc[x radius=\innradx, y radius=\innrady, start angle=180,end angle=360]
		.. controls ++(-0.2,-1) .. (\innradx,-1)
		arc[x radius=\innradx, y radius=\innrady, start angle=360,end angle=180];
	\draw[2facetLine] (-\innradx, 1) .. controls ++(0.2,-1) .. (-\innradx,-1)
		arc[x radius=\innradx, y radius=\innrady, start angle=180,end angle=360]
		.. controls ++(-0.2,1) .. (\innradx, 1);
	\fill[2facetFront](-0.5,-1) .. controls (-0.3,0) .. (-0.5,1)
		arc[x radius=0.5, y radius=0.15, start angle=180, end angle=360]
		.. controls (0.3,0) .. (0.5,-1)
		arc[x radius=0.5, y radius=0.15, start angle=360, end angle=180];
	\draw[2facetLine] (-0.5, 1) .. controls (-0.3,0) .. (-0.5,-1)
		arc[x radius=0.5, y radius=0.15, start angle=180, end angle=360]
		.. controls (0.3,0) .. (0.5,1);
	\draw[2facetLine]
		(0, 1) ellipse[x radius=0.5, y radius=0.15]
		(0, 1) ellipse[x radius=\innradx, y radius=\innrady];
}

\deffoampict{1cap 1cup}(-0.6,-1.2)--(0.6,1.2) M=0 [#1#2]=[.]{%
	\fill[1facetBack](0.5,1) .. controls (0.5,0) and (-0.5,0) .. (-0.5,1)
		arc[x radius=0.5, y radius=0.15, start angle=180, end angle=0];
	\fill[1facetFront](0.5,1) .. controls (0.5,0) and (-0.5,0) .. (-0.5,1)
		arc[x radius=0.5, y radius=0.15, start angle=180, end angle=360];
	\draw[1facetLine]
		(0.5,1) .. controls (0.5,0) and (-0.5,0) .. (-0.5,1)
		(0,1) ellipse[x radius=0.5, y radius=0.15];
	\fill[1facetBack] (0.5,-1) .. controls (0.5,0) and (-0.5,0) .. (-0.5,-1)
		arc[x radius=0.5, y radius=0.15, start angle=180, end angle=0];
	\draw[1facetLine] (-0.5,-1)
		arc[x radius=0.5, y radius=0.15, start angle=180, end angle=0];
	\draw[1facetFront](0.5,-1) .. controls (0.5,0) and (-0.5,0) .. (-0.5,-1)
		arc[x radius=0.5, y radius=0.15, start angle=180, end angle=360];
	\ifx#1.\else
		\ifcase 0#2
			\foamdot(0.2, #10.6);\or
			\foamdot(0.2, #10.6);\or
			\foamdualdot(0.2, #10.6);
		\fi
	\fi
}
\deffoampict{2cap 2cup}(-0.6,-1.2)--(0.6,1.2) M=0 [#1#2]=[.]{%
	\getlength\innradx({$(0.5,0) + (-1.7pt,0)$});
	\getlength\innrady({$(0,0.15) + (0,-1pt)$});
	\fill[2facetInner](0.5,1) .. controls (0.5,0) and (-0.5,0) .. (-0.5,1)
		arc[x radius=0.5, y radius=0.15, start angle=180, end angle=0];
	\fill[2facetBack] (\innradx,1)
		.. controls (\innradx,1.7pt) and (-\innradx,1.7pt) .. (-\innradx,1)
		arc[x radius=\innradx,y radius=\innrady,start angle=180,end angle=-180];
	\fill[2facetInner] (\innradx,1)
		.. controls (\innradx,1.7pt) and (-\innradx,1.7pt) .. (-\innradx,1)
		arc[x radius=\innradx,y radius=\innrady,start angle=180,end angle=360];
	\draw[2facetLine] (\innradx,1)
		.. controls (\innradx,1.7pt) and (-\innradx,1.7pt) .. (-\innradx,1)
		(0,1) ellipse[x radius=\innradx,y radius=\innrady];
	\fill[2facetFront] (0.5,1) .. controls (0.5,0) and (-0.5,0) .. (-0.5,1)
		arc[x radius=0.5, y radius=0.15, start angle=180, end angle=360];
	\draw[2facetLine]
		(0.5,1) .. controls (0.5,0) and (-0.5,0) .. (-0.5,1)
		(0,1) ellipse[x radius=0.5, y radius=0.15];
	\fill[2facetInner] (0.5,-1) .. controls (0.5,0) and (-0.5,0) .. (-0.5,-1)
		arc[x radius=0.5, y radius=0.15, start angle=180, end angle=0];
	\draw[2facetLine] (-0.5,-1)
		arc[x radius=0.5, y radius=0.15, start angle=180, end angle=0];
	\fill[2facetBack] (\innradx,-1)
		.. controls (\innradx,-1.7pt) and (-\innradx,-1.7pt) .. (-\innradx,-1)
		arc[x radius=\innradx,y radius=\innrady,start angle=180,end angle=0];
	\draw[2facetLine] (-\innradx,-1)
		arc[x radius=\innradx,y radius=\innrady,start angle=180,end angle=0];
	\draw[2facetInner] (\innradx,-1)
		.. controls (\innradx,-1.7pt) and (-\innradx,-1.7pt) .. (-\innradx,-1)
		arc[x radius=\innradx,y radius=\innrady,start angle=180,end angle=360];
	\draw[2facetFront] (0.5,-1) .. controls (0.5,0) and (-0.5,0) .. (-0.5,-1)
		arc[x radius=0.5, y radius=0.15, start angle=180, end angle=360];
	\ifx#1.\else
		\ifcase 0#2
			\foamdot(0.2, #10.6);\or
			\foamdot(0.2, #10.6);\or
			\foamdualdot(0.2, #10.6);
		\fi
	\fi
}

\deffoampict{2cyl at 1plane}(-1.8,-0.4)--(1.8,1.25) M=0.4 [#1#2]=[.]{%
	\getlength\innradx({$(0.5,0) + (-1.7pt,0)$});
	\getlength\innrady({$(0,0.15) + (0,-1pt)$});
	\draw[1facetFront] (-1.8,-0.4) -- (-0.6,0.4) -- (1.7,0.4) -- (0.5,-0.4) -- cycle;
	\fill[2facetInner] (0.5,0)
		-- ( 0.5,1) arc[x radius=0.5, y radius=0.15, start angle=0, end angle=180]
		-- (-0.5,0) arc[x radius=0.5, y radius=0.15, start angle=180, end angle=0];
	\draw[seam] (-0.5,0)
		arc[x radius=0.5, y radius=0.15, start angle=180, end angle=0];
	\fill[2facetBack] (\innradx,0) -- (\innradx,1)
		arc[x radius=\innradx, y radius=\innrady, start angle=0, end angle=180]
		-- (-\innradx, 0)
		arc[x radius=\innradx, y radius=\innrady, start angle=180, end angle=0];
	\draw[seam] (-\innradx,0)
		arc[x radius=\innradx, y radius=\innrady, start angle=180, end angle=0];
	\fill[2facetInner] (\innradx,0) -- (\innradx,1)
		arc[x radius=\innradx,y radius=\innrady,start angle=360,end angle=180]
		-- (-\innradx,0)
		arc[x radius=\innradx,y radius=\innrady,start angle=180,end angle=360];
	\draw[2facetLine]
		( \innradx, 1) -- ( \innradx, 0)
		(-\innradx, 1) -- (-\innradx, 0);
	\draw[seam] (\innradx, 0)
		arc[x radius=\innradx,y radius=\innrady,start angle=360,end angle=180];
	\ifx V#2\relax\ifx #1>%
		\draw[2facetLine,line width=1pt,->]
			(0,0.75) ++(220:0.5 and 0.15) -- (0,0.75);
	\fi\fi
	\fill[2facetFront] (0.5,0)
		-- ( 0.5,1) arc[x radius=0.5, y radius=0.15, start angle=360, end angle=180]
		-- (-0.5,0) arc[x radius=0.5, y radius=0.15, start angle=180, end angle=360];
	\draw[2facetLine] (-0.5,1) -- (-0.5,0) (0.5,1) -- (0.5,0)
		(0,1) ellipse[x radius=\innradx, y radius=\innrady]
		(0,1) ellipse[x radius=0.5, y radius=0.15];
	\draw[seam]
		(-0.5,0) arc[x radius=0.5, y radius=0.15, start angle=180, end angle=360];
	\ifx V#2\relax
		\ifx #1<
			\draw[2facetLine,line width=1pt,->]
				(0,0.75) ++(220:0.5 and 0.15) -- ++(220:0.5 and 0.15);
		\fi
	\else\ifx\relax#2\relax\else
		\foamarrow#1[2facetLine]0({$(0,0.85) + (-4pt,1pt)$});
		\foamarrow#1[2facetLine]180({$(0, 1.15) + (4pt,-1pt)$});
	\fi\fi
	\foamarrow#1[seam]7({$(0,-0.15) + (4pt, 1pt)$});
}

\deffoampict{2cup > 1plane}(-1.8,-0.4)--(1.8,1.25) M=0.4 [#1]=[.]{%
	\getlength\innradx({$(0.5,0) + (-1.7pt,0)$});
	\getlength\innrady({$(0,0.15) + (0,-1pt)$});
	\draw[1facetFront] (-1.8,-0.4) -- (-0.6,0.4) -- (1.7,0.4) -- (0.5,-0.4) -- cycle;
	\fill[2facetInner](0.5,1) .. controls (0.5,0) and (-0.5,0) .. (-0.5,1)
		arc[x radius=0.5, y radius=0.15, start angle=180, end angle=0];
	\fill[2facetBack] (\innradx,1)
		.. controls (\innradx,1.7pt) and (-\innradx,1.7pt) .. (-\innradx,1)
		arc[x radius=\innradx,y radius=\innrady,start angle=180,end angle=-180];
	\fill[2facetInner] (\innradx,1)
		.. controls (\innradx,1.7pt) and (-\innradx,1.7pt) .. (-\innradx,1)
		arc[x radius=\innradx,y radius=\innrady,start angle=180,end angle=360];
	\draw[2facetLine] (\innradx,1)
		.. controls (\innradx,1.7pt) and (-\innradx,1.7pt) .. (-\innradx,1)
		(0,1) ellipse[x radius=\innradx,y radius=\innrady];
	\ifx A#1\relax
		\draw[2facetLine,line width=1pt,->] (0,0.25) -- (0,0.6);
	\fi
	\fill[2facetFront] (0.5,1) .. controls (0.5,0) and (-0.5,0) .. (-0.5,1)
		arc[x radius=0.5, y radius=0.15, start angle=180, end angle=360];
	\draw[2facetLine]
		(0.5,1) .. controls (0.5,0) and (-0.5,0) .. (-0.5,1)
		(0,1) ellipse[x radius=0.5, y radius=0.15];
	\ifx V#1\relax
		\draw[2facetLine,line width=1pt,->] (0,0.25) -- (0,-0.1);
	\else
		\foamarrow#1[2facetLine]({$(0,0.85) + (-4pt,1pt)$});
		\foamarrow#1[2facetLine]180({$(0, 1.15) + (4pt,-1pt)$});
	\fi
}

\deffoampict{2planes at 1plane}(-1.8,-0.4)--(1.8,1.1) M=0.35 [#1]=[.]{%
	\draw[1facetFront] (-1.7,-0.4) -- (-0.5,0.4) -- (1.7,0.4) -- (0.5,-0.4) -- cycle;
	\draw[2facetInner] (-0.875, 0.15) rectangle ++(2.2,0.85);
	\draw[seam] (-0.875,0.15) -- ++(2.2,0);
	\draw[2facetBack]  (-0.875, 0.15) ++(-1.4pt,-0.96pt) rectangle ++(2.2,0.85);
	\draw[seam] (-0.875,0.15) ++(-1.4pt,-0.96pt) -- ++(2.2,0);
	\draw[2facetInner] (-1.325,-0.15) ++( 1.4pt, 0.96pt) rectangle ++(2.2,0.85);
	\draw[seam] (-1.325,-0.15) ++(1.4pt,0.96pt) -- ++(2.2,0);
	\foamarrow#1[seam]180({ 0.225, 0.15) ++(-0.7pt,-0.47pt});
	\draw[2facetFront] (-1.325,-0.15) rectangle ++(2.2,0.85);
	\draw[seam] (-1.325,-0.15) -- ++(2.2,0);
	\foamarrow#1[seam]({-0.225,-0.15) ++(0.7pt, 0.47pt});
}

\deffoampict{2saddle at 1plane}(-1.8,-0.4)--(1.8,1.1) M=0.35 [#1]=[.]{%
	\getlength\xradbig({$(0.95,0) + (-3.1pt,0)$});
	\getlength\xradsmall({$(0.5,0) + (-0.3pt,0)$});
	\getlength\yrad({$(0,0.15) + (0,-0.96pt)$});
	\getlength\radmid({$(0,0.375) + (0, 1.7pt)$});
	\draw[1facetFront] (-1.7,-0.4) -- (-0.5,0.4) -- (1.7,0.4) -- (0.5,-0.4) -- cycle;
	\fill[2facetInner] (-0.875, 0.15)
		arc[x radius=0.5, y radius=0.15, start angle=90, end angle=0]
		arc[radius=0.375, start angle=180, end angle=0]
		arc[x radius=0.95, y radius=0.15, start angle=180, end angle=90]
		-- ++(0,0.85) -- ++(-2.2,0) -- cycle;
	\draw[2facetLine] (-0.875,0.15) |- ++(2.2,0.85) -- ++(0,-0.85);
	\draw[seam]
		(-0.875, 0.15) arc[x radius=0.5, y radius=0.15, start angle=90, end angle=0]
		++(0.75, 0) arc[x radius=0.95, y radius=0.15, start angle=180, end angle=90];
	\fill[2facetBack] (-0.875, 0.15) ++(-1.4pt,-0.96pt)
		arc[x radius=\xradsmall, y radius=\yrad, start angle=90, end angle=0]
		arc[radius=\radmid, start angle=180, end angle=0]
		arc[x radius=\xradbig, y radius=\yrad, start angle=180, end angle=90]
		-- ++(0,0.85) -- ++(-2.2,0) -- cycle;
	\draw[2facetLine] (-0.875, 0.15) ++(-1.4pt,-0.96pt) |- ++(2.2,0.85) -- ++(0,-0.85);
	\draw[seam] (-0.875, 0.15) ++(-1.4pt,-0.96pt)
		arc[x radius=\xradsmall, y radius=\yrad, start angle=90, end angle=0]
		++(\radmid, 0) ++(\radmid, 0)
		arc[x radius=\xradbig, y radius=\yrad, start angle=180, end angle=90];
	\draw[2facetInner] (-1.325,-0.15) ++(1.4pt, 0.96pt)
		arc[x radius=\xradbig, y radius=\yrad, start angle=270, end angle=360]
		arc[radius=\radmid, start angle=180, end angle=0]
		arc[x radius=\xradsmall, y radius=\yrad, start angle=180, end angle=270]
		-- ++(0,0.85) -- ++(-2.2,0) -- cycle;
	\draw[seam] (-1.325,-0.15) ++(1.4pt, 0.96pt)
		arc[x radius=\xradbig, y radius=\yrad, start angle=270, end angle=360]
		++(\radmid, 0) ++(\radmid, 0)
		arc[x radius=\xradsmall, y radius=\yrad, start angle=180, end angle=270];
	\draw[2facetFront] (-1.325,-0.15)
		arc[x radius=0.95, y radius=0.15, start angle=270, end angle=360]
		arc[radius=0.375, start angle=180, end angle=0]
		arc[x radius=0.5, y radius=0.15, start angle=180, end angle=270]
		-- ++(0,0.85) -- ++(-2.2,0) -- cycle;
	\draw[seam]
		(-1.325,-0.15) arc[x radius=0.95, y radius=0.15, start angle=270, end angle=360]
		++(0.75,0) arc[x radius=0.5, y radius=0.15, start angle=180, end angle=270];
	\foamarrow#1[seam]8({-0.823,-0.128) ++(0.7pt, 0.47pt});
	\foamarrow#1[seam]-10({0.588,-0.123) ++(0.7pt, 0.47pt});
}

\deffoampict{2bubble at 1plane}(-1.4,-0.3)--(1.4,0.8) M=0 [#1]=[.]{%
	\getlength\xrad(0.5,0);
	\getlength\yrad(0,0.15);
	\getlength\yyrad(0,0.7);
	\draw[1facetFront] (-0.6,0.3) -- (-1.4,-0.3) -- (0.6,-0.3) -- (1.4,0.3) -- cycle;
	\fill[2facetInner] (-0.5,0) ++(-0.85pt,0)
		arc[start angle=180, end angle=0, x radius=\xrad+0.85pt, y radius=\yrad+0.7pt]
		arc[start angle=0, end angle=180, x radius=\xrad+0.85pt, y radius=\yyrad+0.85pt];
	\draw[seam] (-0.5,0) ++(-0.85pt,0)
		arc[start angle=180, end angle=0, x radius=\xrad+0.85pt, y radius=\yrad+0.7pt];
	\fill[2facetBack] (-0.5,0) ++(0.85pt,0)
		arc[start angle=180, end angle=0, x radius=\xrad-0.85pt, y radius=\yrad-0.7pt]
		arc[start angle=0, end angle=180, x radius=\xrad-0.85pt, y radius=\yyrad-0.85pt];
	\draw[seam] (-0.5,0) ++(0.85pt,0)
		arc[start angle=180, end angle=0, x radius=\xrad-0.85pt, y radius=\yrad-0.7pt];
	\fill[2facetInner] (-0.5,0) ++(0.85pt,0)
		arc[start angle=180, end angle=360, x radius=\xrad-0.85pt, y radius=\yrad-0.7pt]
		arc[start angle=0, end angle=180, x radius=\xrad-0.85pt, y radius=\yyrad-0.85pt];
	\draw[seam] (-0.5,0) ++(0.85pt,0)
		arc[start angle=180, end angle=360, x radius=\xrad-0.85pt, y radius=\yrad-0.7pt];
	\draw[2facetLine] (-0.5,0) ++(0.85pt,0)
		arc[start angle=180, end angle=0, x radius=\xrad-0.85pt, y radius=\yyrad-0.85pt];
	\fill[2facetFront] (-0.5,0) ++(-0.85pt,0)
		arc[start angle=180, end angle=360, x radius=\xrad+0.85pt, y radius=\yrad+0.7pt]
		arc[start angle=0, end angle=180, x radius=\xrad+0.85pt, y radius=\yyrad+0.85pt];
	\draw[seam] (-0.5,0) ++(-0.85pt,0)
		arc[start angle=180, end angle=360, x radius=\xrad+0.85pt, y radius=\yrad+0.7pt];
	\draw[2facetLine] (-0.5,0) ++(-0.85pt,0)
		arc[start angle=180, end angle=0, x radius=\xrad+0.85pt, y radius=\yyrad+0.85pt];
	\foamarrow#1[seam](0,-0.15);
}

\deffoampict{no bubble at 1plane}(-1.4,-0.3)--(1.4,0.3) M=0 {%
	\draw[1facetFront] (-0.6,0.3) -- (-1.4,-0.3) -- (0.6,-0.3) -- (1.4,0.3) -- cycle;
}

\deffoampict{2tunel at 1plane}(-1.4,-0.3)--(1.4,0.8) M=0.2 [#1]=[.]{%
	\draw[1facetFront] (-0.6,0.3) -- (-1.4,-0.3) -- (0.6,-0.3) -- (1.4,0.3) -- cycle;
	\getlength\radx(0.2,0);
	\getlength\rady(0,0.6);
	\fill[2facetInner] (1.0,0.75) ++(0,0.85pt)
		arc[x radius=\radx+0.68pt, y radius=\rady+0.34pt,
		    start angle=90, end angle=0] -- ++(-2,0)
		arc[x radius=\radx+0.68pt, y radius=\rady+0.34pt,
		    start angle=0, end angle=90];
	\draw[seam] (1.2,0.15) ++(0.68pt,0.51pt) -- ++(-2,0);
	\draw[2facetLine] (-0.8,0.15) ++(0.68pt,0.51pt)
		arc[x radius=\radx+0.68pt, y radius=\rady+0.34pt,
		    start angle=0, end angle=90]
		++(2,0) arc[x radius=\radx+0.68pt, y radius=\rady+0.34pt,
		    start angle=90, end angle=0];
	\fill[2facetBack] (1.0,0.75) ++(0,-0.85pt)
		arc[x radius=\radx-0.68pt, y radius=\rady-0.34pt,
		    start angle=90, end angle=0] -- ++(-2,0)
		arc[x radius=\radx-0.68pt, y radius=\rady-0.34pt,
		    start angle=0, end angle=90];
	\draw[seam] (1.2,0.15) ++(-0.68pt,-0.51pt) -- ++(-2,0);
	\draw[2facetLine] (-0.8,0.15) ++(-0.68pt,-0.51pt)
		arc[x radius=\radx-0.68pt, y radius=\rady-0.34pt,
		    start angle=0, end angle=90]
		++(2,0) arc[x radius=\radx-0.68pt, y radius=\rady-0.34pt,
		    start angle=90, end angle=0];
	\getlength\rady(0,0.9);
	\fill[2facetInner] (1.0,0.75) ++(0,-0.85pt)
		arc[x radius=\radx-0.68pt, y radius=\rady-1.36pt,
		    start angle=90, end angle=180] -- ++(-2,0)
		arc[x radius=\radx-0.68pt, y radius=\rady-1.36pt,
		    start angle=180, end angle=90] -- cycle;
	\draw[seam] (0.8,-0.15) ++(0.68pt,0.51pt) -- ++(-2,0);
	\draw[2facetLine] (1.0,0.75) ++(0,-0.85pt)
		arc[x radius=\radx-0.68pt, y radius=\rady-1.36pt,
		    start angle=90, end angle=180]
		++(-2,0) arc[x radius=\radx-0.68pt, y radius=\rady-1.36pt,
		    start angle=180, end angle=90] -- ++(2,0);
	\foamarrow#1[seam]{180}(0.4,0.15);
	\fill[2facetFront] (1.0,0.75) ++(0, 0.85pt)
		arc[x radius=\radx+0.68pt, y radius=\rady+1.36pt,
		    start angle=90, end angle=180] -- ++(-2,0)
		arc[x radius=\radx+0.68pt, y radius=\rady+1.36pt,
		    start angle=180, end angle=90] -- cycle;
	\draw[seam] (0.8,-0.15) ++(-0.68pt,-0.51pt) -- ++(-2,0);
	\draw[2facetLine] (1.0,0.75) ++(0, 0.85pt)
		arc[x radius=\radx+0.68pt, y radius=\rady+1.36pt,
		    start angle=90, end angle=180] ++(-2,0)
		arc[x radius=\radx+0.68pt, y radius=\rady+1.36pt,
		    start angle=180, end angle=90] -- ++(2,0);
	\foamarrow#1[seam](-0.4,-0.15);
}

\deffoampict{2shells at 1plane}(-1.4,-0.3)--(1.4,0.8) M=0.2 [#1]=[.]{%
	\getlength\xunit(0.1,0);
	\getlength\yunit(0,0.1);
	\draw[1facetFront] (-0.6,0.3) -- (-1.4,-0.3) -- (0.6,-0.3) -- (1.4,0.3) -- cycle;
	\fill[2facetInner] (-1, 0.75) ++(0, 0.85pt)
		arc[x radius=\xunit*8.0+0.85pt, start angle=90, 
		    y radius=\yunit*7.5+0.85pt,   end angle=0]
		arc[x radius=\xunit*6.0+0.17pt, start angle=0, 
		    y radius=\yunit*1.5+0.51pt,   end angle=90]
		arc[x radius=\xunit*2.0+0.68pt, start angle=0,
		    y radius=\yunit*6.0+0.34pt,   end angle=90];
	\draw[2facetLine] (-1, 0.75) ++(0, 0.85pt)
		arc[x radius=\xunit*2.0+0.68pt, start angle=90,
		    y radius=\yunit*6.0+0.34pt,   end angle=0];
	\draw[seam] (-0.2,0) ++(0.85pt,0)
		arc[x radius=\xunit*6.0+0.17pt, start angle=0,
		    y radius=\yunit*1.5+0.51pt,   end angle=90];
	\fill[2facetBack] (-1, 0.75) ++(0,-0.85pt)
		arc[x radius=\xunit*8.0-0.85pt, start angle=90, 
		    y radius=\yunit*7.5-0.85pt,   end angle=0]
		arc[x radius=\xunit*6.0-0.17pt, start angle=0, 
		    y radius=\yunit*1.5-0.51pt,   end angle=90]
		arc[x radius=\xunit*2.0-0.68pt, start angle=0,
		    y radius=\yunit*6.0-0.34pt,   end angle=90];
	\draw[2facetLine] (-1, 0.75) ++(0,-0.85pt)
		arc[x radius=\xunit*2.0-0.68pt, start angle=90,
		    y radius=\yunit*6.0-0.34pt,   end angle=0];
	\draw[seam] (-0.2,0) ++(-0.85pt,0)
		arc[x radius=\xunit*6.0-0.17pt, start angle=0, 
		    y radius=\yunit*1.5-0.51pt,   end angle=90];
	\fill[2facetInner] (-1, 0.75) ++(0,-0.85pt)
		arc[x radius=\xunit*8.0-0.85pt, start angle=90, 
		    y radius=\yunit*7.5-0.85pt,   end angle=0]
		arc[x radius=\xunit*10 -1.53pt, start angle=360,
		    y radius=\yunit*1.5-0.51pt,   end angle=270]
		arc[x radius=\xunit*2.0-0.68pt, start angle=180,
		    y radius=\yunit*9.0-1.36pt,   end angle=90];
	\draw[2facetLine] (-0.2,0) ++(-0.85pt,0)
		arc[x radius=\xunit*8.0-0.85pt, start angle=0, 
		    y radius=\yunit*7.5-0.85pt,   end angle=90]
		arc[x radius=\xunit*2.0-0.68pt, start angle=90,
		    y radius=\yunit*9.0-1.36pt,   end angle=180];
	\draw[seam] (-0.2,0) ++(-0.85pt,0)
		arc[x radius=\xunit*10 -1.53pt, start angle=360,
		    y radius=\yunit*1.5-0.51pt,   end angle=270];
	\fill[2facetFront] (-1, 0.75) ++(0, 0.85pt)
		arc[x radius=\xunit*8.0+0.85pt, start angle=90, 
		    y radius=\yunit*7.5+0.85pt,   end angle=0]
		arc[x radius=\xunit*10 +1.53pt, start angle=360,
		    y radius=\yunit*1.5+0.51pt,   end angle=270]
		arc[x radius=\xunit*2.0+0.68pt, start angle=180,
		    y radius=\yunit*9.0+1.36pt,   end angle=90];
	\draw[2facetLine] (-0.2,0) ++(0.85pt,0)
		arc[x radius=\xunit*8.0+0.85pt, start angle=0, 
		    y radius=\yunit*7.5+0.85pt,   end angle=90]
		arc[x radius=\xunit*2.0+0.68pt, start angle=90,
		    y radius=\yunit*9.0+1.36pt,   end angle=180];
	\draw[seam] (-0.2,0) ++(0.85pt,0)
		arc[x radius=\xunit*10 +1.53pt, start angle=360,
		    y radius=\yunit*1.5+0.51pt,   end angle=270];
	\fill[2facetInner] (1,0.75) ++(0,0.85pt)
		arc[x radius=\xunit*2.0+0.68pt, start angle=90,  
		    y radius=\yunit*6.0+0.34pt,   end angle=0]
		arc[x radius=\xunit*10 +1.53pt, start angle=90,  
		    y radius=\yunit*1.5+0.51pt,   end angle=180]
		arc[x radius=\xunit*8.0+0.85pt, start angle=180,
		    y radius=\yunit*7.5+0.85pt,   end angle=90];
	\draw[2facetLine] (1,0.75) ++(0,0.85pt)
		arc[x radius=\xunit*2.0+0.68pt, start angle=90,  
		    y radius=\yunit*6.0+0.34pt,   end angle=0];
	\draw[seam] (0.2,0) ++(-0.85pt,0)
		arc[x radius=\xunit*10 +1.53pt, start angle=180,  
		    y radius=\yunit*1.5+0.51pt,   end angle=90];
	\fill[2facetBack] (1,0.75) ++(0,-0.85pt)
		arc[x radius=\xunit*2.0-0.68pt, start angle=90,  
		    y radius=\yunit*6.0-0.34pt,   end angle=0]
		arc[x radius=\xunit*10 -1.53pt, start angle=90,  
		    y radius=\yunit*1.5-0.51pt,   end angle=180]
		arc[x radius=\xunit*8.0-0.85pt, start angle=180,
		    y radius=\yunit*7.5-0.85pt,   end angle=90];
	\draw[2facetLine] (1,0.75) ++(0,-0.85pt)
		arc[x radius=\xunit*2.0-0.68pt, start angle=90,  
		    y radius=\yunit*6.0-0.34pt,   end angle=0];
	\draw[seam] (0.2,0) ++(0.85pt,0)
		arc[x radius=\xunit*10 -1.53pt, start angle=180,  
		    y radius=\yunit*1.5-0.51pt,   end angle=90];
	\fill[2facetInner] (1,0.75) ++(0,-0.85pt)
		arc[x radius=\xunit*8.0-0.85pt, start angle=90,  
		    y radius=\yunit*7.5-0.85pt,   end angle=180]
		arc[x radius=\xunit*6.0-0.17pt, start angle=180,  
		    y radius=\yunit*1.5-0.51pt,   end angle=270]
		arc[x radius=\xunit*2.0-0.68pt, start angle=180,  
		    y radius=\yunit*9.0-1.36pt,   end angle=90];
	\draw[2facetLine] (1,0.75) ++(0,-0.85pt)
		arc[x radius=\xunit*8.0-0.85pt, start angle=90,  
		    y radius=\yunit*7.5-0.85pt,   end angle=180]
		(1,0.75) ++(0,-0.85pt)
		arc[x radius=\xunit*2.0-0.68pt, start angle=90,  
		    y radius=\yunit*9.0-1.36pt,   end angle=180];
	\draw[seam] (0.2,0) ++(0.85pt,0)
		arc[x radius=\xunit*6.0-0.17pt, start angle=180,  
		    y radius=\yunit*1.5-0.51pt,   end angle=270];
	\fill[2facetFront] (1,0.75) ++(0, 0.85pt)
		arc[x radius=\xunit*8.0+0.85pt, start angle=90,  
		    y radius=\yunit*7.5+0.85pt,   end angle=180]
		arc[x radius=\xunit*6.0+0.17pt, start angle=180,  
		    y radius=\yunit*1.5+0.51pt,   end angle=270]
		arc[x radius=\xunit*2.0+0.68pt, start angle=180,  
		    y radius=\yunit*9.0+1.36pt,   end angle=90];
	\draw[2facetLine] (1,0.75) ++(0, 0.85pt)
		arc[x radius=\xunit*8.0+0.85pt, start angle=90,  
		    y radius=\yunit*7.5+0.85pt,   end angle=180]
		(1,0.75) ++(0, 0.85pt)
		arc[x radius=\xunit*2.0+0.68pt, start angle=90,  
		    y radius=\yunit*9.0+1.36pt,   end angle=180];
	\draw[seam] (0.2,0) ++(-0.85pt,0)
		arc[x radius=\xunit*6.0+0.17pt, start angle=180,  
		    y radius=\yunit*1.5+0.51pt,   end angle=270];
	\foamarrow#1[seam]-10(0.414,-0.115);
	\foamarrow#1[seam]8(-0.7,-0.13);
}

\deffoampict{1cup > 2plane}(-1.4,-0.8)--(1.4,1.15) M=0.2 [#1]=[.]{%
	\draw[2facetInner] (-0.6,0.3)
		++(0,-0.85pt) -- ++(-0.8, -0.6) -- ++(2, 0) -- ++(0.8, 0.6) -- cycle;
	\draw[2facetFront] (-0.6,0.3)
		++(0, 0.85pt) -- ++(-0.8, -0.6) -- ++(2, 0) -- ++(0.8, 0.6) -- cycle;
	\foamarrow#1[2facetLine](-0.6,-0.3);
	\foamarrow#1[2facetLine]180(0.6,0.3);
	\fill[1facetBack] (-0.5,1)
			arc[x radius=0.5, y radius=0.15, start angle=180, end angle=0]
			.. controls (0.5,0) and (-0.5,0) .. (-0.5,1);
	\fill[1facetFront] (-0.5,1)
			arc[x radius=0.5, y radius=0.15, start angle=180, end angle=360]
			.. controls (0.5,0) and (-0.5,0) .. (-0.5,1);
	\draw[1facetLine]
		(0,1) ellipse[x radius=0.5, y radius=0.15]
		(0.5,1) .. controls (0.5,0) and (-0.5,0) .. (-0.5,1);
}

\deffoampict{1cup < 2plane}(-1.4,-0.8)--(1.4,1.15) M=0.2 [#1#2]=[>]{%
	\fill[2facetInner] (-0.5,-0.85pt)
			arc[x radius=0.5, y radius=0.15, start angle=180, end angle=0] --
			++(0,0.3) -- ++(-1,0) -- cycle;
	\draw[2facetLine] (-0.5,-0.85pt) ++(0,0.3) -- ++(1,0);
	\fill[2facetFront] (-0.5, 0.85pt)
			arc[x radius=0.5, y radius=0.15, start angle=180, end angle=0] --
			++(0,0.3) -- ++(-1,0) -- cycle;
	\draw[2facetLine] (-0.5,0.85pt) ++(0,0.3) -- ++(1,0);
	\fill[1facetBack]  (-0.5, 1) -- (-0.5,0)
			.. controls (-0.5,-1) and (0.5,-1) .. (0.5,0) -- (0.5,1)
			arc[x radius=0.5, y radius=0.15, start angle=0, end angle=180];
	\draw[seam]
		(-0.5, 0.85pt) arc[x radius=0.5, y radius=0.15, start angle=180, end angle=0]
		(-0.5,-0.85pt) arc[x radius=0.5, y radius=0.15, start angle=180, end angle=0];
	\fill[1facetFront] (-0.5, 1) -- (-0.5,0)
			.. controls (-0.5,-1) and (0.5,-1) .. (0.5,0) -- (0.5,1)
			arc[x radius=0.5, y radius=0.15, start angle=360, end angle=180];
	\draw[1facetLine]
			(0,1) ellipse[x radius=0.5, y radius=0.15]
			(-0.5,0.85pt) .. controls (-0.5,-1) and (0.5,-1) .. (0.5,0.85pt);
	\fill[2facetInner] (0.5,-0.85pt)
			arc[x radius=0.5, y radius=0.15, start angle=360, end angle=180] |-
			++(-0.1,0.3) -- ++(-0.8,-0.6) -- ++(2,0) -- ++(0.8,0.6) -| cycle;
	\draw[seam] (0.5,-0.85pt)
			arc[x radius=0.5, y radius=0.15, start angle=360, end angle=180];
	\draw[2facetLine] (-0.5,0.3) ++(0,-0.85pt) --
			++(-0.1,0) -- ++(-0.8,-0.6) -- ++(2,0) -- ++(0.8,0.6) -- ++(-0.9,0);
	\fill[2facetFront] (0.5, 0.85pt)
			arc[x radius=0.5, y radius=0.15, start angle=360, end angle=180] |-
			++(-0.1,0.3) -- ++(-0.8,-0.6) -- ++(2,0) -- ++(0.8,0.6) -| cycle;
	\draw[seam] (0.5, 0.85pt)
			arc[x radius=0.5, y radius=0.15, start angle=360, end angle=180];
	\draw[2facetLine] (-0.5,0.3) ++(0, 0.85pt) --
			++(-0.1,0) -- ++(-0.8,-0.6) -- ++(2,0) -- ++(0.8,0.6) -- ++(-0.9,0);
	\draw[1facetLine] (-0.5,0.85pt) -- (-0.5,1) (0.5,0.85pt) -- (0.5,1);
	\ifx\relax#2\relax\else
		\foamarrow#1[2facetLine](-0.7,-0.3);
		\foamarrow#1[2facetLine]180(0.7,0.3);
	\fi
	\foamarrow#1[seam](0,-0.15);
}

\deffoampict{1cap > 2plane}(-1.4,-1.15)--(1.4,0.8) M=-0.2 [#1#2]=[>]{%
	\fill[2facetInner] (0.5,-0.85pt)
		arc[x radius=0.5, y radius=0.15, start angle=0, end angle=180] --
		++(-0.9,-0.3) -- ++(0.8,0.6) -- ++(2,0) -- ++(-0.8,-0.6) -- cycle;
	\draw[2facetLine] (0.5,-0.85pt)
		arc[x radius=0.5, y radius=0.15, start angle=0, end angle=180]
		++(-0.9,-0.3) -- ++(0.8,0.6) -- ++(2,0);
	\fill[2facetFront] (0.5, 0.85pt)
		arc[x radius=0.5, y radius=0.15, start angle=0, end angle=180] --
		++(-0.9,-0.3) -- ++(0.8,0.6) -- ++(2,0) -- ++(-0.8,-0.6) -- cycle;
	\draw[2facetLine] (0.5, 0.85pt)
		arc[x radius=0.5, y radius=0.15, start angle=0, end angle=180]
		++(-0.9,-0.3) -- ++(0.8,0.6) -- ++(2,0);
	\fill[1facetBack] (-0.5,-1) -- (-0.5,0)
		.. controls (-0.5,1) and (0.5,1) .. (0.5,0) -- (0.5,-1)
		arc[x radius=0.5, y radius=0.15, start angle=0, end angle=180];
	\draw[1facetLine] (0.5,-1)
		arc[x radius=0.5, y radius=0.15, start angle=0, end angle=180];
	\draw[1facetFront] (-0.5,-1) -- (-0.5,0)
		.. controls (-0.5,1) and (0.5,1) .. (0.5,0) -- (0.5,-1)
		arc[x radius=0.5, y radius=0.15, start angle=360, end angle=180];
	\fill[2facetInner] (-0.5, -0.85pt)
		arc[x radius=0.5, y radius=0.15, start angle=180, end angle=360] --
		++(0.1,-0.3) -- ++(-2,0) -- cycle;
	\draw[2facetLine] (-0.5, -0.85pt)
		arc[x radius=0.5, y radius=0.15, start angle=180, end angle=360]
		++(0.9, 0.3) -- ++(-0.8,-0.6) -- ++(-2,0);
	\fill[2facetFront] (-0.5, 0.85pt)
		arc[x radius=0.5, y radius=0.15, start angle=180, end angle=360] --
		++(0.1,-0.3) -- ++(-2,0) -- cycle;
	\draw[2facetLine] (-0.5, 0.85pt)
		arc[x radius=0.5, y radius=0.15, start angle=180, end angle=360]
		++(0.9, 0.3) -- ++(-0.8,-0.6) -- ++(-2,0);
	\ifx\relax#2\relax\else
		\foamarrow#1[2facetLine](-0.7,-0.3);
		\foamarrow#1[2facetLine]180(0.7,0.3);
	\fi
	\foamarrow#1[seam](0,-0.15);
}

\deffoampict{1cap < 2plane}(-1.4,-1.15)--(1.4,0.8) M=-0.2 [#1]=[.]{%
	\fill[1facetBack] (-0.5,-1)
			arc[x radius=0.5, y radius=0.15, start angle=180, end angle=0]
			.. controls (0.5,0) and (-0.5,0) .. (-0.5,-1);
	\draw[1facetLine]
			(-0.5,-1) arc[x radius=0.5, y radius=0.15, start angle=180, end angle=0];
	\draw[1facetFront] (-0.5,-1)
			arc[x radius=0.5, y radius=0.15, start angle=180, end angle=360]
			.. controls (0.5,0) and (-0.5,0) .. (-0.5,-1);
	\draw[2facetInner] (-0.6, 0.3) ++(0,-0.85pt) -- ++(-0.8,-0.6) --
		++(2,0) -- ++(0.8,0.6) -- ++(-2,0);
	\draw[2facetFront] (-0.6, 0.3) ++(0, 0.85pt) -- ++(-0.8,-0.6) --
		++(2,0) -- ++(0.8,0.6) -- ++(-2,0);
	\foamarrow#1[2facetLine](-0.6,-0.3);
	\foamarrow#1[2facetLine]180(0.6,0.3);
}

\deffoampict{Vsaddle > 2plane}(-1.5,-0.8)--(1.5,1.1) M=0.1 [#1]=[.]{%
	\draw[2facetInner] (-0.7,0.2) ++(0,-0.85pt) --
		++(0.4,0.4) -- ++(1.8,0) -- ++(-0.4,-0.4);
	\draw[2facetFront] (-0.7,0.2) ++(0, 0.85pt) --
		++(0.4,0.4) -- ++(1.8,0) -- ++(-0.4,-0.4);
	\fill[1facetBack] (-0.7,1.1)
			arc[x radius=0.4, y radius=0.2, start angle=90, end angle=0]
			arc[x radius=0.3, y radius=0.5, start angle=180, end angle=360]
			arc[x radius=0.8, y radius=0.2, start angle=180, end angle=90]
			-- (1.1,-0.4) -- (-0.7,-0.4) -- cycle;
	\draw[1facetLine] (1.1, 1.1) |- (-0.7,-0.4) -- (-0.7,1.1);
	\draw[seam]
		(-0.7,0.2) ++(0,-0.85pt) -- ++(1.8,0)
		(-0.7,0.2) ++(0, 0.85pt) -- ++(1.8,0);
	\foamarrow#1[seam]{180}(0.3,0.2);
	\fill[1facetFront] (-1.1,0.7)
			arc[x radius=0.8, y radius=0.2, start angle=270, end angle=360]
			arc[x radius=0.3, y radius=0.5, start angle=180, end angle=360]
			arc[x radius=0.4, y radius=0.2, start angle=180, end angle=270]
			-- (0.7,-0.8) -- (-1.1,-0.8) -- cycle;
	\draw[1facetLine]
			(0.3, 0.9) arc[x radius=0.3, y radius=0.5, start angle=360, end angle=180]
			(-0.7,1.1) arc[x radius=0.4, y radius=0.2, start angle= 90, end angle=0]
			           arc[x radius=0.8, y radius=0.2, start angle=360, end angle=270]
			-- (-1.1,-0.8) -- (0.7,-0.8) -- (0.7,0.7)
			           arc[x radius=0.4, y radius=0.2, start angle=270, end angle=180]
			           arc[x radius=0.8, y radius=0.2, start angle=180, end angle=90];
	\draw[2facetInner] (-1.1,-0.2) ++(0,-0.85pt) --
		++(-0.4,-0.4) -- ++(1.8,0) -- ++(0.4, 0.4);
	\draw[seam] (-1.1,-0.2) ++(0,-0.85pt) -- ++(1.8,0);
	\draw[2facetFront] (-1.1,-0.2) ++(0, 0.85pt) --
		++(-0.4,-0.4) -- ++(1.8,0) -- ++(0.4, 0.4);
	\draw[seam] (-1.1,-0.2) ++(0, 0.85pt) -- ++(1.8,0);
	\foamarrow#1[seam](-0.3,-0.2);
}

\deffoampict{Vsaddle < 2plane}(-1.5,-1.2)--(1.5,0.7) M=-0.3 [#1]=[.]{%
	\fill[2facetInner] (0.3,-0.85pt)
			arc[x radius=0.8, y radius=0.2, start angle=180, end angle=90]
			-- ++(0.4,0.4) -- ++(-1.8,0) -- ++(-0.4,-0.4)
			arc[x radius=0.4, y radius=0.2, start angle=90, end angle=0];
	\draw[2facetLine] (1.1,0.2) ++(0,-0.85pt)
		-- ++(0.4,0.4) -- ++(-1.8,0) -- ++(-0.4,-0.4);
	\fill[2facetFront] (0.3, 0.85pt)
			arc[x radius=0.8, y radius=0.2, start angle=180, end angle=90]
			-- ++(0.4,0.4) -- ++(-1.8,0) -- ++(-0.4,-0.4)
			arc[x radius=0.4, y radius=0.2, start angle=90, end angle=0];
	\draw[2facetLine] (1.1,0.2) ++(0, 0.85pt)
		-- ++(0.4,0.4) -- ++(-1.8,0) -- ++(-0.4,-0.4);
	\fill[1facetBack] (-0.7,0.7)
			arc[x radius=0.4, y radius=0.2, start angle=90, end angle=0] -- (-0.3,0)
			arc[x radius=0.3, y radius=0.4, start angle=180, end angle=360] -- (0.3,0.5)
			arc[x radius=0.8, y radius=0.2, start angle=180, end angle=90]
			-- (1.1,-0.8) -- (-0.7,-0.8) -- cycle;
	\draw[1facetLine] (1.1,0.7) |- (-0.7,-0.8) -- (-0.7,0.7);
	\draw[seam]
		( 0.3,-0.85pt) arc[x radius=0.8, y radius=0.2, start angle=180, end angle=90]
		( 0.3, 0.85pt) arc[x radius=0.8, y radius=0.2, start angle=180, end angle=90]
		(-0.3,-0.85pt) arc[x radius=0.4, y radius=0.2, start angle=0, end angle=90]
		(-0.3, 0.85pt) arc[x radius=0.4, y radius=0.2, start angle=0, end angle=90];
	\fill[1facetFront] (-1.1, 0.3)
			arc[x radius=0.8, y radius=0.2, start angle=270, end angle=360] -- (-0.3,0)
			arc[x radius=0.3, y radius=0.4, start angle=180, end angle=360] -- (0.3,0.5)
			arc[x radius=0.4, y radius=0.2, start angle=180, end angle=270]
			-- (0.7,-1.2) -- (-1.1,-1.2) -- cycle;
	\draw[1facetLine] (-0.3,0.5) -- (-0.3,0)
			arc[x radius=0.3, y radius=0.4, start angle=180, end angle=360] -- (0.3,0.5)
			(-0.7,0.7)
			arc[x radius=0.4, y radius=0.2, start angle=90, end angle=0]
			arc[x radius=0.8, y radius=0.2, start angle=360, end angle=270]
			|- (0.7,-1.2) -- (0.7,0.3)
			arc[x radius=0.4, y radius=0.2, start angle=270, end angle=180]
			arc[x radius=0.8, y radius=0.2, start angle=180, end angle=90];
	\fill[2facetInner] (0.3,-0.85pt)
			arc[x radius=0.4, y radius=0.2, start angle=180, end angle=270]
			-- ++(-0.4,-0.4) -- ++(-1.8,0) -- ++(0.4,0.4)
			arc[x radius=0.8, y radius=0.2, start angle=270, end angle=360];
	\draw[2facetLine] (0.7,-0.2) ++(0,-0.85pt) -- ++(-0.4,-0.4) -- ++(-1.8,0) -- ++(0.4,0.4);
	\draw[seam]
		( 0.3,-0.85pt) arc[x radius=0.4, y radius=0.2, start angle=180, end angle=270]
		(-0.3,-0.85pt) arc[x radius=0.8, y radius=0.2, start angle=360, end angle=270];
	\fill[2facetFront] (0.3, 0.85pt)
			arc[x radius=0.4, y radius=0.2, start angle=180, end angle=270]
			-- ++(-0.4,-0.4) -- ++(-1.8,0) -- ++(0.4,0.4)
			arc[x radius=0.8, y radius=0.2, start angle=270, end angle=360];
	\draw[2facetLine] (0.7,-0.2) ++(0, 0.85pt) -- ++(-0.4,-0.4) -- ++(-1.8,0) -- ++(0.4,0.4);
	\draw[seam]
		( 0.3, 0.85pt) arc[x radius=0.4, y radius=0.2, start angle=180, end angle=270]
		(-0.3, 0.85pt) arc[x radius=0.8, y radius=0.2, start angle=360, end angle=270];
	\foamarrow#1[seam]-19(0.468,-0.163);
	\foamarrow#1[seam]11(-0.782,-0.184);
}

\deffoampict{Asaddle > 2plane}(-1.5,-0.7)--(1.5,1.2) M=0.1 [#1]=[.]{%
	\fill[2facetInner] (0.3,-0.85pt)
			arc[x radius=0.8, y radius=0.2, start angle=180, end angle=90]
			-- ++(0.4,0.4) -- ++(-1.8,0) -- ++(-0.4,-0.4)
			arc[x radius=0.4, y radius=0.2, start angle=90, end angle=0];
	\draw[2facetLine] (1.1,0.2) ++(0,-0.85pt)
		-- ++(0.4,0.4) -- ++(-1.8,0) -- ++(-0.4,-0.4);
	\fill[2facetFront] (0.3, 0.85pt)
			arc[x radius=0.8, y radius=0.2, start angle=180, end angle=90]
			-- ++(0.4,0.4) -- ++(-1.8,0) -- ++(-0.4,-0.4)
			arc[x radius=0.4, y radius=0.2, start angle=90, end angle=0];
	\draw[2facetLine] (1.1,0.2) ++(0, 0.85pt)
		-- ++(0.4,0.4) -- ++(-1.8,0) -- ++(-0.4,-0.4);
	\fill[1facetBack] (-0.7,-0.3)
			arc[x radius=0.4, y radius=0.2, start angle=90, end angle=0] -- (-0.3,0)
			arc[x radius=0.3, y radius=0.4, start angle=180, end angle=0] -- (0.3,-0.5)
			arc[x radius=0.8, y radius=0.2, start angle=180, end angle=90]
			-- (1.1,1.2) -- (-0.7,1.2) -- cycle;
	\draw[1facetLine] (-0.3,-0.5)
			arc[x radius=0.4, y radius=0.2, start angle=0, end angle=90]
			-- (-0.7,1.2) -| (1.1,-0.3)
			arc[x radius=0.8, y radius=0.2, start angle=90, end angle=180];
	\draw[seam]
		( 0.3,-0.85pt) arc[x radius=0.8, y radius=0.2, start angle=180, end angle=90]
		( 0.3, 0.85pt) arc[x radius=0.8, y radius=0.2, start angle=180, end angle=90]
		(-0.3,-0.85pt) arc[x radius=0.4, y radius=0.2, start angle=0, end angle=90]
		(-0.3, 0.85pt) arc[x radius=0.4, y radius=0.2, start angle=0, end angle=90];
	\draw[1facetFront] (-1.1,-0.7)
			arc[x radius=0.8, y radius=0.2, start angle=-90, end angle=0] -- (-0.3,0)
			arc[x radius=0.3, y radius=0.4, start angle=180, end angle=0] -- (0.3,-0.5)
			arc[x radius=0.4, y radius=0.2, start angle=180, end angle=270]
			|- (-1.1,0.8) -- cycle;
	\fill[2facetInner] (0.3,-0.85pt)
			arc[x radius=0.4, y radius=0.2, start angle=180, end angle=270]
			-- ++(-0.4,-0.4) -- ++(-1.8,0) -- ++(0.4,0.4)
			arc[x radius=0.8, y radius=0.2, start angle=270, end angle=360];
	\draw[2facetLine] (0.7,-0.2) ++(0,-0.85pt) -- ++(-0.4,-0.4) -- ++(-1.8,0) -- ++(0.4,0.4);
	\draw[seam]
		( 0.3,-0.85pt) arc[x radius=0.4, y radius=0.2, start angle=180, end angle=270]
		(-0.3,-0.85pt) arc[x radius=0.8, y radius=0.2, start angle=360, end angle=270];
	\fill[2facetFront] (0.3, 0.85pt)
			arc[x radius=0.4, y radius=0.2, start angle=180, end angle=270]
			-- ++(-0.4,-0.4) -- ++(-1.8,0) -- ++(0.4,0.4)
			arc[x radius=0.8, y radius=0.2, start angle=270, end angle=360];
	\draw[2facetLine] (0.7,-0.2) ++(0, 0.85pt) -- ++(-0.4,-0.4) -- ++(-1.8,0) -- ++(0.4,0.4);
	\draw[seam]
		( 0.3, 0.85pt) arc[x radius=0.4, y radius=0.2, start angle=180, end angle=270]
		(-0.3, 0.85pt) arc[x radius=0.8, y radius=0.2, start angle=360, end angle=270];
	\foamarrow#1[seam]-19(0.468,-0.163);
	\foamarrow#1[seam]11(-0.782,-0.184);
}

\deffoampict{Asaddle < 2plane}(-1.5,-1.1)--(1.5,0.8) M=-0.3 [#1]=[.]{%
	\draw[2facetInner] (-0.7,0.2) ++(0,-0.85pt)
		-- ++(0.4,0.4) -- ++(1.8,0) -- ++(-0.4,-0.4);
	\draw[2facetFront] (-0.7,0.2) ++(0, 0.85pt)
		-- ++(0.4,0.4) -- ++(1.8,0) -- ++(-0.4,-0.4);
	\fill[1facetBack] (-0.7,-0.7)
			arc[x radius=0.4, y radius=0.2, start angle=90, end angle=0]
			arc[x radius=0.3, y radius=0.5, start angle=180, end angle=0]
			arc[x radius=0.8, y radius=0.2, start angle=180, end angle=90]
			-- (1.1,0.8) -- (-0.7,0.8) -- cycle;
	\draw[1facetLine] (-0.3,-0.9)
			arc[x radius=0.4, y radius=0.2, start angle=0, end angle=90]
			|- (1.1,0.8) -- (1.1,-0.7)
			arc[x radius=0.8, y radius=0.2, start angle=90, end angle=180];
	\draw[seam]
		(-0.7,0.2) ++(0,-0.85pt) -- ++(1.8,0)
		(-0.7,0.2) ++(0, 0.85pt) -- ++(1.8,0);
	\foamarrow#1[seam]{180}(0.3,0.2);
	\draw[1facetFront] (-1.1,-1.1)
			arc[x radius=0.8, y radius=0.2, start angle=-90, end angle=0]
			arc[x radius=0.3, y radius=0.5, start angle=180, end angle=0]
			arc[x radius=0.4, y radius=0.2, start angle=180, end angle=270]
			|- (-1.1,0.4) -- cycle;
	\draw[2facetInner] (-1.1,-0.2) ++(0,-0.85pt)
		-- ++(-0.4,-0.4)-- ++(1.8,0) -- ++(0.4,0.4);
	\draw[seam] (-1.1,-0.2) ++(0,-0.85pt) -- ++(1.8,0);
	\draw[2facetFront] (-1.1,-0.2) ++(0, 0.85pt)
		-- ++(-0.4,-0.4)-- ++(1.8,0) -- ++(0.4,0.4);
	\draw[seam] (-1.1,-0.2) ++(0, 0.85pt) -- ++(1.8,0);
	\foamarrow#1[seam](-0.3,-0.2);
}

\deffoampict{cup}(-0.7, 0.2)--(0.7,-1.4) [#1]=[1]{%
	\fill[#1facetBack]
		(-0.5,0) .. controls ++(0,-1.0) and ++(0,-1.0) .. (0.5,0)
		arc[x radius=0.5, y radius=0.15, start angle=0, end angle=180];
	\fill[#1facetFront]
		(-0.5,0) .. controls ++(0,-1.0) and ++(0,-1.0) .. (0.5,0)
		arc[x radius=0.5, y radius=0.15, start angle=360, end angle=180];
	\draw[#1facetLine]
		(-0.5,0) .. controls ++(0,-1.0) and ++(0,-1.0) .. (0.5,0)
		(0,0) ellipse[x radius=0.5, y radius=0.15];
}

\deffoampict{cap}(-0.7,-0.2)--(0.7,1.4) [#1]=[1]{%
	\fill[#1facetBack]
		(-0.5,0) .. controls ++(0,1.0) and ++(0,1.0) .. (0.5,0)
		arc[x radius=0.5, y radius=0.15, start angle=0, end angle=180];
	\draw[#1facetLine]
		(0.5,0) arc[x radius=0.5, y radius=0.15, start angle=0, end angle=180];
	\draw[#1facetFront]
		(-0.5,0) .. controls ++(0,1.0) and ++(0,1.0) .. (0.5,0)
		arc[x radius=0.5, y radius=0.15, start angle=360, end angle=180];
}

\deffoampict{saddle}(-1,-1.1)--(1,0.5) [#1]=[1]{%
	\fill[#1facetBack] (-0.6,0.45)
			arc[x radius=0.3, y radius=0.15, start angle=90, end angle=0]
			arc[x radius=0.3, y radius=0.60, start angle=180, end angle=360]
			arc[x radius=0.6, y radius=0.15, start angle=180, end angle=90]
			-- (0.9,-0.7) -- (-0.6,-0.7) -- cycle;
	\draw[#1facetLine] (0.9,0.45) |- (-0.6,-0.7) -- (-0.6,0.45);
	\fill[#1facetFront] (-0.9, 0.15)
			arc[x radius=0.6, y radius=0.15, start angle=270, end angle=360]
			arc[x radius=0.3, y radius=0.60, start angle=180, end angle=360]
			arc[x radius=0.3, y radius=0.15, start angle=180, end angle=270]
			-- (0.6,-1) -- (-0.9,-1) -- cycle;
	\draw[#1facetLine] (-0.3,0.3)
			arc[x radius=0.3, y radius=0.6, start angle=180, end angle=360]
			(-0.6,0.45)
			arc[x radius=0.3, y radius=0.15, start angle=90, end angle=0]
			arc[x radius=0.6, y radius=0.15, start angle=360, end angle=270]
			|- (0.6,-1) -- (0.6,0.15)
			arc[x radius=0.3, y radius=0.15, start angle=270, end angle=180]
			arc[x radius=0.6, y radius=0.15, start angle=180, end angle=90];
}

\deffoampict{dot}(-0.1,-0.1)--(1.3,1.6) [#1]=[1]{%
	\draw[#1facetFront]
		(0,0) .. controls ++(0.7,0.1) and ++(-0.5,-0.3) .. (1.2,0.4) --
		(1.2,1.5) .. controls ++(-0.5,-0.3) and ++(0.7,0.1) .. (0,1.1) -- cycle;
	\foamdot (0.65,0.7);
}

\deffoampict{merge}(-1.6,-0.2)--(1.6,2.2) M=1 [#1]=[1]{%
	\fill[#1facetBack] (-1.5,0)
		arc[x radius=0.5, y radius=0.15, start angle=180, end angle=0]
		.. controls ++(0,0.8) and ++(0,0.8) .. (0.5,0)
		arc[x radius=0.5, y radius=0.15, start angle=180, end angle=0]
		.. controls ++(0,1.2) and ++(0,-1) .. (0.5, 2)
		arc[x radius=0.5, y radius=0.15, start angle=0, end angle=180]
		.. controls ++(0,-1) and ++(0,1.2) .. (-1.5, 0);
	\draw[#1facetLine]
		(-1.5,0) arc[x radius=0.5, y radius=0.15, start angle=180, end angle=0]
		( 0.5,0) arc[x radius=0.5, y radius=0.15, start angle=180, end angle=0];
	\fill[#1facetFront] (-0.5,2)
		.. controls ++(0,-1) and ++(0,1.2) .. (-1.5,0)
		arc[x radius=0.5, y radius=0.15, start angle=180, end angle=360]
		.. controls ++(0,0.8) and ++(0,0.8) .. (0.5,0)
		arc[x radius=0.5, y radius=0.15, start angle=180, end angle=360]
		.. controls ++(0,1.2) and ++(0,-1) .. (0.5, 2)
		arc[x radius=0.5, y radius=0.15, start angle=360, end angle=180];
	\draw[#1facetLine] (-0.5,2)
		.. controls ++(0,-1) and ++(0,1.2) .. (-1.5,0)
		arc[x radius=0.5, y radius=0.15, start angle=180, end angle=360]
		.. controls ++(0,0.8) and ++(0,0.8) .. (0.5,0)
		arc[x radius=0.5, y radius=0.15, start angle=180, end angle=360]
		.. controls ++(0,1.2) and ++(0,-1) .. (0.5, 2)
		(0,2) ellipse[x radius=0.5, y radius=0.15];
}

\deffoampict{split}(-1.6, 0.2)--(1.6,-2.2) M=-1 [#1]=[1]{%
	\fill[#1facetBack] (-1.5,0)
		arc[x radius=0.5, y radius=0.15, start angle=180, end angle=0]
		.. controls ++(0,-0.8) and ++(0,-0.8) .. (0.5,0)
		arc[x radius=0.5, y radius=0.15, start angle=180, end angle=0]
		.. controls ++(0,-1.2) and ++(0,1) .. (0.5,-2)
		arc[x radius=0.5, y radius=0.15, start angle=0, end angle=180]
		.. controls ++(0,1) and ++(0,-1.2) .. (-1.5, 0);
	\draw[#1facetLine]
		(-0.5,-2) arc[x radius=0.5, y radius=0.15, start angle=180, end angle=0];
	\fill[#1facetFront] (-1.5,0)
		arc[x radius=0.5, y radius=0.15, start angle=180, end angle=360]
		.. controls ++(0,-0.8) and ++(0,-0.8) .. (0.5,0)
		arc[x radius=0.5, y radius=0.15, start angle=180, end angle=360]
		.. controls ++(0,-1.2) and ++(0,1) .. (0.5,-2)
		arc[x radius=0.5, y radius=0.15, start angle=360, end angle=180]
		.. controls ++(0,1) and ++(0,-1.2) .. (-1.5, 0);
	\draw[#1facetLine]
		(-0.5,0) .. controls ++(0,-0.8) and ++(0,-0.8) .. (0.5,0)
		( 1.5,0) .. controls ++(0,-1.2) and ++(0,1) .. (0.5,-2)
		arc[x radius=0.5, y radius=0.15, start angle=360, end angle=180]
		.. controls ++(0,1) and ++(0,-1.2) .. (-1.5, 0)
		(-1,0) ellipse[x radius=0.5, y radius=0.15]
		( 1,0) ellipse[x radius=0.5, y radius=0.15];
}

\deffoampict{pocket above}(-0.5,-0.95)--(0.85,0.95) {%
	\fill[2facetBack] (0.15,0.7)
		arc[start angle= 90, end angle=180, x radius=0.60, y radius=0.22]
		arc[start angle=180, end angle=270, x radius=0.45, y radius=0.48]
		arc[start angle=270, end angle=360, x radius=0.15, y radius=0.70];
	\fill[2facetFront] (-0.15,0.3)
		arc[start angle=270, end angle=180, x radius=0.30, y radius=0.18]
		arc[start angle=180, end angle=270, x radius=0.45, y radius=0.48]
		arc[start angle=270, end angle=180, x radius=0.15, y radius=0.30];
	\draw[2facetLine] (0,0)
		arc[start angle=270, end angle=180, x radius=0.45, y radius=0.48]
		(0.15,0.7)
		arc[start angle= 90, end angle=180, x radius=0.60, y radius=0.22]
		arc[start angle=180, end angle=270, x radius=0.30, y radius=0.18];
	\draw[1facetFront]
		(-0.3, 0.1) -- (0.3, 0.9) -- (0.3, -0.1) -- (-0.3, -0.9) -- cycle;
	\fill[2facetBack] (0.15, 0.7) -- (0.8, 0.7) -- (0.8, -0.3)
		arc[start angle=90, end angle=180, x radius=0.5, y radius=0.22]
		arc[start angle=0, end angle=90, x radius=0.3, y radius=0.52]
		arc[start angle=270, end angle=360, x radius=0.15, y radius=0.7];
	\draw[seam] (0,0) arc[start angle=270, end angle=360, x radius=0.15, y radius=0.7];
	\draw[2facetLine] (0.15, 0.7) -| (0.8, -0.3)
		arc[start angle=90, end angle=180, x radius=0.5, y radius=0.22];
	\fill[2facetFront] (-0.15, 0.3) -- (0.5, 0.3) -- (0.5, -0.7)
		arc[start angle=270, end angle=180, x radius=0.2, y radius=0.18]
		arc[start angle=0, end angle=90, x radius=0.3, y radius=0.52]
		arc[start angle=270, end angle=180, x radius=0.15, y radius=0.3];
	\draw[seam] (0,0) arc[start angle=270, end angle=180, x radius=0.15, y radius=0.3];
	\draw[2facetLine] (-0.15, 0.3) -- (0.5, 0.3) -- (0.5, -0.7)
		arc[start angle=270, end angle=180, x radius=0.2, y radius=0.18]
		arc[start angle=0, end angle=90, x radius=0.3, y radius=0.52];
}

\deffoampict{pocket below}(-0.5,-0.95)--(0.9,0.95) {%
	\fill[2facetBack] (0.15,-0.3)
		arc[start angle= 90, end angle=180, x radius=0.60, y radius=0.22]
		arc[start angle=180, end angle= 90, x radius=0.45, y radius=0.52]
		arc[start angle= 90, end angle=  0, x radius=0.15, y radius=0.30];
	\draw[2facetLine] (0.15,-0.3)
		arc[start angle= 90, end angle=180, x radius=0.60, y radius=0.22];
	\fill[2facetFront] (-0.15,-0.7)
		arc[start angle=270, end angle=180, x radius=0.30, y radius=0.18]
		arc[start angle=180, end angle= 90, x radius=0.45, y radius=0.52]
		arc[start angle= 90, end angle=180, x radius=0.15, y radius=0.7];
	\draw[2facetLine] (-0.15,-0.7)
		arc[start angle=270, end angle=180, x radius=0.30, y radius=0.18]
		arc[start angle=180, end angle= 90, x radius=0.45, y radius=0.52];
	\draw[1facetFront]
		(-0.3, 0.1) -- (0.3, 0.9) -- (0.3, -0.1) -- (-0.3, -0.9) -- cycle;
	\fill[2facetBack] (0.15,-0.3) -| (0.85, 0.7)
		arc[start angle=90, end angle=180, x radius=0.5, y radius=0.22]
		arc[start angle=0, end angle=-90, x radius=0.35, y radius=0.48]
		arc[start angle=90, end angle=0, x radius=0.15, y radius=0.3];
	\draw[seam] (0, 0) arc[start angle=90, end angle=0, x radius=0.15, y radius=0.3];
	\draw[2facetLine] (0.15,-0.3) -| (0.85,0.7)
		arc[start angle=90, end angle=180, x radius=0.5, y radius=0.22];
	\fill[2facetFront] (0.35,0.48)
		arc[start angle=0, end angle=-90, x radius=0.35, y radius=0.48]
		arc[start angle=90, end angle=180, x radius=0.15, y radius=0.7] -| (0.55,0.3)
		arc[start angle=270, end angle=180, x radius=0.2, y radius=0.18];
	\draw[seam] (0,0) arc[start angle=90, end angle=180, x radius=0.15, y radius=0.7];
	\draw[2facetLine] (-0.15,-0.7) -| (0.55,0.3)
		arc[start angle=270, end angle=180, x radius=0.2, y radius=0.18]
		(0,0) arc[start angle=270, end angle=360, x radius=0.35, y radius=0.48];
}

\deffoampict{red cup > blue plane}(-1.9,-0.5)--(1.8,1.3) M=0.2 {%
	\draw[1facetFront] (-1.8,-0.4) -- (-0.6,0.4) -- (1.7,0.4) -- (0.5,-0.4) -- cycle;
	\fill[2facetBack] (-0.5,1) .. controls (-0.5,0) and (0.5,0) .. (0.5,1)
		arc[x radius=0.5, y radius=0.15, start angle=0, end angle=180];
	\fill[2facetFront] (-0.5,1) .. controls (-0.5,0) and (0.5,0) .. (0.5,1)
		arc[x radius=0.5, y radius=0.15, start angle=360, end angle=180];
	\draw[2facetLine]
		(-0.5,1) .. controls (-0.5,0) and (0.5,0) .. (0.5,1)
	   (0,1) ellipse[x radius=0.5, y radius=0.15];
}

\deffoampict{red cup < blue plane}(-1.9,-0.8)--(1.8,1.3) M=0.2 {%
	\fill[2facetBack] (-0.5,0) .. controls (-0.5,-1) and (0.5,-1) .. (0.5,0)
		arc[x radius=0.5, y radius=0.15, start angle=0, end angle=180];
	\fill[2facetFront] (-0.5,0) .. controls (-0.5,-1) and (0.5,-1) .. (0.5,0)
		arc[x radius=0.5, y radius=0.15, start angle=360, end angle=180];
	\draw[2facetLine] (-0.5,0) .. controls (-0.5,-1) and (0.5,-1) .. (0.5,0);
	\draw[1facetFront] (-1.8,-0.4) -- (-0.6,0.4) -- (1.7,0.4) -- (0.5,-0.4) -- cycle;
	\fill[2facetBack] (0.5,0)
		-- ( 0.5,1) arc[x radius=0.5, y radius=0.15, start angle=0, end angle=180]
		-- (-0.5,0) arc[x radius=0.5, y radius=0.15, start angle=180, end angle=0];
	\draw[seam]
		(0.5,0) arc[x radius=0.5, y radius=0.15, start angle=0, end angle=180];
	\fill[2facetFront] (0.5,0)
		-- ( 0.5,1) arc[x radius=0.5, y radius=0.15, start angle=360, end angle=180]
		-- (-0.5,0) arc[x radius=0.5, y radius=0.15, start angle=180, end angle=360];
	\draw[2facetLine]
		(0,1) ellipse[x radius=0.5, y radius=0.15]
		( 0.5,1) -- ( 0.5,0)
		(-0.5,1) -- (-0.5,0);
	\draw[seam] (0.5,0)
		arc[x radius=0.5, y radius=0.15, start angle=360, end angle=180];
}

\deffoampict{red saddle > blue plane}(-1.5,-0.7)--(1.5,1.2) M=0.1 {%
	\fill[2facetBack]
		(0.3,-0.5) arc[x radius=0.8, y radius=0.2, start angle=180, end angle=90] --
		(1.1, 0.2) arc[x radius=0.8, y radius=0.2, start angle=90, end angle=180];
	\fill[2facetBack]
		(-0.3,-0.5) arc[x radius=0.4, y radius=0.2, start angle=0, end angle=90] --
		(-0.7, 0.2) arc[x radius=0.4, y radius=0.2, start angle=90, end angle=0];
	\draw[2facetLine]
		( 0.3,-0.5) arc[x radius=0.8, y radius=0.2, start angle=180, end angle=90]
		-- ( 1.1,0.2)
		(-0.3,-0.5) arc[x radius=0.4, y radius=0.2, start angle=0, end angle=90]
		-- (-0.7,0.2);
	\fill[2facetFront]
		( 0.3,-0.5) arc[x radius=0.4, y radius=0.2, start angle=180, end angle=270] --
		( 0.7,-0.2) arc[x radius=0.4, y radius=0.2, start angle=270, end angle=180];
	\fill[2facetFront]
		(-0.3,-0.5) arc[x radius=0.8, y radius=0.2, start angle=0, end angle=-90] --
		(-1.1,-0.2) arc[x radius=0.8, y radius=0.2, start angle=-90, end angle=0];
	\draw[2facetLine] (0.3,0) -- ( 0.3,-0.5)
		arc[x radius=0.4, y radius=0.2, start angle=180, end angle=270] -- (0.7,-0.2)
		(-0.3,0) -- (-0.3,-0.5)
		arc[x radius=0.8, y radius=0.2, start angle=0, end angle=-90] -- (-1.1,-0.2);
	\draw[1facetFront] (1.4,0.5) -- (0.4,-0.5) -- (-1.4,-0.5) -- (-0.4,0.5) -- cycle;
	\fill[2facetBack] (-0.7,0.2)
		arc[x radius=0.4, y radius=0.2, start angle=90, end angle=0]
		arc[x radius=0.3, y radius=0.4, start angle=180, end angle=0]
		arc[x radius=0.8, y radius=0.2, start angle=180, end angle=90] |- (-0.7,1.2);
	\draw[seam] ( 0.3,0) arc[x radius=0.8, y radius=0.2, start angle=180, end angle=90]
	            (-0.3,0) arc[x radius=0.4, y radius=0.2, start angle=  0, end angle=90];
	\draw[2facetLine] (1.1,0.2) |- (-0.7,1.2) -- (-0.7,0.2);
	\fill[2facetFront] (-1.1,0.8) -- (-1.1,-0.2)
		arc[x radius=0.8, y radius=0.2, start angle=270, end angle=360]
		arc[x radius=0.3, y radius=0.4, start angle=180, end angle=0  ]
		arc[x radius=0.4, y radius=0.2, start angle=180, end angle=270] |- cycle;
	\draw[seam] ( 0.3,0) arc[x radius=0.4, y radius=0.2, start angle=180, end angle=270]
	            (-0.3,0) arc[x radius=0.8, y radius=0.2, start angle=360, end angle=270];
	\draw[2facetLine] (-1.1,-0.2) -- (-1.1,0.8) -| (0.7,-0.2);
}

\deffoampict{red saddle < blue plane}(-1.5,-1.1)--(1.5,0.8) M=-0.3 {%
	\fill[2facetBack] (-0.7,-0.7)
		arc[x radius=0.4, y radius=0.2, start angle=90, end angle=0]
		arc[x radius=0.3, y radius=0.5, start angle=180, end angle=0]
		arc[x radius=0.8, y radius=0.2, start angle=180, end angle=90] |- (-0.7,0.2);
	\draw[2facetLine] (0.3,-0.9)
		arc[x radius=0.8, y radius=0.2, start angle=180, end angle=90] -- (1.1,0.2)
		(-0.3,-0.9)
		arc[x radius=0.4, y radius=0.2, start angle=0, end angle=90] -- (-0.7,0.2);
	\draw[2facetFront] (-1.1,-0.2) -- (-1.1,-1.1)
		arc[x radius=0.8, y radius=0.2, start angle=270, end angle=360]
		arc[x radius=0.3, y radius=0.5, start angle=180, end angle=0  ]
		arc[x radius=0.4, y radius=0.2, start angle=180, end angle=270] -- (0.7,-0.2);
	\draw[1facetFront] (1.4,0.5) -- (0.4,-0.5) -- (-1.4,-0.5) -- (-0.4,0.5) -- cycle;
	\draw[2facetBack]  (-0.7, 0.2) -- (-0.7, 0.8) -- (1.1, 0.8) -- (1.1, 0.2);
	\draw[seam] (-0.7, 0.2) -- (1.1, 0.2);
	\draw[2facetFront] (-1.1,-0.2) -- (-1.1, 0.4) -- (0.7, 0.4) -- (0.7,-0.2);
	\draw[seam] (-1.1,-0.2) -- (0.7,-0.2);
}

\deffoampict{inner cylinder}(-0.9,-0.5)--(1.2,2.5) M=1 {%
	\draw[2facetBack] (-0.8,0)
		arc[x radius=0.8, y radius=0.4, start angle=180, end angle=90]
		-- (1.1, 0.4) |- (0,2.4)
		arc[x radius=0.8, y radius=0.4, start angle= 90, end angle=180];
	\fill[1facetBack] (-0.4,0)
		arc[x radius=0.4, y radius=0.15, start angle=180, end angle=0] -- ++(0,2)
		arc[x radius=0.4, y radius=0.15, start angle=0, end angle=180];
	\draw[1facetLine]
		(-0.4,0) arc[x radius=0.4, y radius=0.15, start angle=180, end angle=0];
	\fill[1facetFront] (-0.4,0)
		arc[x radius=0.4, y radius=0.15, start angle=180, end angle=360] -- ++(0,2)
		arc[x radius=0.4, y radius=0.15, start angle=360, end angle=180];
	\draw[1facetLine]
		(0,2) ellipse[x radius=0.4, y radius=0.15]
		(-0.4,2) -- (-0.4,0)
		arc[x radius=0.4, y radius=0.15, start angle=180, end angle=360] -- ++(0,2);
	\draw[2facetFront] (-0.8,2) -- (-0.8,0)
		arc[x radius=0.8, y radius=0.4, start angle=180, end angle=270]
		-| (0.8, 1.6) -- (0,1.6)
		arc[x radius=0.8, y radius=0.4, start angle=270, end angle=180];
}

\deffoampict{inner cylinder cut}(-0.9,-0.5)--(1.2,2.5) M=1 {%
	\begin{scope}
		\clip (-0.6, 1) rectangle (0.6, 2);
		\fill[1facetBack, rotate around={-10:(0,1.5)}]
			(-0.4,1.5) arc[x radius=0.4, y radius=0.12, start angle=180, end angle=0]
			-- ++(-80:0.5) -- ++(190:0.8) -- cycle;
		\fill[1facetFront, rotate around={-10:(0,1.5)}]
			(-0.4,1.5) arc[x radius=0.4, y radius=0.12, start angle=180, end angle=360]
			-- ++(-80:0.5) -- ++(190:0.8);
	\end{scope}
	\draw[1facetLine]
		(-0.4,1.569) -- (-0.4,1)
		( 0.4,1.431) -- ( 0.4,1);
	\draw[2facetBack] (-0.8,0)
		.. controls (-0.8,1.0) and ( 0.6,0.5) .. ( 0.6,1)
		.. controls ( 0.6,1.5) and (-0.8,1.0) .. (-0.8,2)
		arc[x radius=0.8, y radius=0.4, start angle=180, end angle=90]
		-| (1.1, 0.4) -- (0,0.4)
		arc[x radius=0.8, y radius=0.4, start angle= 90, end angle=180];
	\begin{scope}
		\clip (-0.4,2) arc[x radius=0.4, y radius=0.15, start angle=180, end angle=0]
			-| ( 0.6,1) -| (-0.6,2) -- cycle;
		\fill[1facetBack, rotate around={-10:(0,1.5)}]
			(0.4,1.5) arc[x radius=0.4, y radius=0.12, start angle=0, end angle=180]
			-- ++(100:0.7) -- ++(10:0.8) -- cycle;
	\end{scope}
	\begin{scope}
		\clip (-0.4,2) arc[x radius=0.4, y radius=0.15, start angle=180, end angle=360]
			-| ( 0.6,1) -| (-0.6,2) -- cycle;
		\fill[1facetFront, rotate around={-10:(0,1.5)}]
			(0.4,1.5) arc[x radius=0.4, y radius=0.12, start angle=360, end angle=180]
			-- ++(100:0.7) -- ++(10:0.8) -- cycle;
	\end{scope}
	\draw[seam] (0,1.5) ellipse[x radius=0.4, y radius=0.12, rotate=-10];
	\draw[1facetLine]
		(0,2) ellipse[x radius=0.4, y radius=0.15]
		(-0.4,1.569) -- (-0.4,2)
		( 0.4,1.431) -- ( 0.4,2);
	\begin{scope}
		\clip (-0.4,0) arc[x radius=0.4, y radius=0.15, start angle=180, end angle=0]
			-| ( 0.6,1) -| (-0.6,2) -- cycle;
		\fill[1facetBack, rotate around={10:(0,0.5)}]
			(-0.4,0.5) arc[x radius=0.4, y radius=0.12, start angle=180, end angle=0]
			-- ++(260:0.7) -- ++(170:0.8) -- cycle;
		\draw[1facetLine]
			(-0.4,0) arc[x radius=0.4, y radius=0.15, start angle=180, end angle=0];
	\end{scope}
	\begin{scope}
		\clip (-0.4,0) arc[x radius=0.4, y radius=0.15, start angle=180, end angle=360]
			-| ( 0.6,1) -| (-0.6,2) -- cycle;
		\fill[1facetFront, rotate around={10:(0,0.5)}]
			(-0.4,0.5) arc[x radius=0.4, y radius=0.12, start angle=180, end angle=360]
			-- ++(260:0.8) -- ++(170:0.8) -- cycle;
	\end{scope}
	\draw[1facetLine]
		(-0.4,0.431) -- (-0.4,0)
		arc[x radius=0.4, y radius=0.15, start angle=180, end angle=360] -- (0.4,0.569);
	\draw[2facetFront] (-0.8,2)
		arc[x radius=0.8, y radius=0.4, start angle=180, end angle=270]
		-| (0.8,-0.4) -- (0,-0.4)
		arc[x radius=0.8, y radius=0.4, start angle=270, end angle=180]
		.. controls (-0.8,1.0) and ( 0.6,0.5) .. ( 0.6,1)
		.. controls ( 0.6,1.5) and (-0.8,1.0) .. (-0.8,2);
	\draw[seam,shift={(0,0.5)}] [rotate=10] (170:0.4 and 0.12)
		arc[x radius=0.4, y radius=0.12, start angle=170, end angle=-10];
	\begin{scope}
		\clip (-0.6, 1) rectangle (0.6, 0);
		\fill[1facetBack, rotate around={10:(0,0.5)}]
			(-0.4,0.5) arc[x radius=0.4, y radius=0.12, start angle=180, end angle=0]
			-- ++(80:0.5) -- ++(170:0.8) -- cycle;
		\fill[1facetFront, rotate around={10:(0,0.5)}]
			(-0.4,0.5) arc[x radius=0.4, y radius=0.12, start angle=180, end angle=360]
			-- ++(80:0.5) -- ++(170:0.8) -- cycle;
	\draw[seam,shift={(0,0.5)}] [rotate=10] (170:0.4 and 0.12)
		arc[x radius=0.4, y radius=0.12, start angle=170, end angle=370];
	\end{scope}
	\draw[1facetLine]
		( 0.4,0.569) -- ( 0.4,1)
		(-0.4,0.431) -- (-0.4,1);
}

\deffoampict{cutting plane}(-0.9,-0.9)--(1.2,2.9) M=1 {%
	\draw[1facetFront] (0.3,2.4) -- (0.45,2.8) -- (0.45,0.8) -- (0.3,0.4);
	\draw[2facetBack] (-0.8,0)
		arc[x radius=0.8, y radius=0.4, start angle=180, end angle=90]
		-- (1.1, 0.4) |- (0,2.4)
		arc[x radius=0.8, y radius=0.4, start angle= 90, end angle=180];
	\fill[1facetFront] (0, 1.6) -- (0.3,2.4) -- (0.3,0.4) -- (0,-0.4);
	\draw[seam] (0.3,2.4) -- (0.3,0.4);
	\draw[1facetLine] (0.3,0.4) -- (0,-0.4);
	\draw[2facetFront] (-0.8,2) -- (-0.8,0)
		arc[x radius=0.8, y radius=0.4, start angle=180, end angle=270]
		-| (0.8, 1.6) -- (0,1.6)
		arc[x radius=0.8, y radius=0.4, start angle=270, end angle=180];
	\fill[1facetFront] (-0.15,-0.8) -- (0,-0.4) -- (0,1.6) -- (-0.15,1.2);
	\draw[seam] (0,-0.4) -- (0,1.6);
	\draw[1facetLine] (0.45,2.8) -- (-0.15,1.2) -- (-0.15,-0.8) -- (0,-0.4);
}

\deffoampict{cutting plane cut}(-0.9,-0.9)--(1.2,2.9) M=1 {%
	\begin{scope}
		\clip (0.3, 2.4)
			arc[x radius=0.15, y radius=1.12, start angle=360, end angle=270] --
			(0.15,0.72) arc[x radius=0.15, y radius=0.32, start angle= 90, end angle=0]
			|- (-0.8,0) |- cycle;
		\fill[2facetBack] (-0.8,0)
			.. controls (-0.8,1.0) and ( 0.6,0.5) .. ( 0.6,1)
			.. controls ( 0.6,1.5) and (-0.8,1.0) .. (-0.8,2)
			arc[x radius=0.8, y radius=0.4, start angle=180, end angle=90]
			-| (1.1, 0.4) -- (0,0.4)
			arc[x radius=0.8, y radius=0.4, start angle= 90, end angle=180];
	\end{scope}
	\draw[2facetLine]
		(-0.8,0) arc[x radius=0.8, y radius=0.4, start angle=180, end angle=90]
		-- (0.3,0.4);
	\begin{scope}
		\clip (0,1.6)
			arc[x radius=0.15, y radius=0.32, start angle=180, end angle=270] --
			(0.15,0.72) arc[x radius=0.15, y radius=1.12, start angle=90, end angle=180]
			-| (-0.9,2.4) -| cycle;
		\fill[2facetFront] (-0.8,2)
			arc[x radius=0.8, y radius=0.4, start angle=180, end angle=270]
			-| (0.8,-0.4) -- (0,-0.4)
			arc[x radius=0.8, y radius=0.4, start angle=270, end angle=180]
			.. controls (-0.8,1.0) and ( 0.6,0.5) .. ( 0.6,1)
			.. controls ( 0.6,1.5) and (-0.8,1.0) .. (-0.8,2);
		\draw[2facetLine](0,-0.4)
			arc[x radius=0.8, y radius=0.4, start angle=270, end angle=180]
			.. controls (-0.8,1.0) and ( 0.6,0.5) .. ( 0.6,1)
			.. controls ( 0.6,1.5) and (-0.8,1.0) .. (-0.8,2);
	\end{scope}
	\draw[1facetFront]
		(-0.15,-0.8) -- (0.45, 0.8) -- (0.45, 2.8) -- (-0.15,1.2) -- cycle;
	\begin{scope}
		\clip (0.3, 2.4)
			arc[x radius=0.15, y radius=1.12, start angle=360, end angle=270] --
			(0.15,0.72) arc[x radius=0.15, y radius=0.32, start angle= 90, end angle=0]
			-| (1.1,2.4);
		\fill[2facetBack] (-0.8,0)
			.. controls (-0.8,1.0) and ( 0.6,0.5) .. ( 0.6,1)
			.. controls ( 0.6,1.5) and (-0.8,1.0) .. (-0.8,2)
			-- (-0.8,2.4) -| (1.1, 0.4) -| cycle;
	\end{scope}
	\draw[seam]
		(0.15,1.28) arc[x radius=0.15, y radius=1.12, start angle=270, end angle=360]
		(0.3 ,0.4 ) arc[x radius=0.15, y radius=0.32, start angle=0, end angle=90];
	\draw[2facetLine] (0.3,2.4) -| (1.1,0.4) -- (0.3,0.4);
	\begin{scope}
		\clip (0,1.6)
			arc[x radius=0.15, y radius=0.32, start angle=180, end angle=270] --
			(0.15,0.72) arc[x radius=0.15, y radius=1.12, start angle=90, end angle=180]
			-| (0.8,1.6) -- cycle;
		\fill[2facetFront] (-0.8,2)
			arc[x radius=0.8, y radius=0.4, start angle=180, end angle=270]
			-| (0.8,-0.4) -- (0,-0.4)
			arc[x radius=0.8, y radius=0.4, start angle=270, end angle=180]
			.. controls (-0.8,1.0) and ( 0.6,0.5) .. ( 0.6,1)
			.. controls ( 0.6,1.5) and (-0.8,1.0) .. (-0.8,2);
		\draw[2facetLine] (-0.8,0)
			.. controls (-0.8,1.0) and ( 0.6,0.5) .. ( 0.6,1)
			.. controls ( 0.6,1.5) and (-0.8,1.0) .. (-0.8,2);
	\end{scope}
	\draw[seam]
		(0, 1.6) arc[x radius=0.15, y radius=0.32, start angle=180, end angle=270]
		(0,-0.4) arc[x radius=0.15, y radius=1.12, start angle=180, end angle= 90];
	\draw[2facetLine]
		(0,-0.4) -| (0.8,1.6) -- (0,1.6)
		arc[x radius=0.8, y radius=0.4, start angle=270, end angle=90]
		-- (0.3,2.4);
}

\deffoampict{edge zip}(-0.8,-0.5)--(0.8,0.5) {%
	\coordinate (lend) at (-0.3,0.35);
	\coordinate (rend) at ( 0.3,0.35);
	\draw[2facetInner,draw=seamColor] ([yshift=0.6pt,xshift=-0.2pt]lend)
		.. controls ++([yshift=-0.9pt]0,-0.5) and ++([yshift=-0.9pt]0,-0.5)
		.. ([yshift=0.6pt]rend);
	\draw[2facetLine] ([yshift=0.6pt,xshift=-0.2pt]lend) -- ([yshift=0.6pt]rend);
	\fill[1facetBack] (-0.8, 0.5)
		arc[x radius=0.5, y radius=0.15, start angle=90,  end angle=0]
		.. controls ++(0,-0.5) and ++(0,-0.5) .. (0.3,0.35)
		arc[x radius=0.3, y radius=0.15, start angle=180, end angle=90]
		-- (0.6,-0.2) .. controls (0,-0.35) and (-0.2,-0.35) .. (-0.8,-0.2)
		-- cycle;
	\draw[1facetLine] (0.3,0.35)
		arc[x radius=0.3, y radius=0.15, start angle=180, end angle=90]
		-- (0.6,-0.2)
		.. controls (0,-0.35) and (-0.2,-0.35) .. (-0.8,-0.2) -- (-0.8,0.5)
		arc[x radius=0.5, y radius=0.15, start angle=90,  end angle=0];
	\fill[1facetFront] (0.8, 0.2)
		arc[x radius=0.5, y radius=0.15, start angle=270,  end angle=180]
		.. controls ++(0,-0.5) and ++(0,-0.5) .. (-0.3,0.35)
		arc[x radius=0.3, y radius=0.15, start angle=360, end angle=270]
		-- (-0.6,-0.5) .. controls (0,-0.35) and (0.2,-0.35) .. (0.8,-0.5)
		-- cycle;
	\draw[1facetLine] (-0.3,0.35)
		.. controls ++(0,-0.5) and ++(0,-0.5) .. (0.3,0.35);
	\draw[2facetFront,draw=seamColor] ([yshift=-0.6pt]lend)
		.. controls ++([yshift=-0.9pt]0,-0.5) and ++([yshift=-0.9pt]0,-0.5)
		.. ([yshift=-0.6pt,xshift=0.2pt]rend);
	\draw[2facetLine] ([yshift=-0.6pt]lend) -- ([yshift=-0.6pt,xshift=0.2pt]rend);
	\draw[1facetLine] (-0.3,0.35)
		arc[x radius=0.3, y radius=0.15, start angle=360, end angle=270]
		-- (-0.6,-0.5)  .. controls (0,-0.35) and (0.2,-0.35) .. (0.8,-0.5) -- (0.8,0.2)
		arc[x radius=0.5, y radius=0.15, start angle=270,  end angle=180];
}

\deffoampict{edge unzip}(-0.8,-0.5)--(0.8,0.5) {%
	\coordinate (lend) at (-0.3,-0.35);
	\coordinate (rend) at ( 0.3,-0.35);
	\fill[1facetBack] (-0.8,-0.2)
		arc[x radius=0.5, y radius=0.15, start angle=90,  end angle=0]
		.. controls ++(0, 0.5) and ++(0, 0.5) .. (0.3,-0.35)
		arc[x radius=0.3, y radius=0.15, start angle=180, end angle=90]
		-- (0.6, 0.5) .. controls (0, 0.35) and (-0.2, 0.35) .. (-0.8, 0.5)
		-- cycle;
	\draw[1facetLine] (0.3,-0.35)
		arc[x radius=0.3, y radius=0.15, start angle=180, end angle=90]
		-- (0.6, 0.5) .. controls (0, 0.35) and (-0.2, 0.35) .. (-0.8, 0.5)
		-- (-0.8,-0.2)
		arc[x radius=0.5, y radius=0.15, start angle=90,  end angle=0];
	\draw[2facetInner,draw=seamColor] ([yshift=0.6pt,xshift=-0.2pt]lend)
		.. controls ++([yshift=0.9pt]0,0.5) and ++([yshift=0.9pt]0,0.5)
		.. ([yshift=0.6pt]rend);
	\draw[2facetLine] ([yshift=0.6pt,xshift=-0.2pt]lend) -- ([yshift=0.6pt]rend);
	\draw[1facetFront] (0.8,-0.5)
		arc[x radius=0.5, y radius=0.15, start angle=270,  end angle=180]
		.. controls ++(0, 0.5) and ++(0, 0.5) .. (-0.3,-0.35)
		arc[x radius=0.3, y radius=0.15, start angle=360, end angle=270]
		-- (-0.6, 0.2) .. controls (0, 0.35) and (0.2, 0.35) .. (0.8, 0.2)
		-- cycle;
	\draw[2facetFront,draw=seamColor] ([yshift=-0.6pt]lend)
		.. controls ++([yshift=0.9pt]0,0.5) and ++([yshift=0.9pt]0,0.5)
		.. ([yshift=-0.6pt,xshift=0.2pt]rend);
	\draw[2facetLine] ([yshift=-0.6pt,xshift=0.2pt]rend) -- ([yshift=-0.6pt]lend);
}

\DeclarePictureFamily{JW}

\def\pictdrawcoupon(#1)(#2)#3{%
	\draw[V1,color=black,fill=white] (#1) rectangle (#2)
		node[midway,anchor=mid] {\scalebox{0.75}{$#3$}};
}

\defJWpict{5}(-0.75,0)--(0.75,1) M=0.5 [#1]=[n] {%
	\foreach \x in {-0.5,-0.25,0,0.25,0.5}
		\draw[V1] (\x,0) -- (\x,1);
	\pictdrawcoupon(-0.65,0.3)(0.65,0.7){#1};
}

\defJWpict{3}(-0.75,0)--(0.75,1) M=0.5 [#1]=[n-2] {%
	\foreach \x in {-0.25,0,0.25}
		\draw[V1] (\x,0) -- (\x,1);
	\pictdrawcoupon(-0.4,0.3)(0.4,0.7){#1};
}

\defJWpict{3-narrow}(-0.5,0)--(0.5,1) M=0.5 [#1]=[n-2] {%
	\foreach \x in {-0.25,0,0.25}
		\draw[V1] (\x,0) -- (\x,1);
	\pictdrawcoupon(-0.4,0.3)(0.4,0.7){#1};
}

\defJWpict{4+1}(-0.75,0)--(0.75,1) M=0.5 [#1]=[n-1] {%
	\foreach \x in {-0.5,-0.25,0,0.25,0.5}
		\draw[V1] (\x,0) -- (\x,1);
	\pictdrawcoupon(-0.65,0.3)(0.35,0.7){#1};
}

\defJWpict{4+1-narrow}(-0.75,0)--(0.65,1) M=0.5 [#1]=[n-1] {%
	\foreach \x in {-0.5,-0.25,0,0.25,0.5}
		\draw[V1] (\x,0) -- (\x,1);
	\pictdrawcoupon(-0.65,0.3)(0.35,0.7){#1};
}

\defJWpict{3+cap}(-0.75,0)--(0.75,1) M=0.5 [#1]=[n-1] {%
	\foreach \x in {-0.5,-0.25,0}
		\draw[V1] (\x,0) -- (\x,1);
	\draw[V1] (0.25,0) -- (0.25,0.7)
		arc[x radius=0.125, y radius=0.15, start angle=180, end angle=0]
		-- (0.5,0);
	\pictdrawcoupon(-0.65,0.2)(0.35,0.6){#1};
}

\defJWpict{3+cup}(-0.75,0)--(0.75,1) M=0.5 [#1]=[n-1] {%
	\foreach \x in {-0.5,-0.25,0}
		\draw[V1] (\x,0) -- (\x,1);
	\draw[V1] (0.25,1) -- (0.25,0.3)
		arc[x radius=0.125, y radius=0.15, start angle=180, end angle=360]
		-- (0.5,1);
	\pictdrawcoupon(-0.65,0.4)(0.35,0.8){#1};
}

\defJWpict{5-tail}(-0.75,0)--(0.75,2) M=1 [#1]=[n-1] {%
	\foreach \x in {-0.5,-0.25,0}
		\draw[V1] (\x,0) -- (\x,2);
	\draw[V1] (0.25,0) -- (0.25,0.7)
		arc[x radius=0.125, y radius=0.15, start angle=180, end angle=0]
		-- (0.5,0);
	\draw[V1] (0.25,2) -- (0.25,1.3)
		arc[x radius=0.125, y radius=0.15, start angle=180, end angle=360]
		-- (0.5,2);
	\pictdrawcoupon(-0.65,0.2)(0.35,0.6){#1};
	\pictdrawcoupon(-0.65,1.4)(0.35,1.8){#1};
}

\defJWpict{3+loop}(-0.75,0)--(0.75,1) M=0.5 [#1]=[n-2] {%
	\foreach \x in {-0.5,-0.25,0}
		\draw[V1] (\x,0) -- (\x,1);
	\draw[V1] (0.375, 0.5) ellipse[x radius=0.125, y radius=0.35];
	\pictdrawcoupon(-0.65,0.3)(0.1,0.7){#1};
}

\defJWpict{4-tr}(-0.75,0)--(0.75,1) M=0.5 [#1]=[n-1] {%
	\foreach \x in {-0.5,-0.25,0}
		\draw[V1] (\x,0) -- (\x,1);
	\draw[V1] (0.25,0.5) -- (0.25,0.75)
		arc[x radius=0.125, y radius=0.15, start angle=180, end angle=90]
		arc[x radius=0.125, y radius=0.4, start angle=90, end angle=-90]
		arc[x radius=0.125, y radius=0.15, start angle=-90, end angle=-180]
		-- (0.25,0.5);
	\pictdrawcoupon(-0.65,0.3)(0.35,0.7){#1};
}

\defJWpict{cup}(-0.1,-1.5ex)--(2.9,1) M=0.5 [#1#2]=[in] {%
	\foreach \x in {0,0.4,0.8, 2,2.4,2.8}
		\draw[V1] (\x,0.5ex) -- (\x,1);
	\draw[V1] (1.2,1) -- (1.2,0.9) .. controls ++(0,-0.5) and ++(0,-0.5) .. (1.6,0.9) -- (1.6,1);
	\node[anchor=base] at (0,-1.25ex) {$\scriptstyle 1$};
	\node[anchor=base] at (0.4,-1.25ex) {$\scriptstyle\cdots$};
	\node[anchor=base] at (0.8,-1.25ex) {$\scriptstyle#1$};
	\node[anchor=base] at (2.8,-1.25ex) {$\scriptstyle#2$};
}

\defJWpict{4-head}(-0.75,0)--(0.75,2) M=0.5 [#1]=[n-1] {%not currently used
	\foreach \x in {-0.5,-0.25,0,0.25}
		\draw[V1] (\x,0) -- (\x,2);
	\pictdrawcoupon(-0.65,0.2)(0.15,0.6){#1};
	\pictdrawcoupon(-0.65,1.4)(0.15,1.8){#1};
}

\defJWpict{4-tail}(-0.75,0)--(0.75,2) M=0.5 [#1]=[n-1] {%not currently used
	\foreach \x in {-0.5,-0.25}
		\draw[V1] (\x,0) -- (\x,2);
		\draw[V1] (0,0) -- (0,0.7)
			arc[x radius=0.125, y radius=0.15, start angle=180, end angle=0]
			-- (0.25,0);
		\draw[V1] (0,2) -- (0,1.3)
			arc[x radius=0.125, y radius=0.15, start angle=180, end angle=360]
			-- (0.25,2);
	\pictdrawcoupon(-0.65,0.2)(0.15,0.6){#1};
	\pictdrawcoupon(-0.65,1.4)(0.15,1.8){#1};
}

\defJWpict{XYY}(-0.75,0)--(0.75,2) M=0.5 [#1]=[n-1] {%
	\foreach \x in {-0.5,-0.25,0,0.25}
		\draw[V1] (\x,0) -- (\x,2);
	\pictdrawcoupon(-0.65,1.4)(0.15,1.8){#1};
	\pictdrawcoupon(-0.65,0.2)(-0.1,0.4){};
	\pictdrawcoupon(-0.65,0.8)(-0.1,1){};
}

\defJWpict{XaYY}(-0.75,0)--(0.75,2) M=0.5 [#1]=[n-1] {%
	\foreach \x in {-0.5,-0.25}
		\draw[V1] (\x,0) -- (\x,2);
	\draw[V1] (0,0) -- (0,1)
		arc[x radius=0.125, y radius=0.15, start angle=180, end angle=0]
		-- (0.25,0);
	\draw[V1] (0,2) -- (0,1.4)
		arc[x radius=0.125, y radius=0.15, start angle=180, end angle=360]
		-- (0.25,2);
	\pictdrawcoupon(-0.65,1.4)(0.15,1.8){#1};
	\pictdrawcoupon(-0.65,0.2)(-0.1,0.4){};
	\pictdrawcoupon(-0.65,0.8)(-0.1,1){};
}

\defJWpict{XYbY}(-0.75,0)--(0.75,2) M=0.5 [#1]=[n-1] {%
	\foreach \x in {-0.5,0.25}
		\draw[V1] (\x,0) -- (\x,2);
	\draw[V1] (-0.25,0) -- (-0.25,0.4)
		arc[x radius=0.125, y radius=0.15, start angle=180, end angle=0]
		-- (0,0);
	\draw[V1] (-0.25,2) -- (-0.25,0.8)
		arc[x radius=0.125, y radius=0.15, start angle=180, end angle=360]
		-- (0,2);
	\pictdrawcoupon(-0.65,1.4)(0.15,1.8){#1};
	\pictdrawcoupon(-0.65,0.2)(-0.1,0.4){};
	\pictdrawcoupon(-0.65,0.8)(-0.1,1){};
}

\defJWpict{XaYbY}(-0.75,0)--(0.75,2) M=0.5 [#1]=[n-1] {%
	\draw[V1] (-0.5,0) -- (-0.5,2);
	\draw[V1] (0,2) -- (0,1.4)
		arc[x radius=0.125, y radius=0.15, start angle=180, end angle=360]
		-- (0.25,2);
	\draw[V1] (-0.25,0) -- (-0.25,0.4)
		arc[x radius=0.125, y radius=0.15, start angle=180, end angle=0]
		-- (0,0);
	\draw[V1] (-0.25,2) -- (-0.25,0.8)
		arc[x radius=0.125, y radius=0.15, start angle=180, end angle=360]
		-- (0,1)
		arc[x radius=0.125, y radius=0.15, start angle=180, end angle=0]
		-- (0.25,0);
	\pictdrawcoupon(-0.65,1.4)(0.15,1.8){#1};
	\pictdrawcoupon(-0.65,0.2)(-0.1,0.4){};
	\pictdrawcoupon(-0.65,0.8)(-0.1,1){};
}

\defJWpict{XabY}(-0.6,0)--(0.35,1.8) M=0.5 [#1]=[n-1] {%
	\draw[V1] (-0.5,0) -- (-0.5,1.8);
	\draw[V1] (0,1.8) -- (0,1.15)
		arc[x radius=0.125, y radius=0.15, start angle=180, end angle=360]
		-- (0.25,1.8);
	\draw[V1] (-0.25,0) -- (-0.25,0.45)
		arc[x radius=0.125, y radius=0.15, start angle=180, end angle=0]
		-- (0,0);
	\draw[V1] (-0.25,1.8) -- (-0.25,1.05)
		arc[x radius=0.25, y radius=0.25, start angle=180, end angle=270]
		arc[x radius=0.25, y radius=0.25, start angle=90, end angle=0]
		-- (0.25,0);
	\pictdrawcoupon(-0.65,1.2)(0.15,1.6){#1};
	\pictdrawcoupon(-0.65,0.2)(-0.1,0.4){};
}

\defJWpict{XaY}(-0.6,0)--(0.35,1.8) M=0.5 [#1]=[n-1] {%
	\foreach \x in {-0.5,-0.25}
		\draw[V1] (\x,0) -- (\x,1.8);
	\draw[V1] (0,1.8) -- (0,1.15)
		arc[x radius=0.125, y radius=0.15, start angle=180, end angle=360]
		-- (0.25,1.8);
	\draw[V1] (0,0) -- (0,0.65)
		arc[x radius=0.125, y radius=0.15, start angle=180, end angle=0]
		-- (0.25,0);	
	\pictdrawcoupon(-0.65,1.2)(0.15,1.6){#1};
	\pictdrawcoupon(-0.65,0.2)(-0.1,0.4){};
}

\defJWpict{XY}(-0.6,0)--(0.35,1.8) M=0.5 [#1]=[n-1] {%
	\foreach \x in {-0.5,-0.25,0,0.25}
		\draw[V1] (\x,0) -- (\x,1.8);
	\pictdrawcoupon(-0.65,1.2)(0.15,1.6){#1};
	\pictdrawcoupon(-0.65,0.2)(-0.1,0.4){};
}

\defJWpict{TXX}(-0.75,0)--(0.75,1.7) M=0.85 [#1]=[n-1] {%
	\foreach \x in {-0.5,-0.25,0}
		\draw[V1] (\x,0) -- (\x,1.7);
	\draw[V1] (0.25,0.2) -- (0.25,1.5)
		arc[x radius=0.125, y radius=0.15, start angle=180, end angle=0]
		-- (0.5,0.2)
		arc[x radius=0.125, y radius=0.15, start angle=0, end angle=-180];
	\pictdrawcoupon(-0.65,0.3)(0.35,0.7){#1};
	\pictdrawcoupon(-0.65,1)(0.35,1.4){#1};	
}

\defJWpict{nYu}(-0.75,0)--(0.68,1.4) M=0.7 {
	\foreach \x in {-0.5,-0.25,0}
		\draw[V1] (\x,0) -- (\x,1.4);
	\draw[V1] (0.25,0.7)
		arc[x radius=0.16, y radius=0.4, start angle=-180, end angle=180];
	\pictdrawcoupon(-0.65,0.3)(0.15,0.5){};
	\pictdrawcoupon(-0.65,0.9)(0.15,1.1){};
}

\defJWpict{nYbYu}(-0.75,0)--(0.68,1.4) M=0.7 {
	\foreach \x in {-0.5,-0.25}
		\draw[V1] (\x,0) -- (\x,1.4);
	\draw[V1] (0,0) -- (-0,0.5)
		arc[x radius=0.125, y radius=0.15, start angle=180, end angle=0]
		-- (0.25,0.3)
		arc[x radius=0.125, y radius=0.15, start angle=180, end angle=270]
		arc[x radius=0.16, y radius=0.55, start angle=-90, end angle=90]
		arc[x radius=0.125, y radius=0.15, start angle=90, end angle=180]
		-- (0.25,0.9)
		arc[x radius=0.125, y radius=0.15, start angle=360, end angle=180]
		-- (0,1.4);
	\pictdrawcoupon(-0.65,0.3)(0.15,0.5){};
	\pictdrawcoupon(-0.65,0.9)(0.15,1.1){};
}

\defJWpict{cup-cap*}(-0.5,0.2)--(1.25,1.8) M=1 {%
	\foreach \x in {-0.25,0,0.75,1}
		\draw[V1] (\x,0.4) -- (\x,1.6);
	\draw[V1] (0.25,0.4) -- (0.25,0.7)
		arc[x radius=0.125, y radius=0.15, start angle=180, end angle=0]
		-- (0.5,0.4);
	\draw[V1] (0.25,1.6) -- (0.25,1.3)
		arc[x radius=0.125, y radius=0.15, start angle=180, end angle=360]
		-- (0.5,1.6);
	\node[anchor=base] at (-0.25,0.1) {$\scriptstyle 1$};
	\node[anchor=base] at (0,0.1) {$\scriptstyle \cdots$};
	\node[anchor=base] at (0.25,0.1) {$\scriptstyle i$};
	\node[anchor=base] at (0.625,0.1) {$\scriptstyle \cdots$};
	\node[anchor=base] at (1,0.1) {$\scriptstyle n$};
}

\defJWpict{cap*}(-0.5,0.0)--(1.25,1.6) M=0.8 {%
	\foreach \x in {-0.25,0,0.75,1}
		\draw[V1] (\x,0.4) -- (\x,1.2);
	\draw[V1] (0.25,0.4) -- (0.25,0.7)
		arc[x radius=0.125, y radius=0.15, start angle=180, end angle=0]
		-- (0.5,0.4);
	\node[anchor=base] at (-0.25,0.1) {$\scriptstyle 1$};
	\node[anchor=base] at (0,0.1) {$\scriptstyle \cdots$};
	\node[anchor=base] at (0.25,0.1) {$\scriptstyle i$};
	\node[anchor=base] at (0.625,0.1) {$\scriptstyle \cdots$};
	\node[anchor=base] at (1,0.1) {$\scriptstyle n$};
}

\defJWpict{cup*}(-0.5,0.0)--(1.25,1.6) M=0.8 {%
	\foreach \x in {-0.25,0,0.75,1}
		\draw[V1] (\x,0.4) -- (\x,1.2);
	\draw[V1] (0.25,1.2) -- (0.25,0.9)
		arc[x radius=0.125, y radius=0.15, start angle=180, end angle=360]
		-- (0.5,1.2);
	\node[anchor=base] at (-0.25,1.33) {$\scriptstyle 1$};
	\node[anchor=base] at (0,1.33) {$\scriptstyle \cdots$};
	\node[anchor=base] at (0.25,1.33) {$\scriptstyle i$};
	\node[anchor=base] at (0.625,1.33) {$\scriptstyle \cdots$};
	\node[anchor=base] at (1,1.33) {$\scriptstyle n$};
}

\defJWpict{id*}(-0.5,0.0)--(1.25,1.6) M=0.8 {%
	\foreach \x in {-0.25,0,0.25,0.5,0.75,1}
		\draw[V1] (\x,0.4) -- (\x,1.2);
	\node[anchor=base] at (-0.25,0.1) {$\scriptstyle 1$};
	\node[anchor=base] at (0.375,0.1) {$\scriptstyle \cdots$};
	\node[anchor=base] at (0.375,0.1) {$\scriptstyle \cdots$};
	\node[anchor=base] at (1,0.1) {$\scriptstyle n$};
}

\defJWpict{Te}(-0.75,0)--(0.75,2) M=1 {%
	\draw[V1] (-0.5,0.05) -- (-0.5,1.85);
	\draw[V1] (0.25,1.85) -- (0.25,0.3)
		arc[x radius=0.125, y radius=0.15, start angle=180, end angle=360]
		-- (0.5,1.85);
	\draw[V1] (-0.25,1.85) -- (-0.25,1.1)
		arc[x radius=0.125, y radius=0.15, start angle=180, end angle=360]
		-- (0,1.85);
	\node[anchor=base] at (-0.5,2) {$\scriptstyle +$};
	\node[anchor=base] at (-0.25,2) {$\scriptstyle +$};
	\node[anchor=base] at (0,2) {$\scriptstyle -$};
	\node[anchor=base] at (0.25,2) {$\scriptstyle +$};
	\node[anchor=base] at (0.5,2) {$\scriptstyle -$};
}

\defJWpict{Xe}(-0.75,0)--(0.75,2) M=1 {%
	\draw[V1] (-0.5,0.05) -- (-0.5,1.85);
	\draw[V1] (0.25,1.85) -- (0.25,0.3)
		arc[x radius=0.125, y radius=0.15, start angle=180, end angle=360]
		-- (0.5,1.85);
	\draw[V1] (-0.25,1.85) -- (-0.25,1.1)
		arc[x radius=0.125, y radius=0.15, start angle=180, end angle=360]
		-- (0,1.85);
	\pictdrawcoupon(-0.65,0.4)(0.35,0.8){2};
	\pictdrawcoupon(-0.65,1.2)(-0.15,1.6){2};
	\node[anchor=base] at (-0.5,2) {$\scriptstyle +$};
	\node[anchor=base] at (-0.25,2) {$\scriptstyle +$};
	\node[anchor=base] at (0,2) {$\scriptstyle -$};
	\node[anchor=base] at (0.25,2) {$\scriptstyle +$};
	\node[anchor=base] at (0.5,2) {$\scriptstyle -$};
}

\defJWpict{Ye}(-0.75,0)--(0.75,2) M=1 {%
	\draw[V1] (-0.5,1.95) -- (-0.5,0.15);
	\draw[V1] (0.25,0.15) -- (0.25,1.7)
		arc[x radius=0.125, y radius=0.15, start angle=180, end angle=0]
		-- (0.5,0.15);
	\draw[V1] (-0.25,0.15) -- (-0.25,0.9)
		arc[x radius=0.125, y radius=0.15, start angle=180, end angle=0]
		-- (0,0.15);
	\pictdrawcoupon(-0.65,1.6)(0.35,1.2){2};
	\pictdrawcoupon(-0.65,0.8)(-0.15,0.4){2};
	\node[anchor=base] at (-0.5,-0.1) {$\scriptstyle +$};
	\node[anchor=base] at (-0.25,-0.1) {$\scriptstyle +$};
	\node[anchor=base] at (0,-0.1) {$\scriptstyle -$};
	\node[anchor=base] at (0.25,-0.1) {$\scriptstyle +$};
	\node[anchor=base] at (0.5,-0.1) {$\scriptstyle -$};
}

\newcommand{\XYY}{\JWpict{XYY}}
\newcommand{\XaYY}{\JWpict{XaYY}}
\newcommand{\XYbY}{\JWpict{XYbY}}
\newcommand{\XaYbY}{\JWpict{XaYbY}}
\newcommand{\XY}{\JWpict{XY}}
\newcommand{\XaY}{\JWpict{XaY}}
\newcommand{\XabY}{\JWpict{XabY}}
\newcommand{\TXX}{\JWpict{TXX}}

\newcommand{\nYu}{\JWpict{nYu}}
\newcommand{\nYbYu}{\JWpict{nYbYu}}

\newcommand{\JWT}{\JWpict{3-narrow}}

\usepackage[hmargin=3cm,vmargin=3cm]{geometry}

\usepackage{combelow} % For Negut

\usepackage[latin1,utf8]{inputenc} 
\usepackage[english]{babel}
\usepackage[backend=bibtex,doi=false,url=false,isbn=false,style=alphabetic,maxnames=6]{biblatex}
\usepackage{etex}
\AtBeginBibliography{\small}

\usepackage{amsthm,amsmath,amssymb,amscd,amsfonts}
\usepackage{euscript,wasysym,xspace,stmaryrd,latexsym,youngtab,verbatim,mathrsfs}
\usepackage{palatino}
\usepackage{euscript}
\usepackage{mathtools}

\usepackage{graphicx}
\usepackage[all]{xy}
\SelectTips{cm}{}

\usepackage{tikz, tikz-cd}
\usetikzlibrary{matrix, arrows, calc}
\tikzset{frontline/.style={preaction={draw=white,-,line width=6pt}},}  %%% for 3d commutative diagrams

\usepackage[dvips]{epsfig}
\usepackage{pinlabel}
\usepackage{psfrag}

\usepackage{color}
\definecolor{references}{rgb}{0,0,1}

\usepackage[%pdftex,
colorlinks = true,
urlcolor = references, % \href{...}{...} external (URL)
citecolor = references, % \cite{...}
linkcolor = references, % \ref{...} and \pageref{...}
]
 {hyperref}

\newtheorem{thm}{Theorem}[section]

\newtheorem{lemma}[thm]{Lemma}
\newtheorem{theorem}[thm]{Theorem}

\newtheorem{proposition}[thm]{Proposition}
\newtheorem{cor}[thm]{Corollary}
\newtheorem{corollary}[thm]{Corollary}

\theoremstyle{definition}

\newtheorem{definition}[thm]{Definition}

\newtheorem{example}[thm]{Example}

\newtheorem{remark}[thm]{Remark}

\theoremstyle{remark}

\numberwithin{equation}{section}

\def\AS{{\EuScript A}}

\def\CS{{\EuScript C}}

\def\LC{\mathcal{L}}
\def\TC{\mathcal{T}}
\def\VC{\mathcal{V}}

\def\a{\alpha}
\def\b{\beta}

\def\d{\delta}

\def\e{\varepsilon}

\let\phi=\varphi
\let\tilde=\widetilde

\usepackage{bbm}
\def\C{{\mathbbm C}}

\def\R{{\mathbbm R}}
\def\Z{{\mathbbm Z}}
\def\Q{{\mathbbm Q}}
\def\1{\mathbbm{1}}
\newcommand{\one}{\1}

\newcommand{\tw}{\operatorname{tw}}

\renewcommand{\sl}{\mathfrak{sl}}

\renewcommand{\k}{\mathbbm{k}}
\newcommand{\smMatrix}[1]{\left[\begin{smallmatrix}#1\end{smallmatrix}\right]}

\newcommand{\inv}{^{-1}}
\newcommand{\TLC}{\TC \hskip-2pt\LC}

\newcommand{\Ch}{\operatorname{Ch}}
\newcommand{\Cone}{\operatorname{Cone}}

\newcommand{\qs}{\mathbbm{q}}
\newcommand{\qf}{q}

\newcommand{\ts}{\mathbbm{t}}

\renewcommand{\d}{\delta}
\renewcommand{\b}{\beta}
\renewcommand{\a}{\alpha}
\renewcommand{\top}{\mathrm{top}}
\renewcommand{\bot}{\mathrm{bot}}

\newcommand{\nunlhd}{\, \not \hskip-3pt \unlhd \, }

\usepackage{hhline}
\usepackage{enumitem}

\setlist{noitemsep,topsep=0pt,parsep=0pt,partopsep=0pt}
\setenumerate{label={\arabic*)}}

\makeatletter

\begingroup
	\def\definering#1#2{\protected\gdef#1{\@ifstar{#2^*}{#2}}}%
	\definering\scalars\Bbbk
	\definering\Fld{\mathbb F}
\endgroup

\def\Ann{\mathbb A}

\newcommand{\HH}{\operatorname{HH}}

\usepackage{tikz,  tikz-cd}
\usetikzlibrary{matrix, arrows, calc}
\tikzset{frontline/.style={preaction={draw=white,-,line width=6pt}},}  %%% for 3d commutative diagrams

\theoremstyle{plain}
\newtheorem{itheorem}{Theorem}

\newtheorem{iconjecture}[itheorem]{Conjecture}

\makeatletter
\newcommand{\newreptheorem}[2]{%
	\newtheorem*{rep@#1}{\rep@title}
	\newenvironment{rep#1}[1]{%
		\def\rep@title{#2 \ref*{##1}}%
		\begin{rep@#1}}%
		{\end{rep@#1}}}
\makeatother

\newreptheorem{theorem}{Theorem}

\usepackage{subfiles}

\begin{document}

\title{On unification of colored annular $\sl_2$ knot homology}

\begin{abstract} We show that the Khovanov and Cooper--Krushkal models for colored $\sl_2$ homology are equivalent in the case of the unknot,
when formulated in the
quantum annular
Bar-Natan category.  
Again for the unknot, these two theories are shown to be equivalent to a third colored  homology theory, defined using the action of Jones--Wenzl projectors on the quantum annular homology of cables.
The proof is given by conceptualizing the properties of all three models into a Chebyshev system and by proving its uniqueness.
In addition, we show that the classes of the Cooper--Hogancamp 
projectors in the quantum horizontal trace coincide with those of the
Cooper--Krushkal projectors on the passing through strands. As an application, we compute the full quantum Hochschild homology of Khovanov's arc algebras.
Finally, we state precise conjectures formalizing cabling operations and
extending the above results to all knots.
\end{abstract}

\author{Anna Beliakova}
	\address{Universit\"at Z\"urich, Z\"urich, Switzerland}
	\email{anna@math.uzh.ch}
\author{Matthew Hogancamp}
	\address{Northeastern University, Boston, MA}
	\email{m.hogancamp@northeastern.edu}
\author{Krzysztof Putyra}
	\address{Universit\"at Z\"urich, Z\"urich, Switzerland}
	\email{krzysztof.putyra@math.uzh.ch}
\author{Stephan Wehrli}
	\address{Syracuse University, Syracuse, NY}
	\email{smwehrli@syr.edu}
\maketitle

\setcounter{tocdepth}{2}
\tableofcontents

\renewcommand{\BN}{\ccat{BN}}

\section{Introduction}
\label{s:intro}
\begin{figure}[h]
\includegraphics{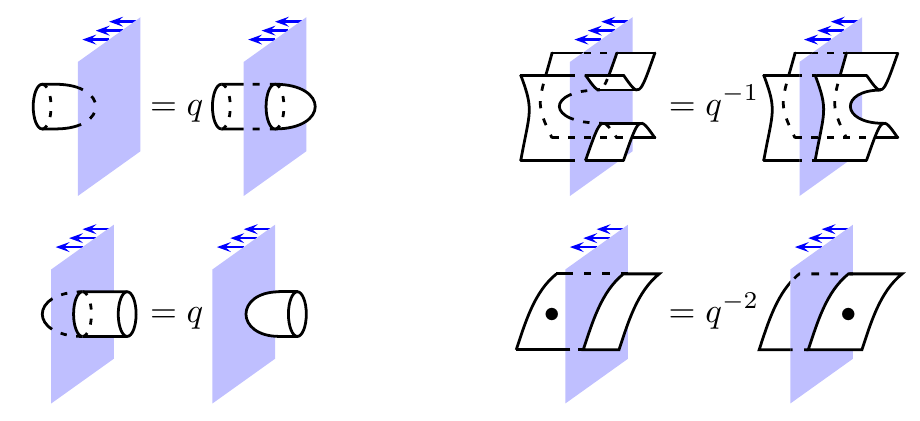}
\caption{Membrane relations.}\label{fig:membrane}
\end{figure}

At the beginning of our century knot theory was revolutionized 
by the Khovanov construction of a chain complex \cite{Kh00} whose graded Euler characteristic is the Jones polynomial, and whose Poincar{\'e} polynomial is a new link invariant detecting the unknot \cite{KronMro11}.
Since then many modifications of this construction were introduced in the literature. The  most relevant  for us is the quantum annular Khovanov homology constructed by  three authors of the present paper in
\cite{BPW19}. This is a triply graded homology theory for
links in a thickened annulus $\Ann$ equipped with a membrane. The differentials
 preserve the quantum and annular gradings and are given by  annular cobordisms  considered up to isotopies fixing the membrane, so that
moving a cobordism $S$ through the membrane
contributes a factor $q^{\pm \chi(S)}$, where $q$ is a 
 quantum parameter and $\chi(S)$ is the Euler characteristic (see Figure~\ref{fig:membrane}).
Such cobordisms   form a category $\cat{BN}_{\!q}(\Ann)$.
Setting $q=1$  we recover 
 the usual annular Khovanov 
homology of Asaeda--Przytycki--Sikora \cite{APS}.

This paper aims to compare different approaches 
extending this construction to colored knots.
Let us recall that given   a framed oriented knot $K\subset S^3$, 
and an $(n+1)$-dimensional representation $V_n$ of 
the quantum group
$U_q(\mathfrak{sl}_2)$, the colored Jones polynomial $J(K, V_n)$
can be defined  by inserting the Jones--Wenzl idempotent into a diagram representing $n$-cable of $K$ and 
by resolving the crossings with the Kauffman bracket skein relations.
Alternatively, the  Reshetikhin--Turaev functor can be used to
define $J(K,V_n)$.

In the categorified  setting there are three different approaches
to define colored annular Khovanov homology, or its
quantization.  To outline these constructions we will need some 
preparation. Let us denote by $\TL$ the
Temperley--Lieb category,
whose objects    are 
 natural numbers and whose morphisms from $n$ to $m$
 are flat $(n,m)$-tangles
(or  crossingless matchings of $n+m$  points).
The unknot considered as $(0,0)$ morphism in $\TL$ takes value $q+q^{-1}$.
For generic $q$
the Temperley--Lieb algebra $\mathrm{TL}_n:=\TL(n,n)$ 
is semisimple  equipped with a family of idempotents $\{p_\e\}$,  $\e\in \{-1,1\}^n$ projecting onto simple components. The most
 famous among them  $p_n$ corresponding to $\e=(1,\dots, 1)$ is called the Jones--Wenzl projector.
There is  a  functor $\HH_0:\TL \to S(\Ann)$, to the skein algebra
of the annulus, sending $p_n$ to its braid closure $[p_n]$ satisfying the following Chebyshev recursion relation:
$$[p_n][p_1]=[p_{n+1}]+[p_{n-1}]$$
This relation also holds for classes  $[V_n]$ in the~Grothendieck 
ring of the~representation category  $\cat{Rep}(U_q(\mathfrak{sl}_2))$ 
or in the Hochschild homology of the latter, thus leading to an isomorphism between $\HH_0(\TL)$  and
$\mathrm{K}_0(\cat{Rep}(U_q(\mathfrak{sl}_2)))$.

The Bar-Natan 2-category $\ccat{BN}$
categorifies $\TL$   by adding as
2-morphisms  surfaces bounded by flat tangles. We will denote 
  by $\cat{BN}_n$ the category $\ccat{BN}(n,n)$.
 There is also a natural categorification of 
 the $0$th Hochschild homology, called the \emph{horizontal trace} of a 2-category.  We will work with the $q$-twisted, or \emph{quantum}, version of the horizontal trace developed  in \cite{BPW19}, where it was shown that the quantum horizontal trace  $\hTr_q(\ccat{BN}^\oplus)$ is equivalent to 
 $\cat{BN}_{\!q}(\Ann)^\oplus$. Here the symbol $\oplus$ indicates the additive closure. 
 
Our first model for colored knot homology 
uses the Cooper--Krushkal categorification
$P_n$ of the Jones--Wenzl idempotent $p_n\in \mathrm{TL}_n$
living in an appropriate completion of the category
of bounded complexes over the Bar-Natan category $\BN_n
$. Let us denote  by $[P_n]$ the class of $P_n$ in 
$\hTr_q(\ccat{BN}^\oplus)$.
In the Coooper--Krushkal model  the colored annular complex 
is constructed by composing $[P_n]$ with the quantum annular complex
for the $n$ cable of $K$.
%. 

The second approach  is due to Khovanov and  inspired by the following decomposition of 
$[V_n]\in \mathrm{K}_0(\cat{Rep}(U_q(\mathfrak{sl}_2)))$
\begin{equation}
\label{eq:Vn}
	[V_n]=\sum_{k=0}^{\lfloor\frac{n}{2}\rfloor}(-1)^k \binom{n-k}{k} \left[V_1^{\otimes(n-2k)}\right].
\end{equation} 
In \cite{Kh05}  Khovanov constructed a family of finite complexes $\VC_n$  over the Temperley--Lieb category $\TLC$ whose Euler characteristics satisfies \eqref{eq:Vn}
where the binomial coefficients are interpreted as  grading shifts.
Combining these complexes with  those for the $n$ cable of $K$
we get the $n$-colored quantum annular complex.

The third construction 
makes use of a  functor
\[
S^1 \times (-): \TL\longrightarrow \cat{BN}_{\!q}(\Ann)
\]
that sends an object $n$ to a collection of $n$ essential circles in $\Ann$ and a cap or cup to a band between neighboring circles. More generally, this functor sends a flat tangle $T$ to the surface $S^1\times T\subset\Ann\times I$. In \cite{BPW19} it is shown that this functor induces an equivalence between 
$\TLC^\oplus$ and $\cat{BN}_{\!q}(\Ann)^\oplus$.
 This allows us to define the colored complex for a knot $K$ by composing 
  $S^1 \times p_n\in \cat{BN}_{\!q}(\Ann)^\oplus$ with the
quantum annular complex of the $n$ cable of $K$. 

To conclude,
all three models  associate with a pair $(K,V_n)$
a (bounded above) chain complex 
over the Karoubi envelope of 
 $\TL^\oplus$. We denote the category of such chain complexes
 by $\Ch^-(\Kar(\TL)^\oplus)$.
Let us also fix a coefficient ring 
$\k$, such that  $q$ is a fixed invertible element of $\k$ and
$1-q^d$ is invertible for all integers $d$ in $(0,2N]$ for a  big enough $N$.

Our main result establishes an equivalence of these
three models for the $n$-colored unknot.

\begin{itheorem}\label{thm:intro mainthm1}
The  quantum annular complexes for the $n$-colored unknot
in all three models are  homotopy equivalent,  meaning that
$$[P_n]\simeq \im p_n\simeq \mathcal{V}_n \in \Ch^-(\Kar(\TL)^\oplus).$$
 where $\im p_n$  denotes the idempotent $p_n$ viewed
as an object of $\Kar(\TL)$.
\end{itheorem}
We remark that the above theorem fails in the setting of the usual (untwisted) horizontal trace, and in general the Cooper--Krushkal  and Khovanov colored homologies are non-isomorphic.

To prove Theorem \ref{thm:intro mainthm1} we formalize
 properties of all three models by introducing a
{\it Chebyshev system} defined for any monoidal category $(\cat C,\otimes, \one)$ by
the following data:
\begin{itemize}
		\item a self-dual object $V$ in $\cat{C}$,
		\item a~family of complexes $V^{(n)}\in \Ch^-(\cat C)$ with $V^{(0)} = \one$ and $V^{(1)}=V$,
		\item chain maps $\pi^{(n)} \colon V^{(n-1)}\otimes V \to V^{(n)}$ for $n\geqslant 1$, 
	\end{itemize}
	such that $\pi^{(1)}\colon \one\otimes V\to V$ is the~canonical isomorphism and there is a~distinguished triangle
	\begin{equation}\label{model:triangle}
		V^{(n-2)} \longrightarrow
		V^{(n-1)}\otimes V \xrightarrow{\ \pi^{(n)}\ }
		V^{(n)} \longrightarrow
		V^{(n-2)}[1]
	\end{equation}
	in the~homotopy category of  bounded complexes over $\cat{C}$ 
	(see
Definition \ref{def:homological-model}	for more details).
We  prove that for a given $V$, a
Chebyshev system in $\cat{C}$ is unique up to homotopy and show that each of our there models defines a 
Chebyshev system in $\Ch^-(\Kar(\TL)^\oplus)$, thus proving the claim.

\subsection{Generalized projectors}
Notice that as an abelian group, the Hochschild homology $\HH_0(\mathrm{TL}_n)$ is free on the classes $[p_n]$.  To express a general endomorphism $x\in \mathrm{TL}_n$ in this basis, 
we use a complete collection of primitive orthogonal idempotents 
$\{p_{\e}\in \mathrm{TL}\}_\e$, 
labelled by certain length $n$ sequences of $\pm1$. In particular,
for any $x\in\mathrm{TL}_n$
we decompose \[
x = \sum_\e x p_\e \quad\text{or}\quad  [x]=\left[\sum_\e p_\e x p_\e\right]  \in \HH_0(\TL).\]
Since each $p_\e$ is a primitive idempotent, $p_\e \mathrm{TL}_n p_\e$ is one-dimensional, hence $ p_\e x p_\e = \Tr_\e(x) p_\e$ for some scalars $\Tr_\e(x)$ and 
\[
[x]  = \sum_\e \Tr_\e(x) [p_\e]. 
\]
Moreover, $[p_\e]$ depends only on the integer $|\e|$ 
(the maximal number of coming through strands) and is
characterized by $\im p_\e \cong \im p_{|\e|}$ or equivalently, by $[p_\e] = [p_{|\e|}]$ in $\HH_0(\mathrm{TL}_n)$, 
where $p_{|\e|}$ is the Jones--Wenzl projector.
   We then have
\[
[x]  =\sum_\e \Tr_\e(x) [p_{|\e|}].
\]
We conjecture that a categorified analogue of this statement holds as well. 

\begin{iconjecture}\label{conj3}
For a given $X\in \Ch^-(\cat{BN}_n)$, there exist complexes $\hTr_{q,\e}(X)\in \Ch^-(\k\mathrm{-mod})$ such that
\[
[X] \ \simeq \ \bigoplus_\e \hTr_{q,\e}(X)\otimes_\k [P_{|\e|}].
\]
where $[P_\e]$ is the class of $P_\e$ in the quantum horizontal trace 
and $\hTr_{q,\e}(X)$ is characterized by
\[
[X\star P_\e] = \hTr_{q,\e}(X)\otimes [P_\e]
\]
where  $P_\e$ is the categorified Temperley--Lieb algebra idempotent constructed by the second author and B. Cooper {\rm \cite{CH12}},
and $\star$ denotes the composition of 1-morphisms in $\BN$.
\end{iconjecture}

We will not pursue this conjecture in this paper, but we will prove the following special case.

\begin{itheorem}\label{thm:intro Pe}
For each $\e$ with $|\e|=k$ we have
\[
[P_\e] \ \simeq \ [P_k] \ \in \Ch^-(\TL^{\oplus}).
\]
\end{itheorem}
  The proof of this theorem and our construction of $P_\e$ is based on homological perturbation theory and 
  a general form of the {\it combing hairs} lemma, that may be of independent
  interest for homotopy theorists.
  
 As an application of Theorem \ref{thm:intro Pe} we show that the higher
quantum Hochschild homology groups of the Khovanov arc algebra $H^n$ vanishes and the rank of the zeroth homology is given by the $n$th Catalan number.

The results and tools developed in this paper form a 
 part of our program aiming to
investigate the behavior of  Khovanov homology 
under cabling operations. Let us state two further conjectures that we plan to address in the sequel.

\subsection{Colored knot homology as a functor}
\label{ss:colors and functors}

Given a knot $K$, a complex semisimple Lie algebra ${\mathfrak{g}}$, and a finite dimensional representation $V$ of $U_q({\mathfrak{g}})$,
the Reshetikhin--Turaev construction produces a link invariant
$P_{\mathfrak{g}}(K,V) \in \Z[q,q^{-1}]$.
The aim of categorified Reshetikhin-Turaev invariants is to replace the Laurent polynomial $P_{\mathfrak{g}}(K,V)$ by a complex of vector spaces $C_{\mathfrak{g}}(K,V)$ whose graded Euler characteristic is $P_{\mathfrak{g}}(K,V)$.  As is well known, in favorable situtations $C_{\mathfrak{g}}(K,V)$ is functorial in $K$, with respect to framed oriented cobordisms.

It is also natural to wonder if such colored knot homology theories can be made functorial in $V$.  That is to say, if we fix $K$, we should expect the assignment $V\mapsto C_{\mathfrak{g}}(K,V)$ to lift to a functor
\[
\cat{Rep}(U_q(\mathfrak{g}))\rightarrow \Ch(\k\mathrm{-mod}).
\]
We plan to investigate the existence of such a functor in the case $\mathfrak{g}=\sl_2$.  In this case, $\cat{Rep}(U_q(\mathfrak{g}))$ is equivalent to the Temperley--Lieb category 
 $\Kar(\TL)^\oplus$.  

 In the sequel our  aim will be to prove the following.

\begin{iconjecture}\label{conj1}
Each framed oriented knot $K$ determines a dg functor 
\[
C(K;-):\Ch^-(\Kar(\TL)^{\oplus})\rightarrow \Ch^-(\k\mathrm{-mod})
\]
with the following properties:
\begin{enumerate}
\item $C(K;V_1^{\otimes n})$ is the quantum annular complex associated to the $n$-cable of $K$.
\item $C(K;V_n)$ is a quantum annular version of Khovanov's colored $\sl_2$-homology.
\end{enumerate}
\end{iconjecture}

\subsection{The space of colors as Hochschild homology}
\label{ss:space of colors}

Let $K\subset S^3$ be a framed oriented knot.  Choose a presentation of $K$ as the closure of a $(1,1)$-tangle diagram $D$ (equipped with the blackboard framing), and consider the following construction.  Suppose we are given an element of the Temperley--Lieb algebra $x\in \mathrm{TL}_n$. 
Let $D^n$ denote the $n$-cable of $D$, and let $\langle D^n\rangle\in \mathrm{TL}_n$ denote the Kauffman bracket.  We can multiply $\langle D^n\rangle x$ and then take the trace
\[
P'(K,x):=\Tr(\langle D^n\rangle x) \in \k.
\]
This assignment does not depend on the choice of $(1,1)$-tangle diagram $D$, and satisfies (by isotopy invariance and the nature of traces)
\[
P'(K,xy) = P'(K,yx).
\]
In fact this identity remains true when $xy\in \mathrm{TL}_n$ and $yx\in \mathrm{TL}_m$ for $m\neq n$.  Thus, we can view the assignment $x\mapsto P'(K,x)$ as a linear map
\[
\HH_0({\TL})\rightarrow \k \, .
\]
Since $\mathrm{TL}_n$ is semisimple, its higher Hochschild homology  groups vanish.  
Replacing the Hochschild homology by its full categorification,
which is a dg version of the horizontal trace defined in \cite{GHW22},
we get the following conjecture.  
\begin{iconjecture}\label{conj2}
Each complex $X\in \Ch^-(\cat{BN}_{n})$ determines a colored complex
\[
C'(K,X)\in \Ch^-(\k\mathrm{-mod}).
\]
The assignment $X\mapsto C'(K,X)$ is functorial in $X$, in the sense that it factors through the dg quantum horizontal trace {\rm \cite{GHW22}}.  Furthermore, if $P_n$ denotes the Cooper--Krushkal categorified Jones--Wenzl idempotent, then $C'(K,P_n)$ is homotopy equivalent to the $V_n$-colored complex $C(K,V_n)$ from Conjecture \ref{conj1}.
\end{iconjecture}

In particular, we expect that the   Khovanov and Cooper--Krushkal models for the
colored  quantum annular homology
coincide for all knots. The results of this paper and the previous one
\cite{BHPW1-pp} will play a crucial role in the proof of
these conjectures.

\subsection{Organization} The paper is organized as follows.
In Section~\ref{sec:preliminaries}, we provide the necessary background on categorical traces and shadows. We also introduce the main categories and bicategories
used in the paper. In Section~\ref{ss:twists},
we discuss
homological perturbation theory for twisted complexes and prove 
a general version of the
combing hairs lemma, needed for our proofs.
In Section~\ref{sec:model-color}, we define the notion of a Chebyshev system in a triangulated category, and
we prove a uniqueness result (Theorem~\ref{thm:uniqueness-of-model}).
As examples of Chebyshev systems, we consider the three models from Theorem~\ref{thm:intro mainthm1}.
In the fourth section, we prove Theorem~\ref{thm:P as model},
which asserts that the Cooper--Krushkal model is indeed a Chebyshev system.
To show this, we establish a new 2-periodic model for the Cooper--Krushkal projector (Theorem~\ref{thm:periodic P}).
The next section is devoted to the categorified idempotents $P_\e$.
We first give a new proof of their existence and then
establish Theorem~\ref{thm:intro Pe} about the quantum horizontal trace of $P_\e$. Finally, we compute the quantum Hochschild homology of the Khovanov arc algebras.

\subsection{Acknowledgements}
AB and KP would like to thank NCCR SwissMAP  
of the Swiss National Science 
Foundation (SNF). In addition, AB was partly supported by the SNF grants 200020$\_$207374 and 200021$\_$178767.
SW was partially supported by a grant from the Simons Foundation ($\#$632059 Stephan Wehrli).

\section{Preliminaries}
\label{sec:preliminaries}

In this section, we will review the homological and categorical constructions
that will be used in the later sections of the paper.
We will also give a self-contained proof of the fact that for generic $q$,
the Temperley--Lieb category is equivalent to the quantum horizontal
trace of the Bar-Natan bicategory (Theorem~\ref{thm:qtrace of BN}).
We will start by fixing some conventions.

\subsection{Conventions}
Throughout this paper,
$\Bbbk$ will be a commutative unital ring and $q\in\Bbbk^\times$ will be a fixed invertible element.
We will  assume
that $1-q^d$ is invertible for all integers $d\in(0,2n]$ and $n$ sufficiently large.

We will write $\star$ for the horizontal composition of 1- and 2-morphisms in a bicategory.
Similarly, we will use $\circ$ to denote the vertical composition of 2-morphisms,
as well as the ordinary composition of morphisms in a category.
Moreover, the symbols $\otimes$ and $\sqcup$ will be used to denote
the monoidal products in the Temperley--Lieb category and in the Bar-Natan
bicategory, respectively. Finally, the superscript $\oplus$ will refer to the additive closure.
This superscript will be omitted in the later sections of the paper, where it is understood
implicitly that categories are completed with respect to finite direct sums.

\subsection{Complexes}

If $\AS$ is a $\k$-linear category, we will let $\Ch(\AS)$ denote category of complexes over $\AS$.  Objects of this category are complexes (with the cohomological convention for differentials)
\[
\cdots \buildrel\d^{k-1}\over\longrightarrow X^k  \buildrel\d^{k}\over\longrightarrow X^{k+1} \buildrel\d^{k+1}\over\longrightarrow \cdots
\]
in which each $X^k \in \AS$.  Morphism spaces between objects $X$ and $Y$ in $\Ch(\AS)$ are complexes
\[
\Hom^k_{\Ch(\AS)}(X,Y)
=  \prod_{i\in \Z} \Hom_\AS(X^i,Y^{i+k})
\]
 equipped with the differential
\[
f \ \mapsto \ \d_Y\circ f - (-1)^{|f|}f\circ \d_X  \ =: \ [\d, f]
\]
where  for $f\in \Hom^k_{\Ch(\AS)}(X,Y)$, $|f|=k$
is  the (cohomological) degree of $f$.
In this way, the category $\Ch(\AS)$ is a \emph{dg category}, i.e.~a category in which hom spaces are complexes of $\k$-modules, and for which composition of morphisms satisfies the appropriate version of the Leibniz rule:
\[
[\d,f\circ g] = [\d,f]\circ g +(-1)^{|f|}f\circ [\d,g].
\]

We use superscripts $-,+,b$ to denote complexes which are bounded from the right, respectively left, respectively both right and left.

The category $\Ch(\AS)$ is symmetric monoidal with the tensor product defined on objects
by
$$(X\otimes Y)^k=\bigoplus_{i+j=k} X^i \otimes Y^j,
\quad \d_{X\otimes Y}=\d_X \otimes \id_Y+\id_X\otimes Y$$ and on morphisms $f,g$ by
$(f \otimes g)(x \otimes y)=(-1)^{|x||g|}f(x)\otimes g(y)$ with braidings
$$c_{X,Y}: X\otimes Y\to Y\otimes X\quad
\text{given by}\quad
 c_{X,Y}(x\otimes y)=(-1)^{|x||y|}y \otimes x .$$

For any $\ell$, the translation  $\ts^\ell$
shifts the complex $X$ by $\ell$ steps to the right, i.e.
$$(\ts^\ell X)^k=X^{k-\ell} \quad\text{and}\quad
\d_{\ts^\ell X}=(-1)^{\ell}\d_X .$$

A morphism $f\in \Hom_{\Ch(\AS)}(X,Y)$ is said to be {\it closed} if $[\d, f]=0$ and {\it null-homotopic} if
$f=[\d,h]$ for some $h\in \Hom_{\Ch(\AS)}(X,Y)$.
In the later case $h$ is called a null-homotopy for $f$.
For $f,g\in \Hom_{\Ch(\AS)}(X,Y)$, we write
$f\simeq g$ and say $f$ and $g$ are {\it homotopic} if $
f-g$ is null-homotopic.

We denote by $\mathcal K(\AS)$ the {\it homotopy} category 
of $\Ch(\AS)$ which has the same objects as  $\Ch(\AS)$ 
and whose morphisms are homotopy classes of degree zero closed
morphisms in $\Ch(\AS)$, i.e.
$$\Hom_{\mathcal K(\AS)}(X,Y):=\frac{\{f\in \Hom^0_{\Ch(\AS)}(X,Y)\, |\,
[\d,f]=0\}}{\d\left(\Hom^{-1}_{\Ch(\AS)}(X,Y)\right)} . $$

An isomorphism $f:X\to Y$
in $\mathcal K(\AS)$ is called a
{\it homotopy equivalence} and in this case we write
$X\simeq Y$. If $X\simeq 0$, we say that $X$ is
{\it contractible}.

If $f\in \Hom^1_{\Ch(\AS)}(X,Y)$ is a degree 1 closed morphism,  the
 {\it cone} 
 of $f$ is the object $\Cone(f):=X\oplus Y$ in $\Ch(\AS)$
equipped with the differential
$\left[\begin{array}{ll}\d_X &0\\ f&\d_Y
\end{array}\right]$. If instead $f:X\to Y$ is a degree $0$ closed morphism, then we replace $f$ by a degree 1
closed morphism $\ts^{-1}X\to Y$ and then apply the previous construction.

The construction of cones allows to translate statements about morphisms
into statements about objects. For example, a closed map
$f:X\to Y$ is a homotopy equivalence if and only if
$\Cone(f)\simeq 0$.

\subsection{Gradings on categories and bicategories}

We will say that a $\k$-linear is {\it pregraded} if its morphism sets are endowed with $\Z$-gradings which are additive under composition. A pregraded category is called {\it graded} if it comes with a ``grading shift'' automorphism $\qs$ along with a natural degree $1$ isomorphism $\mathfrak{s}\colon\one\rightarrow\qs$ from the identity functor to $\qs$.

A bicategory is called {\it $\k$-linear} if its 2-morphism sets are equipped with a $\k$-module structures such that both compositions are bilinear. Such a bicategory is {\it pregraded} if its morphism categories are pregraded in such a way that gradings are additive under horizontal (and vertical) composition. In the non-strict setting, it is also required that the associator and the unitors have degree $0$.

Finally, a pregraded bicategory is called {\it graded} its morphism categories are graded in such a way that the grading shift functors and the natural degree $1$ isomorphisms $\mathfrak{s}\colon\one\rightarrow\qs$ behave nicely under horizontal composition. By this, we mean that $(\qs f)\star g = \qs(f\star g) = f\star (\qs g)$ and correspondingly $\mathfrak{s}\star\one=\mathfrak{s}=\one\star\mathfrak{s}$ for composable 1-morphisms $f$ and $g$.

Note that while the grading shift functors in a graded bicategory take a 1-morphism $f$ to a ``shifted'' 1-morphism $\qs f$, they do not change the degrees of 2-morphisms. Indeed, $|\qs\alpha|=|\mathfrak{s}\circ\alpha\circ\mathfrak{s}^{-1}|=|\alpha|$, where we write $|\alpha|\in\Z$ for the degree of the 2-morphism $\alpha$. In general, any homogeneous 2-morphism $\alpha\colon f\dblto g$ in a graded bicategory corresponds to a homogenous 2-morphism $\mathfrak{s}^m\circ\alpha\circ\mathfrak{s}^{-n}\colon\qs^nf\dblto\qs^m g$ of degree $|\alpha|+m-n$ for $m,n\in\Z$. In particular, the identity 2-morphism $\one\colon f\rightarrow f$ corresponds to $\mathfrak{s}_f\colon f\rightarrow\qs f$.

Given a bicategory which is only pregraded but not already graded, we can extend it to a graded bicategory with the same objects. The 1-morphisms in this graded bicategory are given by pairs of the form $\qs^nf:=(f,n)$, where $f$ is a 1-morphism from the original bicategory, and $n$ is an integer, to be viewed as a formal grading shift. If $\qs^nf=(f,n)$ and $\qs^mg=(g,m)$ are two such pairs, then the 2-morphisms $\qs^nf\dblto\qs^mg$ are given by copies of 2-morphisms $f\dblto g$ from the original bicategory, but with degrees raised by $m-n$. 

In a similar way, a category which is only pregraded but not already graded can be extended to a graded category with objects of the form $\qs^nx:=(x,n)$.

If $\ccat C$ is a pregraded (possibly graded) bicategory, then we denote by $\ccat{C}_0$ the bicategory which contains the same objects and 1-morphisms, but only the 2-morphisms of $\ccat C$ that have degree 0. Likewise, for a pregraded (possibly graded) category $\mathcal{C}$, we denote by $\mathcal{C}_0$ the subcategory which has the same objects but only the morphisms of degree 0.

Note that if $\ccat C$ is graded, then every 2-morphism in $\ccat C$ can be obtained from a 2-morphism in $\ccat C_0$ by pre- and post-composing with suitable powers of $\mathfrak{s}$. A similar remark applies to a graded category $\mathcal{C}$ and its subcategory $\mathcal{C}_0$.

%==============================================
\subsection{Temperley--Lieb category}
\label{ss:tanglecategories}
%==============================================

As before, let $\k$ be a commutative ring and let $q\in\Bbbk^\times$ be a fixed invertible element. The {\it Temperley--Lieb category} is the $\k$-linear category $\TL$ whose objects are nonnegative integers and whose morphisms sets are generated by
isotopy classes of flat tangles modulo the relation
\[
T\cup \bigcirc=(q+q^{-1})T
\]
for any circular component $\bigcirc$. By a {\it flat tangle}, we here mean an unoriented compact
1-manifold $T\subset I^2$ with $n$ bottom endpoints
on $I\times\{0\}$ and $m$ top endpoints on $I\times\{1\}$ for some
$n,m\geq 0$.

The composition in $\TL$ is given by vertical stacking of tangles, so that
$T\circ T'$ denotes the result of placing the tangle $T$ on top of the tangle $T'$ and rescaling vertically.
There is also a horizontal composition of tangles, given by placing tangles side-by-side, which
induces a monoidal structure on $\TL$.
For technical reasons, we will sometimes assume that the objects
in $\TL$ come with formal grading shifts, defined as in~\cite[Section~6]{B-N05}.

The endomorphism algebra $\mathrm{TL}_n:=\TL(n,n)$ is known as the \emph{Temperley--Lieb algebra}.
If the quantum integers $[k]:=q^{k-1}+q^{k-3}+\ldots+q^{1-k}$ are invertible for $1<k\leq n$, then $\mathrm{TL}_n$
contains a distinguished idempotent $p_n$ called the {\it Jones--Wenzl idempotent}. The idempotents $p_n$ can be
defined recursively
by $p_0:=\id_0$, $p_1:=\id_1$, and
\begin{equation}\label{eqn:JWrecursion}
 \JWpict{5}[n]=\JWpict{4+1}[n-1]-\frac{[n-1]}{[n]}\JWpict{5-tail}[n-1]
\end{equation}
for $n>1$, where each box represents a Jones--Wenzl idempotent.

\subsection{Khovanov homology and Bar-Natan's bicategory}
\label{ss:BN cat}

Given an oriented link diagram, we can resolve each of its crossings by applying the Kauffman bracket skein relation
\[%
\begin{tikzpicture}[scale=0.29,baseline = 0,line width = 0.7pt]
	\draw (1,0) -- (0,1);
	\draw[preaction={draw=white,-,line width=4pt}] (0,0) -- (1,1);
\end{tikzpicture}=%
q^{(3\epsilon-1)/2}\bigl(\,
\begin{tikzpicture}[scale=0.29,baseline = 0,line width = 0.7pt]
	\draw (0,0)
	arc[x radius=1.7, y radius=1, start angle=-30, end angle=30];
	\draw (1,1)
	arc[x radius=1.7, y radius=1, start angle=150, end angle=210];
\end{tikzpicture}%
\,-\,q\,%
\begin{tikzpicture}[scale=0.29,baseline = 0,line width = 0.7pt]
	\draw (0,0)
	arc[x radius=1, y radius=1.7, start angle=120, end angle=60];
	\draw (1,1)
	arc[x radius=1, y radius=1.7, start angle=-60, end angle=-120];
\end{tikzpicture}%
\,\bigr)%
\]
where $\epsilon\in\{\pm 1\}$ is the sign of the crossing.
Replacing circles in the resulting terms by factors of $q+q^{-1}$, we
obtain the {\it Jones polynomial} of the link. The same algorithm associates with a tangle a morphism
in $\TL$.

In \cite{Kh00} Khovanov categorified this construction by adding
morphisms between resolutions
in form of saddle cobordisms. Arranging the crossing resolutions
as the vertices of a hypercube and the saddle cobordisms as its edges,
Khovanov assigned to a link diagram a commutative diagram in the category of 2-dimensional cobordisms.
By applying a TQFT functor corresponding to the Frobenius algebra $\k[x]/x^2$, he then defined a chain complex whose homology is a link invariant,
called {\it Khovanov homology}.  The Jones polynomial is recovered as the graded Euler characteristic of Khovanov homology. This construction  extends to tangles in surfaces \cite{APS} and
 is (projectively) functorial with respect to
 tangle cobordisms.

It was observed by Bar-Natan that most of Khovanov's
construction can be performed
formally 
on the level of a~complex $\KhBracket{T}$ called the~\emph{formal Khovanov
bracket} of $T\subset F$ \cite{B-N05}. 
Here $F$ can be any  smooth oriented surface with boundary and it is assumed that $\partial T=B$ for a set $B\subset\partial F$.
The complex $\KhBracket{T}$
is constructed in the additive closure of the {\it Bar-Natan category} $\cat{BN}(F,B)$, the graded $\k$-linear category whose objects 
are flat tangles in $F$ with boundary $B$ and with grading shifts, and whose morphism sets are generated by 
compact surfaces in $F\times I$ decorated by dots.
Such surfaces are required to have boundaries of the form $\partial S=(T\times\{0\})\cup(B\times I)\cup (T'\times\{1\})$ and
are considered up to isotopy and up to the following local relations:
\begin{itemize}
	\item \textit{sphere evaluations:}
	\begin{equation}\label{rel:sphere-eval0}
	\tikzset{x=8mm,y=8mm}%
		\foampict{1sphere}    = 0\hskip 0.1\textwidth
		\foampict{1sphere}[1] = 1
	\end{equation}
	
	\item \textit{neck cutting relation:}
	\begin{equation}\label{rel:neck-cutting0}
	\tikzset{x=7mm,y=7mm}%
		\mathclap{\foampict{1neck}\, = \,\foampict{1cap 1cup}[+1]
		                          \, + \,\foampict{1cap 1cup}[-1]
			                       %\, - h \foampict{1cap 1cup}
			                       }
	\end{equation}
	
	\item \textit{dot reduction:}
	\begin{equation}\label{rel:dots0}
	\tikzset{x=1cm,y=1cm}%
		\foampict{1plane}[2] = 0 % h\!\foampict{1plane}[1] + t\!\foampict{1plane}
	\end{equation}
\end{itemize}
%Here $h$ and $t$ are fixed elements of the~ring of scalars $\scalars$.
We define the (unshifted) degree of a cobordism $S\subset F\times I$ by
\[
|S|:=-\chi(S)+\frac{\#\{\it corners \,in\, S\}}{2}+2\#\{\it dots\,on\, S\}
\]
Note that the relations are degree preserving, and that
 the~neck cutting relation evaluates a~handle attached to
a~plane as a~dot scaled by 2. Because of that, it is common to think of a~dot
as ,,half'' of a~handle, even when 2 is not an~invertible scalar.

We will denote by $\cat{BN}(\Ann):=\cat{BN}(\Ann,\emptyset)$
the Bar-Natan category of the annulus. Moreover, we will denote by $\BN$ the graded bicategory whose objects correspond to nonnegative integers, and whose morphism categories are given by the Bar-Natan categories $\cat{BN}(I^2,B_{n,m})$ where $B_{n,m}\subset\partial I^2$ consists of $n$ points on the bottom boundary of the square $I^2$ and of $m$ points on the top boundary. Note that the horizontal composition in $\BN$ is induced by the composition of tangles, and that $\BN$ has a monoidal structure corresponding to the monoidal structure on the category of flat tangles.

The~formal Khovanov bracket is projectively functorial \cite{B-N05}. Indeed,
there is a~way of associating a chain map with each Reidemeister move
and with any cobordism containing a~unique critical point. One constructs
a chain map for any smooth tangle cobordism by decomposing the~cobordism
into a~sequence of the~above elementary pieces and composing the~associated maps;
choosing a~different decomposition may at most changes the map within its homotopy class, or changes the~global sign of the~map. This global sign can be fixed 
by using the explicit equivalence between $\Uqgl$-foams 
and $\ccat{BN}$ constructed in \cite{BHPW1-pp}.
Hence, 
there is a well-defined bifunctor
\[
\KhBracket{-} \colon \ccat{Tan} \to \mathcal{K}(\ccat{BN}^\oplus),
\]
where $\ccat{Tan}$ denotes the bicategory whose 1-morphisms are oriented tangle diagrams in $I^2$ (possibly with crossings) and whose 2-morphisms are smooth tangle cobordisms in $\mathbb{R}^4$, up to isotopy.

\subsection{Quantum horizontal trace}
\label{ss:qtrace}
There is a functorial construction, called the twisted horizontal trace, sending
suitable tangles in $F\times I$
to links in a surface bundle over a circle with fiber $F$ and monodromy $\sigma \in \Diff(F)$. This construction works more generally for any bicategory $\ccat{C}$ equipped with an endofunctor $\Sigma\colon\ccat{C}\rightarrow\ccat{C}$. Let us recall its definition from \cite{BPW19}.

\begin{definition} The~\emph{$\endofun$-twisted horizontal trace} of $\ccat C$
is the~category $\hTr(\ccat C, {\Sigma})$, whose objects are 1-morphisms $f\colon  { \Sigma}x\to x$ in $\ccat C$ and morphisms from $f\in \ccat C( { \Sigma}x, x)$ to $g\in \ccat C( { \Sigma}y, y)$ are equivalence classes $[p, \alpha]$ of squares
	\begin{equation*}\label{diag:morphism-in-shTr}
		\begin{diagram}
			\Node (X1) at (0,1) {\Sigma x}; \Node (X2) at (1,1) {x};
			\Node (Y1) at (0,0) {\Sigma y}; \Node (Y2) at (1,0) {y};
			\Arrow[->] (X1) -- (X2) node[midway] {$f$};
			\Arrow[->] (Y1) -- (Y2) node[midway,below] {$g$};
			\Arrow[->] (X1) -- (Y1) node[left,midway] {$\Sigma p$};
			\Arrow[->] (X2) -- (Y2) node[midway] {$p$};
			\Arrow[double equal sign distance, double, -{Implies[]}]
				($(X2)!1.5em!(Y1)$) -- ($(Y1)!1.5em!(X2)$)
				node[midway,anchor=-45]{$\alpha$};
		\end{diagram}
	\end{equation*}
	where $\alpha\colon  p\star f\dblto g\star { \Sigma}p$ is a~2-morphism in $\ccat C$, modulo the~relation
	\begin{equation}\label{rel:shTr-for-morphisms}
		\begin{diagram}[baseline={2.1em-height("$\vcenter{}$")}]  % in diagrams.sty, the default unit for the y-coordinate seems to be 4.2em
			\Node (X1) at (0,1) {\Sigma x}; \Node (X2) at (1,1) {x};
			\Node (Y1) at (0,0) {\Sigma y}; \Node (Y2) at (1,0) {y};
			\Arrow[->] (X1) -- (X2) node[midway] {$f$};
			\Arrow[->] (Y1) -- (Y2) node[midway,below] {$g$};
			\Arrow[->] (X1) .. controls ++(235:0.5) and ++(125:0.5) .. (Y1) node[left,midway] {$\Sigma p$};
			\Arrow[->] (X1) .. controls ++(-55:0.5) and ++( 55:0.5) .. (Y1) node[right,pos=0.25] {$\Sigma p'$};
			\Arrow[->] (X2) -- (Y2) node[midway] {$p$};
			\Arrow[double equal sign distance, double, -{Implies[]}]
				($(X2)!1.5em!(0.1,0)$) -- ($(0.1,0)!1.5em!(X2)$)
				node[midway,anchor=135]{$\alpha$};
			\Arrow[double equal sign distance, double, -{Implies[]}]
				(0.15,0.5) -- (-0.15,0.5)
				node[midway,above]{$\Sigma \tau$};
		\end{diagram}
		\quad\sim\quad
		\begin{diagram}[baseline={2.1em-height("$\vcenter{}$")}]
			\Node (X1) at (0,1) {\Sigma x}; \Node (X2) at (1,1) {x};
			\Node (Y1) at (0,0) {\Sigma y}; \Node (Y2) at (1,0) {y};
			\Arrow[->] (X1) -- (X2) node[midway] {$f$};
			\Arrow[->] (Y1) -- (Y2) node[midway,below] {$g$};
			\Arrow[->] (X2) .. controls ++(235:0.5) and ++(125:0.5) .. (Y2) node[left,pos=0.7] {$p$};
			\Arrow[->] (X2) .. controls ++(-55:0.5) and ++( 55:0.5) .. (Y2) node[right,midway] {$p'$};
			\Arrow[->] (X1) -- (Y1) node[midway,left] {$\Sigma p'$};
			\Arrow[double equal sign distance, double, -{Implies[]}]
				($(0.9,1)!1.5em!(Y1)$) -- ($(Y1)!1.5em!(0.9,1)$)
				node[midway,anchor=-45]{$\alpha$};
			\Arrow[double equal sign distance, double, -{Implies[]}]
				(1.15,0.5) -- (0.85,0.5)
				node[midway,above]{$\tau$};
		\end{diagram}
	\end{equation}
	for 1-morphisms $p,p'\colon x\to y$, and 2-morphisms $\alpha\colon p\star f\dblto g\star \Sigma p'$ and $\tau\colon p'\dblto p$. 
	Here $[\one_x, \one_f]$ the~identity on $f$ and the~composition
	\begin{equation}
		[q, \beta]\circ [p, \alpha] := [
			q\star p,
			(\beta \star \one_{\endofun p}) \circ
			(\one_q \star \alpha)
		]
	\end{equation}
	can be visualized as stacking squares one on top of the~other:
	\begin{equation}
		\begin{diagram}[baseline={2.1em-height("$\vcenter{}$")}]
			\Node (X1) at (0,1) {\endofun y}; \Node (X2) at (1,1) {y};
			\Node (Y1) at (0,0) {\endofun z}; \Node (Y2) at (1,0) {z};
			\Arrow[->] (X1) -- (X2) node[midway] {$g$};
			\Arrow[->] (Y1) -- (Y2) node[midway,below] {$h$};
			\Arrow[->] (X1) -- (Y1) node[left,midway] {$\endofun q$};
			\Arrow[->] (X2) -- (Y2) node[midway] {$q$};
			\Arrow[double equal sign distance, double, -{Implies[]}]
				($(X2)!1.5em!(Y1)$) -- ($(Y1)!1.5em!(X2)$)
				node[midway,anchor=-45]{$\beta$};
		\end{diagram}
	\quad\circ\quad
		\begin{diagram}[baseline={2.1em-height("$\vcenter{}$")}]
			\Node (X1) at (0,1) {\endofun x}; \Node (X2) at (1,1) {x};
			\Node (Y1) at (0,0) {\endofun y}; \Node (Y2) at (1,0) {y};
			\Arrow[->] (X1) -- (X2) node[midway] {$f$};
			\Arrow[->] (Y1) -- (Y2) node[midway,below] {$g$};
			\Arrow[->] (X1) -- (Y1) node[left,midway] {$\endofun p$};
			\Arrow[->] (X2) -- (Y2) node[midway] {$p$};
			\Arrow[double equal sign distance, double, -{Implies[]}]
				($(X2)!1.5em!(Y1)$) -- ($(Y1)!1.5em!(X2)$)
				node[midway,anchor=-45]{$\alpha$};
		\end{diagram}
	\qquad:=\qquad
		\begin{diagram}[y=3em,baseline={3em-height("$\vcenter{}$")}]
			\Node (X1) at (0,2) {\endofun x}; \Node (X2) at (1,2) {x};
			\Node (Y1) at (0,1) {\endofun y}; \Node (Y2) at (1,1) {y};
			\Node (Z1) at (0,0) {\endofun z}; \Node (Z2) at (1,0) {z};
			\Arrow[->] (X1) -- (X2) node[midway] {$f$};
			\Arrow[->] (Y1) -- (Y2) node[midway] {$g$};
			\Arrow[->] (Z1) -- (Z2) node[midway,below] {$h$};
			\Arrow[->] (X1) -- (Y1) node[left,midway] {$\endofun p$};
			\Arrow[->] (X2) -- (Y2) node[midway] {$p$};
			\Arrow[->] (Y1) -- (Z1) node[left,midway] {$\endofun q$};
			\Arrow[->] (Y2) -- (Z2) node[midway] {$q$};
			\Arrow[double equal sign distance, double, -{Implies[]}]
				($(X2)!1.5em!(Y1)$) -- ($(Y1)!1.5em!(X2)$)
				node[midway,anchor=-45]{$\alpha$};
			\Arrow[double equal sign distance, double, -{Implies[]}]
				($(Y2)!1.5em!(Z1)$) -- ($(Z1)!1.5em!(Y2)$)
				node[midway,anchor=-45]{$\beta$};
		\end{diagram}
	\end{equation}
	Unitarity and associativity follows from \eqref{rel:shTr-for-morphisms}
	with an~appropriate composition of associators and unitors as $\tau$.
	\end{definition}

When $\ccat C$ is $\k$-linear and $\Sigma\colon\ccat{C}\rightarrow\ccat{C}$ is given by linear maps on 2-morphism sets, we modify the definition of the twisted horizontal trace by allowing formal $\k$-linear combinations of squares as morphisms. Such linear combinations are taken modulo relation \eqref{rel:shTr-for-morphisms} and modulo the relation $[p,a\alpha+\beta]=a[p,\alpha]+[p,\beta]$ for $a\in\k$. This definition ensures that $\hTr(\ccat C,\Sigma)$ is itself a $\k$-linear category.

If $\Sigma=\one$, we recover the horizontal trace from \cite{BHLZ}, and in this case 
we will omit $\Sigma$ from the notation.
For example,
\[\hTr(\ccat{BN})\simeq \cat{BN}(\Ann),\]
where the equivalence sends a flat tangle $T\subset I^2$ to its annular closure $\widehat{T}\subset\Ann$.

In case $\ccat C$ is pregraded, we can define
the {\it quantum endofunctor}, which acts as the identity on objects and 1-morphisms and which satisfies $\Sigma \alpha=q^{-|\alpha|} \alpha$
for any homogeneous 2-morphism $\alpha$ and for a fixed invertible element $q\in\k^\times$. The horizontal trace twisted by this endofunctor is called
{\it quantum horizontal trace} and will be denoted by $\hTr_q(\ccat C)$.

We will denote the quantum horizontal trace of $\BN$ by
\[
\cat{BN}_{\!q}(\Ann):=\hTr_q(\BN)
\]
and call it the {\it quantum annular Bar-Natan category}. This category admits the following graphical description. Choose a radial arc $\mu\subset\Ann$ connecting the two components of $\partial\Ann$, and let $M\subset\Ann\times I$ be the {\it membrane} $M:=\mu\times I$. The category $\cat{BN}_{\!q}(\Ann)$ can be defined in the same way as $\cat{BN}(\Ann)$, except that objects are required to be transverse to $\mu$ and cobordisms in $\Ann\times I$ are required to be transverse to $M$. Isotopic cobordisms are identified if the ambient isotopy fixes the membrane setwise. Otherwise, we rescale the target cobordism as shown in Figure~\ref{fig:membrane}, where the coorientation of the membrane is induced by the orientation of the core of the annulus.

There is a natural functor
\[
S^1 \times (-): \TL\longrightarrow \cat{BN}_{\!q}(\Ann)
\]
that sends an object $n$ to a collection of $n$ essential circles in $\Ann$ and a cap or cup to a band between neighboring circles. More generally, this functor sends a flat tangle $T$ to the surface $S^1\times T\subset\Ann\times I$, which corresponds to the morphism $[T,\one_T]$ in the quantum horizontal trace.

The following result first appeared in \cite{BPW19}.

\begin{theorem}\label{thm:qtrace of BN}
Suppose $1-q^d$ is invertible for all $d\neq 0$. Then the functor $S^1\times (-)$ induces an equivalence of categories 
\[\TL^\oplus\simeq\hTr_q(\ccat{BN})^\oplus\cong\cat{BN}_{\!q}(\Ann)^\oplus,\] where objects in the involved categories are assumed to
come with formal grading shifts.
\end{theorem}

A proof of this proposition will be given in Subsection~\ref{ss:vertical-trace}. Note that the assumption on $q$ implies that in $\cat{BN}_{\!q}(\Ann)$, dots on essential components of a cobordism $S\subset\Ann\times I$ are equal to zero. Indeed, moving such a dot through the membrane, we see that $(1-q^2)$ times the dot has to be equal to zero. This is a special case of the following result, in which $\alpha\colon p\dblto p$ in a 2-endomorphism in a pregraded (and possibly graded) bicategory $\ccat{C}$:

\begin{lemma}\label{lem:annihilate} Suppose $1-q^d$ is invertible in $\k$. If $|\alpha|=d$, then the square
	\begin{equation*}
		\begin{diagram}
			\Node (X1) at (0,1) {x}; \Node (X2) at (1,1) {x};
			\Node (Y1) at (0,0) {y}; \Node (Y2) at (1,0) {y};
			\Arrow[->] (X1) -- (X2) node[midway] {$\one_x$};
			\Arrow[->] (Y1) -- (Y2) node[midway,below] {$\one_y$};
			\Arrow[->] (X1) -- (Y1) node[left,midway] {$p$};
			\Arrow[->] (X2) -- (Y2) node[midway] {$p$};
			\Arrow[double equal sign distance, double, -{Implies[]}]
				($(X2)!1.5em!(Y1)$) -- ($(Y1)!1.5em!(X2)$)
				node[midway,anchor=-45]{$\alpha$};
		\end{diagram}
	\end{equation*}
is zero in $\hTr_q(\ccat C)$.
\end{lemma}

\begin{proof} We have
\begin{equation*}
	\begin{diagram}[baseline={2.1em-height("$\vcenter{}$")}]
			\Node (X1) at (0,1) {x}; \Node (X2) at (1,1) {x};
			\Node (Y1) at (0,0) {y}; \Node (Y2) at (1,0) {y};
			\Arrow[->] (X1) -- (X2) node[midway] {$\one_x$};
			\Arrow[->] (Y1) -- (Y2) node[midway,below] {$\one_y$};
			\Arrow[->] (X1) -- (Y1) node[midway,left] {$p$};
			\Arrow[->] (X2) -- (Y2) node[midway]{$p$};
			\Arrow[double equal sign distance, double, -{Implies[]}]
				($(0.9,1)!1.5em!(Y1)$) -- ($(Y1)!1.5em!(0.9,1)$)
				node[midway,anchor=-45]{$\alpha$};
	\end{diagram}
	\,\,=\,\,
	\begin{diagram}[baseline={2.1em-height("$\vcenter{}$")}]  % in diagrams.sty, the default unit for the y-coordinate seems to be 4.2em
			\Node (X1) at (0,1) {x}; \Node (X2) at (1,1) {x};
			\Node (Y1) at (0,0) {y}; \Node (Y2) at (1,0) {y};
			\Arrow[->] (X1) -- (X2) node[midway] {$\one_x$};
			\Arrow[->] (Y1) -- (Y2) node[midway,below] {$\one_y$};
			\Arrow[->] (X1) .. controls ++(235:0.5) and ++(125:0.5) .. (Y1) node[left,midway] {$p$};
			\Arrow[->] (X1) .. controls ++(-55:0.5) and ++( 55:0.5) .. (Y1) node[right,pos=0.25] {$p$};
			\Arrow[->] (X2) -- (Y2) node[midway] {$p$};
			\Arrow[double equal sign distance, double, -{Implies[]}]
				($(X2)!1.5em!(0.1,0)$) -- ($(0.1,0)!1.5em!(X2)$)
				node[midway,anchor=135]{$\one_p$};
			\Arrow[double equal sign distance, double, -{Implies[]}]
				(0.15,0.5) -- (-0.15,0.5)
				node[midway,above]{$\alpha$};
	\end{diagram}
	\,\,\sim\,\,q^{d}\cdot
	\begin{diagram}[baseline={2.1em-height("$\vcenter{}$")}]
			\Node (X1) at (0,1) {x}; \Node (X2) at (1,1) {x};
			\Node (Y1) at (0,0) {y}; \Node (Y2) at (1,0) {y};
			\Arrow[->] (X1) -- (X2) node[midway] {$\one_x$};
			\Arrow[->] (Y1) -- (Y2) node[midway,below] {$\one_y$};
			\Arrow[->] (X2) .. controls ++(235:0.5) and ++(125:0.5) .. (Y2) node[left,pos=0.7] {$p$};
			\Arrow[->] (X2) .. controls ++(-55:0.5) and ++( 55:0.5) .. (Y2) node[right,midway] {$p$};
			\Arrow[->] (X1) -- (Y1) node[midway,left] {$p$};
			\Arrow[double equal sign distance, double, -{Implies[]}]
				($(0.9,1)!1.5em!(Y1)$) -- ($(Y1)!1.5em!(0.9,1)$)
				node[midway,anchor=-45]{$\one_p$};
			\Arrow[double equal sign distance, double, -{Implies[]}]
				(1.15,0.5) -- (0.85,0.5)
				node[midway,above]{$\alpha$};
	\end{diagram}
	\,\,=\,\,q^{d}\cdot
	\begin{diagram}[baseline={2.1em-height("$\vcenter{}$")}]
			\Node (X1) at (0,1) {x}; \Node (X2) at (1,1) {x};
			\Node (Y1) at (0,0) {y}; \Node (Y2) at (1,0) {y};
			\Arrow[->] (X1) -- (X2) node[midway] {$\one_x$};
			\Arrow[->] (Y1) -- (Y2) node[midway,below] {$\one_y$};
			\Arrow[->] (X1) -- (Y1) node[midway,left] {$p$};
			\Arrow[->] (X2) -- (Y2) node[midway]{$p$};
			\Arrow[double equal sign distance, double, -{Implies[]}]
				($(0.9,1)!1.5em!(Y1)$) -- ($(Y1)!1.5em!(0.9,1)$)
				node[midway,anchor=-45]{$\alpha$};
	\end{diagram}
\end{equation*}
which implies that the square from the lemma is annihilated by $1-q^d$. The result follows because $1-q^d$ is invertible in $\k$
\end{proof}

In the graded setting, the lemma admits a generalization in which $\one_x$ and $\one_y$ are replaced by shifted copies $\qs^n\one_x$ and $\qs^m\one_y$. For the lemma to remain true, $d$ must be the degree of $\alpha$ viewed as a 2-morphism between the {\it unshifted} 1-morphisms $p\star\one_x=p$ and $\one_y\star p=p$. To see this, consider a 2-morphism $\alpha\colon p\star (\qs^n\one_x)\dblto(\qs^m\one_y)\star p$, and apply the original version of the lemma to the 2-morphism $\mathfrak{s}^{-m}\circ\alpha\circ \mathfrak{s}^n\colon p\star\one_x\dblto\one_y\star p$.

In Subsection~\ref{ss:vertical-trace}, we will reinterpret the lemma in terms of the vertical quantum trace.

\begin{remark}
An argument similar to the one used in the proof of the lemma shows that if the commutative ring $\k$ is equipped with a nontrivial grading and if $a\in\k$ has degree $|a|=d$, then $(1-q^d)a$ annihilates all morphisms in the quantum horizontal trace.
\end{remark}

\subsection{Cyclicity isomorphism} For a bicategory with left and right duals, the twisted horizontal trace can be seen as the target category of a {\it universal twisted shadow} (see \cite{BPW19}). In particular, this means that $\hTr(\ccat C,\Sigma)$ comes with the following extra structure: for each object $x\in\ccat C$, there is a functor $[-]_x\colon\ccat C(\Sigma x,x)\to\hTr(\ccat C,\Sigma)$
given by
\[
[f]_x:=f,\qquad [\alpha]_x:=[\one_x,\alpha].
\]
Likewise, for each pair of 1-morphisms $\Sigma x\buildrel f\over\to y\buildrel g\over\to x$ in $\ccat C$, there is a {\it cyclicity morphism}
$\theta_{g,f}\colon[f\star\Sigma g]_y\rightarrow[g\star f]_x$
given by an identity 2-morphism (or by an associator in the non- strict setting):
\begin{equation*}
		\theta_{g,f}:=\quad
		\begin{diagram}[baseline={2.1em-height("$\vcenter{}$")}]
			\Node (X1) at (0,1) {\Sigma y}; \Node (X2) at (1,1) {\Sigma x}; \Node (X3) at (2,1) {y};
			\Node (Y1) at (0,0) {\Sigma x}; \Node (Y2) at (1,0) {y}; \Node (Y3) at (2,0) {x};
			\Arrow[->] (X1) -- (X2) node[midway] {$\Sigma g$};
			\Arrow[->] (Y1) -- (Y2) node[midway,below] {$f$};
			\Arrow[->] (X2) -- (X3) node[midway] {$f$};
			\Arrow[->] (Y2) -- (Y3) node[midway,below] {$g$};
			\Arrow[->] (X1) -- (Y1) node[left,midway] {$\Sigma g$};
			\Arrow[->] (X3) -- (Y3) node[midway] {$g$};
			\Arrow[->] (X2) -- (Y1) node[midway,anchor=center,fill=white] {$\one_{\Sigma x}$};
			\Arrow[->] (X3) -- (Y2) node[midway,anchor=center,fill=white] {$\one_{y}$};
		\end{diagram}
\end{equation*}
When $\ccat{C}$ has appropriate duals (which is the case for the bicategories considered in this paper), then the morphism $\theta_{g,f}$ turns out to be an isomorphism. Moreover, $\theta_{g,f}$ is natural in the sense that the following diagram commutes
\begin{equation*}
		\begin{diagram}[x=7em,y=5em]
			\Node (X1) at (0,1) {[f\star\Sigma g]_y}; \Node (X2) at (1,1) {[g\star f]_x};
			\Node (Y1) at (0,0) {[f'\star\Sigma g']_y}; \Node (Y2) at (1,0) {[g'\star f']_x};
			\Arrow[->] (X1) -- (X2) node[midway] {$\theta_{g,f}$};
			\Arrow[->] (Y1) -- (Y2) node[midway,below] {$\theta_{g',f'}$};
			\Arrow[->] (X1) -- (Y1) node[left,midway] {$[\alpha\star\Sigma\beta]_y$};
			\Arrow[->] (X2) -- (Y2) node[midway] {$[\beta\star\alpha]_x$};
		\end{diagram}
\end{equation*}
for 2-morphisms $\alpha\colon f\dblto f'$ and $\beta\colon g\dblto g'$.

When $\Sigma$ is the quantum endofunctor of a pregraded bicategory $\ccat{C}$, we drop the subscript $[f]_x$ referring to the object $x$ and instead write $[f]$ for $[f]_x\in \hTr_q(\ccat{C})=\hTr(\ccat{C},\Sigma)$. In this case, the naturality property of $\theta$ can be stated as saying that
\begin{equation*}
		\begin{diagram}[x=7em,y=5em]
			\Node (X1) at (0,1) {[f\star g]}; \Node (X2) at (1,1) {[g\star f]};
			\Node (Y1) at (0,0) {[f'\star g']}; \Node (Y2) at (1,0) {[g'\star f']};
			\Arrow[->] (X1) -- (X2) node[midway] {$\theta_{g,f}$};
			\Arrow[->] (Y1) -- (Y2) node[midway,below] {$\theta_{g',f'}$};
			\Arrow[->] (X1) -- (Y1) node[left,midway] {$[\alpha\star\beta]$};
			\Arrow[->] (X2) -- (Y2) node[midway] {$q^{|\beta|}[\beta\star\alpha]$};
		\end{diagram}
\end{equation*}
commutes in $\hTr_q(\ccat C)$ for any 2-morphism $\alpha\colon f\dblto f'$ and any homogeneous 2-morphism $\beta\colon g\dblto g'$. More generally, we have a similar commutative diagram in which $[\alpha\star\beta]$ is replaced by $[\alpha_1\star\beta_1]+\ldots+[\alpha_n\star\beta_n]$ and $q^{|\beta|}[\beta\star\alpha]$ is replaced by
\[
q^{|\beta_1|}[\beta_1\star\alpha_1]+\ldots+q^{|\beta_n|}[\beta_n\star\alpha_n]
\]
for any 2-morphisms $\alpha_i\colon f\dblto f'$ and homogeneous 2-morphisms $\beta_i\colon g\dblto g'$.

The cyclicity morphism $\theta_{g,f}$ depends not just on the 1-morphism $g\star f$, but on the factorization of this 1-morphism into $f$ and $g$. This point is important when working in a graded setting because a 1-morphism of the form $\qs^d(f\star g)$ can be factored as $(\qs^{d-n}f)\star (\qs^ng)$ for any integer $n\in\Z$. The corresponding cyclicity morphism in the quantum horizontal trace depends on the chosen $n$, but using relation~\eqref{rel:grading-shift-trace} from the next subsection, one can see that
\begin{equation}\label{eqn:nchange}
\theta_{\qs^ng,\qs^{d-n}f}=q^{n-m}\theta_{\qs^mg,\qs^{d-m}f}
\end{equation}
for all $m,n\in\Z$.
Thus, changing the chosen $n$ only rescales the cyclicity morphism by an overall power of $q$.

We conclude this subsection by noting that the functor
\[
[-]\colon \ccat{C}(x,x)\longrightarrow\hTr_q(\ccat{C})
\]
is linear and thus induces a functor
\[
[-]\colon \Ch(\ccat{C}(x,x))\longrightarrow\Ch(\hTr_q(\ccat{C}))
\]
between the corresponding categories of chain complexes and their associated homotopy categories.
Applying this functor to the formal Bar-Natan bracket of an $(n,n)$ tangle diagram, we get a chain complex in the quantum annular Bar-Natan category. For generic $q$, this complex can be viewed
as a complex in $\TLC^\oplus\subset\Rep(\Uqsl)$ or, after forgetting the quantum group action, as a complex of $\k$-modules. As such, it generalizes the {\it annular Khovanov complex} defined in [APS].

\subsection{Quantum horizontal trace of $\ccat C_0$}\label{subs:C0}
In this paper, we will usually work with the full Bar-Natan bicategory $\BN$. However, most of our results could also be formulated using the bicategory $\BN_0\subset\BN$ in which 1-morphisms still have formal grading shifts, but 2-morphisms are required to have degree 0 (after taking into account these formal grading shifts). For our results to remain true for $\BN_0$, we have to make a minor modification in our definition of the quantum horizontal trace.

Thus, let $\ccat{C}$ be a graded bicategory and let $\ccat{C}_0\subset\ccat{C}$ be the bicategory which has the same objects and 1-morphisms but only 2-morphisms of degree 0. In the quantum horizontal trace of $\ccat{C}$, the relation
\begin{equation}\label{rel:grading-shift-trace}
		\begin{diagram}[baseline={2.1em-height("$\vcenter{}$")}]  % in diagrams.sty, the default unit for the y-coordinate seems to be 4.2em
			\Node (X1) at (0,1) {x}; \Node (X2) at (1,1) {x};
			\Node (Y1) at (0,0) {y}; \Node (Y2) at (1,0) {y};
			\Arrow[->] (X1) -- (X2) node[midway] {$f$};
			\Arrow[->] (Y1) -- (Y2) node[midway,below] {$g$};
			\Arrow[->] (X1) -- (Y1) node[midway,left] {$\qs p$};
			\Arrow[->] (X2) -- (Y2) node[midway] {$\qs p$};
			\Arrow[double equal sign distance, double, -{Implies[]}]
				($(X2)!1.5em!(0.1,0)$) -- ($(0.1,0)!1.5em!(X2)$)
				node[midway,anchor=135]{$\qs\alpha$};
		\end{diagram}
		\quad\sim\quad q\cdot
		\begin{diagram}[baseline={2.1em-height("$\vcenter{}$")}]
			\Node (X1) at (0,1) {x}; \Node (X2) at (1,1) {x};
			\Node (Y1) at (0,0) {y}; \Node (Y2) at (1,0) {y};
			\Arrow[->] (X1) -- (X2) node[midway] {$f$};
			\Arrow[->] (Y1) -- (Y2) node[midway,below] {$g$};
			\Arrow[->] (X1) -- (Y1) node[midway,left] {$p$};
			\Arrow[->] (X2) -- (Y2) node[midway] {$p$};
			\Arrow[double equal sign distance, double, -{Implies[]}]
				($(0.9,1)!1.5em!(Y1)$) -- ($(Y1)!1.5em!(0.9,1)$)
				node[midway,anchor=-45]{$\alpha$};
		\end{diagram}
\end{equation}
is satisfied because
\begin{equation}\label{rel:grading-shift-explanation}
		\begin{diagram}[baseline={2.1em-height("$\vcenter{}$")}]  % in diagrams.sty, the default unit for the y-coordinate seems to be 4.2em
			\Node (X1) at (0,1) {x}; \Node (X2) at (1,1) {x};
			\Node (Y1) at (0,0) {y}; \Node (Y2) at (1,0) {y};
			\Arrow[->] (X1) -- (X2) node[midway] {$f$};
			\Arrow[->] (Y1) -- (Y2) node[midway,below] {$g$};
			\Arrow[->] (X1) .. controls ++(235:0.5) and ++(125:0.5) .. (Y1) node[left,midway] {$\qs p$};
			\Arrow[->] (X1) .. controls ++(-55:0.5) and ++( 55:0.5) .. (Y1) node[right,pos=0.25] {$p$};
			\Arrow[->] (X2) -- (Y2) node[midway] {$\qs p$};
			\Arrow[double equal sign distance, double, -{Implies[]}]
				($(X2)!1.5em!(0.1,0)$) -- ($(0.1,0)!1.5em!(X2)$)
				node[midway,anchor=135]{$\alpha\circ\mathfrak{s}_p^{-1}$};
			\Arrow[double equal sign distance, double, -{Implies[]}]
				(0.15,0.5) -- (-0.15,0.5)
				node[midway,above]{$\Sigma\mathfrak{s}_p$};
		\end{diagram}
		\quad\sim\quad
		\begin{diagram}[baseline={2.1em-height("$\vcenter{}$")}]
			\Node (X1) at (0,1) {x}; \Node (X2) at (1,1) {x};
			\Node (Y1) at (0,0) {y}; \Node (Y2) at (1,0) {y};
			\Arrow[->] (X1) -- (X2) node[midway] {$f$};
			\Arrow[->] (Y1) -- (Y2) node[midway,below] {$g$};
			\Arrow[->] (X2) .. controls ++(235:0.5) and ++(125:0.5) .. (Y2) node[left,pos=0.7] {$\qs p$};
			\Arrow[->] (X2) .. controls ++(-55:0.5) and ++( 55:0.5) .. (Y2) node[right,midway] {$p$};
			\Arrow[->] (X1) -- (Y1) node[midway,left] {$p$};
			\Arrow[double equal sign distance, double, -{Implies[]}]
				($(0.9,1)!1.5em!(Y1)$) -- ($(Y1)!1.5em!(0.9,1)$)
				node[midway,anchor=-45]{$\alpha\circ\mathfrak{s}_p^{-1}$};
			\Arrow[double equal sign distance, double, -{Implies[]}]
				(1.15,0.5) -- (0.85,0.5)
				node[midway,above]{$\mathfrak{s}_p$};
		\end{diagram}
\end{equation}
where $\Sigma$ denotes the quantum endofunctor, and $\mathfrak{s}$ is the natural degree 1 isomorphism between the identity functor and the grading shift functor $\qs$. 

However, relation \eqref{rel:grading-shift-trace} is a priori not longer true in $\hTr_q(\ccat C_0)$ because the 2-morphism $\mathfrak{s}_p$ in~\eqref{rel:grading-shift-explanation} has nonzero degree.

Since~\eqref{rel:grading-shift-trace} is important for our purposes, we therefore {\it redefine $\hTr_q(\ccat C_0)$ by imposing~\eqref{rel:grading-shift-trace} as an extra relation}, which is required to hold in addition to the usual horizontal trace relation. It is then easy to see that
\begin{equation}\label{eqn:C0iso}
\hTr_q(\ccat C_0)\cong\hTr_q(\ccat C)_0
\end{equation}
where $\hTr_q(\ccat C)_0$ denotes the subcategory of $\hTr_q(\ccat C)$ that contains only degree 0 morphisms, and where the degree of a morphism in $\hTr_q(\ccat C)$ is defined by $|[p,\alpha]|:=|\alpha|$.

The extra relation~\eqref{rel:grading-shift-trace} (or the resulting isomorphism~\eqref{eqn:C0iso}) also ensures that equation~\eqref{eqn:nchange} still holds in $\hTr_q(\ccat C_0)$. A further consequence of~\eqref{eqn:C0iso} is that the naturality condition
\[
q^{|\beta|}[\beta\star\alpha]=\theta_{g',f'}\circ[\alpha\star\beta]\circ\theta_{g,f}^{-1}
\]
still holds in $\hTr_q(\ccat C_0)$ whenever $|\alpha\star\beta|=0$. This can also be seen directly as follows: suppose $\alpha$ has degree $|\alpha|=:n$ (so that $|\beta|=-n$), and let $\tilde{\alpha}\colon f\dblto\qs^{-n} f'$ and $\tilde{\beta}\colon g\dblto\qs^n g'$ be the 2-morphisms $\tilde{\alpha}:=\mathfrak{s}^{-n}\circ\alpha$ and $\tilde{\beta}:=\mathfrak{s}^n\circ\beta$, so that $\tilde{\alpha}\star\tilde{\beta}=\alpha\star\beta$ and $\tilde{\beta}\star\tilde{\alpha}=\beta\star\alpha$. Since $\tilde{\alpha}$ and $\tilde{\beta}$ have degree 0, we then have
\[
[\tilde{\beta}\star\tilde{\alpha}]=\theta_{\qs^ng',\qs^{-n}f'}\circ[\tilde{\alpha}\star\tilde{\beta}]\circ\theta_{g,f}^{-1}
\]
in $\hTr_q(\ccat C_0)$, and since~\eqref{eqn:nchange} still holds in $\hTr_q(\ccat C_0)$, this implies the desired result.

Returning to the quantum horizontal trace of the full bicategory $\ccat C$, we note that $\hTr_q(\ccat C)$ is pregraded via the definition $|[p,\alpha]|:=|\alpha|$ mentioned earlier. We can turn $\hTr_q(\ccat C)$ into a graded category by setting
\begin{equation}\label{eqn:grading-shift-on-qhTr}
\qs[f]:=[\qs f],\qquad\qs[p,\alpha]:=[p,\qs\alpha],\qquad\mathfrak{s}_{[f]}:=[\one,\mathfrak{s}_f].
\end{equation}
Note that the automorphism $\qs\colon\hTr_q(\ccat C)\rightarrow\hTr_q(\ccat C)$ sends the subcategory $\hTr_q(\ccat C)_0\subset\hTr_q(\ccat C)$ to itself, and that every morphism in $\hTr_q(\ccat C)$ can be obtained from a morphism in $\hTr_q(\ccat C)_0$ by pre- and post-composing with suitable powers of $\mathfrak{s}$.

\subsection{Connection with the quantum vertical trace}\label{ss:vertical-trace} Let $\ccat C$ be a pregraded bicategory and $q\in\k^\times$ be a fixed invertible element. The {\it quantum vertical trace} of $\ccat C$ is the category obtained from $\ccat C$ by replacing each morphism category $\ccat C(x,y)$ by its quantum trace $\Tr_q(\ccat C(x,y))$ (see~\cite{BPW19}). Equivalently:

\begin{definition} The~\emph{quantum vertical trace} of $\ccat C$
is the $\k$-linear category $\vTr_q(\ccat C)$ with the same objects as $\ccat C$, and whose morphism are given by formal $\k$-linear combinations of 2-endomorphisms in $\ccat C$, subject to the relation $\alpha\circ\beta\sim q^{|\alpha|}\beta\circ\alpha$ for 2-morphisms $\alpha\colon p\dblto p'$ and $\beta\colon p'\dblto p$. Composition in $\vTr_q(\ccat C)$ is induced by the horizontal composition in $\ccat C$.
\end{definition}

If $\ccat C$ is graded, then we modify the definition of $\vTr_q(\ccat C_0)$ by also imposing the relation $\qs\alpha\sim q\alpha$ for all 2-endomorphisms $\alpha$ in $\ccat C_0$. This is to ensure that
\begin{equation}\label{eqn:vTrdegreezero1}
\vTr_q(\ccat C_0)\cong\vTr_q(\ccat C)_0.
\end{equation}
In general, the quantum vertical trace admits a fully faithful embedding $\vTr_q(\ccat C)\rightarrow\hTr_q(\ccat C)$ given by sending each object $x$ to its identity 1-morphism $\one_x$ and each $\alpha\colon p\dblto p$ to the square from Lemma~\ref{lem:annihilate}. This embedding can be used to reestablish the result from Lemma~\ref{lem:annihilate}. Indeed, in $\vTr_q(\ccat C)$, we have
\[
\alpha\sim\alpha\circ\one\sim q^d(\one\circ\alpha)\sim q^d\alpha
\]
for any homogeneous 2-endomorphism $\alpha$ of degree $|\alpha|=d$, and thus $\alpha$ is annihilated by $1-q^d$. When $1-q^d$ is invertible, this implies the lemma. Moreover, if $1-q^d$ is invertible for all $d\neq 0$, then it follows that
\begin{equation}\label{eqn:vTrdegreezero2}
\vTr_q(\ccat C)=\vTr_q(\ccat C)_0
\end{equation}
because in this case there are no nonzero morphisms of nonzero degree. In particular, the quantum vertical trace is usually only pregraded, with all morphisms supported in degree 0.

We are now ready to prove Theorem~\ref{thm:qtrace of BN}.

\begin{proof}[Proof of Theorem~\ref{thm:qtrace of BN}] Ignoring grading shifts, every object in $\cat{BN}_{\!q}(\Ann)$ can be viewed as a closed 1-manifold $C$ embedded in the annuluar $\Ann$. By applying an ambient isotopy to $C$, we can move trivial components in $C$ away from the arc $\mu\subset\Ann$ and straighten essential components so that they become concentric circles. The resulting 1-manifold $C'\subset\Ann$ is a closure of an identity tangle with some trivial circles. Moreover, the isotopy from $C$ to $C'$ can be realized as a composition of cyclicity isomorphisms in $\hTr_q(\BN)$ and isotopies in $I^2$.

In $\hTr_q(\BN)^\oplus$, we also have a ``delooping isomorphism'' $\bigcirc\cong\qs\emptyset\oplus\qs^{-1}\emptyset$ given by the matrix $\begin{bmatrix} S_\bullet & S \end{bmatrix}^t$ where $S$ is a cap cobordism and $S_\bullet$ denotes $S$ with a dot. Applying this isomorphism to each trivial circle in $C'$, we can reduce $C'$ to a direct sum of shifted copies of a closure of an identity tangle.

We conclude that $\hTr_q(\BN)^\oplus$ is obtained from image of $\vTr_q(\BN)\rightarrow\hTr_q(\BN)$ by allowing grading shifts and direct sums. Since $1-q^d$ is invertible for all $d\neq 0$, we also have $\vTr_q(\BN)=\vTr_q(\BN)_0\cong\vTr_q(\BN_0)$ by~\eqref{eqn:vTrdegreezero2} and \eqref{eqn:vTrdegreezero1}, where the latter isomorphism uses that we have imposed the additional relation $\qs\alpha\sim q\alpha$ in the definition of $\vTr_q(\BN_0)$.

It remains to show that $\vTr_q(\BN_0)\simeq\TLC$. This was essentially shown by Bar-Natan in~\cite[Prop.~10.10]{B-N05}, but we will
prove the result using a slightly different argument. For this, let
\[
\phi\colon \TL\longrightarrow\vTr_q(\BN_0)
\]
denote the functor which sends each object of $\TL$ to the corresponding object of $\BN_0$, and each flat tangle $T$ to the class of its identity 2-endomorphism in $\vTr_q(\BN_0)$.

To see that $\phi$ is compatible with isotopy of tangles and with the relation $\bigcirc=q+q^{-1}$, recall from \cite{BPW19} that the assignment $T\mapsto\one_T$ induces a well-defined map from the split Grothendieck group of $\BN_0(n,m)^\oplus$ to the quantum trace
\[
\Tr_q(\BN_0(n,m)^\oplus)=\Tr_q(\BN_0(n,m)^\oplus)=\vTr_q(\BN_0)(n,m).
\]
Since isotopic tangles are isomorphic as objects of $\BN_0(n,m)$, it follows that their identity 2-morphisms are equal in the vertical trace. Likewise, the ``delooping isomorphism'' implies $\one_\bigcirc\sim\qs\one_\emptyset+\qs^{-1}\one_\emptyset\sim q\one_\emptyset+q^{-1}\one_\emptyset$ in $\vTr_q(\BN_0)$, which corresponds to $\bigcirc=q+q^{-1}$.

It is clear that $\phi$ is bijective on objects, so the theorem will follow if we can prove that $\phi$ is full and faithful. For this, we can focus on a single morphism category $\BN_0(n,m)$. Using the ``delooping isomorphism'' and the invariance of the quantum trace under graded equivalences and under additive closures, we can further restrict our attention to the full subcategory $\mathcal{T}\subset\BN_0(n,m)$ generated by tangles without closed components. That is, we only need to consider 2-morphisms of the form $\alpha\colon \qs^iT\dblto \qs^jT'$ where $T$ and $T'$ are flat tangles without closed components.

Let $S\colon T\dblto T'$ be a dotted cobordism representing such an $\alpha$. We first claim that the unshifted degree of $S$ must satisfy $|S|\geq 0$, and $|S|=0$ if and only if $S$ is isotopic to a scalar multiple of an identity cobordism. To see this, note that the surface $S\subset I^3$ can be obtained from a cobordism $S'\subset I^3$ with only disk components by adding tubes, dots, and closed components. Using that undotted spheres are equal to to zero and that $T$ and $T'$ have no closed components, it is easy to see that $|S|\geq |S'|\geq 0$. The second part of the claim follows because closed components evaluate to scalars.

We conclude that any 2-endomorphism $\alpha\colon \qs^i T\dblto \qs^i T$ in $\mathcal{T}\subset\BN_0(n,m)$ must be a scalar multiple of an identity 2-morphism, which proves that $\phi$ is full.

To see that $\phi$ is faithful, consider a relation of the form $\alpha\circ\beta\sim q^{|\alpha|}\beta\circ\alpha$ for 2-morphisms $\alpha\colon \qs^iT\dblto \qs^jT'$ and $\beta\colon \qs^jT'\dblto \qs^iT$ in $\mathcal{T}$. We can drop the factor of $q^{|\alpha|}$ because $|\alpha|=0$, and we can further assume that the 2-morphisms $\alpha$ and $\beta$ are represented by dotted cobordisms $S\colon T\dblto T'$ and $S'\colon T'\dblto T$. Since $|S|+|S'|=|S\circ S'|=|\alpha\circ\beta|=0$ and since $|S|$ and $|S'|$ are nonnegative, it follows that $|S|=|S'|=0$. However, this implies that $S$ and $S'$ are isotopic to scalar multiples of identity cobordisms, and that $i=j$. 

In summary, the relation $\alpha\circ\beta\sim q^{|\alpha|}\beta\circ\alpha$ only has the effect of identifying the identity 2-morphisms $\qs^i\one_{T'}$ and $\qs^i\one_{T}$ of isotopic tangles in $\mathcal{T}$. By applying Khovanov's universal TQFT, one can further see that $c\one_T$ is nonzero in $\mathcal{T}$ for any $c\neq 0$ (for this to make sense, regard $\one_T\subset I^3$ as a cobordism from the empty 1-manifold to $\partial\one_T$).

We conclude that $\phi$ is faithful, which completes the proof of the theorem.
\end{proof}

\section{Homological toolkit}
\label{ss:twists}
In this section we discuss the dg notions of twisting,  splicing, and homological perturbation.  The tools introduced here will be indispensable when discussing models for the Cooper--Krushkal projector.

We discuss the notion of twisting, and we do so first in an abstract setting.  So let $\AS$ be a $\k$-linear category, and consider the dg category $\Ch(\AS)$ of complexes.

If $(X,\d_X)$ is a complex over $\AS$, then any complex of the form $(X,\d_X+\a)$ will be referred to as a \emph{twist} of $X$, written $\tw_\a(X)$.  If $X$ is equipped with a direct sum decomposition $X = \bigoplus_{i\in I} X_i$  (in which the differential on $X$ is the direct sum of differentials on $X_i$) then $\a$ can be expressed in terms of its components $\a_{ij}\in \Hom^{1}_{\Ch(\AS)}(X_j,X_i)$, and the condition that $\d_X+\a$ be a differential on $X$ can be expressed componentwise as
\[
(\d_{X_i}\circ \a_{ij}+\a_{ij}\circ \d_{X_j}) + \sum_k \a_{ik}\circ \a_{kj} = 0.
\]
for all $i,j\in I$.  The twist $\tw_\a(X)$ is \emph{one-sided} if there is a direct sum decomposition $X\cong \bigoplus_{i\in I} X_i$ where the indexing set $I$ is equipped with a partial order with respect to which the matrix representing $\a$ is strictly lower triangular (i.e.,~$\a_{ij}=0$ unless $i>j$).   Expressions of the form $\tw_\a(\bigoplus_i X_i)$ are also called \emph{one-sided twisted complexes}.

We will denote one-sided twisted complexes of the form $\tw_\a (\bigoplus_{i\in I} X_i)$ diagrammatically, as in the following examples:

\begin{subequations}
\begin{equation}\label{eq:tw diagram 1}
\tw_{\smMatrix{ 0 & 0\\ \a & 0 }}\left(X\oplus Y\right)   \ \ =: \ \ \left(\begin{tikzpicture}[baseline=-.2em]
\tikzstyle{every node}=[font=\small]
\node (a) at (0,0) {$X$};
\node (b) at (2.5,0) {$Y$};
\path[->,>=stealth,shorten >=1pt,auto,node distance=1.8cm,
  thick]
(a) edge node {$\a$} (b);
\end{tikzpicture}\right)
\end{equation}
\begin{equation}\label{eq:tw diagram 2}
\tw_{\smMatrix{0 & 0 & 0\\ \a & 0 & 0 \\ \gamma & \b& 0}}\left(X\oplus Y\oplus Z\right)   \ \ =: \ \ \left(\begin{tikzpicture}[baseline=0em]
\tikzstyle{every node}=[font=\small]
\node (a) at (0,0) {$X$};
\node (b) at (2.5,0) {$Y$};
\node (c) at (5,0) {$Z$};
\path[->,>=stealth,shorten >=1pt,auto,node distance=1.8cm,
  thick]
(a) edge node {$\a$} (b)
(b) edge node {$\b$} (c)
(a) edge [in=150,out=30] node {$\gamma $} (c);
\end{tikzpicture}\right)
\end{equation}
\end{subequations}

In particular the mapping cone of a degree zero chain map $f:X\rightarrow Y$ will be expressed as
\[
\Cone(f)  \ = \ \Big(\ts\inv X\buildrel f\over\longrightarrow Y\Big).
\]
Here, as always, we abuse notation by regarding the degree zero morphism $f:X\rightarrow Y$ also as a degree $1$ morphism $\ts\inv X\rightarrow Y$.

\subsection{Homological perturbation}
\label{sss:hpl}
If $X\simeq Y$ then a given twist $\tw_\a(X)$ may or may not be homotopy equivalent to a corresponding twist of $Y$.  For instance if $X$ is any nonzero but contractible complex and $Y=0$, then there certainly exist twists of $X$ which are not contractible (and therefore not homotopy equivalent to any twist of $Y$).

Homological perturbation concerns the problem twisting a homotopy equivalence $X\simeq Y$ to a homotopy equivalence $\tw_\a(X)\simeq \tw_\b(Y)$, provided that certain conditions are met.

\begin{lemma}\label{lemma:hpl}
Let $I$ be a partially ordered set satisfying the following ascending chain condition: for each $i\in I$ there are finitely many $j\in I$ with $j\geq i$.  Suppose we have a one-sided twist $\tw_\a(\bigoplus_i X_i)$, and complexes $Y_i\simeq X_i$ for all $i$.  Then there exists a one-sided twist $\tw_\b(\bigoplus_{i\in I} Y_i)$ such that $\tw_\b(\bigoplus_{i\in I} Y_i)\simeq \tw_\a(\bigoplus_i X_i)$.

In more details, if $f_i\in \Hom^0(X_i,Y_i)$ and $g_i\in \Hom^0(Y_i,X_i)$ are inverse homotopy equivalences with $h_i\in \End^{-1}(X_i)$ satisfying $[\d_{X_i},h_i] = \Id_{X_i}- g_i\circ f_i$, then the twist $\b$ is given in terms of components by
\begin{eqnarray*}
\b_{kl} &=& f_k\circ \a_{kl}\circ g_l  - \sum_{k>i_1>l} f_k \circ \a_{k,i_1}\circ h_{i_1}\circ \a_{i_1,l}\circ g_l \\
& & + \sum_{k>i_1>i_2>l}  f_k \circ \a_{k,i_1}\circ h_{i_1}\circ \a_{i_1,i_2}\circ h_{i_2}\circ \a_{i_2,l}\circ g_l + \cdots
\end{eqnarray*}
if $k>l$ and $\b_{kl}=0$ otherwise, and the homotopy equivalence $F:\tw_\a(X)\rightarrow \tw_\b(Y)$ is given in terms of components by  $F_{kk}=f_k$ and
\begin{eqnarray*}
F_{kl} &=&  - f_k\circ \a_{k,l}\circ h_l + \sum_{k>i_1>l} f_k \circ \a_{k,i_1}\circ h_{i_1}\circ \a_{i_1,l}\circ h_l \\
& & - \sum_{k>i_1>i_2>l}  f_k \circ \a_{k,i_1}\circ h_{i_1}\circ \a_{i_1,i_2}\circ h_{i_2}\circ \a_{i_2,l}\circ h_l + \cdots
\end{eqnarray*}
if $k>l$, and $F_{kl}=0$ otherwise.
\end{lemma}

\noindent

\begin{proof}
See \cite[Theorems 2 and 3, and Subs. 8.6]{MarklIdeal}.
\end{proof}

\begin{remark}\label{rmk:adjacent component}
In particular, if indices $k<l\in I$ are ``adjacent'', meaning there are no indices $i\in I$ with $k<i<l$, then the component $\b_{kl}$ equals $g_k\circ \a_{kl}\circ f_l$.
\end{remark}

\begin{remark}
The assumption that $I$ satisfies the ascending chain condition can be weakened as follows.  For each $k\in \Z$, let $I_k\subset I$ be the subset consisting of those $i\in I$ for which the chain group $(X_i)^k$ is nonzero.   The statement of Lemma \ref{lemma:hpl} is true under the assumption that $I_k$ satisfies the ascending chain condition for all $k$.   We may refer to this as the \emph{local ascending chain condition relative to $X$}.
\end{remark}
\begin{remark}\label{rmk:contractible one-sided twist}
As a special case of Lemma \ref{lemma:hpl}, if each complex $X_i$ is contractible then any one-sided twist $\tw_\a(\bigoplus_{i\in I}X_i)$ is also contractible, provided that the partially ordered set $I$ satisfies the ascending chain condition.
\end{remark}

\begin{remark}\label{rem:descending}
There is a version of Lemma \ref{lemma:hpl} in which the ascending chain condition is replaced by the \emph{descending chain condition} and the direct sum is replaced by direct product.
\end{remark}

\subsubsection{Splicing}
\label{sss:splicing}
If  we have twisted complexes (mapping cones)
\[
X_1 \ = \ \left(\begin{tikzpicture}[baseline=0em]
\tikzstyle{every node}=[font=\small]
\node (a) at (0,0) {$A$};
\node (b) at (2.5,0) {$E$};
\path[->,>=stealth,shorten >=1pt,auto,node distance=1.8cm,
  thick]
(a) edge node {$\a$} (b);
\end{tikzpicture}\right),\qquad \qquad X_2 \ = \ \left(\begin{tikzpicture}[baseline=0em]
\tikzstyle{every node}=[font=\small]
\node (a) at (0,0) {$\ts\inv E$};
\node (b) at (2.5,0) {$B$};
\path[->,>=stealth,shorten >=1pt,auto,node distance=1.8cm,
  thick]
(a) edge node {$\b$} (b);
\end{tikzpicture}\right),
\]
then the result of \emph{splicing} $X_1$ and $X_2$ is the twisted complex
\[
Z \ \:= \
 \left(\begin{tikzpicture}[baseline=0em]
\tikzstyle{every node}=[font=\small]
\node (a) at (0,0) {$A$};
\node (b) at (2.5,0) {$B$};
\path[->,>=stealth,shorten >=1pt,auto,node distance=1.8cm,
  thick]
(a) edge node {$\b\circ \a$} (b);
\end{tikzpicture}\right).
\]
The spliced complex $Z$ is homotopy equivalent to the mapping cone:
\[
Z \ \simeq \ \left(\begin{tikzpicture}[baseline=0em]
\tikzstyle{every node}=[font=\small]
\node (a) at (0,0) {$X_1$};
\node (b) at (2.5,0) {$X_2$};
\path[->,>=stealth,shorten >=1pt,auto,node distance=1.8cm,
  thick]
(a) edge node {$-\pi$} (b);
\end{tikzpicture}\right)
\]
where $\pi$ is the projection-followed-by-inclusion of the $E$ summand; this has homological degree $1$, given that $E$ appears in $X_2$ with shift $\ts\inv$.

The operation of splicing can be performed for ``longer'' one-sided twisted complexes.  For instance the result of splicing together twisted complexes
\[
X_1 \ \:= \ \left(\begin{tikzpicture}[baseline=0em]
\tikzstyle{every node}=[font=\small]
\node (a) at (0,0) {$A$};
\node (b) at (2.5,0) {$B$};
\node (c) at (5,0) {$E$};
\path[->,>=stealth,shorten >=1pt,auto,node distance=1.8cm,
  thick]
(a) edge node {$\d_{BA}$} (b)
(b) edge node {$\a$} (c)
(a) edge [in=150,out=30] node {$\a'$} (c);
\end{tikzpicture}\right)
\]
and
\[
X_2  \ \:= \ \left(\begin{tikzpicture}[baseline=0em]
\tikzstyle{every node}=[font=\small]
\node (a) at (0,0) {$\ts\inv E$};
\node (b) at (2.5,0) {$C$};
\node (c) at (5,0) {$D$};
\path[->,>=stealth,shorten >=1pt,auto,node distance=1.8cm,
  thick]
(a) edge node {$\b$} (b)
(b) edge node {$\d_{DC}$} (c)
(a) edge [in=150,out=30] node {$\b'$} (c);
\end{tikzpicture}\right),
\]
is the twisted complex
\begin{equation}\label{eq:ABCD}
Z \ \:= \ \left(\begin{minipage}{2.59in}
\begin{tikzpicture}
\node (a) at (0,0){$A$};
\node (b) at (2,0){$B$};
\node (c) at (4,0){$C$};
\node (d) at (6,0){$D$};
\path[->,>=stealth',shorten >=1pt,auto,node distance=1.8cm,font=\small]
(a) edge node[above,xshift=4pt] {$\d_{BA}$} (b)
(b) edge node {$\d_{CB}$} (c)
(c) edge node {$\d_{DC}$} (d)
(a) edge[bend left=30] node {$\d_{DA}$} (d);
\draw[frontline,->,>=stealth',shorten >=1pt,auto,node distance=1.8cm]
(a) to [bend right] node[below] {$\d_{CA}$} (c);
\draw[frontline,->,>=stealth',shorten >=1pt,auto,node distance=1.8cm]
(b) to [bend right] node[below] {$\d_{DB}$} (d);
\end{tikzpicture}
\end{minipage}\right),
\end{equation}
where 
\[
\d_{CB}=\b\circ \a,\qquad \d_{CA} = \b\circ \a',\qquad \d_{DB}:=\b'\circ \a,\qquad \d_{DA}=\b'\circ \a'.
\]

We may iterate the operation of splicing, even countably infinitely many times.  For instance, suppose we have a collection of twisted complexes of the form
\[
X_i \ \:= \ \ \left(\begin{tikzpicture}[baseline=0em]
\tikzstyle{every node}=[font=\small]
\node (a) at (0,0) {$\ts\inv E_{i+1}$};
\node (b) at (2.5,0) {$C_i$};
\node (c) at (5,0) {$E_i$};
\path[->,>=stealth,shorten >=1pt,auto,node distance=1.8cm,
  thick]
(a) edge node {$\a_i$} (b)
(b) edge node {$\b_i$} (c)
(a) edge [in=150,out=30] node {$\gamma_i$} (c);
\end{tikzpicture}\right),
\]
where $C_i$ and $E_i$ are some complexes ($i=0,1,2,\ldots$).  Let $\pi_i\in \Hom^1(X_i,X_{i+1})$ be the degree 1 chain map which projects onto then includes the $E_{i+1}$ term.  Then $\pi_{i+1}\circ \pi_i=0$ and we have a complex of the form
\[
Y = \left(
\begin{tikzpicture}[baseline=0em]
\tikzstyle{every node}=[font=\small]
\node (a) at (0,0) {$X_0$};
\node (b) at (2.5,0) {$X_1$};
\node (c) at (5,0) {$X_2$};
\node (d) at (7.5,0) {$\cdots $};
\path[->,>=stealth,shorten >=1pt,auto,node distance=1.8cm,
  thick]
(a) edge node {$-\pi_0$} (b)
(b) edge node {$-\pi_1$} (c)
(c) edge node {$-\pi_2$} (d);
\end{tikzpicture}
\right).
\]
Strictly speaking we can express this complex as $Y=\tw_\a(\bigoplus_{i\geq 0}  X_i)$ where the twist $\a$ is given in terms of components by $-\pi_i$.  One can also express $Y$ as the total complex of a bicomplex.

We will call $Y$ the \emph{pre-spliced} complex constructed from the complexes $X_i$.
The following lemma shows that $Y$ is homotopy equivalent to a \emph{spliced} complex.

\begin{lemma}\label{lemma:infinite splicing}
Retain notation as in the paragraph above.  Then $Y$ is homotopy equivalent to a twisted complex of the form
\[
Y \simeq \left(
\begin{tikzpicture}[baseline=2em]
\tikzstyle{every node}=[font=\small]
\node (dd) at (-7.5,1.5) {};
\node (dc) at (-7.5,1) {};
\node (db) at (-7.5,.5) {};
\node (da) at (-7.5,0) {$\cdots$};
\node (c) at (-5,0) {$C_{1}$};
\node (b) at (-2.5,0) {$C_{0}$};
\node (a) at (0,0) {$E_0$};
\path[->,>=stealth,shorten >=1pt,auto,node distance=1.8cm,
  thick]
(b) edge node {} (a)
(c) edge node {} (b)
(da) edge node {} (c)
(c) edge [in=150,out=30] node {} (a)
(db) edge [in=150,out=20] node {} (b)
(dc) edge [in=130,out=30] node {} (a)
(dd) edge [in=130,out=10] node {} (b)
;
\end{tikzpicture}
\right),
\]
in which we have not drawn all of the (countably many) components of the differential, but all of them point rightward.
\end{lemma}
\begin{proof}
We can write $Y$ as a twisted complex of the form
\[
\left(
\begin{tikzpicture}[baseline=-4em]
\tikzstyle{every node}=[font=\small]
\node (aa) at (-4,0) {$\ts\inv E_{1}$};
\node (ab) at (-2,0) {$C_0$};
\node (ac) at (0,0) {$E_0$};
\node (ba) at (-6,-2) {$\ts\inv E_{2}$};
\node (bb) at (-4,-2) {$C_1$};
\node (bc) at (-2,-2) {$E_1$};
\node (cb) at (-6,-4) {$\cdots$};
\node (cc) at (-4,-4) {$E_2$};
\path[->,>=stealth,shorten >=1pt,auto,node distance=1.8cm,
  thick]
(aa) edge node {$\a_0$} (ab)
(ab) edge node {$\b_0$} (ac)
(aa) edge [in=150,out=30] node {$\gamma_0$} (ac)
(ba) edge node {$\a_1$} (bb)
(bb) edge node {$\b_1$} (bc)
(ba) edge [in=150,out=30] node {$\gamma_1$} (bc)
(cb) edge node {$\b_2$} (cc)
(aa) edge node {$-\id$} (bc)
(ba) edge node {$-\id$} (cc);
\end{tikzpicture}
\right)
\]
This can be rearranged into a twisted complex of the form
\[
\left(
\begin{tikzpicture}[baseline=0em]
\tikzstyle{every node}=[font=\small]
\node (g) at (-13,0) {$\cdots$};
\node (f) at (-11.5,0) {$C_2$};
\node (e) at (-9,0) {$(\ts\inv E_2\rightarrow E_2)$};
\node (d) at (-6.5,0) {$C_1$};
\node (c) at (-4,0) {$(\ts\inv E_1\rightarrow E_1)$};
\node (b) at (-1.5,0) {$C_0$};
\node (a) at (0,0) {$E_0$};
\path[->,>=stealth,shorten >=1pt,auto,node distance=1.8cm,
  thick]
(g) edge node {} (f)
(f) edge node {} (e)
(e) edge node {} (d)
(d) edge node {} (c)
(c) edge node {} (b)
(b) edge node {} (a)
(g) edge [in=150,out=30] node {} (e)
(e) edge [in=150,out=30] node {} (c)
(c) edge [in=150,out=30] node {} (a);
\end{tikzpicture}
\right)
\]
Apply Lemma \ref{lemma:hpl} to cancel the contractible terms $(\ts\inv E_i \buildrel-\id\over \longrightarrow E_i)$.  Even though there are infinitely many such terms, this cancelation is valid since the differential respects a partial order (satisfying the accending chain condition) on the terms.
\end{proof}

\subsection{Op-well-ordered twisted complexes}
\label{ss:well ordered}

\begin{definition}
A partially ordered set $I$ is \emph{op-well-ordered} if every non-empty subset $S\subset I$ has a unique maximum element.
\end{definition}
Note that if $I$ is op-well-ordered, then the partial order is a total order: every pair of elements $i,j$ is comparable since $\{i,j\}$ has a maximum.  Moreover, every op-well-ordered set satisfies the \emph{ascending chain condition}: any increasing sequence $i_0\leq i_1\leq \cdots$ eventually stabilizes to the maximum element of $\{i_0,i_1,\ldots \}$.

\begin{example}\label{ex:E is well-ordered}
Suppose $r\colon E\rightarrow B$ is a map of sets so that the base $B$ is op-well-ordered and the fibers $r\inv(b)$ are op-well-ordered for all $b$.  Then $E$ has an op-well-ordering given by $e\leq e'$ if $r(e)<r(e')$ in $B$ or $r(e)=r(e')=:b$ and $e\leq e'$ in $r\inv(b)$.
\end{example}

\begin{corollary}
If $I$ and $J$ are op-well-ordered then so is $I\times J$ with the lexicographic ordering $(i,j)\leq (i',j')$ if $i<i'$ or $i=i'$ and $j\leq j'$.\qed
\end{corollary}

Op-well-ordered sets arise as sets which index iterated one-sided twisted complexes.
\begin{definition}
Suppose $\AS$ is an additive category.  An \emph{op-well-ordered twisted complex} in $\Ch^-(\AS)$ is a one-sided twisted complex $\tw_\d(\bigoplus_{i\in I} \ts^{a_i}X_i)$ where $I$ is op-well-ordered, and the shifts $a_i$ are bounded above.

If $\mathcal{X}$ is a collection of complexes in $\Ch^-(\AS)$, we let $\langle \mathcal{X} \rangle$ denote full subcategory of $\Ch^-(\AS)$ consisting of op-well ordered twisted complexes whose terms are finite direct sums of shifts of direct summands of objects in $\mathcal{X}$. 
\end{definition}

\begin{definition}
If $\mathcal{X}$ is a collection of complexes in $\Ch^-(\AS)$, then we let $\overline{\mathcal{X}}$ denote the closure of $\mathcal{X}$ with respect to homotopy equivalences.
\end{definition}

\begin{lemma}\label{lemma:S in T}
Let $\mathcal{X},\mathcal{Y}\subset \Ch^-(\AS)$ be some collection of complexes such that each $N\in \mathcal{X}$ is homotopy equivalent to a complex in $\langle \mathcal{Y}\rangle$.  Then $\langle \mathcal{X} \rangle\subset \overline{\langle \mathcal{Y}\rangle}$.
\end{lemma}
\begin{proof}
Suppose $X=\tw_{\d}(\bigoplus_{i\in I} X_i)$ with $X_i\in \mathcal{X}$ and $I$ op-well-ordered, and suppose each $X_i$ is homotopy equivalent to $\tw_{\e(i)}(\bigoplus_{j\in J_i} Y_{i,j})$ with $J_i$ op-well-ordered for all $i\in I$, and $Y_{i,j}\in \mathcal{Y}$.  Let $E$ denote the set of pairs $(i,j)$ with $j\in J_i$, and let $r\colon E\rightarrow I$ be the projection map $(i,j)\mapsto i$.  Then $E$ has the structure of an op-well-ordered set via $(i,j)\leq (i',j')$ if $i<i'$ or $i=i'$ and $j<j'$.

Since $I$ satisfies the ascending chain condition, the homological perturbation lemma applies, so
\[
\tw_{\d}\left(\bigoplus_{i\in I} X_i\right) 
\simeq  \tw_{\d'}\left(\bigoplus_{i\in I} \tw_{\e_i} \left(\bigoplus_{j\in J_i}Y_{i,j}\right)\right) 
\cong \tw_{\d'+\e}\left(\bigoplus_{(i,j)\in E} Y_{i,j}\right).
\]
where $\d'$ strictly increases the index $i$ and $\e=\sum_i \e(i)$ is a sum of terms which fix the index $i$ but strictly increase $j$. Thus, the twist $\d'+\e$ is strictly increasing with respect to the partial order on $E$, and this last complex is in $\langle \mathcal{Y}\rangle$.
\end{proof}

\begin{cor}\label{cor:closure of <S>}
The category $\langle \mathcal{X}\rangle \subset \Ch^-(\AS)$ is closed under taking op-well-ordered twisted complexes.\qed
\end{cor}
\begin{remark}
This corollary would be false if we had only allowed $\Z_{\leq 0}$-indexed one-sided twisted complexes in the definition of $\langle \mathcal{X} \rangle$. 
\end{remark}

\subsection{Combing hairs}
\label{ss:combing}

We use the phrase ``combing hairs'' to refer to the construction of isomorphisms of the form $\tw_\d(\bigoplus_{i\in I} X_i)\cong \tw_{\d'}(\bigoplus_{i\in I'} X_i)$ in which $I$ is an op-well-ordered-set and $I'=I$ with a coarser partial order $\leq'$.   Being coarser means that the condition $i\leq' j$ is more restrictive than $i\leq j$. Here are some instances where we wish to apply such ideas (we will expand on these momentarily):

\begin{enumerate}
\item Suppose $\mathcal{X}$ is a collection of complexes such that  $\Hom(N,M)$ has no cohomology in positive degrees for all $N,M\in \mathcal{X}$, and $X=\tw_{\d}(\bigoplus_{i\in I} \ts^{a_i} N_i)$ is an op-well-ordered twisted complex with $N_i\in \mathcal{X}$. Let $I'=I$ with the partial order given by $i<' j$ if $i< j$ and $a_i<a_j$.  Combing hairs in such situations amounts to the elimination of components of the twist $\ts^{a_i}N_i\rightarrow \ts^{a_j} N_j$  which have positive cohomological degree when regarded as elements of $\Hom(N_i,N_j)$.
\item Suppose $(\Omega,\unlhd)$ is a partially ordered set and $\mathcal{X}_\omega\subset \Ch^-(\AS)$ is a collection of complexes indexed by $\omega\in \Omega$ such that if $N\in \mathcal{X}_\omega$ and $M\in \mathcal{X}_{\nu}$ then $\Hom(N,M)\simeq 0$ unless $\omega\unlhd \nu$.   Suppose $X=\tw_{\d}(\bigoplus_{i\in I} X_{i})$ is an op-well-ordered twisted complex with $X_i\in \mathcal{X}_{\omega(i)}$. Let $I'=I$ with the partial order given by $i\leq ' j$ if $i\leq j$ and $\omega(i)\unlhd \omega(j)$.  Combing hairs in such situations amounts to the elimination of components of the twist $X_i\rightarrow X_j$ in which $\omega(j)$ is lower or incomparable to $\omega(i)$.
\end{enumerate}

Our main lemma will combine both of the above notions.  For technical reasons we will need some boundedness and compactness assumptions on our complexes in $\mathcal{X}_\omega$.  To this end, fix an additive category $\AS$, suppose $(\Omega,\unlhd)$ is a partially ordered set, and for each $\omega\in \Omega$ suppose we have a collection of complexes $\mathcal{X}_\omega$ such that:
\begin{enumerate}
\item Each complex in $\mathcal{X}_\omega$ is supported in non-positive cohomological degrees.
\item If $N\in \mathcal{X}_\omega$ and $M\in \mathcal{X}_\nu$ then the chain objects $N^k$ (zero unless $k\leq 0$) satisfy
\[
\Hom(N^k, M)\simeq 0
\]
if $\omega{\nunlhd} \nu$ or $\omega=\nu$ and $k<0$.
\end{enumerate}

\begin{lemma}\label{lemma:combing}
Retain the setup from the previous paragraph.  Suppose $I$ is an op-well-ordered set and we have a mapping $I\rightarrow \Omega\times \Z$, ($i\mapsto (\omega(i),a(i))$) and a twisted complex $X=\tw_\d (\bigoplus_{i\in I} \ts^{a(i)} X_i)$ with $X_i\in \mathcal{X}_{\omega(i)}$.  Assume that the mapping $a\colon I\rightarrow \Z$ has finite preimages and $a\inv(k)=0$ for $k\gg 0$.  Let $I'=I$ with the partial order $i<'j$ if $i<j$ and $(\omega(i),a(i))<(\omega(j),a(j))$ in the lexicographic order ($\omega(i)\lhd \omega(j)$ or $\omega(i)=\omega(j)$ and $a(i)<a(j)$).  Then $X\cong \tw_{\d'}(\bigoplus_{i\in I'} \ts^{a(i)} X_i)$ where the twist $\d'$ respects the partial order $<'$.
\end{lemma}

\begin{proof}
Let's introduce some terminology and notation for the proof. We denote by $0$ the unique maximum element of $I$.  We may assume that $I$ has a unique minimum element $-\infty\in I$ (if not, then we may replace $I$ by $\{-\infty\}\sqcup I$ and set $X_{-\infty}=0$). Recall, a subset $K\subset I$ is \emph{convex} if $k\leq i\leq k'$ and $k,k'\in K$ implies $i\in K$.  Given $i\leq j$ in $I$ we have the interval $[i,j]:=\{k\in I\:|\: i\leq k\leq j\}$, and $(i,j]=[i,j]\setminus \{i\}$ (and so on).  Note that intervals are convex.

If $K\subset I$ is convex then  $X_K:=\tw_{\d_K}(\bigoplus_{i\in K}\ts^{a(i)} X_i)$ is a well-defined one-sided twisted complex, where the twist $\d_K$ is obtained from $\d$ by restriction.  A \emph{combing} of $X_K$ is a pair $(Y_K,\phi_K)$ where $Y_K$ is a one-sided twisted complex of the form in the statement (with respect to the coarse partial order) and $\phi_K$ is an isomorphism $X_K\buildrel \cong \over\rightarrow Y_K$ which preserves the fine partial order on $K$.

We wish to show that $X=X_{[-\infty,0]}$ can be combed.  Our plan is to show that each $X_{[i,0]}$ can be combed, by transfinite induction.   In the base case, we observe that $(X_{\{0\}},\id)$ is a combing of $X_{\{0\}}$, trivially.  

Now, fix $i$ and assume by induction that we have constructed combings $(Y_{[j,0]},\phi_{[j,0]})$ of $X_{[j,0]}$, for all $j>i$, and assume that these combings are compatible in the sense that if $j<j'$ then $\phi_{[j',0]}$ is obtained from  $\phi_{[j,0]}$ by restriction.  Then we have a combing $(Y_{(i,0]},\phi_{(i,0]})$ of $X_{(i,0]}$ by taking the union of combings $(Y_{[j,0]},\phi_{[j,0]})$ for $j>i$.

To complete the proof we must construct a combing of $X_{[i,0]}$.  Note that
\begin{equation}\label{eq:X [i,0]}
X_{[i,0]} \cong (\ts^{a(i)}X_{i}\buildrel \partial\over \rightarrow X_{(i,0]}) \cong (\ts^{a(i)} X_i\buildrel \partial'\over \rightarrow Y_{(i,0]}),
\end{equation}
where the second isomorphism is $\smMatrix{\id &0\\0& \phi_{(i,0]}}$. Let us decompose $(i,0]$ into subsets $J_1\sqcup J_2\sqcup J_3$ as follows.  Let $J_3$ consist of those $j>i$ with $\omega(i)\lhd \omega(j)$ or $\omega(i)=\omega(j)$ and $a(i)<a(j)$.  Let $J_2$ consist of those $j>i$ with $\omega(i)=\omega(j)$ and $a(i)\geq a(j)$.  Let $J_1$ consist of those $j>i$ with $\omega(i)\nunlhd \omega(j)$. 

Note that each $J_l$ is convex with respect to the coarse partial order so the combed $Y_{(i,0]}$ can be written as
\[
Y_{(i,0]} = (Y_{J_1}\rightarrow Y_{J_2}\rightarrow Y_{J_3})
\]
(with a long arrow connecting $Y_{J_1}$ to $Y_{J_3}$) where $Y_{J_l}$ contains those direct summands $\ts^{a(j)} X_j$ with $j\in J_l$.

We claim that $\Hom(X_i,Y_{J_1})\simeq 0$.  To prove this claim, consider the chain objects $(X_i)^k$ of $X_i$.  We have $\Hom((X_i)^k, X_j)\simeq 0$ for all $j\in J_1$ since $\omega(i)\nunlhd \omega(j)$.  Since $(X_i)^k$ is an object of $\AS$, it is compact when regarded as an object of $\Ch^-(\AS)$. Thus
\[
\Hom((X_i)^k, \bigoplus_{j\in J_1} X_j)= \bigoplus_{j\in J_1} \Hom((X_i)^k, X_j)\simeq 0.
\]
Now, homological perturbation tells us $\Hom((X_i)^k,Y_{J_1})\simeq 0$; this applies since the indexing set $J_1$ satisfies the ascending chain condition.  Finally, $\Hom(X_i,Y_{J_1})$ is a one-sided twist (indexed by $\Z_{\geq 0}$ with the opposite partial order) of $\prod_{k\geq 0} \Hom((X_i)^{-k},Y_{J_1})$, which is contractible by another application of homological perturbation.

Next, we claim that $\Hom(\ts^{a(i)} X_i,Y_{J_2})$ has no homology in degree 1.  Indeed, for each $k< 0$ the chain group $(X_i)^k$ satisfies $\Hom((X_i)^k, X_j)\simeq 0$ by hypothesis.  It follows that 
\[
\Hom(\ts^{a(i)} X_i,\ts^{a(j)} X_j) \cong \ts^{a(j)-a(i)} \Hom(X_i,X_j) \simeq \ts^{a(j)-a(i)} \Hom((X_i)^0,X_j),
\]
which is supported in cohomological degrees $\leq 0$ since $a(i)\leq a(j)$ from the way $J_2$ is defined. Two applications of homological perturbation tell us that $\Hom(\ts^{a(i)} X_i,Y_{J_2})$ is homotopy equivalent to a complex supported in cohomological degrees $\leq 0$.  In particular, there is no homology in degree 1.

Combining the above two paragraphs, we obtain that the inclusion $Y_{J_3}\hookrightarrow Y_{(i,0]}$ induces an isomorphism $H_1(\Hom(\ts^{a(i)}X_i,Y_{J_3}))\rightarrow H_1(\Hom(\ts^{a(i)}X_i, Y_{(i,0]}))$.  It follows that the connecting differential $\partial'$ in \eqref{eq:X [i,0]} is homotopic to an element $\partial''=\partial'+d(h)$ which has components $X_i\rightarrow X_j$ only if $j\in J_3$ (i.e.~$j$ is greater than $i$ in the coarse partial order). Therefore $\phi_{[i,0]}:=\smMatrix{\id & 0\\ h & \id}\circ \smMatrix{\id &0\\0& \phi_{(i,0]}}$  gives an isomorphism
\[
X_{[i,0]} = (\ts^{a(i)}X_{i}\buildrel \partial\over \rightarrow X_{[i,0]})\rightarrow (\ts^{a(i)}X_{i}\buildrel \partial''\over \rightarrow Y_{(i,0]})=:Y_{[i,0]}.
\]
This defines the desired combing of $X_{[i,0]}$ and completes the inductive proof that $X$ can be combed.
\end{proof}

\section{Chebyshev systems}
\label{sec:model-color}

\begin{definition}
Let $R$ be a ring.  We say that a family of elements $v_k\in R$ ($k\in \Z_{\geq 0}$) \emph{satisfies the Chebyshev II recursion} if $v_0=1$ and $v_{n-1}v_1 = v_n + v_{n-2}$ for $n\geq 2$.
\end{definition}
\begin{remark}
In such a family, the  element $v_1$ determines all of the elements $v_2,v_3,\ldots$. 
\end{remark}

The Chebyshev II recursion is intimately tied with the representation theory of quantum $\LieSL$, as illustrated in the following examples.
\begin{example}
If $R=\Z[q,q\inv]$ then $v_n:=[n]:= q^{n-1}+q^{n-3}+\cdots+q^{1-n}$ satisfies the Chebyshev II recursion.
\end{example}
\begin{example}
If $R=K_0(\TLC)$ then the classes of simple objects $v_n:=[V_n]$ satisfy the Chebyshev II recursion because of the Clebsch-Gordan rule:
\[
V_{n-1}\otimes V_1 \cong V_{n}\oplus V_{n-2}.
\]
\end{example}

We now turn our attention to a categorical analogue of the Chebyshev II recursion.  So we let $(\cat C,\otimes, I)$ be a monoidal category.  Naively, we might say that a collection of objects $V_n\in \cat$ ($n\in \Z_{\geq 0}$) \emph{satisfies the categorified Chebyshev II recursion} if $V_0\cong \one$ and $V_{n-1}\otimes V\cong V_n \oplus V_{n-2}$ for $n\geq 2$.  However this is much too naive; it is much better to have knowledge of the maps involved in such an isomorphism, and include them as part of the structure.  This is accomplished in Definition \ref{def:homological-model} below. 

\begin{definition}\label{def:selfdual}
Let $(\cat C,\otimes, I)$ be a monoidal category.  A \emph{self-dual object} (or \emph{self-biadjoint}) in $V$ is a triple $(V,\ev,\coev)$ where  $V\in \cat{C}$ is an object and $\coev\colon \one\leftrightarrow V\otimes V\colon \ev$ are maps satisfying
\[
(\id_V\otimes \ev)\circ(\coev\otimes \id_V) = \id_V = (\ev\otimes \id_V)\circ(\id_V\otimes \coev).
\]
\end{definition}

\begin{definition}\label{def:homological-model}
	Let $(\cat C,\otimes, I)$ be a monoidal category.   A~\emph{Chebyshev system} in $\cat{C}$ is the data of:
	\begin{itemize}
		\item a self-dual object $(V,\ev,\coev)$ in $\cat{C}$,
		\item a~family of complexes $V^{(n)}\in \Ch^-(\cat C)$ with $V^{(0)} = \one$ and $V^{(1)}=V$,
		\item chain maps $\pi^{(n)} \colon V^{(n-1)}\otimes V \to V^{(n)}$ for $n\geqslant 1$, 
	\end{itemize}
	such that $\pi^{(1)}\colon \one\otimes V\to V$ is the~canonical isomorphism and there is a~distinguished triangle
	\begin{equation}\label{model:triangle}
		V^{(n-2)} \xrightarrow{\ \iota^{(n-2)}\ }
		V^{(n-1)}\otimes V \xrightarrow{\ \pi^{(n)}\ }
		V^{(n)} \longrightarrow
		V^{(n-2)}[1]
	\end{equation}
	in the~homotopy category $\mathcal{K}^-(\cat C)$, where $\iota^{(n-2)}$ is defined by commutativity of the triangle
	\begin{equation}\label{model:iota-vs-pi}
		\begin{tikzpicture}[x=8em,y=12ex]
			\node (V1) at (0,1) {$V^{(n-2)}$};
			\node (V3) at (1,1) {$V^{(n-2)}\otimes V\otimes V$};
			\node (V2) at (1,0) {$V^{(n-1)}\otimes V$};
			\path[->,>=stealth',shorten >=1pt,auto,node distance=1.8cm]
				(V1) edge node {$\scriptstyle\id\otimes\coev$} (V3)
				(V3) edge node {$\scriptstyle\pi^{(n-1)}\otimes\id$} (V2)
				(V1) edge node {$\scriptstyle\iota^{(n-2)}$} (V2);
		\end{tikzpicture}
	\end{equation}
\end{definition}

\begin{remark}
Our interest in Chebyshev systems stems from their expected relation to colored  $\LieSL$ link homology.  Specifically we regard Chebyshev systems as models for colored $\LieSL$ link homology of the unknot.  Under an appropriate cabling functor, one hopes to obtain models for colored $\LieSL$ link homology of other knots as well.
\end{remark}

The~following result shows Chebyshev systems are determined up to homotopy equivalence by $V^{(1)}$.

\begin{theorem}\label{thm:uniqueness-of-model}
	Let $(V^{(n)},\pi^{(n)})$ and $(V^{\prime\,(n)}, \pi^{\prime\,(n)})$
	be two Chebyshev systems in a~category $\cat C$, with $V^{(1)}\cong V^{\prime\, (1)}$ as self-dual objects.  Then there are homotopy equivalences $\theta^{(n)}\colon V^{(n)} \xrightarrow{\ \simeq\ }
	V^{\prime\,(n)}$ that form isomorphisms of distinguished triangles
	\begin{equation}\label{eq:isom-of-triangles}
		\begin{tikzpicture}[x=8em,y=12ex]
			\node (V1) at (0,1) {$V^{(n-2)}$};
			\node (V2) at (1,1) {$V^{(n-1)}\otimes V$};
			\node (V3) at (2,1) {$V^{(n)}$};
			\node (V4) at (3,1) {$V^{(n-2)}[1]$};
			\node (W1) at (0,0) {$V^{\prime\,(n-2)}$};
			\node (W2) at (1,0) {$V^{\prime\,(n-1)}\otimes V$};
			\node (W3) at (2,0) {$V^{\prime\,(n)}$};
			\node (W4) at (3,0) {$V^{\prime\,(n-2)}[1]$};
			\path[->,>=stealth',shorten >=1pt,auto,node distance=1.8cm]
				(V1) edge node {$\scriptstyle\iota^{(n-2)}$} (V2)
				(V2) edge node {$\scriptstyle\pi^{(n)}$} (V3)
				(V3) edge (V4)
				(W1) edge node {$\scriptstyle\iota^{\prime\,(n-2)}$} (W2)
				(W2) edge node {$\scriptstyle\pi^{\prime\,(n)}$} (W3)
				(W3) edge (W4)
				(V1) edge node {$\scriptstyle\theta^{(n-2)}$} (W1)
				(V2) edge node {$\scriptstyle\theta^{(n-1)}\otimes\id$} (W2)
				(V3) edge node {$\scriptstyle\theta^{(n)}$} (W3)
				(V4) edge node {$\scriptstyle\theta^{(n-2)}[1]$} (W4);
		\end{tikzpicture}
	\end{equation}
\end{theorem}
\begin{proof}
	We prove the~thesis by induction on $n$.  Since $V^{(0)}=V^{\prime\, (0)}=\one$, we take $\theta^{(0)}=\id_{\one}$.  Set $V:=V^{(1)}$ and $V':=V^{\prime\, (1)}$.  We denote by $\coev$ and $\coev'$ the coevaluation morphisms for $V$ and $V'$.  Let $\theta^{(1)}\colon V\to V'$ be an isomorphism of self-dual objects, so in particular
	\[
	(\theta^{(1)}\otimes \theta^{(1)})\circ\coev = \coev'.
	\]
	Thus,  for $n=2$ the left square in \eqref{eq:isom-of-triangles} commutes (in this case, $\iota^{(0)}=\coev$ and $\iota^{\prime\,(0)}=\coev'$.  It follows (by properties of distinguished triangles) that there is a map $\theta^{(2)}$ which completes $\theta^{(0)}$ and $\theta^{(1)}\otimes\theta^{(1)}$ to a morphism of triangles.  By the five lemma for triangulated categories, $\theta^{(2)}$ is a homotopy equivalence. 
	
	Now,  let $n\geq 3$ be given and assume by induction that $\theta^{(k)}$ is constructed for $k<n$.  Note that if there exists $\theta^{(n)}$ which fits into a morphism of distinguished triangles  \eqref{eq:isom-of-triangles}, then $\theta^{(n)}$ is automatically a homotopy equivalence, again by the five lemma.  To prove the~existence
of $\theta^{(n)}$ we have to show that the~left square of \eqref{eq:isom-of-triangles}  commutes up to homotopy. To see this, consider the following diagram:
\[
\begin{tikzpicture}[x=8em,y=12ex]
\node (A1) at (0,1) {$V^{(n-2)}$};
\node (A2) at (1.3,1) {$V^{(n-2)}\otimes V\otimes V$};
\node (A3) at (2.6,1) {$V^{(n-1)}\otimes V$};
\node (B1) at (0,0) {$V^{\prime\,(n-2)}$};
\node (B2) at (1.3,0) {$V^{\prime\, (n-2)}\otimes V'\otimes V'$};
\node (B3) at (2.6,0) {$V^{\prime\,(n-1)}\otimes V'$};
\path[->,>=stealth',shorten >=1pt,auto,node distance=1.8cm]
(A1) edge node {$\scriptstyle\id\otimes \coev$} (A2)
(A2) edge node {$\scriptstyle\pi^{(n-1)}\otimes \id$} (A3)
(B1) edge node {$\scriptstyle\id\otimes \coev'$} (B2)
(B2) edge node {$\scriptstyle\pi^{\prime\,(n-1)}\otimes \id$} (B3)
(A1) edge node {$\scriptstyle\theta^{(n-2)}$} (B1)
(A2) edge node {$\scriptstyle\theta^{(n-2)}\otimes \theta^{(1)}\otimes \theta^{(1)}$} (B2)
(A3) edge node {$\scriptstyle\theta^{(n-1)}\otimes\theta^{(1)}$} (B3);
\end{tikzpicture}
\]
The left square commutes since $\theta^{(1)}$ is a morphism of self-dual objects, and the right square commutes by commutativity of the middle square in
\eqref{eq:isom-of-triangles} with $n$ lowered by 1.  The composition of the horizontal arrows are $\iota^{(n-2)}$ and $\iota^{\prime\, (n-2)}$.  This proves that the left square in \eqref{eq:isom-of-triangles} commutes, which gives us existence of $\theta^{(n)}$.
\end{proof}

\begin{remark}
There is a~notion of a~dual Chebyshev system,
obtained from the one in Definition~\ref{def:homological-model} by reversing
the~direction of all arrows and by replacing $\Ch^-$ with $\Ch^+$.  Note that if $\{V^{(n)},\pi^{(n)}\}$ is a Chebyshev system 
	in a~rigid category then $\{(V^{(n)})^\vee,(\pi^{(n)})^\vee\}$ is a dual Chebyshev system , where
	${(-)}^{\vee}$ denotes the~dualization functor. Explicitly,
	$(-)^{\vee}\colon\Ch^-(\cat C)\rightarrow\Ch^+(\cat C)$
	negates the~homological grading and transforms differentials into
	their duals.
\end{remark}

In what follows we discuss several Chebyshev systems, in $\TLC$ and related categories.  By Theorem \ref{thm:uniqueness-of-model} these systems are all equivalent in $\Ch^-(\Kar(\TLC))$.

%% ============================================================================
%%  Jones--Wenzl model
%% ============================================================================
\subsection{The~Jones--Wenzl system}
\label{ssec:JW-model}
In this section we take $\cat{C}:=\Kar(\TLC)$, the idempotent completion of the Temperley--Lieb category $\TLC$.
Recall that $\TLC$ contains the~Jones--Wenzl projectors
$p_n$, each representing a projection of $V_1^{\otimes n}$ onto its simple $(n+1)$-dimensional
quotient $V_n$.
We shall abuse notation and denote by $V$ also
the object $1\in\TL$. Because the~Jones--Wenzl projectors satisfy the identity
\[
	p_{n} \circ (p_{n-1} \otimes \id_V) = p_{n} = (p_{n-1}\otimes\id_V) \circ p_{n},
\]
we have well-defined morphisms in $\Kar(\TLC)$
\begin{gather*}
	\begin{tikzpicture}[x=7em]
		\node[anchor=mid east] (n0) at (0,0) {$\im p_{n-1}\otimes V$};
		\node[anchor=mid west] (n1) at (1,0) {$\im p_{n}$};
		\path[->,>=stealth',node distance=1.8cm,transform canvas={yshift= 2pt}] (n0) edge (n1);
		\path[->,>=stealth',node distance=1.8cm,transform canvas={yshift=-2pt}] (n1) edge (n0);
		\begin{scope}[x=7mm,y=7mm]
			\path (n0) -- (n1)
				node[midway,above] {\scriptsize $\pi^{(n)} = \JWpict[baseline=3mm]{5}[n]$}
				node[midway,below] {\scriptsize $\rho^{(n)} = \JWpict[baseline=3mm]{5}[n]$};
		\end{scope}
	\end{tikzpicture}
	\\[2ex]
	\begin{tikzpicture}[x=7em]
		\node[anchor=mid east] (n0) at (0,0) {$\im p_{n-2}$};
		\node[anchor=mid west] (n1) at (1,0) {$\im p_{n-1}\otimes V$};
		\path[->,>=stealth',node distance=1.8cm,transform canvas={yshift= 2pt}] (n0) edge (n1);
		\path[->,>=stealth',node distance=1.8cm,transform canvas={yshift=-2pt}] (n1) edge (n0);
		\begin{scope}[x=8mm,y=8mm]
			\path (n0) -- (n1)
				node[midway,above] {\scriptsize $\iota^{(n-2)} = \JWpict[baseline=3mm]{3+cup}[n-1]$}
				node[midway,below] {\scriptsize $\kappa^{(n-2)} =\frac{[n-1]}{[n]} \JWpict[baseline=3mm]{3+cap}[n-1]$};
		\end{scope}
	\end{tikzpicture}
\end{gather*}

\begin{theorem}
	The~data $\{\im p_n, \pi^{(n)}\}_{n\geq 0}$ is a~Chebyshev system.
	in $\Kar(\TLC)$.
\end{theorem}
\begin{proof}
	We must construct a~homotopy equivalence $\varphi^{(n)} \colon \im p_{n} \to \Cone(\iota^{(n-2)})$,
	such that $\varphi^{(n)}\circ\pi^{(n)}$ is homotopic to the~canonical inclusion $\im p_{n-1}\otimes V\rightarrow \Cone(\iota^{(n-2)})$.
	For that consider Figure \ref{fig:JW-model}.
	\begin{figure}[t]
		\begin{tikzpicture}[x=10em,y=12ex]
			\node   (im) at (0,1) {$\im p_{n}$};
			\node (cone) at (0,0) {$\Cone(\iota^{(n-2)})$};
			\node   (0) at (1,1) {$\mathllap{\smash{\Big(}\quad}0$};
			\node   (n) at (2,1) {$\im p_{n}\mathrlap{\quad\smash{\Big)}}$};
			\node  (n2) at (1,0) {$\mathllap{\smash{\Big(}}\im p_{n-2}$};
			\node  (n1) at (2,0) {$\im p_{n-1} \otimes V\mathrlap{\smash{\Big)}}$};
			\path (im.center) --  (0.center) node[midway] {$=$}
			    (cone.center) -- (n2.center) node[midway] {$=$};
			\path[->,>=stealth',node distance=1.8cm,transform canvas={xshift=-2pt}]
				(im) edge node[left] {$\scriptstyle\varphi^{(n)}$} (cone)
				 (n) edge node[left] {$\scriptstyle \rho^{(n)}$} (n1);
			\path[->,>=stealth',node distance=1.8cm,transform canvas={xshift=2pt}]
				(cone) edge node[right] {$\scriptstyle\bar\varphi^{(n)}$} (im)
				  (n1) edge node[right] {$\scriptstyle \pi^{(n)}$} (n);
			\path[->,>=stealth',node distance=1.8cm]
				(0) edge[<->] node[left] {$\scriptstyle 0$} (n2)
				(0) edge node[above] {$\scriptstyle 0$} (n)
				(n2) edge node[below] {$\scriptstyle d=\iota^{(n-2)}$} (n1)
				(n1) edge[dashed,out=165,in=15] node[above] {$\scriptstyle h= \kappa^{(n-2)}$}(n2);
		\end{tikzpicture}
		\caption{The~homotopy equivalences between $\im p_{n+2}$ and $\Cone(\iota^{(n)})$.}
		\label{fig:JW-model}
	\end{figure}
	Clearly, the~composition $\bar\varphi^{(n)}\circ\varphi^{(n)}$ is the~identity on $\im p_{n+2}$.
We compute
	\begin{align*}
		dh &= \frac{[n-1]}{[n]}\JWpict{5-tail}[n-1]
		    = \JWpict{4+1}[n-1] - \JWpict{5}[n]
		    = \id_{\im p_{n-1}\otimes V} - \rho^{(n)}\circ\pi^{(n)},
	\\[2ex]
		hd &= \frac{[n-1]}{[n]} \JWpict{4-tr}[n-1]
		    = \frac{[n-1]}{[n]} \left( \JWpict{3+loop}[n-2] - \frac{[n-2]}{[n-1]} \JWpict{3}[n-2] \right) \\
		   &\qquad\qquad= \frac{[n-1][2] - [n-2]}{[n]} \JWpict{3}[n-2]
		    = \JWpict{3}[n-2] = \id_{\im p_{n-2}}.
	\end{align*}
	Hence, $h$ is a homotopy between $\varphi^{(n)}\circ\bar\varphi^{(n)}$ and the~identity on the~cone complex.
	Note that the~above computation shows also that $h$ is a homotopy between $\varphi^{(n)}\circ\pi^{(n)}$
	and the~canonical inclusion of $\im p_{n-1}\otimes V$ into the~cone complex.
\end{proof}

%% ============================================================================
%%  Khovanov model
%% ============================================================================
\subsection{The Khovanov system}
\label{ssec:kh-model}

In this section we take $\cat{C}:=\mathcal{K}^b(\TLC)$, the bounded homotopy category of complexes over the Temperley--Lieb category.  Inside $\cat{C}$ we shall construct a Chebyshev system which is suggested by Khovanov's categorification of the~colored Jones
polynomial.

For that let us review the~complex $\VC_n$ assigned to the~unknot colored by
the~simple $(n+1)$-dimensional representation $V_n$, which we shall interpret
as an~object of $\Ch^-(\TLC)$.

Recall that
\[
	[V_n]=\sum_{k=0}^{\lfloor\frac{n}{2}\rfloor}(-1)^k \binom{n-k}{k} \left[V_1^{\otimes(n-2k)}\right],
\]
in the~Grothendieck group of the~representation category of quantum $\LieSL$.  Following Khovanov, we will construct a family of complexes $\VC_n\in \Ch^-(\TLC)$ whose Euler characteristics satisfy the same relation as above.

The~binomial coefficient $\binom{n-k}{k}$ counts the~number of ways of selecting $k$ disjoint
pairs of neighbored dots from $n$ dots placed on a~horizontal line.
We call such choices \emph{$k$-pairings} and we denote by $I_{k,n}$ the~set of all such $k$-pairings.
Each $s\in I_{k,n}$ represents an~object $\VC_s \in \TLC$ that consists of $n-2k$ dots.

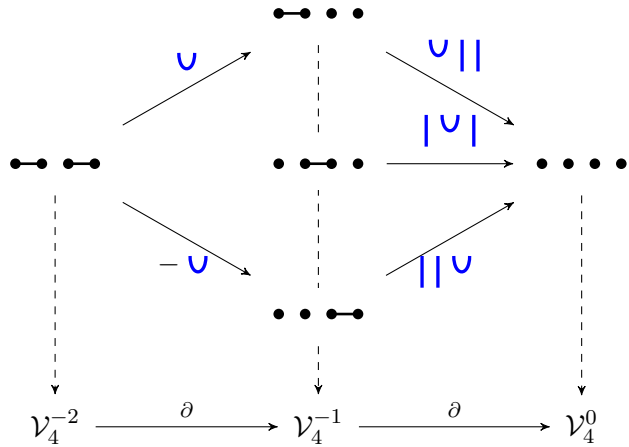
\begin{figure}[h]
	\begingroup
	% #1 = position of the left point
	% #2 = the number of points
	% #3 = a list of left endpoints for pairings
	\def\pairing#1#2{\hbox{\begin{tikzpicture}[x=10pt,y=10pt]
		\fill[white] (#1,1) ++(0.5,0) rectangle (0.5,-1);
		\foreach \x in {1,...,#1} \fill (\x,0) circle[radius=2pt];
		\foreach \x in {#2} \draw[line width=1pt] (\x,0) -- ++(1,0);
	\end{tikzpicture}}}%
	\def\morphism#1#2{\hbox{\begin{tikzpicture}[x=6pt,y=8pt,baseline=(ref.base)]
		\ifnum #1>0\relax
			\foreach \x in {-#1,...,-1} \draw[V1,-] (\x,0) -- (\x,2);
		\fi
		\ifnum #2>0\relax
			\foreach \x in {1,...,#2} \draw[V1,-] (\x,0) ++(1,0) -- ++(0,2);
		\fi
		\draw[V1,-] (0,2) .. controls ++(0,-1.25) and ++(0,-1.25) .. (1,2);
		\node[anchor=base] (ref) at (0.5,1) {$\phantom x$};
	\end{tikzpicture}}}%
	\begin{tikzpicture}[x=35mm,y=2cm]
		%% chain objects
		\node (V2) at (0,-1.75) {$\VC_4^{-2}$};
		\node (V1) at (1,-1.75) {$\VC_4^{-1}$};
		\node (V0) at (2,-1.75) {$\VC_4^0$};
		%% vertical grouping lines
		\path[->,>=stealth',shorten >=2pt,dashed]
			(0,0) edge (V2)
			(1,1) edge (V1)
			(2,0) edge (V0);
		%% pairing
		\node[anchor=center] (p13) at (0,0) {\pairing4{1,3}};
		\node[anchor=center] (p1) at (1, 1) {\pairing4{1}};
		\node[anchor=center] (p2) at (1, 0) {\pairing4{2}};
		\node[anchor=center] (p3) at (1,-1) {\pairing4{3}};
		\node[anchor=center] (p) at (2, 0) {\pairing4{}};
		%% edge maps and differential
		\path[->,>=stealth',shorten >=2pt,shorten <=2pt]
			(p13) edge node[above] {$\morphism00$} (p1)
			(p13) edge node[below] {$-\morphism00$} (p3)
			(p1) edge node[above] {\morphism02} (p)
			(p2) edge node[above] {\morphism11} (p)
			(p3) edge node[below] {\morphism20} (p)
			(V2) edge node[above] {$\scriptstyle\partial$} (V1)
			(V1) edge node[above] {$\scriptstyle\partial$} (V0);
	\end{tikzpicture}
\endgroup
	\caption{The~graph $\Gamma_4$ and the complex $\VC_4$.}
	\label{fig:Gn}
\end{figure}

Let $\Gamma_n$ be a~graph, the~vertices of which are $k$-pairings of $n$ dots for all $k\geq 0$.
Two vertices of $\Gamma_n$ are connected by a~directed edge if $s=s'\cup\{p\}$
for a~single pair $p$.
We decorate a~vertex $s$ with the~object $\VC_s$ and an~edge $s\rightarrow s'$
with the~morphism
\[
	\cup_{s',s} := \JWpict{cup}[{i}],
\]
where $i$ counts unpaired dots in $s$ that lie to the left of the~pairing that does not belong
to $s'$. An~example for $n=4$ is shown in Figure~\ref{fig:Gn}.
Note that $\Gamma_n$ is a~commuting diagram in $\TLC$.
In order to make all squares anticommute, scale $\cup_{s',s}$ by $(-1)^{(s,s')}$,
where $(s,s')$ is the~number of pairs in $s$ that lie to the~right of the unique pair in
$s\setminus s'$.
The~complex $\VC_n$ is the~result of collapsing the~graph
by taking direct sums of $k$-pairings with the same $k$:
\[
	\VC_n^{-k} := \bigoplus_{s\in I_{k,n}} \VC_s
\]
and the~differential $\delta^{-k}\colon\VC_n^{-k}\rightarrow\VC_n^{-k+1}$ is given by the matrix with entries
\[
	\delta^{-k}_{s', s}:=\begin{cases}
		(-1)^{(s,s')}\cup_{s',s} & \textrm{if there is an~edge $s\rightarrow s'$, and}\\
		0& \textrm{otherwise}.
	\end{cases}
\]

We may identify $\VC_n\otimes V$ with the~subcomplex of $\VC_{n+1}$ spanned by those $k$-pairings in which the~right most dot does not belong to a~pair.

\begin{theorem}
	The data $\{\VC_n,\pi^{(n)}\}_{n\geq 0}$ is a Chebyshev system in $\mathcal{K}^b(\TLC)$, where $\pi^{(n)}\colon \VC_{n-1}\otimes V\rightarrow \VC_n$ is the inclusion of a subcomplex. 
\end{theorem}
\begin{proof}
	The~graph $\Gamma_n$ has two types of vertices, corresponding to two types of pairings:
	those in which the~rightmost dot belongs to a~pair, and those for which it does not.
	Vertices of the~second type generate the~subcomplex $\VC_{n-1}\otimes V$,
	whereas vertices of the~first type generate a~quotient complex that is isomorphic to $\ts^{-1}\VC_{n-2}$.
	We can thus write $\VC_n$ as a mapping cone
	\[
		\VC_n = \Cone\left(\VC_{n-2}\rightarrow\VC_{n-1}\otimes V\right),
	\]
	where the chain map in the~cone is induced by all morphisms $\cup_{s',s}$ that correspond
	to edges in $\Gamma_n$ connecting vertices of the~two different types---each is given
	by the~diagram with a~cup on its right edge.
	Explicitly, this chain map can be written as $\id_{\VC_{n-2}} \otimes \coev_V$ (with codomain restricted to $\VC_{n-1}\otimes V$) and it coincides
	with $(\pi^{(n-1)}\otimes \id_V)\circ (\id_{\VC_{n-2}} \otimes \coev_V)$.  This shows that $\VC_{n-2}$, $\VC_n$, and $\VC_{n-1}\otimes V$ are related by a distinguished triangle as in Definition \ref{def:homological-model}.
\end{proof}

%% ============================================================================
%%  Cooper--Krushkal model
%% ============================================================================
\subsection{The Cooper--Krushkal system}
\label{ss:CK model}
One of the main goals of this paper is to show how the Cooper--Krushkal categorified Jones--Wenzl idempotents yield a Chebyshev system.

\begin{definition}
Recall that the~\emph{width} of a~flat tangle $T$ is the~minimal number of points on
a horizontal section of $T$.  In other words, the width of an $(m,n)$-tangle $T$ is the minimal $k$ such that $T$ can be written $T= T'\circ T''$ with $T'$ an $(m,k)$-tangle and $T''$ a $(k, n)$-tangle. 

If $X\in \Ch^-(\ccat{BN}_{m,n})$ then we say that $X$ \emph{has width $\leq k$ up to homotopy} if $X$ is homotopy equivalent to a complex $Y\in \Ch^-(\ccat{BN}_{m,n})$ whose chain groups consist of sums of tangles with width $\leq k$.
\end{definition}

Hereafter fix a~natural number $n$ and the~following tangles
\[
	\one_n = \JWpict{id*}
		\qquad\,
	\cup_i = \JWpict{cup*}
		\qquad\,
	\cap_i = \JWpict{cap*}
		\qquad\,
	B_i = \JWpict{cup-cap*}
\]
for $i=1,\dots,n-1$.   We say that a complex $X\in \Ch^-(\ccat{BN}_{n,n})$ \emph{kills turnbacks} if $\cap_i\star X\simeq 0\simeq X\star \cup_i$ for all $1\leq i\leq n-1$.

\begin{theorem}[\cite{CK12a}]\label{thm:P}
There exists a complex $P_n\in \Ch^-(\BN_{n,n})$ and a chain map $\eta\colon\one_n\rightarrow P_n$ satisfying the following properties:
	\begin{enumerate}[label={(P\arabic*)},ref={P\arabic*}]
		\item $P_n$ kills turnbacks.
		\item $\Cone(\eta)$ has width $<n$ up to homotopy.
	\end{enumerate}
	These properties determine the~pair $(P_n, \eta)$ uniquely up to homotopy.
\end{theorem}

\begin{remark}  In fact, it can be shown \cite[Cor.~5.13]{H12a} that for any $X\in \Ch^-(\ccat{BN})$ we have $X\star \cup_i\simeq 0$ for all $1\leq i\leq n-1$ if and only if $\cap_i\star X\simeq 0$ for all $1\leq i\leq n-1$.  Thus, to show $P_n$ kills turnbacks it suffices to show that $P_n$ kills turnbacks from the right ($P_n\star \cup_i\simeq 0$ for $1\leq i\leq n-1$) or from the left ($\cap_i\star P_n\simeq 0$ for $1\leq i\leq n-1$).
\end{remark}

\begin{definition}\label{def:CK proj}
A \emph{Cooper--Krushkal projector on $n$ strands} is any pair $(P_n,\eta)$ satisfying axioms (P1) and (P2) of Theorem \ref{thm:P}.
\end{definition}
\begin{proposition}\label{prop:CKwidth}
If $X\in \Ch^-(\BN_{m,n})$ and $Y\in \Ch^-(\BN_{n,m})$ have width $<n$ up to homotopy then
\[
X\star P_n\simeq 0\simeq P_n\star Y.
\]
\end{proposition}

\begin{proof} This follows immediately from the definitions
and from homological perturbation theory.
\end{proof}

\begin{remark}
In practice, axiom (P2) is typically realized in the following way.  Suppose $P_n$ can be expressed as a chain complex of the form
\[
	\dots \rightarrow P_n^{-3}
	      \rightarrow P_n^{-2}
	      \rightarrow P_n^{-1}
			\rightarrow \one_n,
\]
where each $P_n^k$ for $k<0$ is spanned by tangles of width smaller than $n$, and let $\eta\colon \one_n\rightarrow P_n$ be the inclusion of the degree zero chain object.  Then
\[
\Cone(\eta) \  \ \simeq \ \ \dots \rightarrow P_n^{-3}
	      \rightarrow P_n^{-2}
	      \rightarrow P_n^{-1}
			\rightarrow 0
\]
and axiom (P2) is satisfied.
\end{remark}

The axioms imply that $P_n$ absorbs $P_{n-1}$ in the following sense.

\begin{lemma}
Let $(P_{n-1},\eta_{n-1})$ be a Cooper--Krushkal projector on $n-1$ strands.  Then the map
\[
P_n\cong \one_n\star P_n \xrightarrow{(\eta_{n-1}\sqcup \id_{\one_1})\star \id_{P_n}} (P_{n-1}\sqcup \one_1)\star P_n
\]
is a homotopy equivalence.
\end{lemma}
\begin{proof}
Recall that a chain map $f$ is a homotopy equivalence if and only if $\Cone(f)$ is contractible.  The cone of the map in the statement is isomorphic to $(\Cone(\eta_{n-1})\sqcup \one_1)\star P_n$ which is contractible since $\Cone(\eta_{n-1})$ has width $<n$ up to homotopy and $P_n$ kills turnbacks.  See \cite[Prop.~3.3]{CK12a} for details.
\end{proof}

\begin{proposition}\label{prop:rel unit}
There is a unique map $\nu\colon P_{n-1}\sqcup \one_1\rightarrow P_n$ up to homotopy with the property that $\nu\circ \eta_{n-1} \simeq \eta_n$.  
\end{proposition}
\begin{proof}
For existence we may define $\nu$ to be the composition of chain maps
\[
P_{n-1}\sqcup \one_1\cong (P_{n-1}\sqcup \one_1)\star \one_n \xrightarrow{\id\star \eta_n}(P_{n-1}\sqcup \one_1)\star P_n\simeq P_n
\]
(we leave it as an exercise to show $\nu\circ \eta_{n-1}\simeq \eta_n$).   Uniqueness and other properties of $\nu$ follow from  general facts about categorical idempotents \cite[Thm.~1.6]{Hog17a}.
\end{proof}

\begin{theorem}\label{thm:P as model}
Given $X\in \Ch^-(\ccat{BN}_{n,n})$, let $[X]\in \Ch^-(\hTr_q(\ccat{BN}))$ denote the quantum horizontal trace of $X$. 	There is a Chebyshev system $\{[P_n], \pi^{(n)}\}_{n\geq 0}$,
where
\[
\pi^{(n)} = [\nu] \colon [P_{n-1}\sqcup \one_1] \longrightarrow [P_{n}]
\]
and $\nu$ is as in Proposition \ref{prop:rel unit}.
\end{theorem}

In this theorem, it is assumed that $1-q^d$ is invertible for all integers $d>0$.
We will prove Theorem~\ref{thm:P as model} in the next section by using the following lemma:

\begin{lemma}\label{lemma:zeroendo} Let $f\colon P_n\rightarrow P_n$ be an endomorphism of quantum degree $d$,
and suppose $1-q^d$ is invertible. Then
the induced endomorphism $[f]\colon [P_n]\rightarrow[P_n]$ is null-homotopic.
\end{lemma}

\begin{proof} Let \[\theta\colon [P_n\star P_n]\longrightarrow [P_n\star P_n]\] be the chain isomorphism
whose components are given by sums of cyclicity isomorphisms $\theta_{P_n^j,P_n^i}\colon [P_n^i\star P_n^j]\rightarrow [P_n^j\star P_n^i]$. Then
the naturality of the cyclicity isomorphisms implies
\[
[\id\star f]\ =\ \theta\circ[f\star\id]\circ\theta^{-1}\ \simeq\ \theta\circ[\id\star f]\circ\theta^{-1}\ =\ q^d[f\star\id]\ \simeq\ q^d[\id\star f]
\]
where we have used that $\id\star f\simeq f\star\id$ for any endomorphism $f\colon P_n\rightarrow P_n$
(this is a general fact about unital categorical idempotents, see \cite[Lemma~3.18]{BD14}). Hence
$(1-q^d)[\id\star f]\simeq 0$, and so $[\id\star f]\simeq 0$. This also implies $[f]\simeq 0$ because
$(\eta_n\star\id)\circ f=(\id\star f)\circ (\eta_n\star\id)$ and $\eta_n\star\id$ is a homotopy equivalence.
\end{proof}

\section{The categorical idempotent model}
\label{sec:categorical-model}

The purpose of this section is to prove Theorem \ref{thm:P as model}.  To this,  we must carefully choose our preferred model for $P_n$.  We begin with some setup in \S \ref{ss:setup}.  Then in \S \ref{ss:2 periodic model} we establish our model for $P_n$ which manifests a sort of 2-periodicity similar to that of $P_2$.  
\subsection{Setup}
\label{ss:setup}
 It will be beneficial to set up some notation which will be used throughout.  Firstly, let
\[
s\colon \JWpict{5-tail}[n-1] \to \JWpict{4+1}[n-1]^{\star 2}\, \qquad s^\ast\colon  \JWpict{4+1}[n-1]^{\star 2}\to \JWpict{5-tail}[n-1]
\]
denote degree $\qs$ chain maps induced by saddle cobordisms.  By abuse we will also regard $s,s^\ast$ as morphisms
\[
\JWpict{5-tail}[n-1] \leftrightarrow \JWpict{4+1}[n-1]
\]
(via the homotopy equivalence $(P_{n-1}\sqcup\one_1)^{\star 2}\simeq P_{n-1}\sqcup\one_1$).  Note that since 
\[
s^\ast\circ s\colon \JWpict{5-tail}[n-1] \to \JWpict{5-tail}[n-1]
\]
is a sum of dots on the ``cup'' and ``cap''.

Recall also that in \cite{H14a}
distinguished endomorphisms $u_1,\ldots,u_n \in \End_{\Ch^-(\BN_n)}(P_n)$ were constructed.   The endomorphism $u_k$ generates the $\k$-module of endomorphisms of $P_n$ of degree $\ts^{2-2k}\qs^{2k}$ modulo homotopy, and is uniquely characterized up to homotopy and invertible scalar by this fact.

\noindent

Finally, let
\[
u_{n-1}^{\top}\, , \, u_{n-1}^{\bot} \colon  \JWpict{5-tail}[n-1] \to  \JWpict{5-tail}[n-1]
\]
denote the degree $\ts^{4-2n}\qs^{2n-2}$ morphisms given by the action of applying $u_{n-1}$ on the top or bottom factors of $P_{n-1}$, respectively.

\subsection{The 2-periodic model for $P_n$}
\label{ss:2 periodic model}

Our proof of Theorem \ref{thm:P as model}  uses the following description of $P_n$.

\begin{theorem}\label{thm:periodic P}
There is a model for $P_n$ of the form
\begin{equation}\label{eq:periodic P}
\left(\cdots \xrightarrow{u_{n-1}^{\top}-u_{n-1}^{\bot}}   \ts^{-3}\mathbbm{a}_3\JWpict{5-tail}[n-1] \xrightarrow{s^\ast\circ s} \ts^{-2}\mathbbm{a}_2 \JWpict{5-tail}[n-1] \xrightarrow{u_{n-1}^{\top}-u_{n-1}^{\bot}} \ts\inv \mathbbm{a}_1 \JWpict{5-tail}[n-1] \xrightarrow{s} \JWpict{4+1}[n-1]\right)
\end{equation}
(where we are omitting longer arrows which point to the right) in which the shifts satisfy:
\[
\mathbbm{a}_1:=\qs\, \qquad \mathbbm{a}_{2k}=\ts^{4-2n}\qs^{2n-2}\mathbbm{a}_{2k-1}\, \qquad \mathbbm{a}_{2k+1} = \qs^2\mathbbm{a}_{2k}
\]
for all $k\geq 1$.
Moreover, with respect to this expression, the endomorphism $u_n$ acts as
\begin{equation}\label{eq:periodic U}
\begin{tikzpicture}[baseline=0em]
\tikzstyle{every node}=[font=\small]
\node (A0) at (0,0) {$\JWpict{4+1}[n-1]$};
\node (A1) at (-2.5,0) {$\JWpict{5-tail}[n-1]$};
\node (A2) at (-5,0) {$\JWpict{5-tail}[n-1]$};
\node (A3) at (-7.5,0) {$\cdots $};
\node (B0) at (5,-3) {$\JWpict{4+1}[n-1]$};
\node (B1) at (2.5,-3) {$\JWpict{5-tail}[n-1]$};
\node (B2) at (0,-3) {$\JWpict{5-tail}[n-1]$};
\node (B3) at (-2.5,-3) {$\JWpict{5-tail}[n-1]$};
\node (B4) at (-5,-3) {$\JWpict{5-tail}[n-1]$};
\node (B5) at (-7.5,-3) {$\cdots$};
\path[->,>=stealth,shorten >=1pt,auto,node distance=1.8cm,
  thick]
(A3) edge node{} (A2)
(A2) edge node{} (A1)
(A1) edge node{$s$} (A0)
(B5) edge node{} (B4)
(B4) edge node{} (B3)
(B3) edge node{} (B2)
(B2) edge node{} (B1)
(B1) edge node{$s$} (B0)
(A0) edge node{$s$} (B2)
(A1) edge node{$\id$} (B3)
(A2) edge node{$\id$} (B4);
\end{tikzpicture}
\end{equation}
together with arrows pointing strictly to the right.
\end{theorem}

\begin{remark}
The expression \eqref{eq:periodic P} categorifies the right-hand side of the recursion~\eqref{eqn:JWrecursion} for the Jones--Wenzl idempotent.
\end{remark}

This theorem will be a consequence of the following.

\begin{lemma}\label{lemma:Qn four terms}
There is a complex of the form
\begin{equation}\label{eq:Qn four terms}
Q_n  = 
\left(\begin{tikzpicture}[baseline=0em]
\node (A0) at (-.3,0) {$\JWpict{4+1}[n-1]$};
\node (A1) at (-3,0) {$\ts\inv\mathbbm{a}_1\JWpict{5-tail}[n-1]$};
\node (A2) at (-7.3,0) {$\ts^{-2}\mathbbm{a}_2\JWpict{5-tail}[n-1]$};
\node (A3) at (-10.7,0) {$\ts^{-3}\mathbbm{a}_3\JWpict{4+1}[n-1]$};
\path[->,>=stealth',shorten >=1pt,auto,node distance=1.8cm,font=\small]
(A3) edge node[above,xshift=4pt] {$s^\ast$} (A2)
(A2) edge node {$u_{n-1}^{\top}-u_{n-1}^{\bot}$} (A1)
(A1) edge node {$s$} (A0)
(A3) edge[bend left=30] node {$\gamma$} (A0);
\draw[frontline,->,>=stealth',shorten >=1pt,auto,node distance=1.8cm]
(A2) to [bend right] node[below] {$h$} (A0);
\draw[frontline,->,>=stealth',shorten >=1pt,auto,node distance=1.8cm]
(A3) to [bend right] node[below] {$k$} (A1);
\end{tikzpicture}
\right)
\end{equation}
for the appropriate shifts $\mathbbm{a}_i$.  Specifically $\mathbbm{a}_1=\qs$, $\mathbbm{a}_2 = \ts^{4-2n}\qs^{2n-1}$, $\mathbbm{a}_3=\ts^{4-2n}\qs^{2n}$.  This complex kills turnbacks.
\end{lemma}
Compare \eqref{eq:Qn four terms} with the expression (4.10) in \cite{HogSym-GT}.

\begin{proof}
The existence of $h$ and $k$ will follow from the vanishing of $s\circ (u_{n-1}^{\top}-u_{n-1}^{\bot})$ and  $(u_{n-1}^{\top}-u_{n-1}^{\bot})\circ s^\ast$ up to homotopy.   We prove only the first of these, since the second is similar.  Consider the diagram
\[
\begin{tikzpicture}
\node (A1) at (0,0) {$\JWpict{5-tail}[n-1]$};
\node (A2) at (3,0) {$\JWpict{4+1}[n-1]^{\star 2}$};
\node (A3) at (6,0) {$\JWpict{4+1}[n-1]^{\star 2}$};
\node (B1) at (0,4) {$\JWpict{5-tail}[n-1]$};
\node (B2) at (3,4) {$\JWpict{4+1}[n-1]^{\star 2}$};
\node (B3) at (6,4) {$\JWpict{4+1}[n-1]^{\star 2}$};
\path[->,>=stealth',shorten >=1pt,auto,node distance=1.8cm,font=\small]
(A1) edge node[left] {$u_{n-1}^{\top}-u_{n-1}^{\bot}$} (B1)
(A2) edge node[left] {$u_{n-1}^{\top}-u_{n-1}^{\bot}$} (B2)
(A3) edge node {$0$} (B3)
(A1) edge node {$s$} (A2)
(A2) edge node {$\id$} (A3)
(B1) edge node {$s$} (B2)
(B2) edge node {$\id$} (B3)
;
\end{tikzpicture}
\]
The first square obviously commutes, and the second square commutes since $f\star \id\simeq \id \star f$ as endomorphisms of $P_{n-1}\star P_{n-1}$ for any endomorphism $f$ of $P_{n-1}$ (this is a general fact about unital categorical idempotents, see \cite[Lemma~3.18]{BD14}). This implies that $s\circ (u_{n-1}^{\top}-u_{n-1}^{\bot})\simeq 0$, which gives us the existence of $h$.  A similar argument gives us the existence of $k$.  Finally, the obstruction to constructing $\gamma$ is an endomorphism of $P_{n-1}\sqcup \one_1$ of degree $\ts^{3-2n}\qs^{2n}$.   Any such obstruction vanishes up to homotopy by an argument identical to that in \cite[Lemma~4.34]{HogSym-GT}.  This completes the construction of $Q_n$.
\end{proof}

Using this $Q_n$ we construct our 2-periodic model for $P_n$.
\begin{proof}[Proof of Theorem \ref{thm:periodic P}]
It is clear that splicing together infinitely many copies of $Q_n$ will result in a complex of the form \eqref{eq:periodic P}.  In more details, $Q_n$ comes equipped with an endomorphism $\partial$ of degree $\ts^{1-2n}\qs^{2n}$ which projects onto the left-most term of \eqref{eq:Qn four terms} and includes this as the right-most term. The pre-spliced complex is given by

\begin{equation}\label{eq:Qn prespliced}
\left(
\begin{tikzpicture}[baseline=0em]
\tikzstyle{every node}=[font=\small]
\node (a) at (0,0) {$Q_n$};
\node (b) at (2.5,0) {$\ts^{2-2n}\qs^{2n}Q_n$};
\node (c) at (6,0) {$\ts^{4-4n}\qs^{4n} Q_n$};
\node (d) at (9,0) {$\cdots $};
\path[->,>=stealth,shorten >=1pt,auto,node distance=1.8cm,
  thick]
(a) edge node {} (b)
(b) edge node {} (c)
(c) edge node {} (d);
\end{tikzpicture}
\right)
\end{equation}
which obviously carries an endomorphism $U'$ of degree $\ts^{2-2n}\qs^{2n}$ which just shifts the above expression one unit to the right.  By construction, $\Cone(U')\simeq Q_n$.  After cancelling contractible terms we obtain the honest spliced complex \eqref{eq:periodic P} with a chain endomorphism $U$ of the form \eqref{eq:periodic U} (this is an easy exercise in Gaussian elimination, but also note that the homotopy equivalence $\Cone(U)\simeq \Cone(U')\simeq Q_n$ puts serious restrictions on $U$).  Let us denote this spliced complex by $P_n'$.   Once we show that $P_n'\simeq P_n$, the fact that $U$ is a model for $\pm u_n$ will follow for degree reasons.

It is clear that $P_n'$ satisfies axiom (P2) from Theorem \ref{thm:P}, so to complete the proof we must show that $P_n'$ kills turnbacks, say, from the right.   The one-sided twisted complex \eqref{eq:Qn prespliced} satisfies a local ascending chain condition because $Q_n$ is bounded from above, so by homologial perturbation, $P_n'$ kills turnbacks if and only if $Q_n$ does.  We assume by induction that we have already established $P_{n-1}'\simeq P_{n-1}$.   

Consider the expression \eqref{eq:Qn four terms} for $Q_n$.
This complex annihilates $B_k$ for $1\leq k\leq n-2$ because $P_{n-1}\sqcup \one_1$ does.  So we must show that $B_{n-1}$ annihilates $Q_n$.  For this it will be useful to find a ``six-term'' expression for $Q_n$.  We expand the lower projector $P_{n-1}$ on the right in \eqref{eq:Qn four terms},  according to the recursion \eqref{eq:periodic P} for $P_{n-1}$, 
obtaining the %following
expression for $Q_n$
shown in Figure~\ref{fig:Qnexpanded}.
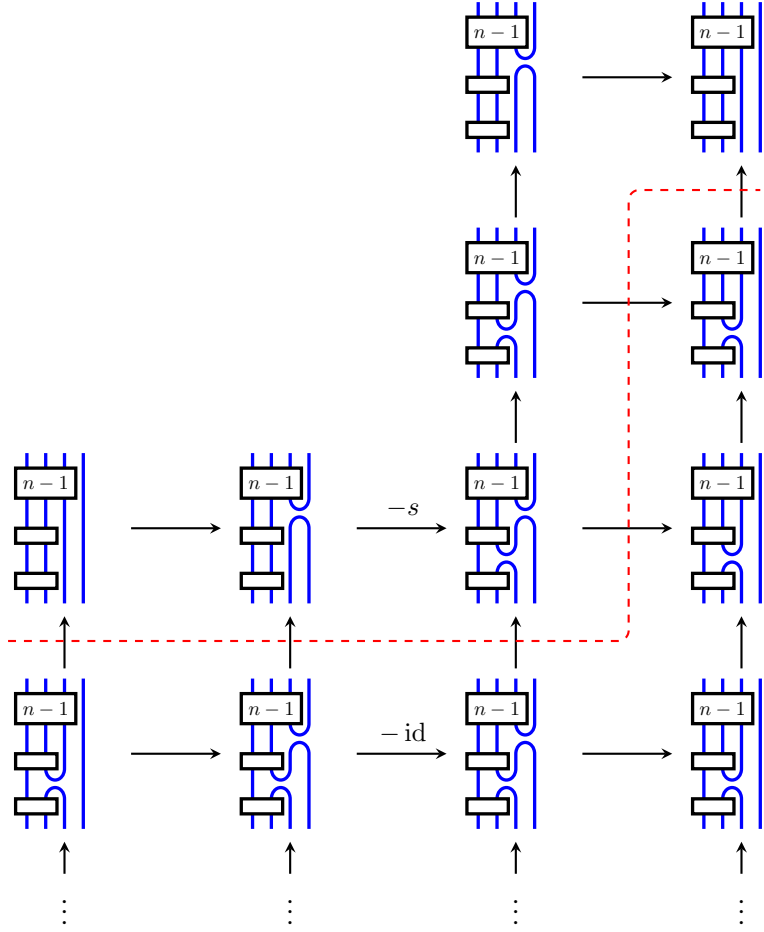
\begin{figure}
\begin{tikzpicture}[baseline=0em]
\tikzstyle{every node}=[font=\small]
\node (a0) at (0,0) {$\XYY$};
\node (a1) at (0,-3) {$\XYbY$};
\node (a2) at (0,-6) {$\XYbY$};
\node (a3) at (0,-9) {$\XYbY$};
\node (a4) at (0,-11) {$\vdots$};
\node (b0) at (-3,0) {$\XaYY$};
\node (b1) at (-3,-3) {$\XaYbY$};
\node (b2) at (-3,-6) {$\XaYbY$};
\node (b3) at (-3,-9) {$\XaYbY$};
\node (b4) at (-3,-11) {$\vdots$};
\node (c2) at (-6,-6) {$\XaYY$};
\node (c3) at (-6,-9) {$\XaYbY$};
\node (c4) at (-6,-11) {$\vdots$};
\node (d2) at (-9,-6) {$\XYY$};
\node (d3) at (-9,-9) {$\XYbY$};
\node (d4) at (-9,-11) {$\vdots$};
\path[->,>=stealth,shorten >=1pt,auto,node distance=1.8cm,
  thick]
(a1) edge node {} (a0)
(a2) edge node {} (a1)
(a3) edge node {} (a2)
(a4) edge node {} (a3)
(b1) edge node {} (b0)
(b2) edge node {} (b1)
(b3) edge node {} (b2)
(b4) edge node {} (b3)
(c3) edge node {} (c2)
(c4) edge node {} (c3)
(d3) edge node {} (d2)
(d4) edge node {} (d3)
(d3) edge node {} (c3)
(c3) edge node {$-\id$} (b3)
(b3) edge node {} (a3)
(d2) edge node {} (c2)
(c2) edge node {$-s$} (b2)
(b2) edge node {} (a2)
(b1) edge node {} (a1)
(b0) edge node {} (a0);
\draw[rounded corners,dashed,thick,red] (-9.75,-7.5) -- (-1.5,-7.5) -- (-1.5,-1.5) -- (0.25,-1.5);
\end{tikzpicture}
\caption{The complex $Q_n$ after expanding the lower projector.}\label{fig:Qnexpanded}
\end{figure}

In this figure, we have omitted longer arrows including arrows pointing strictly to the top right. The arrows labeled $-s$ and $-\id$ come from the expansion of $-u_{n-1}^{bot}$, while the arrows corresponding to $u_{n-1}^{top}$ would point to the top right.

For the purpose of this proof, the arrows coming form the maps $h$ and $\gamma$ in $Q_n$ are irrelevant,
but it will be important to make sure that all arrows in the expansion of the map $k$ point (non-strictly) to the top right.
This can be achieved by combing hairs, or directly as follows. From the proof of Lemma~\ref{lemma:Qn four terms},
we know that $k=s^*\circ k'$ where $s^*$ is a saddle cobordism and $k'$ is
a null-homotopy for the endomorphism $u:=u_{n-1}^{\top}-u_{n-1}^{\bot}$ of $P_{n-1}^{\star 2}\sqcup\one_1$.
If we now expand the lower projector in $P_{n-1}^{\star 2}\sqcup\one_1$, we obtain
the one-sided twisted complex shown in the leftmost column of Figure~\ref{fig:Qnexpanded}.
In this complex, all terms beyond the first one are contractible. Using homological
perturbation theory, we thus get a homotopy equivalence, denoted $F$, from this twisted complex to its leading term.
The inverse homotopy equivalence, denoted $G$, is given by inclusion, and there is a homotopy $H$
between the identity map of $P_{n-1}^{\star 2}\sqcup\one_1$ and $GF$ such that no arrows in $H$ 
point downwards. In the definition of $k$, we can now replace $k'$ by a pair of homotopies
$u\simeq GFu\simeq 0$ given by $Hu$ and $GFk'$. This yields a new map $k$,
whose components point (non-strictly) to the top right.

With this in mind, we can now cancel the contractible terms in Figure~\ref{fig:Qnexpanded} that
lie below the dashed line. After absorbing
$(P_{n-1}\sqcup \one_1)\star (P_{n-2}\sqcup\one_2)\simeq P_{n-1}\sqcup \one_1$, this yields a homotopy equivalence
\[
Q_n \ \simeq \ \left(
\begin{tikzpicture}[baseline=0em]
\tikzstyle{every node}=[font=\small]
\node (a) at (0,0) {\XY};
\node (b) at (2.5,0) {\XaY};
\node (c) at (5,0) {\XabY};
\node (d) at (7.5,0) {\XabY};
\node (e) at (10,0) {\XaY};
\node (f) at (12.5,0) {\XY};
\path[->,>=stealth,shorten >=1pt,auto,node distance=1.8cm,
  thick]
(a) edge node {} (b)
(b) edge node {} (c)
(c) edge node {} (d)
(d) edge node {} (e)
(e) edge node {} (f);
\end{tikzpicture}
 \right)
\]
in which all length one maps are saddle maps except for the map in the middle which is $u_{n-2}^{\top} - u_{n-2}^{\bot}$
(here we use that the components of $k$ point to the top right). 
The left half and right half each are killed by $-\star B_{n-1}$ (see for example \cite[Subs.~7.1.1]{CK12a}, so $Q_n$ is killed by $-\star B_{n-1}$ as well. This completes the proof.
\end{proof}

\subsection{The quantum trace of $P_n$}
\label{ss:trace of Pn}
In this section, we use the complex $[Q_n]$
to show that $[P_n]$ is homotopy equivalent to $V_n:=\im p_n$.
In particular, this will lead to a proof of 
Theorem~\ref{thm:P as model}.

\begin{lemma}\label{lemma:trace of Qn}
We have
\[
[Q_n] \ \simeq \ 
\left(\begin{tikzpicture}[baseline=0em]
\node (a) at (0,0){$\left[\JWpict{4+1}[n-1]\right] $ };
\node (b) at (-3,0){$\left[\JWpict{4-tr}[n-1]\right]$ };
\node (c) at (-6,0){$\left[\JWpict{4-tr}[n-1]\right]$};
\node (d) at (-9,0){ $\left[\JWpict{4+1}[n-1]\right]$};
\path[->,>=stealth',shorten >=1pt,auto,node distance=1.8cm,font=\small]
(b) edge node {$a$} (a)
(c) edge node {$[u_{n-1}]$} (b)
(d) edge node {$c$} (c);
\end{tikzpicture}
\right)
\]
where
we are omitting longer arrows, and where $a$ and $c$
are both induced by cyclicity isomorphisms and saddles.
\end{lemma}
\begin{proof}
We first observe that we have the following diagram which commutes up to homotopy:
\[
\begin{tikzpicture}
\node (a1) at (0,0){$\left[\JWpict{5-tail}[n-1]\right]$ };
\node (a2) at (0,-4){$\left[\JWpict{5-tail}[n-1]\right] $ };
\node (b1) at (4,0){$\left[\TXX\right]$};
\node (b2) at (4,-4){$\left[\TXX\right]$};
\node (c1) at (8,0){$\left[\TXX\right]$};
\node (c2) at (8,-4){$\left[\TXX\right]$};
\node (d1) at (12,0){$\left[\JWpict{4-tr}[n-1]\right]$};
\node (d2) at (12,-4){ $\left[\JWpict{4-tr}[n-1]\right]$};
\path[->,>=stealth',shorten >=1pt,auto,node distance=1.8cm,font=\small]
(a1) edge node[left] {$[u^{\top}_{n-1}-u^{\bot}_{n-1}]$}(a2)
(b1) edge node[left] {$[u^{\bot}_{n-1}- \qf^{2n-2}u^{\top}_{n-1}]$}(b2)
(c1) edge node[left] {$[u^{\top}_{n-1}- \qf^{2n-2}u^{\top}_{n-1}]$}(c2)
(d1) edge node[left] {$[(1- \qf^{2n-2})u_{n-1}]$}(d2)
(a1) edge node[label=below:$\theta$] {$\cong$} (b1)
(b1) edge node {$\id$} (c1)
(c1) edge node {$\simeq$} (d1)
(a2) edge node[label=below:$\theta$]  {$\cong$} (b2)
(b2) edge node {$\id$} (c2)
(c2) edge node {$\simeq$} (d2);
\end{tikzpicture}
\]
The composition of horizontal maps gives a homotopy equivalence \[\left[\JWpict{5-tail}[n-1]\right]\simeq \left[\JWpict{4-tr}[n-1]\right]\] 
Applying this homotopy equivalence
to the two middle terms in $[Q_n]$, we get 
\[
[Q_n]  \ \simeq \ 
\left(\begin{tikzpicture}[baseline=0em]
\node (a) at (0,0){$\left[\JWpict{4+1}[n-1]\right] $ };
\node (b) at (-3,0){$\left[\JWpict{4-tr}[n-1]\right]$ };
\node (c) at (-6,0){$\left[\JWpict{4-tr}[n-1]\right]$};
\node (d) at (-9,0){ $\left[\JWpict{4+1}[n-1]\right]$};
\path[->,>=stealth',shorten >=1pt,auto,node distance=1.8cm,font=\small]
(b) edge node {$a$} (a)
(c) edge node {$f$} (b)
(d) edge node {$c$} (c);
\end{tikzpicture}
\right)
\]
where $f=(1-\qf^{2n-2})[u_{n-1}]$. The lemma
now follows because $1-q^{2n-2}$ is assumed to be invertible.
\end{proof}

We are now ready to prove:

\begin{lemma}\label{lemma:PnVn}
$[P_n]\simeq V_n$ \end{lemma}

\begin{proof}
By Theorem~\ref{thm:periodic P}, we have an expansion
\[
\JWpict{4-tr}[n-1]  \simeq 
\left(\begin{tikzpicture}[baseline={0em-height("$\vcenter{}$")}]
\node (a) at (0,0){$\nYu$ };
\node (b) at (-3,0){$\nYbYu$ };
\node (c) at (-6,0){$\nYbYu$};
\node (d) at (-9,0){$\nYbYu$};
\node (e) at (-11,0){ $\cdots$};
\path[->,>=stealth',shorten >=1pt,auto,node distance=1.8cm,font=\small]
(b) edge node {$s$} (a)
(c) edge node {} (b)
(d) edge node {} (c)
(e) edge node{} (d);
\end{tikzpicture}
\right)
\]
in which the right-hand side is homotopic to
\begin{equation}\label{eqn:tracedstrandexpansion}
\begin{tikzpicture}[baseline={0em-height("$\vcenter{}$")}]
\node (ar) at (1.8,0){$\JWT$};
\node (a) at (0,0){$\JWT$ };
\node (b) at (-2.5,0){$\JWT$ };
\node (c) at (-5,0){$\JWT$};
\node (d) at (-7.5,0){$\JWT$};
\node (e) at (-9.5,0){ $\cdots$};
\path[->,>=stealth',shorten >=1pt,auto,node distance=1.8cm,font=\small]
(b) edge node {$\id$} (a)
(b) edge[bend right=35] node {} (ar)
(c) edge node {} (b)
(d) edge node {} (c)
(e) edge node{} (d);
\path
(a) edge[draw=none] node {$\oplus$} (ar);
\end{tikzpicture}
\end{equation}
We now replace 
the middle terms in the complex from Lemma~\ref{lemma:trace of Qn} by the above expansion.
This yields the
expression for $[Q_n]$
shown in Figure~\ref{fig:traceQnexpanded}.

As indicated by the dashed lines, the complex from Figure~\ref{fig:traceQnexpanded}
can be viewed as a one-sided twisted complex of the form
\[
\begin{tikzpicture}[baseline=0em]
\tikzstyle{every node}=[font=\small]
\node (a) at (0,0) {$C$};
\node (b) at (1.5,0) {$B$};
\node (c) at (3,0) {$A$};
\path[->,>=stealth,shorten >=1pt,auto,node distance=1.8cm,
  thick]
(a) edge node {} (b)
(b) edge node {} (c)
(a) edge [in=150,out=30] node {} (c);
\end{tikzpicture}
 \] 
 where $B$ and the bottom part of $C$ are contractible.
 After collapsing these contractible pieces, this twisted complex becomes
\begin{equation}\label{eqn:traceQncollapsed}
\begin{tikzpicture}[baseline=0em]
\node (a) at (0,0){$\left[\JWpict{4+1}[n-1]\right] $ };
\node (b) at (-3,0){$\left[\JWT\right]$ };
\node (c) at (-6,0){$\left[\JWT\right]$};
\node (d) at (-9,0){ $\left[\JWpict{4+1}[n-1]\right]$};
\path[->,>=stealth',shorten >=1pt,auto,node distance=1.8cm,font=\small]
(b) edge node {$a'$} (a)
(c) edge node {} (b)
(d) edge node {} (c);
\end{tikzpicture}
\end{equation}
We now use induction on $n$ to replace the remaining copies of
$[P_{n-1}\sqcup\one_1]$ and $[P_{n-2}]$ by $V_{n-1}\otimes V_1$ and $V_{n-2}$.
Taking into account the shifts of the cohomological grading (but not of the quantum grading), this yields
\begin{equation}\label{eqn:Qnfirst}
[Q_n]\ \simeq\ \left(
\begin{tikzpicture}[baseline=0em]
\node (a) at (0,0){$V_{n-1}\otimes V_1$ };
\node (b) at (-3,0){$\ts^{-1}V_{n-2}$ };
\node (c) at (-6,0){$\ts^{2-2n}V_{n-2}$};
\node (d) at (-9,0){ $\ts^{1-2n}V_{n-1}\otimes V_1$};
\path[->,>=stealth',shorten >=1pt,auto,node distance=1.8cm,font=\small]
(b) edge node {} (a)
(d) edge node {} (c);
\end{tikzpicture}
\right)
\end{equation}
On the other hand, the proof of Theorem~\ref{thm:periodic P} shows that
\begin{equation}\label{eqn:Qnsecond}
[Q_n]\ \simeq\
\left(
\begin{tikzpicture}[baseline=0em]
\node (a) at (0,0){$[P_n]$ };
\node (b) at (-2.2,0){$[\ts^{1-2n}P_n]$ };
\path[->,>=stealth',shorten >=1pt,auto,node distance=1.8cm,font=\small]
(b) edge node {$[u_n]$} (a);
\end{tikzpicture}
\right)
\ \simeq \
\left(
\begin{tikzpicture}[baseline=0em]
\node (a) at (0,0){$[P_n]$ };
\node (b) at (-2.2,0){$[\ts^{1-2n}P_n]$ };
\path[->,>=stealth',shorten >=1pt,auto,node distance=1.8cm,font=\small]
(b) edge node {$0$} (a);
\end{tikzpicture}
\right)
\ = \
[\ts^{1-2n}P_n]\ \oplus\ [P_n]
\end{equation}
where the second homotopy equivalence follows from Lemma~\ref{lemma:zeroendo} because $u_n$ has nonzero
quantum degree.

Taking into account grading shifts, it is now easy to see that \eqref{eqn:Qnsecond} can only hold
if both arrows in~\eqref{eqn:Qnfirst} are nonzero
and $[P_n]\simeq V_n$. Indeed, we have $V_{n-1}\otimes V_1\cong V_n\oplus V_{n-2}$,
and the two middle terms in~\eqref{eqn:Qnfirst} cannot be connected by an arrow for
grading reasons (here we can assume $n>2$ because the $n=2$ case is straightforward).
\end{proof}

\begin{figure}
\begin{center}
\begin{tikzpicture}[baseline=0em]
\tikzstyle{every node}=[font=\small]
\node (a0) at (-0.15,0) {$\left[\JWpict{4+1-narrow}\right]$};
\node (b0) at (-3,0) {$\left[\JWT\right]$};
\node (b0l) at (-4.8,0) {$\left[\JWT\right]$};
\node (b1) at (-3,-3) {$\left[\JWT\right]$};
\node (b2) at (-3,-6) {$\left[\JWT\right]$};
\node (b3) at (-3,-9) {$\left[\JWT\right]$};
\node (b4) at (-3,-12) {$\left[\JWT\right]$};
\node (b5) at (-3,-14) {$\vdots$};
\node (c2l) at (-6.9,-6) {$\left[\JWT\right]$};
\node (c2r) at (-5.1,-6) {$\left[\JWT\right]$};
\node (c3) at (-6,-9) {$\left[\JWT\right]$};
\node (c4) at (-6,-12) {$\left[\JWT\right]$};
\node (c5) at (-6,-14) {$\vdots$};
\node (d2) at (-9.5,-6) {$\left[\JWpict{4+1-narrow}\right]$};

\path
(c2l) edge[draw=none] node {$\oplus$} (c2r)
(b0l) edge[draw=none] node {$\oplus$} (b0);
\path[->,>=stealth,shorten >=1pt,auto,node distance=1.8cm,
  thick]
(b1) edge node {} (b0)
(b1) edge node {$\id$} (b0l)
(b2) edge node {} (b1)
(b3) edge node {} (b2)
(b4) edge node {} (b3)
(b5) edge node {} (b4)
(c3) edge node {$\id$} (c2l)
(c3) edge node {} (c2r)
(c4) edge node {} (c3)
(c5) edge node {} (c4)
(c4) edge node {$\id$} (b4)
(c3) edge node {$\id$} (b3)
(d2) edge node {} (c2l)
(c2r) edge node {$\id$} (b2)
(b0) edge node {$a'$} (a0);
\draw[rounded corners,dashed,thick,red] (-3.9,2) -- (-3.9,-1) -- (-1.5,-1) -- (-1.5,-14.5);
\draw[rounded corners,dashed,thick,red] (-6,2) -- (-6,-7) -- (-1.7,-7) -- (-1.75,-14.5);
\node at (-1.8,1.3) {\resizebox{0.4cm}{!}{${\color{red} A}$}};
\node at (-5,1.3) {\resizebox{0.4cm}{!}{${\color{red} B}$}};
\node at (-8,1.3) {\resizebox{0.4cm}{!}{${\color{red} C}$}};
\end{tikzpicture}
\end{center}
\caption{Expansion of $[Q_n]$ without longer arrows.}\label{fig:traceQnexpanded}
\end{figure}
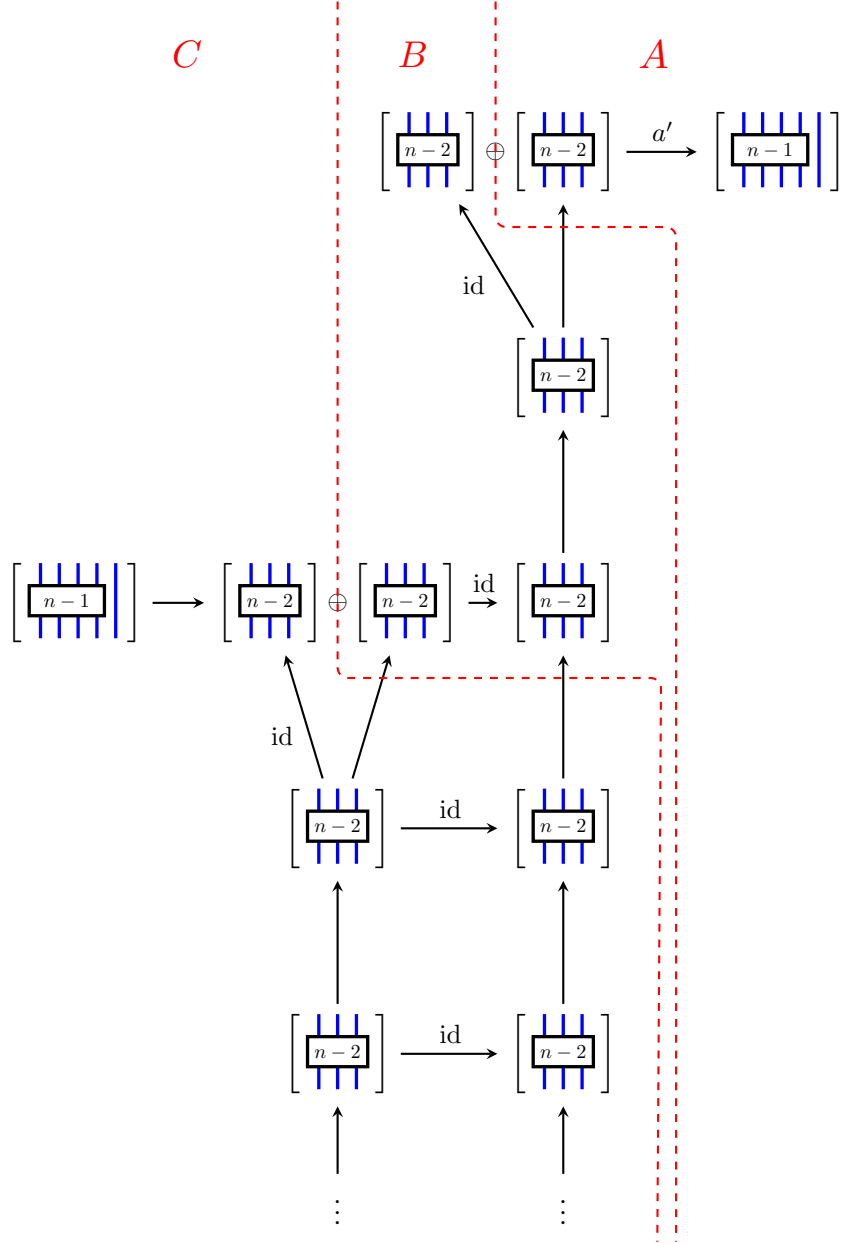

\begin{remark}\label{rem:mapaprime}
The map $a$ in Lemma~\ref{lemma:trace of Qn} is induced by the composition $s\circ\theta^{-1}$ where $s$ is a saddle and $\theta$ is
a chain isomorphism given by cyclicity isomorphisms. From this, one can easily see that the map $a'$
in \eqref{eqn:traceQncollapsed} and in Figure~\ref{fig:traceQnexpanded} is induced by an annulus $[P_{n-2}]\rightarrow[P_{n-2}\sqcup\one_2]$.
On the other hand, the leftmost map in \eqref{eqn:traceQncollapsed} is more difficult to describe because it potentially depends on longer arrows
connecting the leftmost term in Figure~\ref{fig:traceQnexpanded} to the rightmost column of $C$.
As a consequence, it is not clear whether this map is nonzero, although this follows indirectly
from the proof of Lemma~\ref{lemma:PnVn}.
\end{remark}

Although Theorem~\ref{thm:P as model} follows implicitly from
Lemma~\ref{lemma:PnVn} and its proof, we will now give a concrete
proof of this theorem by showing explicitly that $[P_n]$ deformation
retracts to the complex $A$ from Figure~\ref{fig:traceQnexpanded}.
Note, however, that a part of this proof relies on Lemma~\ref{lemma:PnVn}.

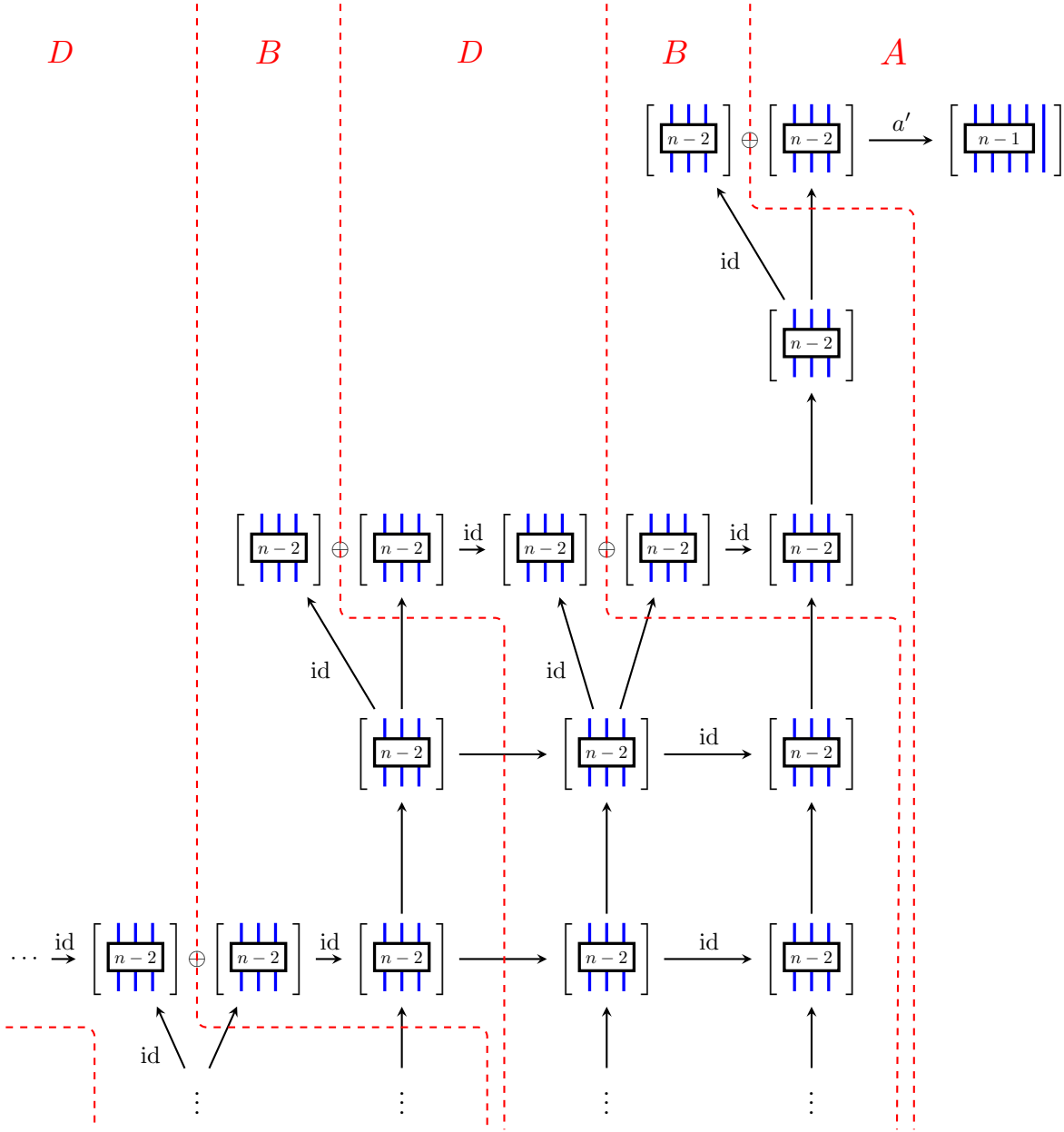
\begin{figure}
\begin{center}
\begin{tikzpicture}[baseline=0em]
\tikzstyle{every node}=[font=\small]
\node (a0) at (-0.15,0) {$\left[\JWpict{4+1-narrow}\right]$};
\node (b0) at (-3,0) {$\left[\JWT\right]$};
\node (b0l) at (-4.8,0) {$\left[\JWT\right]$};
\node (b1) at (-3,-3) {$\left[\JWT\right]$};
\node (b2) at (-3,-6) {$\left[\JWT\right]$};
\node (b3) at (-3,-9) {$\left[\JWT\right]$};
\node (b4) at (-3,-12) {$\left[\JWT\right]$};
\node (b5) at (-3,-14) {$\vdots$};
\node (c2l) at (-6.9,-6) {$\left[\JWT\right]$};
\node (c2r) at (-5.1,-6) {$\left[\JWT\right]$};
\node (c3) at (-6,-9) {$\left[\JWT\right]$};
\node (c4) at (-6,-12) {$\left[\JWT\right]$};
\node (c5) at (-6,-14) {$\vdots$};
\node (d2) at (-9,-6) {$\left[\JWT\right]$};
\node (d2l) at (-10.8,-6) {$\left[\JWT\right]$};
\node (d3) at (-9,-9) {$\left[\JWT\right]$};
\node (d4) at (-9,-12) {$\left[\JWT\right]$};
\node (d5) at (-9,-14) {$\vdots$};
\node (e4l) at (-12.9,-12) {$\left[\JWT\right]$};
\node (e4r) at (-11.1,-12) {$\left[\JWT\right]$};
\node (e5) at (-12,-14) {$\vdots$};
\node (f4) at (-14.5,-12) {$\cdots$};
\path
(e4l) edge[draw=none] node {$\oplus$} (e4r)
(d2l) edge[draw=none] node {$\oplus$} (d2)
(c2l) edge[draw=none] node {$\oplus$} (c2r)
(b0l) edge[draw=none] node {$\oplus$} (b0);
\path[->,>=stealth,shorten >=1pt,auto,node distance=1.8cm,
  thick]
(b1) edge node {} (b0)
(b1) edge node {$\id$} (b0l)
(b2) edge node {} (b1)
(b3) edge node {} (b2)
(b4) edge node {} (b3)
(b5) edge node {} (b4)
(c3) edge node {$\id$} (c2l)
(c3) edge node {} (c2r)
(c4) edge node {} (c3)
(c5) edge node {} (c4)
(d3) edge node {} (d2)
(d3) edge node {$\id$} (d2l)
(d4) edge node {} (d3)
(d5) edge node {} (d4)
(e5) edge node {$\id$} (e4l)
(e5) edge node {} (e4r)
(f4) edge node {$\id$} (e4l)
(e4r) edge node {$\id$} (d4)
(d4) edge node {} (c4)
(c4) edge node {$\id$} (b4)
(d3) edge node {} (c3)
(c3) edge node {$\id$} (b3)
(d2) edge node {$\id$} (c2l)
(c2r) edge node {$\id$} (b2)
(b0) edge node {$a'$} (a0);
\draw[rounded corners,dashed,thick,red] (-3.9,2) -- (-3.9,-1) -- (-1.5,-1) -- (-1.5,-14.5);
\draw[rounded corners,dashed,thick,red] (-6,2) -- (-6,-7) -- (-1.7,-7) -- (-1.75,-14.5);
\draw[rounded corners,dashed,thick,red] (-9.9,2) -- (-9.9,-7) -- (-7.5,-7) -- (-7.5,-14.5);
\draw[rounded corners,dashed,thick,red] (-12,2) -- (-12,-13) -- (-7.75,-13) -- (-7.75,-14.5);
\draw[rounded corners,dashed,thick,red] (-14.8,-13) -- (-13.5,-13) -- (-13.5,-14.5);
\node at (-1.8,1.3) {\resizebox{0.4cm}{!}{${\color{red} A}$}};
\node at (-5,1.3) {\resizebox{0.4cm}{!}{${\color{red} B}$}};
\node at (-8,1.3) {\resizebox{0.4cm}{!}{${\color{red} D}$}};
\node at (-10.95,1.3) {\resizebox{0.4cm}{!}{${\color{red} B}$}};
\node at (-14,1.3) {\resizebox{0.4cm}{!}{${\color{red} D}$}};
\end{tikzpicture}
\end{center}
\caption{Expansion of $[P_n]$ without longer arrows.}\label{fig:tracePnexpanded}
\end{figure}

\begin{proof}[Proof of Theorem~\ref{thm:P as model}]
By using
an analogue of Lemma~\ref{lemma:trace of Qn}, we can write the trace of $P_n$ as
\[
[P_n] \ \simeq \ 
\left(\begin{tikzpicture}[baseline=0em]
\node (a) at (0,0){$\left[\JWpict{4+1}[n-1]\right] $ };
\node (b) at (-3,0){$\left[\JWpict{4-tr}[n-1]\right]$ };
\node (c) at (-6,0){$\left[\JWpict{4-tr}[n-1]\right]$};
\node (d) at (-9,0){ $\left[\JWpict{4-tr}[n-1]\right]$};
\node (e) at (-12,0){ $\cdots$};
\path[->,>=stealth',shorten >=1pt,auto,node distance=1.8cm,font=\small]
(b) edge node {$a$} (a)
(c) edge node {$[u_{n-1}]$} (b)
(d) edge node {$b$} (c)
(e) edge node  {$[u_{n-1}]$} (d);
\end{tikzpicture}
\right)
\]
where $b$
is given by a dot placed on the looped strand.
In this expansion of $[P_n]$, we now replace each term (beyond the first) by
the expansion from~\eqref{eqn:tracedstrandexpansion}.
This yields the expression for $[P_n]$ shown in Figure~\ref{fig:tracePnexpanded}.

Note that the complexes $A$ and $B$ in Figure~\ref{fig:tracePnexpanded}
are identical with the ones from Figure~\ref{fig:traceQnexpanded}.
While the complex $B$ is obviously contractible, it is
not clear whether $D$ is contractible because
there could be longer arrows connecting the leftmost term in $D$
to its rightmost column.

However, a simple grading argument shows that $A$ becomes a direct
summand of $[P_n]$ if one replaces $[P_{n-2}]$
and $[P_{n-1}\sqcup\one_1]$ by $V_{n-2}$ and $V_{n-1}\otimes V_1$.
Because the resulting complex $V_{n-2}\rightarrow V_{n-1}\otimes V_1$ necessarily contains $V_n$
as a direct summand, and because $[P_n]\simeq V_n$ by
Lemma~\ref{lemma:PnVn}, it then follows that the complex
$\ldots\rightarrow D\rightarrow B\rightarrow D\rightarrow B$ in the expansion
of $[P_n]$ must be contractible. Hence $[P_n]$ deformation
retracts to the complex
\[
A\ = \ \left(
\begin{tikzpicture}[baseline=0em]
\node (a) at (0,0){$\left[\JWpict{4+1}[n-1]\right] $ };
\node (b) at (-3,0){$\left[\JWT\right]$ };
\path[->,>=stealth',shorten >=1pt,auto,node distance=1.8cm,font=\small]
(b) edge node {$a'$} (a);
\end{tikzpicture}
\right)
\]
where $a'$ is as described in Remark~\ref{rem:mapaprime}.
The theorem now follows from the definition of a model.
\end{proof}

\section{Quantum traces of generalized projectors}
\label{sec:generalized-projectors}

In~\cite{CH12}, B.~Cooper and the second author introduced a family of complexes $P_\e\in\Ch^-(\cat{BN}_n)$ which
categorify the indecomposable idempotents in $\TL$. In this section, we give a new and streamlined proof of the existence of these complexes.
We then show that the class of $P_\e$ in the quantum horizontal trace
depends only on the through-degree of $P_\e$ (Theorem~\ref{thm:intro Pe}).
As in the previous sections, we need to assume that $1-q^d$ is invertible for $d>0$.

\subsection{Idempotents in $\TL$}
\label{ss:TL idempts}

Recall the Karoubi envelope of a $\k$-linear category $\CS$ is the category $\Kar(\CS)$ whose objects are formal expressions $\im e$, where $e$ is an idempotent endomorphism in $\CS$.  If $e,e'$ are idempotent endomorphisms in $\CS$, then any morphism in $\CS$ of the form $e\circ f\circ e'$ is regarded as a morphism $\im e\leftarrow \im e'$ in $\Kar(\CS)$.

Recall the Temperley--Lieb category $\TL$.  This category has identity maps $\id_n$ indexed by $n\in \Z_{\geq 0}$, and the space of morphisms $\im \id_n \leftarrow \im \id_m$ is the Temperley--Lieb space of $(n,m)$ tangles, denoted $\TL_{n,m}$.  We also denote $\TL_n:=\TL_{n,n}$.

We also have the Jones--Wenzl idempotents $p_n = \id_n\circ p_n \circ \id_n$.

\begin{proposition}
There is a monoidal functor from $\Kar(\TL)$ to the category of finite dimensional $U_q(\sl_2)$ representations sending
\[
\im \id_n \ \mapsto \ V_1^{\otimes n},\qquad \im p_n \mapsto V_n,
\]
where $V_n$ denotes the $n+1$ dimensional simple representation of $U_q(\sl_2)$.  This functor is an equivalence of categories.
\end{proposition}

As is well known
\begin{equation}\label{eq:Vn V1 decomp}
V_n\otimes V_1 \cong V_{n+1}\oplus V_{n-1}
\end{equation}
for $n\geq 1$, and
\[
V_0\otimes V_1\cong V_1.
\]
Since the summands on the right-hand side above are distinct simples there is a unique pair of complementary idempotent endomorphisms of $V_n\otimes V_1$ which project onto summands isomorphic to $V_{n\pm 1}$. 

Iterating \eqref{eq:Vn V1 decomp} allows us to obtain a preferred collection of idempotent endomorphisms which project $V_1^{\otimes n}$ onto its simple direct summands.  The goal of this section is to arrive at such idempotent decompositions inside $\TL$, by means which are amenable to categorification.  We begin by introducing the central idempotents $p_{n,k}\in \TL_n$, which under the functor to $U_q(\sl_2)$ representations become the projections of $V^{\otimes n}$ onto its $V_k$-isotypic summand.

\begin{proposition}\label{prop:central idempotents}
There exist a unique family of elements $p_{n,k}\in \TL_n$ indexed by integers $0\leq k\leq n$ with $k\equiv n$ (mod 2) such that
\begin{enumerate}
\item $\id_n = \sum_k p_{n,k}$.
\item $p_{n,k}$ is a sum of elements factoring through $p_k\in \TL_k$.
\end{enumerate}
\end{proposition}
\begin{proof}
The proof of existence is by induction on $n$.  In the base cases we have $p_{0,0}=\id_0$ and $p_{1,1}=\id_1$. So assume that $p_{n,k}$ have been constructed as in the statement. By hypothesis, $p_{n,k}$ can be written as
\[
p_{n,k} = \sum_i x_{n,k,i} p_k y_{n,k,i}
\]
for some elements $x_{n,k,i}\in \TL_{n,k}$ and $y_{n,k,i}\in \TL_{k,n}$. Adding a strand yields
\[
p_{n,k}\sqcup \id_1  = \sum_k (x_{n,k,i}\sqcup \id_1) (p_k\sqcup \id_1)(y_{n,k,i}\sqcup \id_1).
\]
The recursion satisfied by the Jones--Wenzl idempotents tells us that $p_k\sqcup \id_1$ is a sum of $p_{k+1}$ plus terms which factor through $p_{k-1}$.  It follows that $p_{n,k}\sqcup \id_1$ can be written as a sum of terms factoring through $p_{k+1}$ or $p_{k-1}$, hence by induction that $\id_{n+1}$ can be written as a sum of terms factoring through $p_l$ for various $l$, which completes the inductive proof on the existence of $\{p_{n,k}\}_{n,k}$.

We make some observations.  If $f\in \TL_{n,m}$ factors through $p_k$, then $l\neq k$ implies $f p_{n,l}=0$, hence
\[
f = f\sum_l p_{n,l} = fp_{n,k}.
\]
A similar argument shows that $f = p_{m,k}f$.  Now, suppose that $f\in \TL_n$ is a sum of terms factoring through $p_k$, and $\id_n-f$ is a sum of terms factoring through $p_l$ for various $l\neq k$.  Then $p_{n,k} f = f$ and $p_{n,k}(\id_n - f) =0$.  In other words,
\[
f = p_{n,k}f = p_{n,k}.
\]
This proves the uniqueness of $p_{n,k}$.
\end{proof}

Additional properties of the $p_{n,k}$ are recorded below.

\begin{corollary}\label{cor:pnk are idemp}
The $p_{n,k}$ are a complete collection of mutually orthogonal idempotents in $\TL_n$, with $p_{n,n}=p_n$.\qed
\end{corollary}

\begin{corollary}\label{cor:pnk are central}
The $p_{n,k}$ are central, in the sense that $p_{n,k} f = fp_{m,k}$ for all $f\in \TL_{n,m}$ (interpreted as zero if $n<k$ or $m<k$).
\end{corollary}
\begin{proof}
If $f$ factors through $p_k$ for some $k$ then $f = p_{n,k}f = fp_{m,k}$.  Any $f\in \TL_{n,m}$ is a sum of terms factoring through some $p_k$, which completes the proof.
\end{proof}

\begin{corollary}\label{cor:p steps}
We have $(p_{n-1,k}\sqcup \id_1)p_{n,l} =0$ unless $l=k+1$ or $l=k-1$.
\end{corollary}
\begin{proof}
We have seen that $p_{n-1,k}\sqcup \id_1$ is a sum of terms factoring through $p_{k+1}$ or $p_{k-l}$; the corollary follows.
\end{proof}

We remark that the $p_{n,k}$ are not primitive idempotents.  To construct a complete collection of primitive idempotents in $\TL_n$ we use the following.

\begin{definition}\label{def:pe}
A sequence $\e\in \{1,-1\}^n$ will be called \emph{admissible} if $\e_1+\cdots+\e_i\geq 0$ for all $1\leq i\leq n$.  For each admissible $\e\in \{1,-1\}$ define an element $p_\e\in \TL_n$ by the formula
\[
p_\e = \prod_{i=1}^n (p_{i, \e_1+\cdots+\e_i}\sqcup \id_{n-i})
\]
\end{definition}

\begin{proposition}\label{prop:pe are complete idemp}
The set $\{p_\e\:|\: \e\in \{1,-1\}^n \text{ is admissible}\}$ is a complete collection of orthogonal idempotents in $\TL_n$.  Morover, $\im p_\e \cong \im p_{|\e|}$ in $\Kar(\TL)$; in particular the idempotents $p_\e$ are primitive.
\end{proposition}
\begin{proof}
We multiply together our idempotent decompositions of $\id_n$:
\[
\id_n = \prod_{m=1}^n \left(\sum_k p_{m,k}\sqcup \id_{n-m}\right) = \sum_{k_1,\ldots,k_n} p_{1,k_1}' p_{2,k_2}'\cdots p_{n,k_n}',
\]
where we are abbreviating by writing $p_{m,k}':=p_{m,k}\sqcup \one_{n-m}$, and the sum on the right is indexed by sequences $(k_1,\ldots,k_n)$ with $0\leq k_m\leq m$ and $k_m\equiv m$ (mod 2). Observe that all of the factors $p_{m,k}'$ commute with one another, hence the   $p_{1,k_1}' p_{2,k_2}'\cdots p_{n,k_n}'$ are orthogonal idempotents.   Corollary \ref{cor:p steps} says that $p_{1,k_1}' p_{2,k_2}'\cdots p_{n,k_n}'$ is zero unless $k_m = k_{m-1}\pm 1$ for all $2\leq m\leq n$.  We obtain an idempotent decomposition
\[
\id_n = \sum_{\e} p_\e
\]
indexed by admissible sequences $\e\in \{1,-1\}^n$ as claimed. 

 To see that $\im p_\e \cong \im p_{|\e|}$ in $\Kar(\TL)$, we argue by induction on $n$.  The statement is trivial for $n=0$ or $n=1$, since in these cases $p_\e=p_n=\id_n$.  Assume by induction that $\im p_\e\cong \im p_k$, where $k=|\e|$. The decomposition $p_\e \sqcup \id_1 = p_{(\e, 1)} + p_{(\e, -1)}$ yields an isomorphism
 \[
\im(p_\e)\otimes \im(\id_1) \cong \im(p_{(\e, 1)})\oplus \im(p_{(\e, -1)}).
 \]
On the other hand,
\[
\im(p_\e)\otimes \im(\id_1)\cong \im(p_k)\otimes \im(\id_1)\cong \im(p_{k+1})\oplus \im(p_{k-1})
\]
Since $p_{(\e, \pm 1)}$ factors through $p_{k\pm 1}$, the only possibility is that $\im(p_{(\e, \pm 1)})\cong \im(p_{k\pm 1})$.  This completes the proof.
\end{proof}

\subsection{The categorified central idempotents}
\label{ss:cat Pnk}

In this section we describe the categorical analogues of the idempotent decompositions from \S \ref{ss:TL idempts}. As discussed in \S \ref{ss:CK model}, the Cooper--Krushkal projectors $P_n\in \Ch^-(\BN_n)$ are categorical analogues of the Jones--Wenzl idempotents $p_n\in \TL_n$. 

\begin{definition}
Given integers $0\leq k\leq n$ with $k\equiv n$ (mod 2), we let $\mathcal{X}_{n,k}$ denote the collection of complexes which are direct sums of $q$-shifts of complexes of the form $T\star P_k \star T'$, where $T,T'$ are objects of $\BN_{n,k}$ and $\BN_{k,n}$, respectively.  We say that $X\in \Ch^-(\BN_n)$ \emph{factors through} $P_k$ if $X\simeq X'$ with $X'\in  \langle \mathcal{X}_{n,k}\rangle$.
\end{definition}

Note that we do not allow cohomological shifts in the definition, so each complex in $\mathcal{X}_{n,k}$ is supported in cohomological degrees $\leq 0$.

\begin{lemma}\label{lemma:cat branching rule}
If $X\in \Ch^-(\BN_n)$ factors through $P_k$ then $X\sqcup \one_1\in \Ch^-(\BN_{n+1})$  is homotopy equivalent to a complex factoring through $P_{k+1}$ and $P_{k-1}$. 
\end{lemma}
\begin{proof}
The recursion satisfied by Cooper-Krushkal projectors ensures that if $X=T\star P_k\star T'$ then $X\sqcup \one_1 = (T\sqcup \one_1)\star (P_k\sqcup \one_1)\star (T'\sqcup \one_1)$ is homotopy equivalent to a complex in $\langle\mathcal{X}_{n+1,k+1}\cup \mathcal{X}_{n+1,k-1}\rangle$.  Then Lemma \ref{lemma:S in T} completes the proof.
\end{proof}

\begin{lemma}\label{lemma:homs between pks}
Suppose we have $X\in \mathcal{X}_{n,k}$ and $Y\in \mathcal{X}_{n,l}$.  We have:
\begin{enumerate}
\item if $k<l$ then each chain object $X^m$ satisfies $\Hom(X^m,Y)\simeq 0$.
\item if $k=l$ and $m<0$, then  $\Hom(X^m,Y)\simeq 0$.
\end{enumerate}
\end{lemma}
\begin{proof}
We use the duality isomorphism $\Hom(X^m,Y)\cong \Hom(\one_n, (X^m)^\vee\star Y)$ together with the fact that the through degree of $X^m$ is $\leq k$, with equality if and only if $m=0$. Since $Y$ factors through $P_l$, $X^m\star Y$ is contractible if the through degree of $X^m$ is strictly less than $l$, which will occur if $k<l$ or $k=l$ and $m<0$.
\end{proof}

\begin{theorem}\label{thm:Pnk}
There exist a unique family of complexes $P_{n,k}\in \Ch^-(\BN_n)$ such that
\begin{enumerate}
\item $\one_n \simeq \left(P_{n,n}\rightarrow P_{n,n-2}\rightarrow \cdots \rightarrow P_{n,\text{0 or 1}}\right)$.
\item $P_{n,k}$ factors through $P_k$.
\end{enumerate}
\end{theorem}
\begin{proof}
Fix $n$, and let $\Omega = \{k\:|\: 0\leq k\leq n, k\equiv n \text{ (mod 2)}\}$ with order $k\unlhd l$ if $k>l$ (note the reversal of order). Iterating Lemma \ref{lemma:cat branching rule} we see that $\one_n$ is homotopy equivalent to a complex $E\in \langle \mathcal{X}_{n,k} \:|\: k\in \Omega\rangle$. The rest of the construction is achieved by combing hairs.

To be precise,  note that the collections of complexes $\mathcal{X}_{n,k}$ (indexed by $k\in \Omega$)  satisfy the hypotheses of Lemma \ref{lemma:combing} by Lemma \ref{lemma:homs between pks}. Combing hairs results in a complex $E'\simeq \one_n$ which can be written as $E=\tw_{\d'}(\bigoplus_{i\in I'} \ts^{a(i)} X_i)$ in which $I$ is a partially ordered set, $X_i\in \mathcal{X}_{\omega(i)}$, such $i\leq j$ implies $\omega(i)\lhd \omega(j)$ or $\omega(i)=\omega(j)$ and $a(i)<a(j)$.  We may collect terms according to the through degree $k$ and the cohomological shift $a$,  obtaining
\[
\one_n\simeq \tw_{\d+\e}\left(\bigoplus_{k\in \Omega}\bigoplus_{a\leq 0} E'_{k,a} \right)
\]
where the twist $\d$ maps $E'_{k,a}$ to $E'_{l,b}$ only if $k\lhd l$ and the twist $\e$ maps $E'_{k,a}$ to $E'_{l,b}$ only if $k=l$ and $a\leq b$.  Then the complexes $P_{n,k}=\tw_\e(\bigoplus_{a\leq 0} E'_{k,a})$ are as in the statement.
\end{proof}

\begin{proposition}
Given $X\in \Ch^-(\BN_{n,m})$ we have $P_{n,k}\star X\simeq X\star P_{m,k}$, naturally in $X$ up to homotopy.
\end{proposition}
\begin{proof}
Define complexes
\[
P_{n,\leq k} = (P_{n,k}\rightarrow P_{n,k-2}\rightarrow \cdots \rightarrow P_{n,\text{0 or 1}}),
\]
\[
P_{n,\geq k} =(P_{n,n}\rightarrow P_{n,n-2}\rightarrow \cdots \rightarrow P_{n,k}).
\]
Then $P_{n,k}\simeq P_{n,\geq k}\star P_{n,\leq k}$.  So it suffices to show that $P_{n,\geq k}$ and $P_{n,\leq k}$ are central, as in the statement.  The centrality of these complexes follows from 
\cite[Theorem~4.24]{Hog17a}.
\end{proof}

The following is straightforward.
\begin{lemma}\label{lemma:Pnk killing}
Suppose $X\in \Ch^-(\BN_n)$ factors through $P_k$.  Then $P_{n,l}\star X\simeq 0$ if $l\neq k$, and $P_{n,k}\star X\simeq X$. with $l\neq k$, then $P_{n,k}\star X\simeq 0$.
\end{lemma}
\begin{proof}
Suppose $X,Y\in \Ch^-(\BN_n)$ are complexes such that $X$ factors through $P_k$ and $Y$ factors through $P_l$.  Then $X\star Y\simeq 0$.  In particular, $P_{n,l}\star X\simeq 0$.  Using the idempotent decomposition of $\one_n$ into central idempotents, we obtain:
\[
X\simeq \tw_\d\left(\bigoplus_k P_{n,k}\right)\star X \cong \tw_\d\left(\bigoplus_k P_{n,k}\star X\right) \simeq P_{n,k}\star X,
\]
where in the last step we contracted the contractible complexes $P_{n,l}\star X$ ($l\neq k$) using homological perturbation.
\end{proof}

\subsection{The primitive categorified idempotents}
\label{ss:cat Pe}

The following defines the categorical analogues of the idempotent elements $p_\e\in \TL_n$ from Definition \ref{def:pe}.

\begin{definition}\label{def:cat Pe}
For each admissible $\e\in \{1,-1\}$ define $P_\e \in \Ch^-(\BN_n)$ by the formula
\[
P_\e = \star_{i=1}^n (P_{i, \e_1+\cdots+\e_i}\sqcup \one_{n-i}).
\]
\end{definition}

\begin{proposition}\label{prop:cat pe are complete idemp}
We have $P_\e\star P_\nu\simeq 0$ if $\e\neq \nu$, and $\one_n\simeq \tw_\d(\bigoplus_\e P_\e)$ in which the twist $\d$ strictly decreases the partial order $(\e_1,\ldots,\e_n)T\geq (\nu_1,\ldots,\nu_n)$ if $\e_1+\cdots+\e_i\geq \nu_1+\cdots+\nu_i$ for all $i=1,\ldots,n$.
\end{proposition}
\begin{proof}
This is analagous to Proposition \ref{prop:pe are complete idemp}, and we leave many details to the reader.   We tensor together our idempotent decompositions of $\one_n=\one_m\sqcup \one_{n-m}$ for $1\leq m\leq n$:
\[
\id_n \simeq \star_{m=1}^n \left(\sum_k P_{m,k}\sqcup \id_{n-m}\right) = \bigoplus_{k_1,\ldots,k_n} P_{1,k_1}' P_{2,k_2}'\cdots P_{n,k_n}',
\]
with a twist that decreases one or more indices $k_i$, 
where we are abbreviating by writing $P_{m,k}':=P_{m,k}\sqcup \one_{n-m}$, and the sum on the right is indexed by sequences $(k_1,\ldots,k_n)$ with $0\leq k_m\leq m$ and $k_m\equiv m$ (mod 2). Observe that all of the factors $P_{m,k}'$ commute with one another up to homotopy, hence the   $P_{1,k_1}' P_{2,k_2}'\cdots P_{n,k_n}'$ are orthogonal idempotents up to homotopy.   Lemma \ref{lemma:Pnk killing} and Lemma \ref{lemma:cat branching rule} tell us that $P_{1,k_1}' P_{2,k_2}'\cdots P_{n,k_n}'$ is contractible unless $k_m = k_{m-1}\pm 1$ for all $2\leq m\leq n$.  Homological perturbation gives an idempotent decomposition
\[
\one_n \simeq \tw_\d\left(\bigoplus_\e P_\e\right)
\]
in which the twist decreases one or more indices $k_i$.
\end{proof}

\subsection{The quantum annular trace of categorified idempotents}
\label{ss:trace of Pe}

\begin{reptheorem}{thm:intro Pe}
For every admissible $\e\in \{1,-1\}^n$ we have $[P_\e]\simeq V_{|\e|}$ in $\Ch^-(\Kar(\cat{BN}_{\!q}(\Ann)))$.
\end{reptheorem}
\begin{proof}
Induction on $n$.  In the base case $n=0$ there is nothing to prove.  Assume by induction that $[P_\e]\simeq V_{k}$, where $k=|\e|$.  Then $[P_\e\sqcup \one_1]\simeq V_{k+1}\oplus V_{k-1}$.  But also $[P_\e\sqcup \one_1]\simeq ([P_{(\e, 1)}]\rightarrow [P_{(\e, -1)}])$, by construction.  Computations in the trace tell us that $P_{(\e, \pm 1)}$ is homotopy equivalent to a complex built from copies of $[P_{k\pm 1}]\simeq V_{k\pm 1}$.  The only possibility is $[P_{(\e,\pm 1)}]\simeq V_{k\pm 1}$.
\end{proof}

\begin{corollary}\label{cor:tracePnk} We have $[P_{n,k}]\simeq V_k^{\oplus C_{n,k}}\cong\im(p_{n,k})$ where $C_{n,k}$ is the number of admissible
$\e\in\{1,-1\}^n$ with $|\e|=k$. Explicitly,
$C_{n,k}=\frac{k+1}{m+1}\binom{n}{m}$
where $m=(n+k)/2$.
\end{corollary}
\begin{proof} It follows from Proposition~\ref{prop:cat pe are complete idemp} that $P_{n,k}$ is homotopy equivalent
to a one-sided twisted complex built from idempotents $P_\e$ with $|\e|=k$. Because of Theorem~\ref{thm:intro Pe},
this implies that $[P_{n,k}]$ is a one-sided twisted complex built from copies of $V_k$, and since all of these
copies are in cohomological degree zero, there can be no differentials between them.

The explicit formula for $C_{n,k}$ can obtained from the hook length formula by identifying admissible sequences $\e\in\{1,-1\}^n$ 
with $|\e|=k$ with standard Young tableaux of shape $\lambda=(m,n-m)$. Specifically, an admissible sequence $\e\in\{1,-1\}^n$
corresponds to a standard two-row Young tableau whose first row contains the indices $i$ with $\e_i=1$,
and whose second row contains the indices $i$ with $\e_i=-1$.
\end{proof}

\section{Quantum Hochschild homology of $H^n$}
As an application of the results from the previous section,
we will now compute the full quantum Hochschild homology
of Khovanov's arc ring $H^n$ \cite{Kh02}.
We will first recall the definition of $H^n$.

Let $B^n$ denote the set of all isotopy classes
of flat tangles without closed components,
and with $2n$ upper endpoints and no lower endpoints. Given $a,b\in B^n$,
let $a^\vee$ denote the reflection of $a$ along a horizontal line, and
$a^\vee b:=a^\vee\star b$ denote the vertical composition of $a^\vee$ and $b$.
Note that $a^\vee b$ is a disjoint collection of closed components in the
plane $\R^2$.

By applying Khovanov's functor from~\cite[Subs.~2.3]{Kh00} (for $c=0$), we can assign
to each pair of elements $a,b\in B^n$
a graded $\Bbbk$-module
\[
{}_aH^n_b:=\qs^n\mathcal{F}(a^\vee b).
\]
Since $\mathcal{F}$ is naturally isomorphic to $\Hom_{\mathbf{BN}_{0,0}}(\emptyset,-)$ \cite[Subs~9.1]{B-N05}, this $\Bbbk$-module can
also be viewed as the morphism set
\[
\qs^n\Hom_{\mathbf{BN}_{0,0}}(\emptyset,a^\vee b)\cong \Hom_{\mathbf{BN}_{2n,2n}}(a,b).
\]
As a graded $\Bbbk$-module, Khovanov's arc algebra $H^n$ can now be defined as the direct sum
\[
H^n:=\bigoplus_{a,b}{}_aH^n_b
\]
for $a,b\in B^n$. The algebra multiplication is induced by cobordism maps
\[
{}_aH^n_b\otimes_{\Bbbk}{}_bH^n_c\longrightarrow{}_aH^n_c,
\]
or equivalently, by the composition in
$\mathbf{BN}_{2n,2n}$.

The algebra $H^n$ comes with a distinguished collection of orthogonal idempotents $(a)\in H^n$, one for
each $a\in B^n$, which are given by labeling each circle in $a^\vee\star a$ by the generator $1$
of Khovanov's Frobenius algebra $\Bbbk[x]/(x^2)$.
Correspondingly, there are projective left modules
\[
H^n_a:=H^n(a)=\bigoplus_b{}_bH^n_a
\]
and projective right modules
\[
{}_aH^n:=(a)H^n=\bigoplus_b{}_aH^n_b.
\]
Khovanov's construction further gives rise to a functor
\[
\mathcal{F}_n\colon\mathbf{BN}_{2n,2n}\longrightarrow H^n\mathrm{-}\mathbf{Bimod}
\]
which sends a flat $(2n,2n)$ tangle $T$ to a sweet graded $H^n$-bimodule $\mathcal{F}_n(T)$, in the sense of~\cite[Def.~1]{Kh02}.
As a $\Bbbk$-module, this bimodule decomposes as
\[
\mathcal{F}_n(T):=\bigoplus_{a,b}{}_a\mathcal{F}(T)_b
\]
for $a,b\in B^n$, where
\[{}_a\mathcal{F}(T)_b:=\qs^n\mathcal{F}(a^\vee\star T\star b).\]
Moreover, Khovanov proved that the functors $\mathcal{F}_n$ extend to a bifunctor
\[
\mathbf{F}\colon\mathbf{BN}_{even}\longrightarrow\mathbf{Bimod},
\]
where $\mathbf{Bimod}$ denotes the bicategory of sweet graded bimodules,
and $\mathbf{BN}_{even}\subset\mathbf{BN}$ denotes the bicategory whose
objects are given by even integers $2n\geq 0$, and whose $1$-morphism categories
are identical with the ones in $\mathbf{BN}$.

Before proceeding,
let us also recall that the \emph{$0$th quantum Hochschild homology}
of a graded bimodule $M$ over a graded algebra $A=\bigoplus_dA_d$ is given by the \emph{quantum space of coinvariants}
\[
\operatorname{coInv}_q(M):= M/[A,M]_q,
\]
where \[[A,M]_q: = \operatorname{Span}_{\Bbbk}\{am- q^dma\ |\ a\in A_d, m\in M \}.\]
More generally, the \emph{full quantum Hochschild homology} of $M$ can be computed as the homology of
the complex
$\operatorname{coInv}_q(\mathcal{P})$
where $\mathcal{P}\rightarrow M$ is a projective resolution of the graded bimodule $M$.

It was observed in~\cite{BPW19} that
the assignment $M\mapsto\operatorname{coInv}_q(M)$ extends to a \emph{quantum} \emph{shadow} on the bicategory
$\mathbf{Bimod}$. Since this bicategory has left and right duals,  Theorem~3.5 from \cite{BPW19} implies that this quantum shadow factors
through a functor
\[
\overline{\operatorname{coInv}_q}\colon\operatorname{hTr}_q(\mathbf{Bimod})\longrightarrow\Bbbk\mathrm{-mod}.
\]

We can now prove 
the following theorem, which holds
under the assumption that $1-q^d$ is invertible for all integers $d$ in $(0,4n]$.

\begin{theorem}\label{thm:qhh} The higher quantum Hochschild homology of $H^n$ vanishes
and the $0$th quantum Hochschild homology of $H^n$ is isomorphic to $\Bbbk^{\oplus C_{2n,0}}$ where
$C_{2n,0}=|B^n|$ is the $n$th Catalan number.
\end{theorem}

\begin{proof}
By Proposition~\ref{prop:cat pe are complete idemp}, we have a homotopy equivalence
\[
\one_{2n}\ \simeq\ \Bigl(P_{2n,2n}\longrightarrow P_{2n,2n-2}\longrightarrow\ldots\longrightarrow P_{2n,0}\Bigr)
\]
and applying $\mathcal{F}_n$ gives
\begin{equation}\label{eqn:Hnexpansion}
H^n\ \simeq\ \Bigl(\mathcal{P}_{2n,2n}\longrightarrow \mathcal{P}_{2n,2n-2}\longrightarrow\ldots\longrightarrow \mathcal{P}_{2n,0}\Bigr)
\end{equation}
where $\mathcal{P}_{2n,k}:=\mathcal{F}_n(P_{2n,k})$.

Now observe that, as a complex of $\Bbbk$-modules, $\mathcal{P}_{2n,k}$
is built from complexes of the form
\[{}_a\mathcal{F}_n(T\star P_k\star T')_b=\mathcal{F}(a^\vee\star T\star P_k\star T'\star b)\]
for $a,b\in B^n$. In particular, $\mathcal{P}_{2n,k}$ is acyclic for $k>0$
because $P_k$ kills turnbacks.

On the other hand, $P_{2n,0}$ is built from flat tangles with through-degree $0$,
and after removing closed components, every such tangle can be written
as $a\star(b^\vee)$ for $a,b\in B^n$.
Since closed components in a flat tangle $T$ don't affect the bimodule multiplication
on $\mathcal{F}_n(T)$,
this shows that $\mathcal{P}_{2n,0}$ is built from projective bimodules of the form
$\mathcal{F}_n(a\star (b^\vee))=H^n_a\otimes_\Bbbk {}_bH^n$
for $a,b\in B^n$.

In summary, we see that the inclusion of the subcomplex $\mathcal{P}_{2n,0}$
into the right-hand side of \eqref{eqn:Hnexpansion} defines a quasi-isomorphism
$\mathcal{P}_{2n,0}\rightarrow H^n$, and hence $\mathcal{P}_{2n,0}$
is a projective resolution of $H^n$.
By the remarks preceding this proof, we can therefore compute
the quantum Hochschild homology of $H^n$ from the complex
\[
\operatorname{coInv}_q(\mathcal{P}_{2n,0}) =
\overline{\operatorname{coInv}_q}([\mathcal{P}_{2n,0}])
=\overline{\operatorname{coInv}_q\circ\mathbf{F}}([P_{2n,0}])
\]
The theorem now follows from Corollary~\ref{cor:tracePnk} for $k=0$.
\end{proof}

\begin{remark}
An analogous result for the Chen--Khovanov algebras $A^n$ was shown in \cite[Prop.~6.6]{BPW19}.
Unlike Theorem~\ref{thm:qhh}, the latter result holds without any restrictions on $q$.
\end{remark}

\begin{remark}
Explicitly, the $0$th quantum Hochschild homology of $H^n$ is spanned
by the idempotents $(a)$ for $a\in B^n$. Indeed, these idempotents
remain linearly independent
in $\operatorname{coInv}_q(H^n)$ because they are mutually orthogonal and
span the degree $0$ part of the positively graded algebra $H^n$.
\end{remark}

\printbibliography

@article{BHLZ,
	author = {Beliakova, A. and Habiro, K. and Lauda, A. and {\v{Z}}ivkovi{\'c}, M.},
	journal = {Math. Ann.},
	pages = {397--440},
	title = {Trace decategorification of categorified quantum $\mathfrak{sl}_2$},
	volume = {367},
	year = {2017}}

@article{GHW22,
	author = {Gorsky, E. and Hogancamp, M. and Wedrich, P.},
	journal = {Int. Math. Res. Not.},
	number = {15},
	pages = {11304--11400},
	title = {Derived Traces of Soergel Categories},
	volume = {2022},
	year = {2022}}

@article{APS,
	author = {Asaeda, M.~M.~ and Przytycki, J.~H.~ and Sikora, A.~S.},
	journal = {Algebr. Geom. Topol.},
	pages = {1177--1210},
	title = {Categorification of the Kauffman bracket skein module of $I$-bundles over surfaces},
	volume = {4},
	year = {2004}}

@article{BPW19,
	author = {Beliakova, A.~ and Putyra, K.~ and Wehrli, S.},
	journal = {Inventiones mathematicae},
	pages = {383--492},
	title = {Quantum link homology via trace functor I},
	volume = {215},
	year = {2019}}

@article{KronMro11,
	author = {Kronheimer, P.~B. and Mrowka, T.~S.},
	journal = {Publications math{\'e}matiques de l'IH{\'E}S},
	pages = {97--208},
	title = {Khovanov homology is an unknot-detector},
	volume = {113},
	year = {2011}}

@article{K}

@article{MarklIdeal,
	author = {Markl, M.~},
	doi = {10.1081/AGB-100106814},
	eprint = {https://www.tandfonline.com/doi/pdf/10.1081/AGB-100106814},
	fjournal = {Communications in Algebra},
	issn = {0092-7872},
	journal = {Comm. Algebra},
	mrclass = {16E05 (16E45 18D50)},
	mrnumber = {1856940},
	mrreviewer = {Michel Van den Bergh},
	number = {11},
	pages = {5209--5232},
	title = {Ideal perturbation lemma},
	url = {https://doi.org/10.1081/AGB-100106814},
	volume = {29},
	year = {2001}}

@online{Hog17a,
	author = {Hogancamp, M.~},
	eprint = {1703.01001},
	eprinttype = {arxiv},
	title = {Idempotents in triangulated monoidal categories},
	year = {2017}}

@article{B-N05,
	author = {Bar-Natan, D.},
	doi = {10.2140/gt.2005.9.1443},
	fjournal = {Geometry and Topology},
	issn = {1465-3060},
	journal = {Geom. Topol.},
	mrclass = {57M27 (57M25 57R56)},
	mrnumber = {2174270 (2006g:57017)},
	mrreviewer = {Justin Sawon},
	pages = {1443--1499},
	title = {Khovanov's homology for tangles and cobordisms},
	url = {http://dx.doi.org/10.2140/gt.2005.9.1443},
	volume = {9},
	year = {2005}}

@article{CK12a,
	author = {Cooper, B.~ and Krushkal, V.~},
	doi = {10.4171/QT/27},
	fjournal = {Quantum Topology},
	issn = {1663-487X},
	journal = {Quantum Topol.},
	mrclass = {57R56 (57Mxx)},
	mrnumber = {2901969},
	number = {2},
	pages = {139--180},
	title = {Categorification of the {J}ones-{W}enzl projectors},
	url = {http://dx.doi.org/10.4171/QT/27},
	volume = {3},
	year = {2012}}

@article{CH12,
	author = {Cooper, B. and Hogancamp, M.},
	doi = {10.2140/agt.2015.15.2657},
	fjournal = {Algebraic \& Geometric Topology},
	journal = {Algebr. Geom. Topol.},
	number = {5},
	pages = {2657--2705},
	title = {An Exceptional Collection For Khovanov Homology},
	url = {http://dx.doi.org/10.2140/agt.2015.15.2657},
	volume = {15},
	year = {2015}}

@online{H12a,
	author = {Hogancamp, M.},
	eprint = {1209.2732},
	eprinttype = {arxiv},
	title = {{Morphisms between categorified spin networks}},
	year = {2012}}

@online{H14a,
	author = {Hogancamp, M.},
	eprint = {1405.2574},
	eprinttype = {arxiv},
	title = {{A polynomial action on colored sl(2) link homology}},
	year = {2014}}

@online{BHPW1-pp,
	author = {Beliakova, A.~ and Hogancamp, M.~ and Putyra, K.~ and Wehrli, S.},
	eprint = {1903.12194},
	eprinttype = {arxiv},
	title = {On the functoriality of sl(2) tangle homology},
	year = {2019}}

@article{HogSym-GT,
	author = {Hogancamp, M.},
	doi = {10.2140/gt.2018.22.2943},
	fjournal = {Geometry \& Topology},
	issn = {1465-3060},
	journal = {Geom. Topol.},
	mrclass = {18G60 (57M27)},
	mrnumber = {3811775},
	number = {5},
	pages = {2943--3002},
	title = {Categorified {Y}oung symmetrizers and stable homology of torus links},
	url = {https://doi.org/10.2140/gt.2018.22.2943},
	volume = {22},
	year = {2018}}

@article{Kh05,
	author = {Khovanov, M.},
	doi = {10.1142/S0218216505003750},
	fjournal = {Journal of Knot Theory and its Ramifications},
	issn = {0218-2165},
	journal = {J. Knot Theory Ramifications},
	mrclass = {57M27},
	mrnumber = {2124557 (2006a:57016)},
	mrreviewer = {Marta M. Asaeda},
	number = {1},
	pages = {111--130},
	title = {Categorifications of the colored {J}ones polynomial},
	url = {http://dx.doi.org/10.1142/S0218216505003750},
	volume = {14},
	year = {2005}}

@article{Kh02,
	author = {Khovanov, M.},
	doi = {10.2140/agt.2002.2.665},
	fjournal = {Algebraic \& Geometric Topology},
	issn = {1472-2747},
	journal = {Algebr. Geom. Topol.},
	mrclass = {57M27 (57R56)},
	mrnumber = {1928174 (2004d:57016)},
	mrreviewer = {Jacob Andrew Rasmussen},
	pages = {665--741},
	title = {A functor-valued invariant of tangles},
	url = {http://dx.doi.org/10.2140/agt.2002.2.665},
	volume = {2},
	year = {2002}}

@article{Kh00,
	author = {Khovanov, M.},
	coden = {DUMJAO},
	doi = {10.1215/S0012-7094-00-10131-7},
	fjournal = {Duke Mathematical Journal},
	issn = {0012-7094},
	journal = {Duke Math. J.},
	mrclass = {57M27 (57R56)},
	mrnumber = {1740682 (2002j:57025)},
	number = {3},
	pages = {359--426},
	title = {A categorification of the {J}ones polynomial},
	url = {http://dx.doi.org/10.1215/S0012-7094-00-10131-7},
	volume = {101},
	year = {2000}}

@article{BD14,
	author = {Boyarchenko, M. and Drinfeld, V.},
	doi = {10.1007/s00029-013-0133-7},
	fjournal = {Selecta Mathematica. New Series},
	issn = {1022-1824},
	journal = {Selecta Math. (N.S.)},
	mrclass = {20Gxx (20C15)},
	mrnumber = {3147415},
	number = {1},
	pages = {125--235},
	title = {Character sheaves on unipotent groups in positive characteristic: foundations},
	url = {http://dx.doi.org/10.1007/s00029-013-0133-7},
	volume = {20},
	year = {2014}}

\end{document}